\documentclass[11pt]{amsart}

\usepackage{amsmath}
\usepackage{amsxtra}
\usepackage{amscd}
\usepackage{amsthm}
\usepackage{amsfonts}
\usepackage{amssymb}
\usepackage{mathtools}
\usepackage{eucal}
\usepackage[all]{xy}
\usepackage{graphicx}
\usepackage{tikz-cd}
\usepackage{mathrsfs}
\usepackage{euler}
\DeclareFontShape{OT1}{cmr}{m}{up}{<->ssub*cmr/m/n}{}
\DeclareFontShape{OT1}{cmr}{b}{up}{<->ssub*cmr/b/n}{}
\usepackage{color}
\usepackage{longtable}
\usepackage{caption}
\usepackage{enumerate}
\usepackage{faktor}
\usepackage{bbm}
\usepackage{rotating}
\usepackage[T1]{fontenc}
\usepackage[utf8]{inputenc}
\usepackage{tikz}
\usepackage{tikz-cd}
\usetikzlibrary{matrix,arrows,decorations.pathmorphing}
\usepackage{hyperref}
\usepackage{float}
\usepackage{placeins}

\makeatletter
\@ifpackagelater{hyperref}{2023/11/09}{}{%
  \@ifpackagelater{hyperref}{2023/10/21}{%
    \AtBeginDocument{%
      \renewcommand\Hy@theorem@refstepcounter[1]{%
        \H@refstepcounter{#1}%
        \hyper@makecurrent{#1}%
        \global\let\Hy@dth@currentHref\@currentHref
        \AddToHookNext{cmd/deferred@thm@head/before}{%
          \AddToHookNext{para/begin}{%
            \let\Hy@tempa\@currentHref
            \MakeLinkTarget*{\Hy@dth@currentHref}%
            \global\let\@currentHref\Hy@tempa
          }%
        }%
      }%
    }%
  }{}%
}
\makeatother

\theoremstyle{plain}
\newtheorem{maintheorem}{Theorem}

\newtheorem{maincorollary}[maintheorem]{Corollary}
\newtheorem{theorem}{Theorem}[section]

\newtheorem{proposition}[theorem]{Proposition}

\newtheorem{lemma}[theorem]{Lemma}

\newtheorem{corollary}[theorem]{Corollary}

\theoremstyle{definition}
\newtheorem{definition}[theorem]{Definition}

\newtheorem{conjecture}[theorem]{Conjecture}

\newtheorem{example}[theorem]{Example}

\theoremstyle{remark}
\newtheorem{remark}[theorem]{Remark}

\newtheorem*{remark*}{Remark}
\newtheorem{introremark}{Remark}[maintheorem]

\numberwithin{equation}{section}

\DeclareMathAlphabet{\mathcal}{OMS}{cmsy}{m}{n}

\newcommand{\one}{\mathbbm{1}}

\newcommand{\C}{\mathbb{C}}

\newcommand{\F}{\mathbb{F}}

\renewcommand{\P}{\mathbb{P}}
\newcommand{\Q}{\mathbb{Q}}
\newcommand{\R}{\mathbb{R}}

\newcommand{\Z}{\mathbb{Z}}

\newcommand{\scl}{\mathscr{L}}

\newcommand{\cala}{\mathcal{A}}
\newcommand{\calb}{\mathcal{B}}
\newcommand{\calc}{\mathcal{C}}
\newcommand{\cald}{\mathcal{D}}
\newcommand{\cale}{\mathcal{E}}

\newcommand{\calo}{\mathcal{O}}

\newcommand{\calr}{\mathcal{R}}

\newcommand{\calu}{\mathcal{U}}

\newcommand{\mfa}{\mathfrak{a}}
\newcommand{\mfb}{\mathfrak{b}}
\newcommand{\mfc}{\mathfrak{c}}

\newcommand{\mff}{\mathfrak{f}}

\newcommand{\mfp}{\mathfrak{p}}
\newcommand{\mfP}{\mathfrak{P}}
\newcommand{\mfq}{\mathfrak{q}}

\newcommand{\rmn}{\mathrm{N}}

\newcommand{\ab}{{\mathrm{ab}}}

\newcommand{\Char}{\mathrm{char}\,}

\renewcommand{\Im}{{\mathrm{Im}}}
\newcommand{\re}{\mathrm{Re}}
\newcommand{\Ker}{{\mathrm{Ker}}}
\newcommand{\coker}{{\mathrm{coker}}}
\newcommand{\gal}{\mathrm{Gal}}
\newcommand{\Gal}{{\mathrm{Gal}}}

\newcommand{\Hom}{{\mathrm{Hom}}}

\newcommand{\ind}{\mathrm{Ind}}

\newcommand{\rank}{{\mathrm{rank}}}

\newcommand{\SL}{{\mathrm{SL}}}

\newcommand{\disc}{{\mathrm{disc}}}

\newcommand{\Frob}{{\mathrm{Frob}}}

\newcommand{\ord}{{\mathrm{ord}}}

\newcommand{\Tr}{{\mathrm{Tr}}}

\newcommand{\Det}{\mathrm{Det}}

\newcommand{\Cl}{\mathrm{Cl}}

\newcommand{\rs}{\mathrm{RS}}

\newcommand{\ur}{\mathrm{ur}}

\newcommand{\length}{\mathrm{length}}

\newcommand{\paren}[1]{\mathopen{}\left(#1\right)\mathclose{}}
\newcommand{\set}[1]{\mathopen{}\left\{#1\right\}\mathclose{}}
\newcommand{\sbrac}[1]{\mathopen{}\left[#1\right]\mathclose{}}
\newcommand{\abrac}[1]{%
\mathopen{}\left\langle#1\right\rangle\mathclose{}}

\newcommand{\verts}[1]{%
\mathopen{}\left\lvert#1\right\rvert\mathclose{}}

\newcommand\restr[2]{{
  \left.\kern-\nulldelimiterspace
  #1 
  \right|_{#2} 
  }}

\newcommand{\widesim}[2][2]{
  \mathrel{\overset{#2}{\scalebox{#1}[1]{$\sim$}}}
}

\renewcommand{\setminus}{-}

\title[Stark units over complex cubic fields]{Stark units and cyclotomic Iwasawa theory over complex cubic fields}

\author{Robin Zhang}
\address{Institut des Hautes \'{E}tudes Scientifiques}
\email{rzhang@ihes.fr}

\date{November 24, 2026}

\begin{document}


\begin{abstract}
  For cyclotomic extensions of complex cubic fields, we prove the integral rank-$1$ and rational rank-$2$ abelian Stark conjectures. More generally, we prove a relative theorem on the higher-rank Rubin--Stark conjectures for cubic extensions of number fields. We then apply our Stark unit constructions to the Bergeron--Charollois--Garc\'{i}a conjectures on elliptic gamma values. Finally, we prove an Iwasawa main conjecture relating norm-compatible systems of Stark units to ideal class groups in cyclotomic towers over complex cubic fields.
\end{abstract}

\maketitle

\setcounter{tocdepth}{1}
\tableofcontents


\section{Introduction}
\label{sec:introduction}

In a series of four papers, Stark
\cite{stark-1971, stark-1975, stark-1976, stark-1980}
formulated rational conjectures relating
the leading terms of Artin $L$-functions $L(s, \rho)$
at $s = 0$ to regulators of algebraic units.
These conjectures apply for every
rank (i.e. order of vanishing)
$r := \ord_{s = 0} L(s, \rho)$ and,
when $\rho = \chi$ is a character of an
abelian extension $K/k$,
predict that $L^{(r)}(0, \chi)/r! = a_\chi R(\chi)$
where $R(\chi)$ is a nonzero $r \times r$ determinant
formed from logarithms of units
and where the coefficients $a_\chi \in \Q(\chi)^\times$
vary equivariantly under Galois conjugation.
In rank $0$, Stark \cite{stark-1975} deduced
the rational abelian conjecture
from Siegel's \cite{siegel-1970} rationality theorem.

In rank $1$, Stark \cite{stark-1980} formulated an
integral refinement that specifies a unit with
a fixed normalization
and requires a root of the unit to
generate an abelian extension.
Let $K/k$ be a finite abelian extension
with ring of integers $\calo_K$
and $W_K := \lvert \mu(K) \rvert$
roots of unity.
Let $S$ be a finite set of at least three places of $k$,
including every archimedean place
and every place ramified in $K/k$.
Write $L_S$ for the product of the Euler factors
at finite places outside $S$.
Let an archimedean place $v \in S$ split completely in $K$
and fix a place $w$ of $K$ above $v$.
We write $x^\sigma := \sigma(x)$
for $x \in K$ and $\sigma \in \Gal(K/k)$
and use the ordinary absolute value at real places
and its square at complex places.
Stark \cite{stark-1980} conjectures
that there is a unit $\epsilon \in \calo_K^\times$ such that
\begin{enumerate}
  \item for every character $\chi$ of $\Gal(K/k)$,
    \[
      L_S'(0, \chi) = -\frac{1}{W_K}
      \sum_{\sigma \in \Gal(K/k)}\chi(\sigma)
      \log\verts{\epsilon^\sigma}_w;
    \]
  \item $\lvert \epsilon\rvert_{w'} = 1$
    at every place $w'$ not above $v$; and
  \item $K(\epsilon^{1/W_K})/k$ is abelian;
\end{enumerate}
this is the integral rank-$1$
abelian Stark conjecture
(in the form given by
\cite[Chapitre IV, Conjecture 2.2]{tate-stark-1984}).
As Stark \cite[p.~63]{stark-1975} and
Tate \cite[p.~95]{tate-stark-1984} observed,
it gives a conjectural approach to
Hilbert's 12$^{\mathrm{th}}$ problem and
Kronecker's \textit{Jugendtraum}
on explicit constructions of abelian field extensions.
Rubin \cite[Conjecture~B$'$]{rubin-1996} extended this
integral refinement to higher rank using exterior powers
of groups of $S$-units;
this is the rank-$r$ Rubin--Stark conjecture
recalled in Conjecture~\ref{conj:integral-rank-two}.

Stark \cite{stark-1975, stark-1980} proved the integral
rank-$1$ abelian Stark conjecture over $k = \Q$ using
cyclotomic units and over imaginary quadratic fields $k$
using elliptic units.
Tate \cite[Proposition~4.3]{tate-stark-1982},
\cite[Corollaire~IV.3.3]{tate-stark-1984} showed that the
integral rank-$1$ abelian Stark conjecture at an
archimedean place is automatically true when the base
field has at least $2$ complex places.
Hence beyond the classical cases,
the integral rank-$1$ abelian Stark conjecture
has two types of open cases:
\begin{itemize}
  \item \textit{Type I}: when $k$ is totally real with
    $[k : \Q] \geq 2$;
  \item \textit{Type II}: when $k$ is almost totally real
    with $[k : \Q] \geq 3$;
\end{itemize}
where a number field is
called \textit{almost totally real} if it has
exactly one complex place.
The smallest fields in type I are
real quadratic fields.\footnote{In
the final months of writing this paper,
the author learned about
the remarkable results on the
algebraicity of the conjectural Stark units
over real quadratic fields $k$
by his colleagues Radchenko and Wheeler
in \cite{radchenko-wheeler-I}.
The author then learned that an independent researcher
communicated an AI-generated proof of further results
using \cite{radchenko-wheeler-I},
which prompted Radchenko and Wheeler to produce
their own AI-generated proof of assertions~(1) and~(2)
of the integral rank-$1$ abelian Stark conjecture
over real quadratic fields
in \cite{radchenko-wheeler-II}.
In light of these sudden developments,
we note that the next smallest fields in type I
are totally real cubic fields.}

We focus mainly on the smallest fields in type II,
which are complex cubic fields.\footnote{The author feels fortunate that
nobody has yet used AI
to reconstruct this paper's proofs
as outlined in the publicly-available lecture
\cite{zhang-swissmap}
more than one year ago.
The author is in awe of the technical power
of new AI models, but
expresses great concern
about how they will affect communication
in the mathematics community going forward.}
After proving our main results on
the integral abelian Stark conjectures
for cyclotomic extensions $K$ of complex cubic fields $k$
and, more generally,
over cubic extensions $k/k_0$ of number fields,
we give applications
to explicit class field theory,
the Bergeron--Charollois--Garc\'{i}a conjecture
on special values of elliptic gamma functions,
and a cyclotomic Iwasawa main
conjecture over complex cubic fields.

\subsection{The Stark conjectures over complex cubic fields}
\label{sec:intro-stark}
Let $k$ be a complex cubic field.
We call a finite extension $K/k$ \emph{cyclotomic}
if $K \subset k\Q^{\ab}$,
where $\Q^{\ab}$ is the maximal abelian extension of $\Q$.
For cyclotomic extensions $K/k$, Artin formalism reduces the
rational Stark conjectures to the classical cases over
$\Q$ and the quadratic resolvent of $k$.
However, the integral Stark conjectures require more
control over unit groups.
To obtain the normalization in assertion~(1)
of the integral rank-$1$ abelian Stark conjecture over $k$,
for instance,
we use the abelian extension in assertion~(3)
for Stark units over the quadratic resolvent of $k$
to prove divisibility properties in groups of units.

We first prove the integral rank-$1$ abelian Stark conjecture
for cyclotomic extensions of complex cubic fields.

\begin{maintheorem}
  \label{thm:intro-rank-one}
  \begin{samepage}
    Let
    \begin{itemize}
      \item $k$ be a complex cubic field;
      \item $K/k$ be a cyclotomic extension;
      \item $S$ be a finite set of places of $k$ containing
        the archimedean places and every finite place above
        a rational prime ramified in $K/\Q$; and
      \item $w$ be a place of $K$ above the complex place of
        $k$.
    \end{itemize}
  \end{samepage}
  Then the integral rank-$1$ abelian Stark conjecture
  for $(K/k, S)$ holds at $w$.
\end{maintheorem}

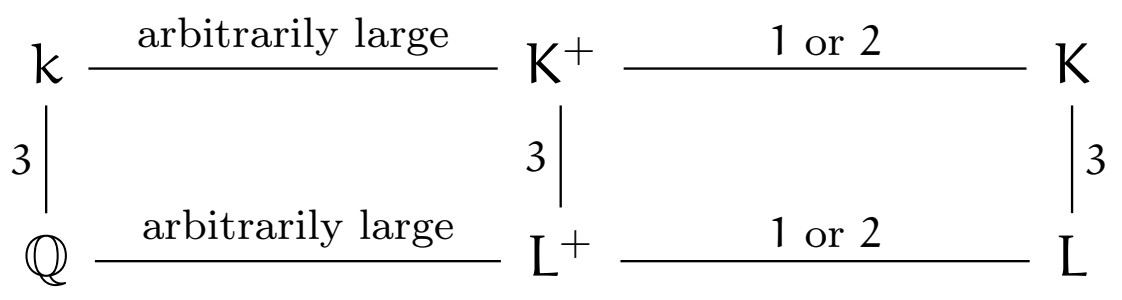
\begin{figure}[ht]
  \[
    \begin{tikzcd}[column sep=7em]
      k \arrow[r, dash, "\text{arbitrarily large}"] \arrow[d,
      dash, "3"']
      & K^+ \arrow[r, dash, "1 \text{ or } 2"] \arrow[d, dash,
      "3"']
      & K \arrow[d, dash, "3"] \\
      \Q \arrow[r, dash, "\text{arbitrarily large}"]
      & L^+ \arrow[r, dash, "1 \text{ or } 2"]
      & L
    \end{tikzcd}
  \]
  \caption{The field extensions in Theorem~\ref{thm:intro-rank-two}.}
  \label{figure:field-extensions-intro}
\end{figure}

For the rank-$2$ abelian Stark conjecture,
consider the maximal abelian subextension
$L := K \cap \Q^{\ab}$
of $K$ over $\Q$.
In particular, $K = kL$
and restriction identifies $\Gal(K/k)$ with $\Gal(L/\Q)$.
Fix complex conjugation $c$ at an embedding
of $K$ that is real on $k$,
and set $L^+ := L^{\abrac{c}}$ and $K^+ := K^{\abrac{c}}$;
then $L^+$ is totally real and $K^+ = kL^+$.
Both archimedean places of $k$ split completely in $K^+/k$.
Figure~\ref{figure:field-extensions-intro} records the field
extensions and their degrees.

\FloatBarrier

We construct a cyclotomic unit in $\calo_{L^+}$
and a unit in $\calo_{K^+}$
with the following regulator.

\begin{maintheorem}
  \label{thm:intro-rank-two}
  \begin{samepage}
    Let
    \begin{itemize}
      \item $k$ be a complex cubic field;
      \item $K/k$ be a cyclotomic extension;
      \item $S$ be a finite set of places of $k$ containing
        the archimedean places, every finite place above
        a rational prime ramified in $K/\Q$, and all places
        above two rational primes; and
      \item $w_\R$ and $w_\C$ be real and complex places of
        $K^+$ above the same real place of $L^+$.
    \end{itemize}
  \end{samepage}
  Then there are units $u_\R \in \calo_{L^+}^\times$ and
  $u_\C \in \calo_{K^+}^\times$ such that
  $\rmn_{K^+/L^+}(u_\C) = 1$ and
  \[
    L_S''(0, \chi)
    = \frac{1}{W_K}
    \paren{\sum_{\alpha \in \Gal(K^+/k)}\chi(\alpha)
    \log\verts{u_\R^\alpha}_{w_\R}}
    \paren{\sum_{\beta \in \Gal(K^+/k)}\chi(\beta)
    \log\verts{u_\C^\beta}_{w_\C}}
  \]
  for every character $\chi$ of $\Gal(K^+/k)$,
  and the rank-$2$ Rubin--Stark conjecture
  for $(K^+/k, S)$ over $\Z\sbrac{\frac{1}{6}}$ holds.
  Furthermore, the rank-$2$ Rubin--Stark conjecture
  for $(K^+/k, S)$ holds over $\Z$
  when $[K^+ : k]$ is odd
  and either $3 \nmid [K^+ : k]$
  or $\disc(k) \notin -3(\Q^\times)^2$.
\end{maintheorem}

\begin{introremark}
  \label{rem:intro-rank-two}
  The higher-rank Stark conjectures are usually
  stated in terms of elements of exterior algebras.
  For the units in Theorem~\ref{thm:intro-rank-two},
  we identify $\frac{1}{6W_K} u_\R \wedge u_\C$,
  restricted to characters whose $L$-functions vanish
  to order exactly $2$, with the rational rank-$2$
  Stark element in
  Corollary~\ref{cor:rational-rank-two-element}.
  We thus prove the rational rank-$2$ abelian Stark
  conjecture of Stark \cite{stark-1975}.
\end{introremark}

\subsection{Relative Stark over cubic extensions of the base field}

In the proofs of Theorems~\ref{thm:intro-rank-one} and
\ref{thm:intro-rank-two},
we construct explicit Stark units
from cyclotomic units and elliptic units.
We next compare the integral rank-$1$ abelian Stark
conjecture and the Rubin--Stark conjecture over cubic
extensions $k/k_0$ of number fields with
the corresponding conjectures over $k_0$ and
the quadratic resolvent of $k/k_0$.
For the abelian extension $K/k$, we replace the
cyclotomic condition $K \subset k\Q^\ab$ by
$K \subset k \cdot k_0^\ab$.
We focus on the case that $k/k_0$ is non-Galois cubic;
the corresponding assertions for Galois cubic extensions
are proved in Proposition~\ref{prop:galois-rank-one}
and Theorem~\ref{thm:Galois-Rubin-Stark}.

\begin{maintheorem}
  \label{thm:intro-relative}
  Let
  \begin{itemize}
    \item $k/k_0$ be a non-Galois cubic extension of number fields;
    \item $K/k$ be a finite abelian extension
      contained in $k \cdot k_0^\ab$;
    \item $S_{k_0}$ be a finite set of places of $k_0$
      containing the archimedean places, every prime
      ramified in $K/k_0$, and at least two finite
      places;
    \item $F$ be the quadratic resolvent of $k/k_0$;
    \item $S$ and $S_F$ be the sets of places of $k$
      and $F$ above $S_{k_0}$, respectively.
  \end{itemize}
  Then:
  \begin{enumerate}[(i)]
    \item The integral rank-$1$ abelian Stark
      conjecture for $(K/k, S)$ holds at an
      archimedean place $v$ of
      $k$ that splits completely in $K$
      if the integral rank-$1$ abelian Stark
      conjecture for $(KF/F, S_F)$ holds at
      any place of $F$ above $\restr{v}{k_0}$.
    \item For integers $r_{k_0}, r_F \geq 0$
      and $L := K \cap k_0^\ab$,
      the integral rank-$(r_{k_0}~+~r_F)$
      Rubin--Stark conjecture for $(K/k, S)$
      holds over $\Z\sbrac{\frac{1}{3W_K}}$
      if the rank-$r_{k_0}$ Rubin--Stark
      conjecture for $(L/k_0, S_{k_0})$ and the
      rank-$r_F$ Rubin--Stark conjecture
      for $(KF/F, S_F)$ hold
      over $\Z\sbrac{\frac{1}{3W_K}}$.
  \end{enumerate}
\end{maintheorem}

\begin{figure}[H]
  \centering
  \usetikzlibrary{arrows}
  \begin{tikzpicture}[
    x = 1cm, y = 0.8cm,
    every node/.style = {font = \normalsize, inner sep = 3pt},
    implication/.style = {double equal sign distance, -implies},
    is-conjecture label/.style = {font = \scriptsize},
    rs-conjecture label/.style = {font = \scriptsize},
    theorem label/.style = {font = \tiny, inner sep = 1pt}
  ]
    \node (Q) at (-3.8, -2) {$\Q$};
    \node (kzero) at (-3.8, 0) {$k_0$};
    \node (realdots) at (-3.8, 3.5) {$k_0^\ab$};
    \node (F) at (0, 1) {$F$};
    \node (Fab) at (0, 4.5) {$F^\ab$};
    \node (k) at (4.5, 1.6) {$k$};
    \node (kcyc) at (4.5, 5.1) {$k \cdot k_0^\ab$};

    \draw (Q) -- (kzero);
    \draw (kzero) -- node[above left, font = \small] {$2$} (F);
    \draw (kzero.east)
      to[out = 8, in = 208]
      node[pos = 0.68, below, font = \small] {$3$} (k.south west);

    \draw[dotted, line cap = round]
      (kzero.north) --
      node[right, blue, rs-conjecture label, yshift=0.5em] (rubinzero)
        {\shortstack{rank-$r_{k_0}$\\Rubin--Stark}}
      (realdots.south);

    \draw[dotted, line cap = round]
      (F.north) --
      node[left, blue, rs-conjecture label, yshift=0.3em] (rubinF)
        {\shortstack{rank-$r_F$\\Rubin--Stark}}
      node[right, red, is-conjecture label, yshift=0.6em] (starkF)
        {\shortstack{rank-$1$\\integral Stark}}
      (Fab.south);

    \draw[dotted, line cap = round]
      (k.north) --
      node[left, red, is-conjecture label, yshift=0.5em] (starkk)
        {\shortstack{rank-$1$\\integral Stark}}
      node[right, blue, rs-conjecture label, yshift=0.3em] (rubink)
        {\shortstack{rank-$(r_{k_0}+r_F)$\\Rubin--Stark}}
      (kcyc.south);

    \draw[red, implication]
      ([yshift = 0mm]starkF.south)
      to[out = -50, in = 180] (1.65, 2.3)
      -- node[below, theorem label, yshift = -2pt]
        {Theorem~\ref*{thm:intro-relative}(i)}
      (2.75, 2.3)
      to[out = 0, in = -120]
      ([yshift = 0mm]starkk.south);

    \coordinate (join) at (-2.1, 4.45);

    \draw[blue, double equal sign distance]
      ([yshift = 1mm]rubinzero.north)
      to[out = 90, in = -90, looseness = 1.1] (join);

    \draw[blue, double equal sign distance]
      ([yshift = 1mm]rubinF.north)
      to[out = 90, in = -90, looseness = 1.1] (join);

    \draw[blue, implication]
      (join)
      to[out = 90, in = 180] (-0.7, 5.8)
      -- node[above, theorem label, yshift = 2pt]
        {Theorem~\ref*{thm:intro-relative}(ii)}
      (5.2, 5.8)
      to[out = 0, in = 90]
      ([yshift = 1mm]rubink.north);

  \end{tikzpicture}
  \caption{The implications in
    Theorem~\ref{thm:intro-relative},
    with \\
    part~(i) in \textcolor{red}{red}
    and part~(ii) over $\Z[1/(3W_K)]$
    in \textcolor{blue}{blue}.}
  \label{figure:relative-Stark-implication}
\end{figure}
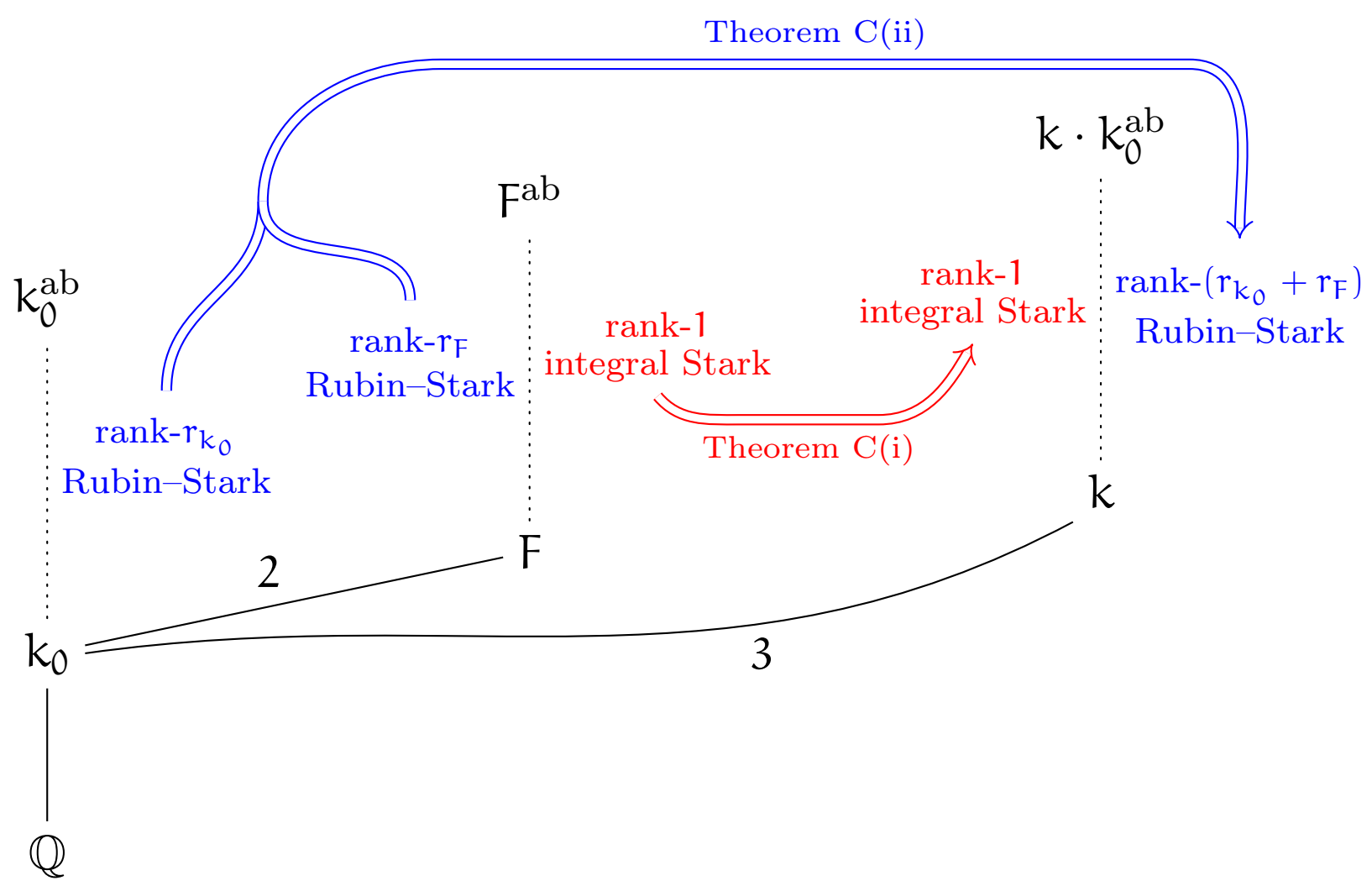

\begin{introremark}
  \label{rem:intro-relative}
  Theorem~\ref{thm:intro-relative}(i)
  generalizes the argument of
  Theorem~\ref{thm:intro-rank-one}
  but with the stronger assumption
  that $S_{k_0}$ contains at least two finite places.
  Theorem~\ref{thm:intro-relative}(ii)
  generalizes Theorem~\ref{thm:intro-rank-two},
  but without the explicit construction of Stark units
  and with weaker integrality statements
  at $2$ and $3$.
\end{introremark}

Theorem~\ref{thm:intro-relative}(ii)
leads to new unconditional cases of the
higher-rank Rubin--Stark conjectures,
such as the rank-$3$ Rubin--Stark conjecture
over $\Z[1/6]$
for a degree-$10$ extension of a
totally imaginary sextic field
in Example~\ref{ex:relative-rank-three}.

\subsection{Explicit class field theory and elliptic gamma functions}
\label{subsec:intro-elliptic-gamma}

Following the vision of
Stark \cite[p.~63]{stark-1975} and
Tate \cite[p.~95]{tate-stark-1984},
we apply
Theorem~\ref{thm:intro-rank-one}
and Theorem~\ref{thm:intro-rank-two}
to problems in explicit class field theory.

Let $K/k$ be a totally imaginary cyclic cyclotomic
extension of a complex cubic field,
with the finite set $S$
and the Stark unit $\epsilon$ as in
Theorem~\ref{thm:intro-rank-one}.
In Corollary~\ref{cor:rank-one-field-generation},
we prove that $K = k(\epsilon)$
if no finite place in $S$
splits completely in $K/k$.
For the units in Theorem~\ref{thm:intro-rank-two},
we obtain Corollary~\ref{cor:rank-two-field-generation}:
$L^+ = \Q(u_\R)$ and $K^+ = k(u_\C)$ if $K^+/k$ is cyclic
and no finite place in $S$ splits completely in $K^+/k$.
In both corollaries, we use $L_S'(0, \chi) \neq 0$ or
$L_S''(0, \chi) \neq 0$ for a faithful character $\chi$.

In statistical mechanics, Ruijsenaars
\cite[Section~III.B]{ruijsenaars} introduced the elliptic
gamma function
\begin{equation}
  \label{eq:elliptic-gamma-function}
  \Gamma(z, \tau_1, \tau_2)
    := \prod_{j, \ell = 0}^{\infty}
    \frac{1 - e^{2\pi i((j + 1)\tau_1 + (\ell + 1)\tau_2 - z)}}
    {1 - e^{2\pi i(j\tau_1 + \ell\tau_2 + z)}},
\end{equation}
a meromorphic function for $z, \tau_1, \tau_2 \in \C$ with
$\Im(\tau_1) > 0$ and $\Im(\tau_2) > 0$.
Felder--Varchenko
\cite[Theorem~4.1 and Section~7]{felder-varchenko} proved
transformation formulas involving $\SL_3(\Z)$, analogous
to the $\SL_2(\Z)$ transformation formulas for the first
Jacobi theta function.
Recently,
Bergeron--Charollois--Garc\'{i}a
\cite{bergeron-charollois-garcia}
propose using elliptic gamma values to construct units
in narrow ray class fields of complex cubic fields,
analogously to how theta values
give the elliptic units used by Stark \cite{stark-1980}
over imaginary quadratic fields.

To state our result on these elliptic gamma values,
let $k$ be a complex cubic field and
$\mff \subset \calo_k$ be a nonzero ideal such that
the narrow ray class field $K := k(\mff)$ is totally
imaginary. For an integral ideal $\mfb \subset \calo_k$
coprime to $\mff$, define $\scl := \mff\mfb^{-1}$.
We will compare the value attached to $\scl$ with
the conjugate of a unit under the Artin symbol
$\sigma_\mfb \in \Gal(K/k)$ of $\mfb$.

Choose an auxiliary ideal $\mfa \subset \calo_k$ coprime
to $6\rmn(\mff)$ with odd prime norm
$\rmn(\mfa) = \verts{\calo_k/\mfa}$.
To specify the arguments of the elliptic gamma functions,
let $q$ be the least positive integer in $\mff \cap \Z$.
We call an element $h = \lambda/q$ \emph{admissible} if
$\lambda \in \scl$, $h - 1 \in \scl$, and $\lambda/\rmn(\mfa)$
belongs to $\mfa^{-1}\scl$ and generates
$\mfa^{-1}\scl/\scl$.

Bergeron--Charollois--Garc\'{i}a
\cite[Definition~2.3]{bergeron-charollois-garcia}
define special values of elliptic gamma functions
using complex parameters $\tau_1$ and $\tau_2$
determined by $\scl$ and a unit of $k$.
Fix complex embeddings
\begin{align*}
  x: k &\longrightarrow \C \\
  \iota : K &\longrightarrow \C
\end{align*}
such that the restriction of $\iota$ to $k$ is $x$,
and let $v$ be the place of $K$ corresponding to $\iota$.
Let $x_\R$ be the real embedding of $k$, and orient
$\R\otimes_\Q k$ by $(x_\R, \re(x), \Im(x))$.
Let $\eta$ generate the group of units $u \in \calo_k^\times$
satisfying $u - 1 \in \mff$ and $x_\R(u) > 0$, choosing
the generator with $0 < x_\R(\eta) < 1$.
We now use $\eta$ and an admissible $h = \lambda/q$ to
determine the parameters $\tau_1$ and $\tau_2$.
Choose a surjective $\Z$-linear form
$a : \scl \twoheadrightarrow \Z$
such that $a(\lambda) = a(\eta^{-1}\lambda) = 0$,
and define $b(t) := a(\eta^{-1}t)$ for $t \in \scl$.
Choose $\alpha, \beta \in \scl$ with
\[
  b(\alpha) = a(\beta) = 0, \;
  a(\alpha) > 0, \,
  \text{ and } \, b(\beta) > 0.
\]
Let $\gamma$ generate $\ker a \cap \ker b$,
with $(\alpha, \beta, \gamma)$ positively oriented.
Choose the sign of $a$ and these elements so that
\begin{align*}
  \tau_1 := -x(\alpha)/x(\gamma) \quad \text{ and } \quad
  \tau_2 := x(\beta)/x(\gamma)
\end{align*}
have positive imaginary parts;
these choices are possible by
\cite[Lemma~2.1 and Equation~2.6]{bergeron-charollois-garcia}.
Let $\cald$ be a set of representatives modulo $\Z\gamma$
for the $\delta \in \scl$ satisfying
$0 \leq a(\delta) < a(\alpha)$ and
$0 \leq b(\delta) < b(\beta)$;
such a set is necessarily finite.
With $z_\delta := x(h + \delta - \alpha)/x(\gamma)$,
we follow Bergeron--Charollois--Garc\'{i}a
\cite[Sections~2.1--2.2]{bergeron-charollois-garcia}
and define
\[
  \Gamma_{\mfa, h}(\scl)
  := \prod_{\delta \in \cald}
  \frac{\Gamma(z_\delta, \tau_1, \tau_2)^{\rmn(\mfa)}}
  {\Gamma(\rmn(\mfa)z_\delta, \rmn(\mfa)\tau_1, \rmn(\mfa)\tau_2)}.
\]
For fixed $\scl$, $\mfa$, $h$, and $x$,
Bergeron--Charollois--Garc\'{i}a
\cite[Sections~2.1--2.2]{bergeron-charollois-garcia}
prove that $\Gamma_{\mfa, h}(\scl)$ is nonzero and
independent of $\alpha$, $\beta$,
and the representatives in $\cald$.

Bergeron--Charollois--Garc\'{i}a
\cite[Section~2.3]{bergeron-charollois-garcia}
conjecture that for each
ideal $\mfa \subset \calo_k$ coprime
to $6\rmn(\mff)$ with odd prime norm
$\rmn(\mfa)$:
\begin{enumerate}[(1)]
  \item $\Gamma_{\mfa, h}(\scl)$ is independent
    of the choice of admissible $h$;
  \item there are units $u_{\scl, \mfa}' \in \calo_K^\times$
    indexed by $\scl$ such that
    \begin{enumerate}[(i)]
      \item $\Gamma_{\mfa, h}(\scl) = \iota(u_{\scl, \mfa}')$,
      \item $\lvert u_{\scl, \mfa}'\rvert_{w'} = 1$
      at every place $w'$ of $K$ above the real place of $k$, and
      \item $\sigma_\mfc(u_{\scl, \mfa}')
        = u_{\mfc^{-1}\scl, \mfa}'$
        for every integral ideal $\mfc$ of $k$
        coprime to $\mff$.
    \end{enumerate}
\end{enumerate}
In Proposition~\ref{prop:elliptic-gamma-independence},
we prove the independence of
$\Gamma_{\mfa, h}(\scl)$ from the choice of admissible $h$
in full generality.
For cyclotomic $K/k$, we use
Theorem~\ref{thm:intro-rank-one} to construct units
with the absolute values and reciprocity law predicted
by the Bergeron--Charollois--Garc\'{i}a conjecture.
We write
$\lvert z \rvert_\C := z\overline{z} = \lvert z \rvert^2$,
where $\lvert z \rvert$ is the usual absolute value on $\C$.

\begin{maincorollary}
  \label{cor:intro-elliptic-gamma}
  Let
  \begin{itemize}
    \item $k$ be a complex cubic field;
    \item $K := k(\mff)$ be a totally imaginary
      narrow ray class field with
      $K \subset k\Q^{\ab}$;
    \item $S_\mff$ be the set of archimedean places
      of $k$ and prime divisors of $\mff$;
  \end{itemize}
    If $S_\mff$ contains every place of $k$ above
    a finite set of at least two rational primes containing
    all primes ramified in $K/\Q$,
    then there is a unit $u_{\mff, \mfa} \in \calo_K^\times$
    whose conjugates
    $u_{\scl, \mfa} := u_{\mff, \mfa}^{\sigma_\mfb}$ satisfy
    \[
      \verts{\Gamma_{\mfa, h}(\scl)}_\C
        = \verts{u_{\scl, \mfa}}_v
    \]
    for every ideal $\mfb$ coprime to $\mff$ and every
    admissible $h$. Furthermore, each $u_{\scl, \mfa}$ has
    absolute value $1$ at every place of $K$
    above the real place of $k$.
\end{maincorollary}

In the proof of Corollary~\ref{cor:intro-elliptic-gamma},
we construct the units $u_{\scl, \mfa}$ explicitly.
To prove the full conjecture under the hypotheses of
Corollary~\ref{cor:intro-elliptic-gamma},
it remains to choose a single $\xi \in \mu(K)$ such that
$\Gamma_{\mfa, h}(\mff\mfb^{-1})$ equals the complex image of
$(\xi u_{\mff, \mfa})^{\sigma_\mfb}$ for every ideal $\mfb$
coprime to $\mff$ and every admissible $h$.

As a demonstration,
we apply our study of Stark units and elliptic gamma values
to the specific cyclic sextic extension
of a complex cubic field
studied by Dasgupta
\cite{dasgupta1999},
in particular verifying the Stark unit
predicted by \cite[Proposition~5.5.3]{dasgupta1999}
and giving an explicit example
for the algebraicity of $\verts{\Gamma_{\mfa, h}(\scl)}_\C$.

\subsection{Cyclotomic Iwasawa theory}

Finally, we choose the units of
Theorems~\ref{thm:intro-rank-one} and \ref{thm:intro-rank-two}
compatibly under norms in a cyclotomic tower.
We compare the submodules they generate
with inverse limits of ideal class groups,
following Iwasawa~\cite{iwasawa-1973} and Rubin~\cite{rubin-1991}.

Let $k$ be a complex cubic field with quadratic resolvent
$F$, and let $L_0/\Q$ be finite abelian with $F \subset L_0$.
Let $K_0 := kL_0$, and let $p > 3$ be unramified in $K_0$,
with $p \nmid [L_0 : \Q]$.
For the degree-$p^n$ layer $\Q_n$ of the cyclotomic
$\Z_p$-extension $\Q_\infty/\Q$,
define $K_n := K_0\Q_n$ and $K_\infty := K_0\Q_\infty$.
Fix a complex embedding of $K_\infty$ that is real on $k$,
and define the fixed field $K_n^+ := K_n^{\abrac{c}}$
of complex conjugation $c$.
Let $\Sigma_\Q$ be a finite set of places of $\Q$,
including $\infty$, $p$, and every prime ramified in $K_0/\Q$;
let $\Sigma_k$ consist of the places of $k$ above $\Sigma_\Q$.

In Theorem~\ref{thm:finite-norm-relations} and
Corollary~\ref{cor:Stark-unit-norm}, we show that the units
$\epsilon_n \in \calo_{K_n}^\times$ of
Theorem~\ref{thm:intro-rank-one} and
the units $u_{\R, n}$ and $u_{\C, n}$ in
$\calo_{K_n^+}^\times$ of
Theorem~\ref{thm:intro-rank-two}, with $S = \Sigma_k$,
can be chosen so that, for every $n \geq 0$,
\begin{align*}
    \rmn_{K_{n + 1}/K_n}(\epsilon_{n + 1})
      &= \epsilon_n, \\
    \rmn_{K_{n + 1}^+/K_n^+}(u_{\R, n + 1})
      &= u_{\R, n}, \text{ and } \\
    \rmn_{K_{n + 1}^+/K_n^+}(u_{\C, n + 1})
      &= u_{\C, n}.
\end{align*}

We will compare these units with ideal classes after
passing to inverse limits under norms.
To describe the Galois action, let
$\Delta_0 := \Gal(L_0/\Q)$ and
$\Gamma := \Gal(\Q_\infty/\Q)$.
We use the restriction isomorphism
$\Gal(K_\infty/k) \cong \Delta_0\times\Gamma$
to view both inverse limits as modules over the same
Iwasawa algebra, which we now define.
Let $\calo$ be the ring of integers of a finite extension
of $\Q_p$ containing the character values of $\Delta_0$
and a primitive third root of unity, and
define the Iwasawa algebra
$\Lambda := \calo[[\Gamma]]$.
Fix a character $\vartheta : \Delta_0 \longrightarrow
\calo^\times$, with idempotent
\[
  e_\vartheta := \lvert\Delta_0\rvert^{-1}
    \sum_{\delta \in \Delta_0}\vartheta(\delta)\delta^{-1}
    \in \calo[\Delta_0].
\]
Viewing complex conjugation $c$ as an element of $\Delta_0$,
define $r_\vartheta := 1$ if $\vartheta(c) = -1$ and
$r_\vartheta := 2$ if $\vartheta(c) = 1$.
Define the inverse limits of units and ideal classes
\begin{align*}
  \cale_\Sigma &:= e_\vartheta\varprojlim_n
    (\calo_{K_n, \Sigma_k}^\times\otimes_\Z\calo), \\
  A_\Sigma &:= e_\vartheta\varprojlim_n
    (\Cl(\calo_{K_n, \Sigma_k})\otimes_\Z\calo),
\end{align*}
with norms as transition maps.
Here tensoring with $\calo$ means $p$-adic completion
followed by extension of scalars from $\Z_p$ to $\calo$,
and $\Cl(\calo_{K_n, \Sigma_k})$ is the ideal class group
modulo the classes of primes above $\Sigma_k$.
To record the omitted finite Euler factors, let $Y_\Sigma$
be the free $\calo$-module on the finite primes of
$K_\infty$ above $\Sigma_\Q$, and set
$B_\Sigma := e_\vartheta\ker(Y_\Sigma \longrightarrow \calo)$,
where each prime maps to $1$.

Regarding $u_{\R, n}$ and $u_{\C, n}$
as units of $K_n$ via $K_n^+ \subset K_n$
and writing units additively after tensoring with $\calo$,
define the submodule of $\cale_\Sigma$
\[
  \calc_\vartheta^\Sigma :=
  \begin{cases}
    \Lambda e_\vartheta(\epsilon_n)_n
    & \text{if $\vartheta(c) = -1$}, \\
      \Lambda e_\vartheta(u_{\R, n})_n
      + \Lambda e_\vartheta(u_{\C, n})_n
    & \text{if $\vartheta(c) = 1$}.
  \end{cases}
\]
For a finitely generated torsion $\Lambda$-module $M$,
write $\Char_\Lambda(M)$ for its characteristic ideal.
The following result is a complex cubic analogue of
the Iwasawa main conjectures proved by
Mazur--Wiles~\cite{mazur-wiles} over $\Q$ and
Rubin~\cite{rubin-1991} over imaginary quadratic fields.

\begin{maintheorem}
  \label{thm:intro-determinant-main}
  The modules $\cale_\Sigma$ and $\calc_\vartheta^\Sigma$
  are free of rank $r_\vartheta$ over $\Lambda$, and
  the modules $A_\Sigma$, $B_\Sigma$, and
  $\cale_\Sigma/\calc_\vartheta^\Sigma$ are torsion.
  Moreover,
  \[
    \Char_\Lambda(\cale_\Sigma/\calc_\vartheta^\Sigma)
    = \Char_\Lambda(A_\Sigma)\Char_\Lambda(B_\Sigma).
  \]
\end{maintheorem}

\subsection{Proof overview and organization}

In Section~\ref{sec:relative-cubic}, we review the basic
theory of non-Galois cubic extensions $k/k_0$ and
Artin $L$-functions over $k$.
In particular, we prove a decomposition of Artin
$L$-functions over $k$ into products of $L$-functions
over $k_0$ and the quadratic resolvent $F$ of $k/k_0$.
We use this decomposition to deduce the
higher-rank rational Stark conjectures over $k$
from those over $k_0$ and $F$
in Corollary~\ref{cor:relative-rational-Stark}.
By comparing Euler factors, we obtain the corresponding
decomposition of partial $L$-functions for any finite
set $S$ of places of $k$ containing the archimedean places
in Proposition~\ref{prop:local-Artin-factorization}.

In Section~\ref{sec:rank-one}, we prove
Theorems~\ref{thm:intro-rank-one} and \ref{thm:intro-relative}(i).
First, we construct a Stark unit for $K/k$
under the assumptions that
$k$ is almost totally real
and $K$ contains the quadratic resolvent $F$
in Proposition~\ref{prop:relative-rank-one}.
We then use the descent argument of
Tate~\cite[Chapter~IV, Proposition~3.5]{tate-stark-1984}
to remove the assumption that $K$ contains $F$
and prove Theorem~\ref{thm:intro-relative}(i).
We separately treat the case when $k/k_0$ is Galois
in Proposition~\ref{prop:galois-rank-one}.
To prove Theorem~\ref{thm:intro-rank-one},
we specialize to $k_0 = \Q$ and apply the theorem of
Stark~\cite[Theorem~2]{stark-1980}
over imaginary quadratic fields.
Finally, we give an application of
Theorem~\ref{thm:intro-rank-one} to explicit class field theory
in Corollary~\ref{cor:rank-one-field-generation}.

Section~\ref{sec:rank-r} proves
Theorems~\ref{thm:intro-relative}(ii) and
\ref{thm:intro-rank-two}.
First, we compare rational Stark elements over $k_0$,
$F$, and $k$ for a non-Galois cubic extension $k/k_0$,
and use this comparison to prove the integrality needed
for Theorem~\ref{thm:intro-relative}(ii).
We separately treat the case when $k/k_0$ is Galois
in Theorem~\ref{thm:Galois-Rubin-Stark}.
We then construct explicit rank-$2$ Stark elements over
almost totally real cubic extensions
of totally real fields
and study their integrality
in Theorem~\ref{thm:rank-two-integral-refinement}.
When $k_0 = \Q$, we use the classical Stark theorems
recalled in Section~\ref{subsec:classical-Stark}
to obtain Theorem~\ref{thm:intro-rank-two},
including its assertions at $2$ and $3$.
Finally, we give an application of
Theorem~\ref{thm:intro-rank-two} to explicit class field
theory in Corollary~\ref{cor:rank-two-field-generation}.

In Section~\ref{sec:applications},
we compare the
Kronecker limit type formula of
Bergeron--Charollois--Garc\'{i}a
\cite[Theorem~3.2]{bergeron-charollois-garcia} with the
regulator of a Stark unit for $K/k$;
this allows us to construct algebraic units
with the absolute values and
reciprocity predicted by Bergeron--Charollois--Garc\'{i}a
\cite[Section~2.3]{bergeron-charollois-garcia}
for elliptic gamma values $\Gamma_{\mfa, h}(\scl)$,
proving Corollary~\ref{cor:intro-elliptic-gamma}.
For cyclotomic $K/k$, uniqueness of these algebraic units
up to roots of unity reduces the equality of complex values
to one common choice of a root of unity in
Proposition~\ref{prop:elliptic-gamma-phase}.
Then for general $K/k$,
we prove the assertion of
the Bergeron--Charollois--Garc\'{i}a conjecture
\cite[Section~2.3]{bergeron-charollois-garcia}
that the elliptic gamma value $\Gamma_{\mfa, h}(\scl)$
is independent from the choice of admissible $h$
in Proposition~\ref{prop:elliptic-gamma-independence}.
We deduce a product identity for elliptic gamma functions from
the work of Morain \cite{morain-gamma-II}
and use some linear algebra over a finite field
to compare elliptic gamma values.

In Section~\ref{sec:Iwasawa}, we prove
Theorem~\ref{thm:intro-determinant-main}.
First, we construct norm-compatible systems
of units and rational rank-$2$ Stark elements
from the Stark units
of Theorems~\ref{thm:intro-rank-one} and
\ref{thm:intro-rank-two} in
Theorem~\ref{thm:finite-norm-relations} and
Corollaries~\ref{cor:rank-two-element-norm} and
\ref{cor:Stark-unit-norm}.
We then decompose $\cale_\Sigma$, $A_\Sigma$, and $B_\Sigma$
as $\Lambda$-modules into components over $\Q$ and $F$
in Proposition~\ref{prop:derived-Artin-decomposition}.
In Lemma~\ref{lem:elliptic-unit-conjugation}, we identify
the submodules generated by our explicit Stark units
in the components over $\Q$ and $F$.
For the component over $\Q$, we use the Iwasawa main
conjecture of Mazur--Wiles \cite{mazur-wiles} in the
formulation of Huber--Kings
\cite[Theorem~4.4.1 and Corollary~5.1.3]{huber-kings}.
For the component over $F$, we use Bullach--Hofer
\cite[Theorem~5.8]{bullach-hofer}, who prove the $p$-part
of the equivariant Tamagawa number conjecture at $s = 0$
for abelian extensions of imaginary quadratic fields
under hypotheses that we verify for every $K_n/F$.
We use these results to prove the equalities of
characteristic ideals in
Propositions~\ref{prop:rational-source-determinant} and
\ref{prop:elliptic-source-identification}.
We then combine these equalities with the decompositions
of the modules and the identities for our Stark units
to prove Theorem~\ref{thm:intro-determinant-main}.

Finally in \hyperref[app:dasgupta]{Appendix 1},
we verify the Stark unit predicted by Dasgupta
\cite[Section~5.5, Proposition~5.5.3]{dasgupta1999}
for an explicit sextic extension of a complex cubic field
by an application of
our decomposition of partial $L$-functions
with the lower bound for Mahler measures of
Mossinghoff--Rhin--Wu~\cite[Theorem~1.1]{mossinghoff-rhin-wu}.

\subsection*{AI usage}
GPT 5.6 Sol and GPT 6 Sol were used
to assist with copyediting,
the \LaTeX \, formatting of Figure
\ref{figure:relative-Stark-implication},
finding the reference
\cite[Tags~00LX and 090V]{stacks-project},
merging what was originally
two separate papers,
and dealing with the notation conflicts from that merge.
GPT 6 Sol was used
to find the field extensions
and finite sets of places in
Example~\ref{ex:relative-rank-three}.

In an earlier version of this paper,
the proof of Theorem~\ref{thm:Dasgupta-exact}
of \hyperref[app:dasgupta]{Appendix 1}
used a height bound from \cite{voutier},
which the author found
using GPT 5.6 Sol.
Following a suggestion of
Wadim Zudilin,
the proof of
Theorem~\ref{thm:Dasgupta-exact}
now states a simpler argument
using a bound from \cite{mossinghoff-rhin-wu}
instead of \cite{voutier}.

\subsection*{Acknowledgements}

The author thanks
Eran Assaf,
Dick Gross,
Cristian Popescu,
and Wadim Zudilin
for many inspiring discussions
and helpful comments on this paper.
This material is partially based on work
that was supported
by the National Science Foundation under Grant
No. DMS-2303280.


\section{Artin \texorpdfstring{$L$}{L}-functions for
non-Galois cubic extensions}
\label{sec:relative-cubic}

Our results on the Stark conjectures
over non-Galois cubic extensions $k/k_0$
begin with the same decomposition of Artin
$L$-functions:
for $F$ the quadratic resolvent of $k/k_0$,
we express an Artin $L$-function over $k$
as a product of one over $k_0$ and
one over $F$. We obtain this decomposition from the
permutation representation of $S_3$.
We first describe the fields and their archimedean places,
then determine how to choose the omitted Euler factors
to get an analogous decomposition for partial $L$-functions.

\subsection{The field extensions and archimedean places}
\label{subsec:Stark-hypotheses}

\subsubsection*{The field extensions}

We follow the classical approach to the arithmetic of
non-Galois cubic fields via their quadratic resolvents,
going back to Hasse \cite{hasse-1930}.
Let $k/k_0$ be a non-Galois cubic extension of number
fields, and let $P \in k_0[x]$ be an irreducible cubic
polynomial defining $k$.
Its quadratic resolvent and Galois closure are
\begin{align*}
  F &:= k_0(\sqrt{\disc(P)}), \\
  E &:= kF.
\end{align*}
The Galois group $\Gal(E/k_0)$ is isomorphic to $S_3$,
and $F = E^{A_3}$.
Thus $E/k$ has degree $2$ and $E/F$ has degree $3$,
while $[E : k_0] = 6$.

We now describe the abelian extensions of $k_0$ and $F$
that will enter the decomposition of an Artin
$L$-function over $k$.
Let $K/k$ be a finite abelian extension contained in
$k \cdot k_0^\ab$, and set $L := K \cap k_0^\ab$.
Since $k/k_0$ is non-Galois cubic,
we have $k \cap k_0^\ab = k_0$.
Restriction therefore identifies
$\Gal(k \cdot k_0^\ab/k)$ with $\Gal(k_0^\ab/k_0)$.
By Galois correspondence, we obtain
\begin{align*}
  K &= kL, \\
  \Gal(K/k) &\cong \Gal(L/k_0).
\end{align*}
Every root of unity in $K$ belongs to $k_0^\ab$,
so $\mu(K) = \mu(L)$.

Later,
we will inflate characters
from $\Gal(K/k)$ to $\Gal(KF/k)$
to reduce the Stark conjectures
for general abelian $K/k$
to the case when $K$ contains $F$;
hence for the rest of this section,
we assume that $F \subset K$.
With this assumption,
then $F \subset L$, $K = EL$,
and the maximal abelian subextension of $E/k_0$ is $F$
so $E \cap L = F$.
Restriction gives
\[
  \Gal(K/F)
  \cong \Gal(E/F) \times \Gal(L/F),
\]
and hence $K/F$ is abelian.
We will fix a generator $\tau$
of the order-$3$ cyclic group
$\Gal(K/L) \cong \Gal(E/F)$.

\subsubsection*{The archimedean places}

We next compare the archimedean places of $k$ and $F$.
For an archimedean place $v$ of $k_0$, write $k_{0,v}$
for the completion of $k_0$ at $v$.
The following table lists the possible algebras
up to isomorphism.
\begin{center}
  \renewcommand{\arraystretch}{1.15}
  \begin{tabular}{ccc}
    \hline
    $k_{0,v}$
      & $k \otimes_{k_0} k_{0,v}$
      & $F \otimes_{k_0} k_{0,v}$ \\
    \hline
    $\R$ & $\R^3$ & $\R^2$ \\
    $\R$ & $\R \times \C$ & $\C$ \\
    $\C$ & $\C^3$ & $\C^2$ \\
    \hline
  \end{tabular}
\end{center}
Each real place of $k_0$ contributes the same number
of complex places to $k$ and $F$.
A complex place of $k_0$ contributes three complex places
to $k$ and two to $F$.
It follows that $k$ has exactly one complex place
if and only if $F$ does.
Either condition forces $k_0$ to be totally real.

For the rest of this subsection, assume that $k$ is
almost totally real, meaning that it has exactly one
complex place.
There is then a unique real place $v_0$ of $k_0$ with
\[
  k \otimes_{k_0, v_0} \R \cong \R \times \C.
\]
We use the real place of $k$ above $v_0$ to choose
complex conjugation fixing $k$.
Fix an embedding of an algebraic closure of $k$ into
$\C$ extending this real place, and let $c$ be the
resulting complex conjugation.
Since $K = EL$ is Galois over $k_0$, the automorphism
$c$ restricts to an automorphism of $K$.
It fixes $k$ and acts nontrivially on $F$.

We use $c$ to define the fixed fields
\begin{align*}
  L^+ &:= L^{\abrac{c}}, \\
  K^+ &:= K^{\abrac{c}} = kL^+.
\end{align*}
Since $c$ acts nontrivially on $F$, we have
$F \cap L^+ = k_0$.
It follows that $L = L^+F$ and $K = K^+F$.
The decomposition group of a place of $L$ above $v_0$
is generated by $c$, so $v_0$ splits completely in
$L^+/k_0$.
The chosen embedding is real on both $L^+$ and $K^+$,
and therefore
\[
  \mu(L^+) = \mu(K^+) = \set{1, -1}.
\]
When $k_0 = \Q$, the field $L^+$ is totally real,
since it is Galois over $\Q$ and has a real embedding.

For the comparison of characters below, we record
the resulting group identifications.
Restriction gives
\begin{align*}
  \Gal(K/kF)
    &\cong \Gal(L/F)
    \cong \Gal(L^+/k_0), \\
  \Gal(K^+/k)
    &\cong \Gal(L^+/k_0).
\end{align*}
Since $c$ has order $2$ and acts nontrivially on $F$,
and $\tau$ generates $\Gal(K/L)$, we have
\begin{align*}
  \Gal(K/k)
    &= \Gal(K/kF) \times \abrac{c}, \\
  \Gal(K/F)
    &= \Gal(K/kF) \times \abrac{\tau}.
\end{align*}
The restrictions of $c$ and $\tau$ to $E$ are a
transposition and a $3$-cycle in $\Gal(E/k_0) \cong S_3$.
Since $\tau$ fixes $L$, the relations in $S_3$ give
$\tau^3 = 1$ and
$c\tau c^{-1} = \tau^{-1}$.

\begin{figure}[!ht]
  \centering
  \begin{tikzcd}[row sep=1.3em, column sep=1.9em]
    & E
    \arrow[rr, dash, "\Gal(K/E)"]
    \arrow[dd, dashed, no head]
    && K
    \arrow[dd, dash,
    "\abrac{\tau}"] \\
    k
    \arrow[ur, dash]
    \arrow[rr, dash]
    \arrow[dd, dash]
    && K^+\arrow[ur, dash, near start,
    "\abrac{c}"]
    \arrow[dd, dash] \\
    & F
    \arrow[rr, dashed, no head]
    && L \\
    k_0
    \arrow[rr, dash]
    \arrow[ur, dash]
    && L^+\arrow[ur, dash]
  \end{tikzcd}
  \caption{The field extensions when $k$ is almost totally
    real and $K$ contains $F$.}
  \label{figure:field-extensions-cube}
\end{figure}
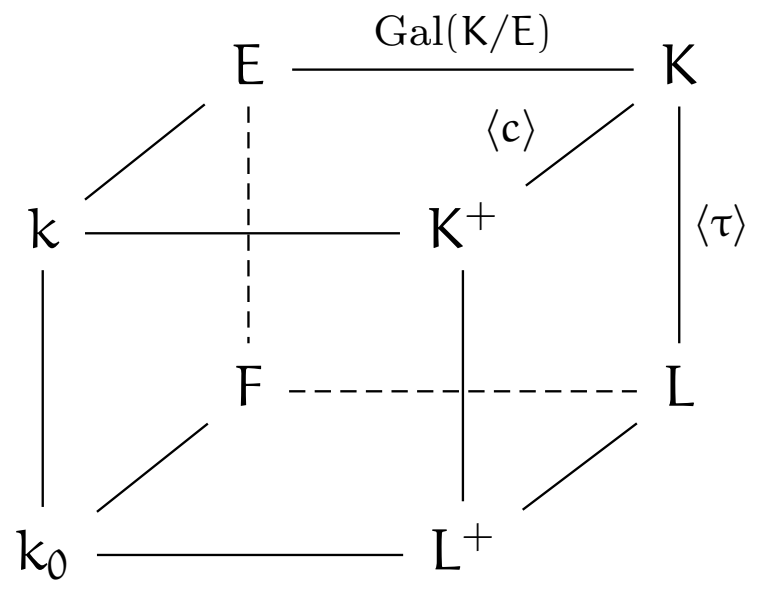

We summarize the field extensions
in Figure~\ref{figure:field-extensions-cube},
in which the vertical extensions have degree $3$,
the diagonal extensions have degree $2$,
and the horizontal extensions have
arbitrarily large Galois group $\Gal(K/E)$.

\FloatBarrier

\subsection{A decomposition of Artin
\texorpdfstring{$L$}{L}-functions}
\label{subsec:decomposition-Artin}

We return to a non-Galois cubic extension $k/k_0$ of
arbitrary number fields, with $K$ containing $F$.
We first fix our conventions for the $L$-functions and
the units of any finite abelian extension $K/k$.

Let $S$ be a finite set of places of $k$
that contains the archimedean places,
and write $W_K := \lvert \mu(K) \rvert$
for the number of roots of unity in $K$.
For a character $\chi$ of $\Gal(K/k)$, let $V_\chi$ be the
one-dimensional complex representation with character $\chi$.
To define its Euler factor at a finite place $v$ of $k$,
we use the inertia subgroup $I_v \subset \Gal(K/k)$ and
the residue field cardinality $\rmn(v)$.
The arithmetic Frobenius element
acts on the residue field at a place
of $K$ above $v$ by $x \longmapsto x^{\rmn(v)}$.
Write $\Frob_v$ for its induced action
on $V_\chi^{I_v}$ and set
\begin{align*}
  L_v(s, \chi) &:= \det(1 - \Frob_v\rmn(v)^{-s} \mid
V_\chi^{I_v})^{-1}, \\
  L_S(s, \chi) &:= \prod_{v \notin S}L_v(s, \chi).
\end{align*}
Tate \cite[Chapter~I, Sections~1 and 3]{tate-stark-1984}
shows that Euler factors are multiplicative in direct sums
and compatible with induction.

For $\sigma \in \Gal(K/k)$, define the partial $L$-function by Fourier
inversion:
\begin{equation}
  \label{eq:partial-zeta}
  L_S(s, \sigma) := \frac{1}{\verts{\Gal(K/k)}}
  \sum_\chi \chi(\sigma)^{-1}L_S(s, \chi).
\end{equation}
Then $L_S(s, \chi) = \sum_{\sigma \in \Gal(K/k)}
\chi(\sigma)L_S(s, \sigma)$ for every character $\chi$ of $\Gal(K/k)$.
When $S$ contains every place ramified in $K/k$,
the partial $L$-function $L_S(s, \sigma)$ agrees
with the partial zeta function
for $\re(s) > 1$:
\begin{equation}
  \label{eq:partial-L-series}
  L_S(s, \sigma)
    = \sum_{\substack{(\mfa, S) = 1 \\ [\mfa, K/k] = \sigma}}
      \rmn(\mfa)^{-s},
\end{equation}
where the sum runs over the nonzero integral ideals
$\mfa$ of $k$ coprime to the finite places in $S$
and the Artin symbol $[\mfa, K/k]$
sends an unramified prime
to the corresponding Frobenius element
and extends multiplicatively to ideals.
Equation~\ref{eq:partial-zeta} also defines
$L_S(s, \sigma)$ when $S$ omits a place ramified in $K/k$.

We next fix the notation for units and their logarithms.
At real places we use the ordinary absolute value,
and at complex places we use its square.
At a finite place $w$ of $K$, set
$\lvert x\rvert_w := \rmn(w)^{-\ord_w(x)}$, where
$\ord_w$ has value group $\Z$.
Write $S_K$ for the places of $K$ above $S$.
The ring $\calo_{K, S}$ of $S$-integers consists of the
$x \in K$ with $\ord_w(x) \geq 0$ at every finite place
$w \notin S_K$. Its unit group $\calo_{K, S}^\times$
consists of the elements of $K^\times$ with zero valuation
at every finite place outside $S_K$.
For an $S$-unit $u$, the product formula gives
$\sum_{w \in S_K}\log\lvert u\rvert_w = 0$.
The $S$-unit theorem
(cf. \cite[Chapter~I, Section~7 and Corollary~11.7]
{neukirch}) identifies
$\R\otimes_\Z\calo_{K, S}^\times$ with the real tuples
$(a_w)_{w \in S_K}$ satisfying $\sum_w a_w = 0$, through
$u \longmapsto (-\log\lvert u\rvert_w)_{w \in S_K}$.
For an automorphism $g$ of $K$, we let $g$ act on places by
$\lvert x\rvert_{g w} := \lvert g^{-1}(x)\rvert_w$.

Our constructions start with units satisfying the integral
rank-$1$ Stark conjecture, including its abelian extension
condition. We record all three requirements here.
Let an archimedean place $v \in S$ split completely in $K$,
and fix a place $w$ of $K$ above $v$.

\begin{definition}
  \label{def:global-rank-one-Stark-unit}
  A unit $\epsilon \in \calo_K^\times$ is a \emph{Stark unit}
  for $(K/k, S, v, w)$ if:
  \begin{enumerate}[(i)]
    \item For every character $\chi$ of $\Gal(K/k)$,
      \[
        L_S'(0, \chi) = -\frac{1}{W_K}
        \sum_{\sigma \in \Gal(K/k)}\chi(\sigma)
        \log\verts{\epsilon^\sigma}_w.
      \]
    \item The unit $\epsilon$ has absolute value $1$ at every
      archimedean place $w'$ of $K$ not above $v$.
    \item The extension $K(\epsilon^{1/W_K})/k$ is abelian.
  \end{enumerate}
\end{definition}

This is the formulation of Stark \cite{stark-1980} for a
global unit.
Applying Fourier inversion
to Equation~\ref{eq:partial-zeta}
shows that part~(i) is equivalent to
$\log\lvert \epsilon^\sigma\rvert_w = -W_K L_S'(0, \sigma)$
for every $\sigma \in \Gal(K/k)$.
Parts~(i) and (ii) determine $\epsilon$ up to a root of unity.
We abbreviate $(K/k, S, v, w)$ to $(K/k, S, w)$ when $v$
is determined by $w$.

We return to the cubic extension $k/k_0$ and $K = kL$.
We use the permutation representation of $S_3$ to express
the Artin $L$-functions over $k$ in terms of $L$-functions
over $k_0$ and $F$.

For each positive integer $m$, write
$\mu_m := \set{\zeta \in \overline{\Q} : \zeta^m = 1}$.
Choose a nontrivial character
$\omega_F : \Gal(E/F) \longrightarrow \mu_3$.
For a character $\chi$ of $\Gal(K/k)$, let $\chi_{k_0}$ be the
corresponding character of $\Gal(L/k_0)$,
and let $\chi_F$ be the restriction of $\chi_{k_0}$ to $\Gal(L/F)$.
Inflate $\chi_F$ and $\omega_F$ to $\Gal(K/F)$ through
restriction to $L$ and $E$, respectively. Their product
$\chi_F\omega_F$ is a character of the abelian extension
$K/F$.
We write $\ind_k^{k_0}$ and $\ind_F^{k_0}$ for induction
of Artin representations to $k_0$, and $\one$ for the
trivial representation. The permutation
representation on the three embeddings of $k$ over $k_0$
is the sum of the trivial representation and the standard
representation $\ind_F^{k_0}(\omega_F)$ of $S_3$. Hence
\begin{equation}
  \label{eq:Artin-decomposition}
  \ind_k^{k_0}(\one) \cong \one_{k_0}
  \oplus \ind_F^{k_0}(\omega_F).
\end{equation}
We twist Equation~\ref{eq:Artin-decomposition} by
$\chi_{k_0}$ to obtain
$\ind_k^{k_0}(\chi) \cong \chi_{k_0} \oplus
\ind_F^{k_0}(\chi_F\omega_F)$.
By the compatibility of Artin $L$-functions with direct
sums and induction, we therefore have
\begin{equation}
  \label{eq:complete-L-factorization}
  L(s, \chi) = L(s, \chi_{k_0})L(s, \chi_F\omega_F).
\end{equation}
The same argument applies to partial $L$-functions when
$S_{k_0}$ is a finite set of places of $k_0$ containing
the archimedean places, and $S$ and $S_F$ consist of the
places above $S_{k_0}$ in $k$ and $F$:
$L_S(s, \chi) = L_{S_{k_0}}(s, \chi_{k_0})
L_{S_F}(s, \chi_F\omega_F)$.
We will use Equation~\ref{eq:complete-L-factorization}
to compare the rational Stark conjectures for $K/k$,
$L/k_0$, and $K/F$.
For the integral results, we will also compare the Euler
factors for general $S$ in
Section~\ref{subsec:decomposition-partial}.

To formulate the higher-rank rational Stark conjectures,
take $S_{k_0}$ to contain all primes ramified in $K/k_0$
and define the rank
$r := \ord_{s = 0} L_S(s, \chi)$.
Choose $S$-units $u_1, \ldots, u_r$ of $K$ and places
$w_1, \ldots, w_r$ above $S$ for which the determinant
below is nonzero for $\chi$ and its Galois conjugates.
The rank-$r$ rational abelian Stark conjecture asserts that
\[
  \frac{L_S^{(r)}(0, \chi)}{r!}
  = a_\chi\det\paren{\sum_{\sigma \in \Gal(K/k)}\chi(\sigma)
  \log\verts{u_j^\sigma}_{w_i}}_{1 \leq i, j \leq r}
\]
for some $a_\chi \in \Q(\chi)^\times$.
For every $\alpha \in \Gal(\Q(\chi)/\Q)$, the units
$u_1, \ldots, u_r$ and places $w_1, \ldots, w_r$ must give
the identity for $\alpha\circ\chi$ with coefficient
$\alpha(a_\chi)$.
The $S$-unit theorem and the order-of-vanishing formula of
Tate \cite[Chapter~I, Proposition~3.4]{tate-stark-1984}
provide choices with nonzero determinant.
For $r = 0$, the determinant is $1$; for $r = 1$, it is
the sum appearing in Theorem~\ref{thm:intro-rank-one}.
Tate \cite[Chapters~I--II]{tate-stark-1984} proves that the
rational Stark conjecture is independent of the choices of units,
places, and the finite set $S$,
and is compatible with direct sums,
induction, and inflation.
We can therefore use the decomposition in
Equation~\ref{eq:Artin-decomposition} to deduce the
rational Stark conjecture over $k$ from the conjectures
over $k_0$ and $F$.

\begin{corollary}
  \label{cor:relative-rational-Stark}
  Let $k/k_0$ be a non-Galois cubic extension
  and let $K/k$ be a finite abelian extension such that
  $F \subset K \subset k \cdot k_0^\ab$.
  The rational Stark conjectures for $L/k_0$ and $K/F$
  imply the rational Stark conjecture for $K/k$.
\end{corollary}

\subsection{A decomposition of partial
  \texorpdfstring{$L$}{L}-functions}
\label{subsec:decomposition-partial}

We have proved the decomposition of partial $L$-functions
when $S$ and $S_F$ consist of all places above $S_{k_0}$.
For a general finite set $S$ of places of $k$ containing
the archimedean places, we now choose finite sets
$S_{k_0}$ and $S_F$ of places of $k_0$ and $F$,
each containing the archimedean places, so that
\[
  L_S(s, \chi)
    = L_{S_{k_0}}(s, \chi_{k_0})
      L_{S_F}(s, \chi_F\omega_F).
\]

The integral rank-$1$ Stark conjecture uses a finite set $S$
containing the archimedean and ramified places and a place
$v \in S$ that splits completely. If $v$ is archimedean and
$\lvert S\rvert \geq 3$, the conjectured Stark unit is a global unit
(Stark \cite{stark-1980}; Tate
\cite[Chapter~IV, Section~2]{tate-stark-1984}).
The explicit constructions in
Sections~\ref{sec:rank-one} and
\ref{subsec:relative-rank-two-rational} use a Stark unit
for $K/F$ in rank $1$ and an additional Stark unit
for $L^+/k_0$ in rank $2$.
For these constructions, we choose sets satisfying
Equation~\ref{eq:local-factorization} below and impose
the following additional conditions on the finite places:
\begin{enumerate}[(i)]
  \item In rank $1$: let $S$ contain
    every finite place ramified in $K/k$,
    let $S_{k_0}$ contain every finite place ramified in $L/k_0$,
    and let $S_F$ contain at least two finite places in $F$
    and every finite place ramified in $K/F$.
  \item In rank $2$: let $S$, $S_{k_0}$, and
    $S_F$ contain the finite places ramified in $K^+/k$,
    $L^+/k_0$, and $K/F$, respectively;
    furthermore, suppose that $S_{k_0}$ and $S_F$
    each contain at least two finite places.
\end{enumerate}
The two finite places ensure that $\lvert S_F\rvert \geq 3$ and,
in rank $2$, that $\lvert S_{k_0}\rvert \geq 3$.
The rank-$1$ Stark units over $F$ and $k_0$ are therefore
global units. In Theorem~\ref{thm:intro-rank-two}, taking
all places above two rational primes ensures two finite
places in both $S_\Q$ and $S_F$. Two finite places of $k$
need not suffice for this purpose, since they may lie above
one rational prime. For Theorem~\ref{thm:intro-rank-one},
Lemma~\ref{lem:rank-one-single-prime} also treats the case
where $S_F$ has only one finite place.

We now return to arbitrary finite sets $S$ containing
the archimedean places. Starting from
Equation~\ref{eq:complete-L-factorization}, we compare
the products of Euler factors at the finite places in
$S$, $S_{k_0}$, and $S_F$. The induction property in
Tate \cite[Chapter~I, Section~3]{tate-stark-1984} lets us
make this comparison through local representations.
The Weil group of a nonarchimedean local field is the
subgroup of its absolute Galois group whose residue action
is an integral power of Frobenius.
We write $W_{k_v}$, $W_{F_\mfq}$, and $W_{(k_0)_\mfp}$
for the Weil groups of the completions
$k_v$, $F_\mfq$, and $(k_0)_\mfp$.
For a set $A$ of places, write $\mathbf{1}_A(v) = 1$
if $v \in A$ and $\mathbf{1}_A(v) = 0$ otherwise.
We require, for every finite prime $\mfp$ of $k_0$,
\begin{equation}
  \label{eq:local-factorization}
    \bigoplus_{\substack{v \in S\\v \mid \mfp}}
    \ind_{W_{k_v}}^{W_{(k_0)_\mfp}}\one
    \cong \one^{\oplus\mathbf{1}_{S_{k_0}}(\mfp)}
    \oplus
    \bigoplus_{\substack{\mfq \in S_F\\\mfq \mid \mfp}}
    \ind_{W_{F_\mfq}}^{W_{(k_0)_\mfp}}
    \paren{\restr{\omega_F}{W_{F_\mfq}}}.
\end{equation}

When $S$ contains every place of $k$ above $\mfp$,
Equation~\ref{eq:local-factorization} requires $\mfp$ in
$S_{k_0}$ and every place of $F$ above $\mfp$ in $S_F$.
We now determine $S_{k_0}$ and $S_F$ for arbitrary $S$.

\begin{proposition}
  \label{prop:local-Artin-factorization}
  For every finite set $S$ containing the archimedean
  places of $k$, there is a unique pair $(S_{k_0}, S_F)$
  satisfying Equation~\ref{eq:local-factorization}
  with $S_F$ invariant under $\Gal(F/k_0)$.
  Furthermore, for every character $\chi$ of $\Gal(K/k)$,
  we have
  \[
    L_S(s, \chi) = L_{S_{k_0}}(s, \chi_{k_0})
    L_{S_F}(s, \chi_F\omega_F).
  \]
\end{proposition}

\begin{proof}
  We determine the local representations attached to the
  places in $S$ for each decomposition group in
  $\Gal(E/k_0)$, then compare their Euler factors.
  Fix a finite prime $\mfp$ of $k_0$ and a decomposition group
  $D_\mfp \subset \Gal(E/k_0) = S_3$.
  The places of $k$ above $\mfp$ correspond to the
  $D_\mfp$-orbits on $S_3/S_2$, where $S_2 = \Gal(E/k)$.
  The places of $F$ above $\mfp$ correspond to the
  $D_\mfp$-orbits on $S_3/A_3$.
  Up to conjugacy, the possibilities for $D_\mfp$ are
  $1$, $S_2$, $A_3$, and $S_3$.

  If $D_\mfp = 1$, there are three places of $k$ and
  two conjugate places of $F$ above $\mfp$, and every
  local summand is trivial. If $S$ contains none of the three
  places, include no
  place above $\mfp$ in $S_{k_0}$ or $S_F$.
  If $S$ contains one place, include only $\mfp$ in
  $S_{k_0}$. If $S$ contains two places, include both places
  of $F$ above $\mfp$ in $S_F$ and omit $\mfp$ from $S_{k_0}$.
  If $S$ contains all three places, include $\mfp$ in
  $S_{k_0}$ and both places of $F$ above $\mfp$ in $S_F$. These
  choices are the only ones
  invariant under $\Gal(F/k_0)$ that have the required
  dimension.

  If $D_\mfp \cong S_2$, its action on
  $S_3/S_2$ has a fixed point and an orbit of
  size $2$. The corresponding completions of $k$ have
  degrees $1$ and $2$ over $(k_0)_\mfp$.
  The degree-$1$ place contributes the trivial
  representation. Include $\mfp$ in $S_{k_0}$ precisely
  when $S$ contains the degree-$1$ place.
  The degree-$2$ place contributes the regular
  representation of $S_2$, induced from the unique place
  of $F$ above $\mfp$. Include that place of $F$ in $S_F$
  precisely when $S$ contains the degree-$2$ place of $k$.
  These choices of $S_{k_0}$ and $S_F$ are independent.

  If $D_\mfp \cong A_3$, there is one place of $k$
  and two conjugate places of $F$ above $\mfp$.
  The representation for the place of $k$ is
  $\one \oplus \omega_F \oplus \omega_F^{-1}$.
  The two places of $F$ contribute $\omega_F$ and
  $\omega_F^{-1}$. If the place of $k$ lies in $S$, include
  $\mfp$ in $S_{k_0}$ and both places of $F$ in $S_F$.
  Otherwise, include none of these places.

  If $D_\mfp = S_3$, there is one place in each of
  $k$ and $F$ above $\mfp$. By Equation
  \ref{eq:Artin-decomposition}, the place of $k$ belongs
  to $S$ exactly when $\mfp$ belongs to $S_{k_0}$ and
  the place of $F$ belongs to $S_F$.

    We have now determined the unique pair $(S_{k_0}, S_F)$.
  To compare the Euler factors for $\chi$, we tensor
  Equation~\ref{eq:local-factorization} by
  $\restr{\chi_{k_0}}{W_{(k_0)_\mfp}}$.
  Since $\chi$ and $\chi_F$ are restrictions of $\chi_{k_0}$,
  we can move this tensor product inside each induction.
  Using the induction and direct-sum properties in Tate
  \cite[Chapter~I, Section~3]{tate-stark-1984}, we obtain
  \begin{align*}
    \prod_{\substack{v \in S \\ v \mid \mfp}}L_v(s, \chi)
      &= L_\mfp\biggl(s,
        \restr{\chi_{k_0}}{W_{(k_0)_\mfp}} \otimes
        \bigoplus_{\substack{v \in S \\ v \mid \mfp}}
        \ind_{W_{k_v}}^{W_{(k_0)_\mfp}}\one\biggr) \\
      &= L_\mfp(s, \chi_{k_0})^{
        \mathbf{1}_{S_{k_0}}(\mfp)}
        \prod_{\substack{\mfq \in S_F \\ \mfq \mid \mfp}}
        L_\mfp\biggl(s,
        \ind_{W_{F_\mfq}}^{W_{(k_0)_\mfp}}
        \paren{\restr{\chi_F\omega_F}{W_{F_\mfq}}}\biggr) \\
      &= L_\mfp(s, \chi_{k_0})^{\mathbf{1}_{S_{k_0}}(\mfp)}
        \prod_{\substack{\mfq \in S_F \\ \mfq \mid \mfp}}
        L_\mfq(s, \chi_F\omega_F).
  \end{align*}
  We multiply this equality over the finite primes $\mfp$
  and remove the resulting Euler factors from
  Equation~\ref{eq:complete-L-factorization} to obtain
  $L_S(s, \chi) = L_{S_{k_0}}(s, \chi_{k_0})
  L_{S_F}(s, \chi_F\omega_F)$.
\end{proof}

The proof of Proposition~\ref{prop:local-Artin-factorization}
establishes the decomposition of partial $L$-functions
for any triple $(S, S_{k_0}, S_F)$ satisfying
Equation~\ref{eq:local-factorization}.
The proofs of Theorems~\ref{thm:intro-rank-one},
\ref{thm:intro-rank-two}, and \ref{thm:intro-relative}
begin with a fixed set of places of $k_0$
rather than a fixed set of places of $k$;
we write down the decomposition of partial $L$-functions
in this case.

\begin{corollary}
  \label{cor:rationally-saturated-factorization}
  Let $S_{k_0}$ be a finite set of places of $k_0$,
  including every archimedean place,
  and let $S$ and $S_F$ be the set of places
  above $S_{k_0}$ in $k$ and $F$,
  respectively.
  Then Equation~\ref{eq:local-factorization} holds
  and for every character $\chi$ of $\Gal(K/k)$,
  \begin{equation*}
    L_S(s, \chi) = L_{S_{k_0}}(s, \chi_{k_0})
    L_{S_F}\paren{s, \chi_F\omega_F}.
  \end{equation*}
\end{corollary}

\begin{proof}
  At each finite prime in $S_{k_0}$, we restrict
  Equation~\ref{eq:Artin-decomposition} to the local Weil
  group to obtain Equation~\ref{eq:local-factorization}.
  Outside $S_{k_0}$, both sides of
  Equation~\ref{eq:local-factorization} are zero.
  Therefore, we can apply
  Proposition~\ref{prop:local-Artin-factorization}
  to deduce the decomposition of partial $L$-functions.
\end{proof}

For other choices of $S_{k_0}'$ and $S_F'$, we obtain
a decomposition of partial $L$-functions with an
explicit quotient of Euler factors.

\begin{proposition}
  \label{prop:local-correction-factors}
  Let $S_{k_0}'$ and $S_F'$ be finite sets of places of
  $k_0$ and $F$ containing the archimedean places.
  For every character $\chi$ of $\Gal(K/k)$, we have
  \[
    L_S(s, \chi)
      = L_{S_{k_0}'}(s, \chi_{k_0})
        L_{S_F'}(s, \chi_F\omega_F)
        \frac{\prod_{\mfp \in S_{k_0}'}L_\mfp(s, \chi_{k_0})
        \prod_{\mfq \in S_F'}L_\mfq(s, \chi_F\omega_F)}
        {\prod_{v \in S}L_v(s, \chi)},
  \]
  where the products range over finite places.
\end{proposition}

\begin{proof}
  To express $L_S(s, \chi)$ in terms of
  $L_{S_{k_0}'}(s, \chi_{k_0})$ and
  $L_{S_F'}(s, \chi_F\omega_F)$, we use
  Equation~\ref{eq:complete-L-factorization}:
  \begin{align*}
    L_S(s, \chi)
      &= \frac{L(s, \chi)}{\prod_{v \in S}L_v(s, \chi)} \\
      &= \frac{L(s, \chi_{k_0})L(s, \chi_F\omega_F)}
        {\prod_{v \in S}L_v(s, \chi)} \\
      &= L_{S_{k_0}'}(s, \chi_{k_0})
        L_{S_F'}(s, \chi_F\omega_F)
        \frac{\prod_{\mfp \in S_{k_0}'}L_\mfp(s, \chi_{k_0})
        \prod_{\mfq \in S_F'}L_\mfq(s, \chi_F\omega_F)}
        {\prod_{v \in S}L_v(s, \chi)}.
  \end{align*}
\end{proof}


\section{The integral rank-\texorpdfstring{$1$}{1}
Stark conjecture}
\label{sec:rank-one}

In this section,
we prove the relative implication in
Theorem~\ref{thm:intro-relative}(i),
specialize to
$k_0 = \Q$ to obtain Theorem~\ref{thm:intro-rank-one},
and then conclude with
an application of Theorem~\ref{thm:intro-rank-one}
to explicit class field theory.

For the construction below, let $k_0$ be totally real,
let $k/k_0$ be an almost totally real cubic extension
with quadratic resolvent $F$, and let $K/k$ be finite
abelian with $K \subset k \cdot k_0^\ab$.
Recall that $L = K \cap k_0^\ab$ and $K = kL$.

\subsection{A relative theorem on the integral
rank-\texorpdfstring{$1$}{1} Stark conjecture}
\label{subsec:relative-rank-one-theorem}

We first construct a Stark unit for $K/k$ from a Stark
unit for $K/F$, assuming that $K$ contains $F$,
or equivalently that $L$ contains $F$.
We then descend from $KF$ to $K$ to remove this assumption
and prove Theorem~\ref{thm:intro-relative}(i)
at the complex place of $k$
(in the other cases, $1$ is a Stark unit).

Recall that $v_0$ is the real place of $k_0$ below the
complex place of $k$,
and $c$ is complex conjugation at
an embedding extending the real place of $k$ above $v_0$.
Choose a generator $\tau$ of $\Gal(K/L)$.

\subsubsection{Stickelberger elements}
\label{subsec:Stickelberger-elements}

We construct units by taking products of Galois conjugates
with exponents determined by partial $L$-values at $s = 0$.
For $x \in K^\times$ and
$a = \sum_g n_g g \in \Z[\Gal(K/k)]$, write
$x^a := \prod_g g(x)^{n_g}$.
For any finite abelian extension $K/k$ and finite set $S$
of places of $k$ containing the archimedean places, set
\[
  \Theta_{K/k, S}(s) := \sum_{\sigma \in \Gal(K/k)}
  L_S(s, \sigma)\sigma^{-1}.
\]
We evaluate characters on the group ring by $\C$-linear
extension. Fourier inversion in
Equation~\ref{eq:partial-zeta} gives
$\chi^{-1}(\Theta_{K/k, S}(s)) = L_S(s, \chi)$.

For the rank-$1$ construction, we use the partial
$L$-values for $L/k_0$.
Let $S_{k_0}$ contain the archimedean places and all primes
ramified in $L/k_0$, and set
$\Theta_{k_0}(0) := \Theta_{L/k_0, S_{k_0}}(0)$.
We regard this element as belonging to $\Q[\Gal(K/k)]$
via the restriction isomorphism
$\Gal(K/k) \cong \Gal(L/k_0)$.
We need two integrality properties.
The first allows us to use $W_L\Theta_{k_0}(0)$ as an
exponent in a product of conjugates:
\begin{equation}
  \label{eq:Stickelberger-integrality}
  W_L\Theta_{k_0}(0) \in \Z[\Gal(K/k)].
\end{equation}
To verify the condition
on $K(\epsilon^{1/W_K})/k$ being abelian in
Definition~\ref{def:global-rank-one-Stark-unit}(iii),
we will also use the following property:
if $g \in \Gal(K/k)$ and $m \in \Z$ satisfy
$g(\zeta) = \zeta^m$ for every $\zeta \in \mu(L)$, then
\begin{equation}
  \label{eq:Stickelberger-action-integrality}
  (g - m)\Theta_{k_0}(0) \in \Z[\Gal(K/k)].
\end{equation}

We first prove Equations~\ref{eq:Stickelberger-integrality}
and \ref{eq:Stickelberger-action-integrality} when $k_0 = \Q$.
Let $f$ be the conductor of $L/\Q$, and choose a primitive
$f$-th root of unity $\zeta_f$.
For an integer $a$ coprime to $f$, let $\sigma_a$ be the
restriction to $L$ of $\zeta_f \longmapsto \zeta_f^a$.
The Hurwitz zeta value $\zeta(0, a/f) = 1/2 - a/f$
for $1 \leq a < f$ gives the sum
\[
  \theta_f := \sum_{\substack{1 \leq a < f \\
    (a, f) = 1}}
  \paren{\frac{1}{2} - \frac{a}{f}}\sigma_a^{-1}.
\]
Removing the Euler factors at the additional finite
primes in $S_\Q$ gives
\[
  \Theta_\Q(0) = \theta_f
  \prod_{\substack{\ell \in S_\Q \text{ finite} \\
    \ell \nmid f}}(1 - \sigma_\ell^{-1}).
\]

For an odd positive integer $b$ coprime to $f$,
multiplication by $b$ permutes the invertible residue
classes modulo $f$.
Reindexing $\sigma_b\theta_f$ using the least positive
residue $ab - f\lfloor ab/f\rfloor$ gives
\[
  (b - \sigma_b)\theta_f
  = \sum_{\substack{1 \leq a < f \\
    (a, f) = 1}}
  \paren{\frac{b - 1}{2} - \left\lfloor\frac{ab}{f}\right\rfloor}
  \sigma_a^{-1}.
\]
Since $b$ is odd, every coefficient in this sum is an
integer. The factors $1 - \sigma_\ell^{-1}$ are also
integral, so
$(b - \sigma_b)\Theta_\Q(0) \in \Z[\Gal(K/k)]$.

To prove Equation~\ref{eq:Stickelberger-integrality},
take $b = \ell$ for an odd prime $\ell$ outside $S_\Q$
that splits completely in $L$.
Then $\sigma_\ell = 1$, so
$(\ell - 1)\Theta_\Q(0) \in \Z[\Gal(K/k)]$.
Stark \cite[p.~220]{stark-1980} observes that the greatest
common divisor of these integers $\ell - 1$ is $W_L$.
Expressing $W_L$ as an integral linear combination of
finitely many such integers proves
Equation~\ref{eq:Stickelberger-integrality}.

To prove Equation~\ref{eq:Stickelberger-action-integrality},
let $g \in \Gal(K/k)$ and $m \in \Z$ satisfy
$g(\zeta) = \zeta^m$ for every $\zeta \in \mu(L)$.
Dirichlet's theorem on primes in arithmetic progressions
provides an odd prime $\ell$ outside $S_\Q$, coprime to
$fW_L$, with $\sigma_\ell = g$.
Comparing the actions of $g$ and $\sigma_\ell$ on $\mu(L)$
gives $\ell \equiv m \pmod{W_L}$.
The integrality of $(\ell - g)\Theta_\Q(0)$, proved with
$b = \ell$, and Equation~\ref{eq:Stickelberger-integrality}
now give
\[
  (g - m)\Theta_\Q(0)
  = -(\ell - g)\Theta_\Q(0)
  + \frac{\ell - m}{W_L}W_L\Theta_\Q(0)
  \in \Z[\Gal(K/k)].
\]

For a general totally real base field $k_0$,
the integrality result proved independently by
Cassou-Nogu\`es \cite{cassou-nogues} and
Deligne--Ribet \cite[Theorem~0.4]{deligne-ribet}
implies that $a\Theta_{k_0}(0) \in \Z[\Gal(K/k)]$
whenever $a \in \Z[\Gal(K/k)]$ annihilates $\mu(L)$
(here $\Gal(K/k)$ acts on $\mu(L)$ by restriction to $L$).
Since $\zeta^{W_L} = 1$ for every $\zeta \in \mu(L)$,
taking $a = W_L$ proves
Equation~\ref{eq:Stickelberger-integrality}.
If $g(\zeta) = \zeta^m$ for every $\zeta \in \mu(L)$,
then $g - m$ also annihilates $\mu(L)$.
Taking $a = g - m$ proves
Equation~\ref{eq:Stickelberger-action-integrality}.

We next determine how $\Theta_{k_0}(0)$ behaves under the
fixed complex conjugation $c$.
We call a character $\chi$ of $\Gal(K/k)$ \emph{odd}
if $\chi(c) = -1$ and \emph{even} if $\chi(c) = 1$.
Under restriction, the corresponding character
$\chi_{k_0}$ of $\Gal(L/k_0)$ has the same value on $c$.
We show that $L_{S_{k_0}}(0, \chi_{k_0}) = 0$
for every even $\chi$.

For a nontrivial even character $\chi$, the character
$\chi_{k_0}$ is trivial on complex conjugation at $v_0$.
Thus $L_{S_{k_0}}(0, \chi_{k_0}) = 0$ by the
order-of-vanishing formula of
\cite[Chapter~I, Proposition~3.4]{tate-stark-1984}.
For the trivial character, if $k_0 \neq \Q$, the Dedekind
zeta function of $k_0$ has positive order at $0$.
If $k_0 = \Q$, the set $S_{k_0}$ contains a prime $\ell$
ramified in $F/\Q$; removing its Euler factor introduces
$1 - \ell^{-s}$, which vanishes at $0$.
Thus the trivial character also gives zero.
For odd $\chi$, we have $\chi^{-1}(1 + c) = 0$.
Evaluating at every character therefore proves
\begin{equation}
  \label{eq:odd-Stickelberger}
  (1 + c)\Theta_{k_0}(0) = 0.
\end{equation}

\subsubsection{Conjugation of the logarithms of a Stark unit}
\label{subsec:rank-one-conjugation}

Before comparing the regulators over $k$ and $F$,
we prove a symmetry of the logarithms of a Stark unit
for $K/F$.
Let $S_F$ contain the archimedean places of $F$,
every finite place ramified in $K/F$, and at least two
finite places. Fix a place $w$ of $K$ above the complex
place of $F$, and let $\epsilon_F \in \calo_K^\times$
be a Stark unit for $(K/F, S_F, w)$.
Assume either that $S_F$ is invariant under $\Gal(F/k_0)$
or that $S_F$ belongs to a triple $(S, S_{k_0}, S_F)$
satisfying Equation~\ref{eq:local-factorization}.
We show that $\log\lvert\epsilon_F^h\rvert_w$ is unchanged
when $h$ is replaced by $c h c^{-1}$.

\begin{lemma}
  \label{lem:elliptic-conjugation-symmetry}
  The unit $\epsilon_F \in \calo_K^\times$
  satisfies
  \begin{equation*}
    \log\verts{\epsilon_F^h}_w
    = \log\verts{\epsilon_F^{c h c^{-1}}}_w
  \end{equation*}
  for every $h \in \Gal(K/F)$. In particular,
  $\log\lvert \epsilon_F^{\beta\tau}\rvert_w =
  \log\lvert \epsilon_F^{\beta\tau^{-1}}\rvert_w$
  for every $\beta \in \Gal(K/kF)$.
\end{lemma}

\begin{proof}
  For a character $\psi$ of $\Gal(K/F)$, set
  $\psi^c(h) := \psi(c h c^{-1})$.
  Conjugation on $F$ permutes the Euler factors and gives
  $L_{S_F}(s, \psi) = L_{cS_F}(s, \psi^c)$.
  If $S_F = cS_F$, we obtain
  $L_{S_F}(s, \psi) = L_{S_F}(s, \psi^c)$.
  Otherwise, the proof of
  Proposition~\ref{prop:local-Artin-factorization} shows that
  $S_F$ can contain just one of two conjugate places above
  $\mfp$ only when $\mfp$ splits completely in $E/k_0$.
  The decomposition groups in $\Gal(K/F)$ at these places
  lie in $\Gal(K/kF)$, which $c$ centralizes.
  The Euler factors at the two places are therefore equal
  for every $\psi$. Replacing $S_F$ by $cS_F$ gives
  $L_{S_F}(s, \psi) = L_{S_F}(s, \psi^c)$ in this case too.

  By Fourier inversion in Equation~\ref{eq:partial-zeta},
  we therefore have
  $L_{S_F}(s, h) = L_{S_F}(s, c h c^{-1})$.
  Since $\epsilon_F$ is a Stark unit, we also have
  $\log\lvert\epsilon_F^h\rvert_w = -W_K L_{S_F}'(0, h)$
  by Definition~\ref{def:global-rank-one-Stark-unit}(i).
  We conclude that
  $\log\lvert\epsilon_F^h\rvert_w
  = \log\lvert\epsilon_F^{c h c^{-1}}\rvert_w$.
\end{proof}

\subsubsection{The rank-\texorpdfstring{$1$}{1} Stark regulator}
\label{subsec:rank-one}

Choose finite sets $S$, $S_{k_0}$, and $S_F$ of places of
$k$, $k_0$, and $F$, respectively, each containing all
archimedean places of its field and satisfying
Equation~\ref{eq:local-factorization}.
Assume also that
\begin{itemize}
  \item $S$ contains every finite place ramified in $K/k$;
  \item $S_{k_0}$ contains every finite place ramified in
    $L/k_0$; and
  \item $S_F$ contains every finite place ramified in $K/F$
    and at least two finite places.
\end{itemize}
These are the hypotheses in
Section~\ref{subsec:decomposition-partial}(i).
Fix a place $w$ of $K$ above the complex place of $k$.
The nonidentity element of the decomposition group at $w$
in $\Gal(K/k_0)$ and the element $c$ restrict to distinct
transpositions
in $\Gal(E/k_0)$. Their product generates $\Gal(K/L)$.
After replacing $\tau$ by $\tau^{-1}$ if necessary, we obtain
\begin{equation}
  \label{eq:place-normalization}
  c w = \tau w.
\end{equation}

Let $\epsilon_F \in \calo_K^\times$ be a Stark unit for
$(K/F, S_F, w)$. The decomposition of partial $L$-functions
in Proposition~\ref{prop:local-Artin-factorization} suggests applying
$\Theta_{k_0}(0)$ to $\epsilon_F$ to incorporate the factor
$L_{S_{k_0}}(0, \chi_{k_0})$.
Equation~\ref{eq:Stickelberger-integrality} allows us to define
\begin{equation}
  \label{eq:rank-one-unit-definition}
  u := \epsilon_F^{W_L\Theta_{k_0}(0)} \in \calo_K^\times.
\end{equation}

To compute the logarithms of $u$, we use character sums.
For a character $\chi$ of $\Gal(K/k)$ and $x \in K^\times$,
write
\begin{align}
  \label{eq:calr}
  \calr_w(x, \chi)
    := \sum_{\sigma \in \Gal(K/k)}\chi(\sigma)
      \log\verts{x^\sigma}_w \in \C.
\end{align}
Reindexing the sum over $\Gal(K/k)$ gives
the following identity for $a \in \Z[\Gal(K/k)]$:
\begin{equation}
  \label{eq:character-logarithm-action}
  \calr_w(x^a, \chi) = \chi^{-1}(a) \calr_w(x, \chi).
\end{equation}
For an odd character $\chi$, we can express
$\calr_w(\epsilon_F, \chi)$ in terms of
$L_{S_F}'(0, \chi_F\omega_F)$.
Recall that $\omega_F$ is a nontrivial character of
$\Gal(E/F)$ and $\chi_F$ is the restriction of
$\chi_{k_0}$ to $\Gal(L/F)$;
in the following lemma,
both are then viewed as characters of $\Gal(K/F)$
via inflation.

\begin{lemma}
  \label{lem:elliptic-regulator-comparison}
  Let $\chi$ be an odd character of $\Gal(K/k)$. Then
  \[
    \calr_w(\epsilon_F, \chi)
      = -W_K L_{S_F}'\paren{0, \chi_F\omega_F}.
  \]
\end{lemma}

\begin{proof}
  We rewrite the regulator over $K/F$ as a sum over
  $\Gal(K/kF)$ and then compare it with
  $\calr_w(\epsilon_F, \chi)$.
  Definition~\ref{def:global-rank-one-Stark-unit}(i) gives
  \[
    -W_K L_{S_F}'(0, \chi_F\omega_F)
    = \sum_{\beta \in \Gal(K/kF)}\chi(\beta)
    \sum_{j = 0}^2\omega_F(\tau^j)
    \log\verts{\epsilon_F^{\beta\tau^j}}_w.
  \]
  By Lemma~\ref{lem:elliptic-conjugation-symmetry},
  the logarithms for $j = 1$ and $j = 2$ agree.
  Since $\omega_F(\tau) + \omega_F(\tau^2) = -1$,
  we obtain
  \[
    -W_K L_{S_F}'(0, \chi_F\omega_F)
    = \sum_{\beta \in \Gal(K/kF)}\chi(\beta)
    \paren{\log\verts{\epsilon_F^\beta}_w
    - \log\verts{\epsilon_F^{\beta\tau}}_w}.
  \]
  Equation~\ref{eq:place-normalization} and the convention
  $\lvert x\rvert_{gw} = \lvert g^{-1}(x)\rvert_w$ give
  \[
    \log\verts{\epsilon_F^{\beta c}}_w
    = \log\verts{\epsilon_F^{\beta\tau^{-1}}}_w
    = \log\verts{\epsilon_F^{\beta\tau}}_w,
  \]
  where the second equality again uses
  Lemma~\ref{lem:elliptic-conjugation-symmetry}.
  Finally, $\chi(c) = -1$ and
  $\Gal(K/k) = \Gal(K/kF)\times\abrac{c}$ imply
  \begin{align*}
    \calr_w(\epsilon_F, \chi)
    &= \sum_{\beta \in \Gal(K/kF)}\chi(\beta)
      \paren{\log\verts{\epsilon_F^\beta}_w
      - \log\verts{\epsilon_F^{\beta c}}_w} \\
    &= -W_K L_{S_F}'(0, \chi_F\omega_F).
  \end{align*}
\end{proof}

We combine Lemma~\ref{lem:elliptic-regulator-comparison}
with the decomposition of partial $L$-functions in
Proposition~\ref{prop:local-Artin-factorization}
to compute the regulator of $u$.

\begin{proposition}
  \label{prop:rank-one-regulator}
  The unit $u = \epsilon_F^{W_L\Theta_{k_0}(0)}
  \in \calo_K^\times$ satisfies
  \[
    L_S'(0, \chi) = -\frac{1}{W_K W_L}
    \sum_{\sigma \in \Gal(K/k)}\chi(\sigma)
    \log\verts{u^\sigma}_w
  \]
  for every character $\chi$ of $\Gal(K/k)$.
  Moreover, $\lvert u\rvert_{w'} = 1$ at every archimedean
  place $w'$ of $K$ above a real place of $k$.
\end{proposition}

\begin{proof}
  Let $\chi$ be a character of $\Gal(K/k)$.
  Since $F$ has a complex place, we have
  $L_{S_F}(0, \chi_F\omega_F) = 0$ by
  \cite[Chapter~I, Proposition~3.4]{tate-stark-1984}.
  We can therefore differentiate the decomposition in
  Proposition~\ref{prop:local-Artin-factorization} to obtain
  \[
    L_S'(0, \chi) = L_{S_{k_0}}(0, \chi_{k_0})
    L_{S_F}'(0, \chi_F\omega_F).
  \]
  For odd $\chi$, we compute the regulator of
  $u = \epsilon_F^{W_L\Theta_{k_0}(0)}$ using
  Equation~\ref{eq:character-logarithm-action} and
  Lemma~\ref{lem:elliptic-regulator-comparison}:
  \begin{align*}
    \calr_w(u, \chi)
      &= W_L L_{S_{k_0}}(0, \chi_{k_0})
    \calr_w(\epsilon_F, \chi) \\
      &= -W_K W_L L_S'(0, \chi).
  \end{align*}
  For even $\chi$, both $\calr_w(u, \chi)$ and
  $L_{S_{k_0}}(0, \chi_{k_0})$ vanish by
  Equation~\ref{eq:odd-Stickelberger}, so the same identity
  holds.
  We have therefore proved that
  \[
    \calr_w(u, \chi) = -W_K W_L L_S'(0, \chi)
  \]
  for every character $\chi$ of $\Gal(K/k)$.
  By Fourier inversion in Equation~\ref{eq:partial-zeta},
  we can also write this identity as
  \[
    L_S'(0, \sigma) = -\frac{1}{W_K W_L}
      \log \verts{u^\sigma}_w
  \]
  for every $\sigma \in \Gal(K/k)$.

  We now check the absolute values above real places of $k$.
  Let $w'$ lie above the real place of $k$ above $v_0$.
  Equations~\ref{eq:Stickelberger-integrality} and
  \ref{eq:odd-Stickelberger} give
  $W_L\Theta_{k_0}(0) \in (1 - c)\Z[\Gal(K/k)]$.
  Since $c w' = w'$, each quotient
  $\epsilon_F^g/\epsilon_F^{cg}$ has absolute value $1$ at
  $w'$. Hence $\lvert u\rvert_{w'} = 1$.
  At any place $w'$ of $K$ above a real place of $k_0$
  different from $v_0$, every conjugate of $\epsilon_F$
  occurring in $u$ has absolute value $1$ by
  Definition~\ref{def:global-rank-one-Stark-unit}(ii).
  Hence $\lvert u\rvert_{w'} = 1$ at every place of $K$ above
  a real place of $k$.
\end{proof}

\subsubsection{Construction of a Stark unit}
\label{subsec:rank-one-construction}

We now construct $\epsilon \in \calo_K^\times$ with
$u = \zeta\epsilon^{W_L}$ for some $\zeta \in \mu(K)$.
By Proposition~\ref{prop:rank-one-regulator}, we then have
$\calr_w(\epsilon, \chi) = -W_K L_S'(0, \chi)$
for every character $\chi$ of $\Gal(K/k)$.
We obtain $\epsilon$ from the abelian extension
$K(\epsilon_F^{1/W_K})/F$, then apply
Equation~\ref{eq:Stickelberger-action-integrality} to prove that
$K(\epsilon^{1/W_K})/k$ is abelian.

\begin{proposition}
  \label{prop:relative-rank-one}
  Let
  \begin{itemize}
    \item $k_0$ be a totally real number field;
    \item $k$ be an almost totally real cubic extension
      of $k_0$ with quadratic resolvent $F$;
    \item $K/k$ be finite abelian such that
      $F \subset K \subset k \cdot k_0^\ab$;
    \item $S$ be a finite set of places of $k$,
      including every archimedean place and every finite
      place ramified in $K/k$;
    \item $S_{k_0}$ be a finite set of places of $k_0$,
      including every archimedean place and every finite
      place ramified in $L/k_0$;
    \item $S_F$ be a finite set of places of $F$,
      including every archimedean place, every finite
      place ramified in $K/F$, and at least two finite
      places; and
    \item $w$ be a place of $K$ above the complex
      place of $k$.
  \end{itemize}
  If $S$, $S_{k_0}$, and $S_F$ satisfy
  Equation~\ref{eq:local-factorization}
  and $\epsilon_F$ is a Stark unit for $(K/F, S_F, w)$,
  then there is a unit $\epsilon \in \calo_K^\times$
  such that
  \begin{enumerate}[(i)]
    \item $\epsilon_F^{W_L\Theta_{k_0}(0)} = \zeta\epsilon^{W_L}$
      for some $\zeta \in \mu(K)$
      and for $W_L\Theta_{k_0}(0)$ from
      \eqref{eq:Stickelberger-integrality};
    \item $K(\epsilon^{1/W_K})/k$ is abelian.
  \end{enumerate}
  Consequently, $\epsilon$ is a Stark unit for
  $(K/k, S, w)$.
\end{proposition}

\begin{proof}
  Let $u := \epsilon_F^{W_L\Theta_{k_0}(0)}$.
  We first construct $\epsilon$ and verify that
  $u = \zeta\epsilon^{W_L}$.
  We then prove that $K(\epsilon^{1/W_K})/k$ is Galois
  and that its automorphisms commute.
  Finally, we verify the regulator and the absolute values
  required of a Stark unit.

  \emph{Step 1. Construction of $\epsilon$.}
  Recall that $\Gal(K/k) = \Gal(K/kF)\times\abrac{c}$
  and that $\tau$ generates $\Gal(K/L)$.
  Every root of unity in $K$ belongs to
  $L = K \cap k_0^{\ab}$,
  so $\mu(K) = \mu(L)$
  and we can define $W := W_K = W_L$.
  In particular, we have that $\tau$ fixes $\mu(K)$.
  Equations~\ref{eq:odd-Stickelberger} and
  \ref{eq:Stickelberger-integrality} give a unique
  $a \in \Z[\Gal(K/kF)]$ such that
  \begin{equation}
    \label{eq:Stark-integral-exponent}
    W\Theta_{k_0}(0) = (1 - c)a.
  \end{equation}
  Indeed, the coefficients of $W\Theta_{k_0}(0)$ on
  $\Gal(K/kF)$ and $c\Gal(K/kF)$ differ only by
  a sign.

  We first show that $c(\epsilon_F)$ and
  $\tau(\epsilon_F)$ differ by a root of unity.
  For this, we compare their absolute values at the
  archimedean places. For every
  $h \in \Gal(K/F)$, the relation $c w = \tau w$ and
  Lemma~\ref{lem:elliptic-conjugation-symmetry} give
  \begin{equation*}
    \log\verts{(c\epsilon_F)^h}_w = \log\verts{
    \epsilon_F^{\tau^{-1}c h c^{-1}}}_w =
    \log\verts{\epsilon_F^{\tau h}}_w.
  \end{equation*}
  As $h$ varies, we find that $c(\epsilon_F)$ and
  $\tau(\epsilon_F)$ have equal absolute values at every
  place above the complex place of $F$.
  At every other archimedean place, both units have absolute
  value $1$ by
  Definition~\ref{def:global-rank-one-Stark-unit}(ii).
  Hence $c(\epsilon_F)/\tau(\epsilon_F)$ is a global unit
  with absolute value $1$ at every archimedean place.
  Kronecker's theorem
  (see \cite[Chapter~I, Proposition~7.1]{neukirch})
  states that a global unit of absolute value $1$ at every
  archimedean place is a root of unity,
  hence there is a $\xi \in \mu(K)$ such that
  \begin{equation}
    \label{eq:finite-elliptic-unit-conjugation}
    c(\epsilon_F) = \xi\tau(\epsilon_F).
  \end{equation}

  We next find a $W$-th root of
  $\epsilon_F^{1 - \tau}$ in $K$; raising this root to
  $a$ will give $\epsilon$.
  Choose $\alpha$ with $\alpha^W = \epsilon_F$.
  Since $\epsilon_F$ is a Stark unit for $(K/F, S_F, w)$,
  the extension $K(\alpha)/F$ is abelian by
  Definition~\ref{def:global-rank-one-Stark-unit}(iii).
  Choose a lift
  $\widetilde\tau \in \Gal(K(\alpha)/F)$ of $\tau$.
  Set $\beta := \alpha/\widetilde\tau(\alpha)$.
  For $g \in \Gal(K(\alpha)/K)$, the quotient
  $g(\alpha)/\alpha$ belongs to $\mu_W$.
  Since $g$ commutes with $\widetilde\tau$ and $\tau$
  fixes $\mu_W$, the numerator and denominator of
  $\beta$ acquire the same factor under $g$.
  Hence $\beta \in K$, with
  $\beta^W = \epsilon_F^{1 - \tau}$.
  In particular, $\beta$ is a global unit.

  Because $\tau$ fixes $\xi$, Equation
  \ref{eq:finite-elliptic-unit-conjugation} implies
  $c(\epsilon_F^{1 - \tau}) = \epsilon_F^{\tau - 1}$.
  Thus $(\beta c(\beta))^W = 1$. At the embedding defining $c$,
  $\beta c(\beta)$ is positive, so $c(\beta) = \beta^{-1}$.
  Set $\epsilon := \beta^a$. Since $c$ commutes with
  $\Gal(K/kF)$, we also have $c(\epsilon) = \epsilon^{-1}$.
  Equations~\ref{eq:Stark-integral-exponent} and
  \ref{eq:finite-elliptic-unit-conjugation} give
  $u = \epsilon_F^{(1 - c)a} = \xi^{-a}\beta^{Wa} =
  \xi^{-a}\epsilon^W$, proving part~(i).

  \emph{Step 2. The Galois property.}
  We show that $K(\epsilon^{1/W})/k$ is Galois by
  constructing a $W$-th root of every conjugate of $\epsilon$
  inside this field.
  For each $h \in \Gal(K/kF)$,
  choose an integer $\kappa_h$ such that
  $h(\xi') = (\xi')^{\kappa_h}$ for every
  $\xi' \in \mu_W$, taking $\kappa_1 = 1$.
  The action $h(\xi') = (\xi')^{\kappa_h}$ and Equation
  \ref{eq:Stickelberger-action-integrality} imply that
  $(h - \kappa_h)\Theta_{k_0}(0) \in \Z[\Gal(K/k)]$.
  Comparing coefficients on $\Gal(K/kF)$ in Equation
  \ref{eq:Stark-integral-exponent} therefore shows that
  $(h - \kappa_h)a/W \in \Z[\Gal(K/kF)]$.

  Choose $y \in \C$ with $y^W = \epsilon$.
  For every $h \in \Gal(K/kF)$, the element
  $y^{\kappa_h}\beta^{(h - \kappa_h)a/W}$ is a
  $W$-th root of $h(\epsilon)$; its inverse is a
  $W$-th root of $ch(\epsilon)$.
  Since $\mu_W \subset K$, the field $K(y)$ contains
  all these roots. Hence $K(y)/k$ is Galois.
  Every lift $\widetilde h$ of $h$ to $\Gal(K(y)/k)$
  satisfies
  $\widetilde h(y) = \eta_h y^{\kappa_h}
  \beta^{(h - \kappa_h)a/W}$ for some $\eta_h \in \mu_W$.

  \emph{Step 3. Commutativity.}
  We first check that the automorphisms over $kF$ commute.
  Let $\widetilde h, \widetilde t \in \Gal(K(y)/k)$ have
  restrictions $h, t \in \Gal(K/kF)$ to $K$. Applying
  $\widetilde h$ to
  $\widetilde t(y) = \eta_t y^{\kappa_t}
  \beta^{(t - \kappa_t)a/W}$ gives
  \begin{equation*}
    \widetilde h\widetilde t(y) =
    \eta_t^{\kappa_h}\eta_h^{\kappa_t}
    y^{\kappa_h\kappa_t}
    \beta^{(ht - \kappa_h\kappa_t)a/W}.
  \end{equation*}
  The exponent $(ht - \kappa_h\kappa_t)a/W$ belongs to
  $\Z[\Gal(K/kF)]$: it is the sum of
  $\kappa_t(h - \kappa_h)a/W$ and $h(t - \kappa_t)a/W$.
  The expression for $\widetilde h\widetilde t(y)$ is
  symmetric in $h$ and $t$.
  The restrictions of $\widetilde h$ and $\widetilde t$ to
  $K$ commute because $\Gal(K/k)$ is abelian. The two composites
  also agree on $y$, so they agree on $K(y)$.
  Thus $\Gal(K(y)/kF)$ is abelian,
  including its subgroup $\Gal(K(y)/K)$.

  To finish the proof that $K(\epsilon^{1/W_K})/k$ is abelian,
  we check that complex conjugation commutes with $\Gal(K(y)/kF)$.
  Continue to write $c$ for complex conjugation on $K(y)$.
  Since $c(\epsilon) = \epsilon^{-1}$,
  $y$ has absolute value $1$ and $c(y) = y^{-1}$.
  The identities $c(\eta_h) = \eta_h^{-1}$ and
  $c(\beta) = \beta^{-1}$ give
  $c\widetilde h(y) = \widetilde h c(y)$.
  The restrictions of $c$ and $\widetilde h$ to $K$ also
  commute. Hence $c$ centralizes $\Gal(K(y)/kF)$.
  Since $\Gal(K(y)/kF)$ and $c$ generate $\Gal(K(y)/k)$,
  the extension $K(y)/k$ is abelian.

  \emph{Step 4. The regulator and absolute values.}
  We verify the two remaining conditions for $\epsilon$
  to be a Stark unit for $(K/k, S, w)$.
  Taking logarithms in $u = \xi^{-a}\epsilon^W$ and
  applying Proposition~\ref{prop:rank-one-regulator} verifies
  Definition~\ref{def:global-rank-one-Stark-unit}(i),
  since $W = W_K$.
  Because $K/k$ is abelian, $c$ fixes every place above
  the real place of $k$ above $v_0$. The identity $c(\epsilon) =
  \epsilon^{-1}$ therefore
  implies that $\lvert\epsilon\rvert_{w'} = 1$ at every
  place $w'$ above the real place of $k$ over $v_0$.
  At places above the other real places of $k_0$,
  we have $\lvert u\rvert_{w'} = 1$ by
  Proposition~\ref{prop:rank-one-regulator}.
  Since $u = \xi^{-a}\epsilon^W$ and $\xi$ is a root of
  unity, we also have $\lvert\epsilon\rvert_{w'} = 1$.
  Thus $\epsilon$ satisfies all three conditions
  of being a Stark unit in
  Definition~\ref{def:global-rank-one-Stark-unit}.
\end{proof}

\subsubsection{Descent from \texorpdfstring{$KF$}{KF}
to \texorpdfstring{$K$}{K}}
\label{subsec:Stark-descent}

We now remove the assumption that $K$ contains $F$.
We first apply Proposition~\ref{prop:relative-rank-one}
to $KF = k(LF)$, after checking its hypotheses.
We then descend the Stark unit from $KF$ to $K$
using the following result of Tate.

\begin{lemma}[{Tate~\cite[Chapter~IV, Proposition~3.5]{tate-stark-1984}}]
  \label{lem:Stark-descent}
  Let
  \begin{itemize}
    \item $K'/k$ be finite abelian;
    \item $k \subset K \subset K'$;
    \item $S$ be a finite set of at least three places of $k$,
      including every archimedean place and every finite
      place ramified in $K'/k$;
    \item $v$ be a fixed archimedean place of $k$
      that splits completely in $K'$;
    \item $v'$ be a place of $K'$ above $v$.
  \end{itemize}
  If $\epsilon' \in \calo_{K'}^\times$ is a Stark
  unit for $(K'/k, S, v, v')$, then there exists a Stark
  unit $\epsilon \in \calo_K^\times$ for
  $(K/k, S, v, \restr{v'}{K})$.
\end{lemma}

To apply Lemma~\ref{lem:Stark-descent} with $K' = KF$,
we must check that the chosen places include those ramified
in $KF/k$. The following lemma allows us to use the same
finite primes of $k_0$ as for $K/k_0$.

\begin{lemma}
  \label{lem:resolvent-ramification}
  Let $k/k_0$ be a non-Galois cubic extension of number
  fields with quadratic resolvent $F$, and let $K/k$ be
  a finite extension. A finite prime of $k_0$ is
  unramified in $KF/k_0$ if and only if it is unramified
  in $K/k_0$.
\end{lemma}

\begin{proof}
  Let $\mfp$ be a finite prime of $k_0$
  unramified in $K/k_0$.
  Fix an algebraic closure $\overline{(k_0)_\mfp}$,
  and let $(k_0)_\mfp^\ur$ be the maximal unramified
  extension of $(k_0)_\mfp$ in $\overline{(k_0)_\mfp}$.
  Since $k \subset K$, each completion $k_v$ with
  $v \mid \mfp$ is unramified over $(k_0)_\mfp$.
  Thus $(k_0)_\mfp\sigma(k) \subset (k_0)_\mfp^\ur$
  for every $k_0$-embedding
  $\sigma : k \longrightarrow \overline{(k_0)_\mfp}$.

  Let $E$ be the Galois closure of $k/k_0$.
  Every completion of $E$ above $\mfp$ is isomorphic to
  the compositum
  \[
    (k_0)_\mfp\bigl(\sigma(k) :
    \sigma \in \Hom_{k_0}(k, \overline{(k_0)_\mfp})\bigr)
    \subset (k_0)_\mfp^\ur.
  \]
  Hence $E/k_0$ is unramified at $\mfp$, and so is
  $F/k_0$ because $F \subset E$.
  For a place $w$ of $KF$ above $\mfp$,
  let $v := \restr{w}{K}$ and $z := \restr{w}{F}$;
  then $(KF)_w = K_vF_z$ inside $(KF)_w$.
  Both $K_v$ and $F_z$ are unramified over $(k_0)_\mfp$,
  so $K_vF_z/(k_0)_\mfp$ is unramified.
  Thus $KF/k_0$ is unramified at $\mfp$.
  Conversely, $KF/k_0$ being unramified at $\mfp$
  implies that its subextension $K/k_0$ is unramified too.
\end{proof}

Finally, we are ready to prove Theorem~\ref{thm:intro-relative}(i).
We construct a Stark unit for $KF/k$
and then descend to $K/k$.

\begin{proof}[Proof of Theorem~\ref{thm:intro-relative}(i)]
  First, we observe that $1$ is a Stark unit
  if $k$ is not almost totally real.
  Recall that $L = K \cap k_0^\ab$ and $K = kL$.
  If $v$ is real, every place of $k$ above
  $\restr{v}{k_0}$ splits completely in $K$, since
  $L/k_0$ is abelian.
  Thus, unless $k$ is almost totally real and $v$ is
  its complex place, at least two archimedean places
  of $k$ split completely in $K$.
  Since $\verts{S} \geq 3$,
  we have $L_S'(0, \chi) = 0$ for every $\chi$ by
  \cite[Chapter~I, Proposition~3.4]{tate-stark-1984}
  and so $1$ is a Stark unit.

  We may therefore assume that $k$ is almost totally real
  and $v$ is its complex place.
  Then $k_0$ is totally real, and the place of $F$ above
  $\restr{v}{k_0}$ is its unique complex place.
  Let $K' := KF = k(LF)$.
  The extension $LF/k_0$ is abelian and contains $F$.
  Lemma~\ref{lem:resolvent-ramification} shows that
  $S_{k_0}$ contains every finite prime ramified in $K'/k_0$.
  Thus $S$, $S_{k_0}$, and $S_F$ contain the finite places
  ramified in $K'/k$, $LF/k_0$, and $K'/F$, respectively,
  and $S_F$ contains at least two finite places.
  The sets $S$ and $S_F$ consist of the places above
  $S_{k_0}$, so they satisfy
  Equation~\ref{eq:local-factorization} for $K'/k$ by
  Corollary~\ref{cor:rationally-saturated-factorization}.

  By hypothesis, the integral rank-$1$ Stark conjecture
  holds for $(K'/F, S_F)$.
  We can therefore apply Proposition~\ref{prop:relative-rank-one}
  with the sets $S$, $S_{k_0}$, and $S_F$ to obtain a Stark
  unit for $(K'/k, S)$.
  We then apply Lemma~\ref{lem:Stark-descent} at the complex
  place of $k$ to obtain a Stark unit for $(K/k, S)$.
\end{proof}

We now treat the case when $k/k_0$ is Galois.
Here, the integral rank-$1$ Stark conjecture
follows from the order of vanishing at $s = 0$
being at least $2$.

\begin{proposition}
  \label{prop:galois-rank-one}
  Let $k/k_0$ be a Galois cubic extension of number fields,
  and let $K/k$ be finite abelian with
  $K \subset k \cdot k_0^\ab$.
  Let $S_{k_0}$ be a finite set of places of $k_0$
  containing the archimedean places, every prime ramified
  in $K/k_0$, and at least two finite places, and let
  $S$ consist of the places of $k$ above $S_{k_0}$.
  Then the integral rank-$1$ abelian Stark conjecture for
  $(K/k, S)$ holds at every archimedean place of $k$
  splitting completely in $K$, with Stark unit $1$.
\end{proposition}

\begin{proof}
  Let $v$ be an archimedean place of $k$ splitting
  completely in $K$; since $k \subset k_0^\ab$, the
  extension $K/k_0$ is abelian, so all three places of
  $k$ above $\restr{v}{k_0}$ split completely in $K$.
  Since $\verts{S} \geq 3$,
  we have $L_S'(0, \chi) = 0$ for every $\chi$ by
  \cite[Chapter~I, Proposition~3.4]{tate-stark-1984}
  and so $1$ is a Stark unit.
\end{proof}

\subsection{Specialization to \texorpdfstring{$k_0 = \Q$}{k0 = Q}}
\label{subsec:rank-one-specialization}

We now take $k_0 = \Q$,
in which case the quadratic resolvent $F$
of a complex cubic field $k/\Q$
is an imaginary quadratic field.
We will use the cyclotomic and elliptic units
of Stark~\cite{stark-1980} to
apply our relative argument over $\Q$ and prove
Theorem~\ref{thm:intro-rank-one}.
We recall the classical theorems,
treat the case of one ramified rational prime,
and then prove
Theorem~\ref{thm:intro-rank-one}.

\subsubsection{The classical Stark theorems}
\label{subsec:classical-Stark}

The next two results concern arbitrary abelian extensions
of $\Q$ and of an imaginary quadratic field, respectively.
We will apply them to $L^+/\Q$ and $KF/F$.

Let $L^+/\Q$ be finite, totally real, and abelian;
note that $W_{L^+} = 2$.
Let $S_\Q$ be a finite set of places of $\Q$ including
$\infty$, every finite prime ramified in $L^+/\Q$,
and at least two finite primes.
Fix a real place $w_{L^+}$ of $L^+$.
Stark~\cite{stark-1975, stark-1980} proves the
integral rank-$1$ abelian Stark conjecture over $\Q$
using Hurwitz zeta functions and cyclotomic units.

\begin{theorem}[{Stark~\cite[Theorem~2]{stark-1980}}]
  \label{thm:rational-Stark}
  There is a cyclotomic unit $\epsilon_\Q \in \calo_{L^+}^\times$
  such that
  \[
    L_{S_\Q}'(0, \chi) = -\frac{1}{2}
    \sum_{\alpha \in \Gal(L^+/\Q)}\chi(\alpha)
    \log\verts{\epsilon_\Q^\alpha}_{w_{L^+}}
  \]
  for every character $\chi$ of $\Gal(L^+/\Q)$.
  Furthermore, $L^+(\epsilon_\Q^{1/2})/\Q$ is abelian.
\end{theorem}

Over an imaginary quadratic field,
we replace cyclotomic units by elliptic units.
Let $F$ be imaginary quadratic and let $K/F$ be finite
abelian. Let $S_F$ be a finite set of places of $F$
including its complex place, every finite place ramified
in $K/F$, and at least one finite place.
Fix a place $w$ of $K$ above the complex place of $F$.
Stark~\cite[Theorem~2]{stark-1980} proves the
integral rank-$1$ abelian Stark conjecture over $F$
using CM theory and elliptic units.

\begin{theorem}[{Stark~\cite[Theorem~2]{stark-1980}}]
  \label{thm:imaginary-quadratic-Stark}
  There is an $S_F$-unit $\epsilon_F \in K^\times$ such that
  \[
    L_{S_F}'(0, \psi) = -\frac{1}{W_K}
    \sum_{h \in \Gal(K/F)}\psi(h)
    \log\verts{\epsilon_F^h}_w
  \]
  for every character $\psi$ of $\Gal(K/F)$.
  Furthermore, $K(\epsilon_F^{1/W_K})/F$ is abelian.
  Moreover, if $S_F$ contains at least two finite places,
  then $\epsilon_F \in \calo_K^\times$.
\end{theorem}

When $S_F$ has only one finite place $\mfq$, we may choose
$\epsilon_F$ so that $\ord_{\mfP}(\epsilon_F)$
is independent of the prime $\mfP$ of $K$ above
$\mfq$, as in Tate's
\cite[Chapter~IV, Sections~2--3]{tate-stark-1984}
formulation of the integral rank-$1$ Stark conjecture.

\subsubsection{One ramified rational prime}
\label{subsec:rank-one-one-prime}

The construction in Proposition~\ref{prop:relative-rank-one}
starts with a global unit $\epsilon_F$.
To prove Theorem~\ref{thm:intro-rank-one}, we must also
treat cyclotomic extensions ramified over $\Q$ at only
one rational prime. Then $S_F$ may contain only one
finite place, and $\epsilon_F$ need not be a global unit.
Recall that $c$ is complex conjugation at the real
place of $k$, and $\tau$ generates $\Gal(K/L)$.
We show that the finite valuations cancel in
$\epsilon_F^{1 - \tau}$ and $c(\epsilon_F)/\tau(\epsilon_F)$,
so we can still construct a global Stark unit for $K/k$.

\begin{lemma}
  \label{lem:rank-one-single-prime}
  Let $k$ be a complex cubic field with quadratic resolvent
  $F$, and let $K/k$ be a finite cyclotomic
  extension that contains $F$.
  If $K/\Q$ is ramified at only one finite prime $p$,
  $S$ is the set of places of $k$ above $\infty$ and $p$,
  and $w$ is any place of $K$ above the complex place of $k$,
  then there is a Stark unit $\epsilon \in \calo_K^\times$
  for $(K/k, S, w)$.
\end{lemma}

\begin{proof}
  Recall that $L = K \cap \Q^{\ab}$.
  Choose a generator $\tau$ of $\gal(K/L)$
  so that $c w = \tau w$.
  Since $F/\Q$ is quadratic and unramified outside $p$,
  the prime $p$ ramifies in $F$.
  Let $\mfq$ be its unique prime in $F$,
  let $S_\Q := \set{\infty, p}$,
  and let $S_F$ consist of
  $\mfq$ and the complex place of $F$.

  Choose $\epsilon_F$ as in
  Theorem~\ref{thm:imaginary-quadratic-Stark}, with equal
  valuations at the primes above $\mfq$.
  The automorphisms $c$ and $\tau$ permute these primes.
  Thus $\epsilon_F^{1 - c}$, $\epsilon_F^{1 - \tau}$, and
  $c(\epsilon_F)/\tau(\epsilon_F)$ are global units.

  Since $S_F$ is invariant under conjugation,
  we can apply the argument of
  Lemma~\ref{lem:elliptic-conjugation-symmetry}
  to the logarithms of $\epsilon_F$.
  Using $c w = \tau w$, we find that
  $c(\epsilon_F)/\tau(\epsilon_F)$ has absolute value $1$
  at every archimedean place.
  Since $c(\epsilon_F)/\tau(\epsilon_F)$ is a global unit,
  $\xi := c(\epsilon_F)/\tau(\epsilon_F)$ belongs to $\mu(K)$.
  Thus $c(\epsilon_F) = \xi\tau(\epsilon_F)$.

  We now find a $W$-th root of $\epsilon_F^{1 - \tau}$ in
  $K$, where $W := W_K = W_L$.
  Choose $\alpha^W = \epsilon_F$.
  The extension $K(\alpha)/F$ is abelian.
  Choose a lift $\widetilde\tau$ of $\tau$ and set
  $\beta := \alpha/\widetilde\tau(\alpha)$.
  The quotient $\beta$ belongs to $K$ because $\tau$
  fixes $\mu(K)$ and $\widetilde\tau$ commutes with
  every automorphism of $K(\alpha)/K$.
  As $\beta^W = \epsilon_F^{1 - \tau}$, the element $\beta$ is
  a global unit. The identity
  $c(\epsilon_F) = \xi\tau(\epsilon_F)$ gives
  $(\beta c(\beta))^W = 1$. Positivity at the embedding
  defining $c$ therefore gives $c(\beta) = \beta^{-1}$.

  Write $W\Theta_\Q(0) = (1 - c)a$ with $a \in \Z[\Gal(K/kF)]$,
  using Equations~\ref{eq:Stickelberger-integrality} and
  \ref{eq:odd-Stickelberger}, and set $\epsilon := \beta^a$.
  Then $\epsilon_F^{W\Theta_\Q(0)} = \xi^{-a}\epsilon^W$.
  We compute the logarithms as in the proofs of
  Lemma~\ref{lem:elliptic-regulator-comparison} and
  Proposition~\ref{prop:rank-one-regulator} to obtain
  $\calr_w(\epsilon_F^{W\Theta_\Q(0)}, \chi)
  = -W^2 L_S'(0, \chi)$
  for every character $\chi$ of $\Gal(K/k)$.
  Since $\xi$ is a root of unity and
  $\epsilon_F^{W\Theta_\Q(0)} = \xi^{-a}\epsilon^W$,
  we conclude that
  $\calr_w(\epsilon, \chi) = -WL_S'(0, \chi)$.
  Also, $c(\epsilon) = \epsilon^{-1}$, so $\epsilon$ has
  absolute value $1$ at every place above the real place of $k$.

  It remains to prove that $K(\epsilon^{1/W})/k$ is abelian.
  For $h \in \Gal(K/kF)$, choose $\kappa_h \in \Z$ with
  $h(\zeta) = \zeta^{\kappa_h}$ for every $\zeta \in \mu(K)$.
  Equation~\ref{eq:Stickelberger-action-integrality} gives
  \[
    (h - \kappa_h)a / W \in \Z[\Gal(K/kF)].
  \]
  The calculation of lifts
  in the proof of Proposition~\ref{prop:relative-rank-one},
  with $\epsilon = \beta^a$ and $c(\beta) = \beta^{-1}$,
  shows that $K(\epsilon^{1/W})/k$ is abelian.
  Thus $\epsilon$ satisfies all three conditions in
  Definition~\ref{def:global-rank-one-Stark-unit}.
\end{proof}

\subsubsection{Cyclotomic extensions of complex cubic fields}
\label{subsec:cyclotomic-extensions}

We now prove Theorem~\ref{thm:intro-rank-one}.
We first construct a Stark unit for the places above
$\infty$ and the rational primes ramified in $K/\Q$.
We then include the remaining finite places in $S$ by
modifying the unit with the corresponding Euler factors.

\begin{proof}[Proof of Theorem~\ref{thm:intro-rank-one}]
  Let $\Sigma_\Q$ consist of $\infty$ and every rational
  prime ramified in $K/\Q$. Let $S_0$ and $S_F$ consist
  of the places above $\Sigma_\Q$ in $k$ and its quadratic
  resolvent $F$, respectively. By hypothesis, $S_0 \subseteq S$.
  Define $K' := KF$ and choose a place $w'$ of $K'$ above $w$.
  Lemma~\ref{lem:resolvent-ramification} shows that
  $\Sigma_\Q$ includes every finite prime ramified in $K'/\Q$.
  Also, $\Sigma_\Q$ includes at least one finite prime,
  since $k \neq \Q$ has a ramified rational prime
  by Minkowski's theorem.

  Suppose first that $\Sigma_\Q$ contains at least two
  finite primes. By Theorem~\ref{thm:imaginary-quadratic-Stark},
  there is a Stark unit $\epsilon_F$ for $(K'/F, S_F, w')$.
  With $L' := K' \cap \Q^{\ab}$, we have $K' = kL'$ and
  $F \subset L'$. The sets $S_0$, $\Sigma_\Q$, and $S_F$
  satisfy Equation~\ref{eq:local-factorization} by
  Corollary~\ref{cor:rationally-saturated-factorization}.
  We can therefore apply
  Proposition~\ref{prop:relative-rank-one}
  to construct a Stark unit for $(K'/k, S_0, w')$.
  If $\Sigma_\Q$ contains only one finite prime, we obtain
  a Stark unit for $(K'/k, S_0, w')$ from
  Lemma~\ref{lem:rank-one-single-prime} instead.
  In either case, we descend this unit using
  Lemma~\ref{lem:Stark-descent} to obtain a Stark unit
  $\epsilon_0 \in \calo_K^\times$ for $(K/k, S_0, w)$.

  We now enlarge $S_0$ to $S$ as in Tate
  \cite[Chapter~IV, Proposition~3.4]{tate-stark-1984},
  keeping track of the resulting unit. Define
  \[
    a := \prod_{\mfp \in S \setminus S_0}
      (1 - [\mfp, K/k]^{-1}) \in \Z[\Gal(K/k)]
  \]
  and $\epsilon := \epsilon_0^a$.
  Every prime in $S \setminus S_0$ is finite and unramified
  in $K/k$. Every $L_{S_0}(s, \chi)$ vanishes at $0$,
  since the complex place of $k$ splits completely.
  Removing the Euler factors at $S \setminus S_0$ therefore
  multiplies $L_{S_0}'(0, \chi)$ by $\chi^{-1}(a)$.
  We use Equation~\ref{eq:character-logarithm-action}
  to check that $\epsilon = \epsilon_0^a$ has the required
  regulator:
  \begin{align*}
    \calr_w(\epsilon, \chi)
      &= \chi^{-1}(a)\calr_w(\epsilon_0, \chi) \\
      &= -W_K L_S'(0, \chi).
  \end{align*}
  Since every conjugate of $\epsilon_0$ has absolute value
  $1$ at every place of $K$ above the real place of $k$,
  the same is true of $\epsilon$.

  Finally, the abelian extension
  $K(\epsilon_0^{1/W_K})/k$ contains a $W_K$-th root of
  $\epsilon_0^a$, obtained by multiplying conjugates of
  a $W_K$-th root of $\epsilon_0$. Thus
  $K(\epsilon^{1/W_K})/k$ is abelian.
  Together with the regulator and absolute values already
  proved, this gives all the assertions of the theorem.
\end{proof}

\subsection{Explicit class field theory over complex cubic fields}
\label{subsec:rank-one-explicit-cft}

We now ask when the Stark unit of
Theorem~\ref{thm:intro-rank-one} generates $K/k$.
For totally imaginary cyclic extensions, we answer this
question using the splitting of the finite places in $S$.
This realizes an instance of the approach of
Stark \cite[p.~63]{stark-1975} and
Tate \cite[p.~95]{tate-stark-1984} to explicit class field theory.

\begin{corollary}
  \label{cor:rank-one-field-generation}
  Let $K/k$ be a totally imaginary cyclic cyclotomic extension
  of a complex cubic field, and let $S$ and $w$ be as in
  Theorem~\ref{thm:intro-rank-one}.
  If no finite place in $S$ splits completely in $K/k$,
  then $K = k(\epsilon)$,
  where $\epsilon \in \calo_K^\times$
  is the unit constructed in the proof of
  Theorem~\ref{thm:intro-rank-one}.
\end{corollary}

\begin{proof}
  Choose a faithful character $\chi$ of $\Gal(K/k)$.
  Since the complex place of $k$ is the only place of $S$
  that splits completely in $K/k$, we have
  \mbox{$L_S'(0, \chi) \neq 0$} by Tate
  \cite[Chapter~I, Proposition~3.4]{tate-stark-1984}.
  If $\alpha \in \Gal(K/k)$ fixes $\epsilon$, we can
  reindex the sum in Theorem~\ref{thm:intro-rank-one}
  by $\alpha$ to obtain
  $L_S'(0, \chi) = \chi(\alpha)^{-1}L_S'(0, \chi)$.
  Thus $\chi(\alpha) = 1$, and so $\alpha = 1$ because
  $\chi$ is faithful. We conclude that $K = k(\epsilon)$.
\end{proof}

\begin{remark}
  The order-of-vanishing formula in Rubin
  \cite[Section~2.1, Equation~6]{rubin-1996}
  explains the hypotheses in
  Corollary~\ref{cor:rank-one-field-generation}:
  if the real place or a finite place in $S$ splits
  completely in $K/k$, then $L_S'(0, \chi) = 0$
  for every faithful character $\chi$.
  Under the hypotheses of the corollary, we also have
  $K = k(u)$ when $F \subset K$, where $u$ is the unit of
  Proposition~\ref{prop:rank-one-regulator}. Indeed,
  $u = \zeta\epsilon^{W_L}$ for some $\zeta \in \mu(K)$,
  so $\calr_w(u, \chi) = W_L \calr_w(\epsilon, \chi)$
  for every character $\chi$.
\end{remark}


\section{The integral rank-\texorpdfstring{$r$}{r}
Stark conjecture}
\label{sec:rank-r}

In this section,
we prove Theorems~\ref{thm:intro-relative}(ii) and
\ref{thm:intro-rank-two}.
We first compare rational Stark elements over $k$, $k_0$,
and the quadratic resolvent $F$ to prove
Theorem~\ref{thm:intro-relative}(ii) in arbitrary rank.
We then treat the case when $k/k_0$ is Galois in
Theorem~\ref{thm:Galois-Rubin-Stark},
and thereafter assume that $k/k_0$ is not Galois.
When $k/k_0$ is an almost totally real cubic extension
and $r = 2$,
we construct explicit units
and obtain stronger integrality results
than are available to us in the general rank $r$ case.
We conclude by proving Theorem~\ref{thm:intro-rank-two}
and giving an application to explicit class field theory.

\subsection{Stark elements over cubic extensions}
\label{subsec:cubic-Rubin-Stark}

We will use the decomposition of Artin $L$-functions to
construct a rational Stark element over $k$ from those
over $k_0$ and $F$. We begin by recalling the definition
in arbitrary rank.

Let $K/k$ be a finite abelian extension, and let $S$
be a finite set of places of $k$ containing the
archimedean places and every finite place ramified in
$K/k$.
Choose $r \geq 0$ distinct archimedean places of $k$
splitting completely in $K$, with ordered places
$w_1, \ldots, w_r$ of $K$ above them, and assume
$\lvert S\rvert > r + 1$.
The rank-$r$ regulator is a determinant of logarithms
of $S$-units at $w_1, \ldots, w_r$.
We work in $U_S := \Q\otimes_\Z\calo_{K, S}^\times$.
For an $S$-unit $a$, define
\[
  \ell_w(a) := -\sum_{\sigma \in \Gal(K/k)}
  \log\verts{a^\sigma}_w\sigma^{-1},
\]
and extend $\ell_w$ $\Q$-linearly to $U_S$.
The rank-$r$ regulator is given by
\[
  \calr_{w_1, \ldots, w_r}(a_1\wedge\cdots\wedge a_r)
  := \det\paren{\ell_{w_i}(a_j)}_{i, j}.
\]
For $r = 0$, the exterior power is the coefficient ring
and the empty determinant is $1$.

We compare the regulator with the leading terms of the
$L$-functions that vanish to order exactly $r$.
For a character $\chi$ of $\Gal(K/k)$, the projection to
its component is
$e_\chi := [K : k]^{-1}\sum_{\sigma \in \Gal(K/k)}
\chi(\sigma)\sigma^{-1} \in \Q(\chi)[\Gal(K/k)]$.
We therefore define
\begin{align*}
  e_{S, r}
  &:= \sum_{\substack{
    \chi : \Gal(K/k) \longrightarrow \C^\times \\
    \ord_{s = 0}L_S(s, \chi) = r}}e_{\chi^{-1}}, \\
  \Theta_S^{(r)}(0)
  &:= \sum_{\substack{
    \chi : \Gal(K/k) \longrightarrow \C^\times \\
    \ord_{s = 0}L_S(s, \chi) = r}}
    \frac{L_S^{(r)}(0, \chi)}{r!}e_{\chi^{-1}}.
\end{align*}
The element $e_{S, r}$ has rational coefficients:
by Tate \cite[Chapter~I, Proposition~3.4]{tate-stark-1984},
the order of vanishing is preserved by Galois conjugation
of the character values.
The element $\Theta_S^{(r)}(0)$ has real coefficients,
since complex conjugation interchanges the summands for
$\chi$ and $\overline\chi$.

The rational Stark conjecture predicts an element with
regulator $\Theta_S^{(r)}(0)$; see Stark
\cite[Section~1]{stark-1975} and Tate
\cite[Chapter~I, Sections~5 and~9]{tate-stark-1984}.
We require its components at characters whose $L$-functions
vanish to order greater than $r$ to be zero.

\begin{definition}[Rational rank-$r$ Stark element]
  \label{def:rational-rank-two-element}
  A rational rank-$r$ Stark element for $(K/k, S)$
  and the ordered places $(w_1, \ldots, w_r)$ is an element
  $\epsilon_{S, \Q}^{(r)} \in
  e_{S, r}\bigwedge_{\Q[\Gal(K/k)]}^r U_S$ such that
  \[
    \calr_{w_1, \ldots, w_r}(\epsilon_{S, \Q}^{(r)})
    = \Theta_S^{(r)}(0).
  \]
\end{definition}

The regulator is injective on
$e_{S, r}\bigwedge_{\Q[\Gal(K/k)]}^r U_S$ by Rubin
\cite[Lemma~2.7(iii)]{rubin-1996}.
Thus $\epsilon_{S, \Q}^{(r)}$ is unique when it exists.

To formulate integrality, we evaluate exterior products
against integral linear maps.
For a module $M$ over a commutative ring $R$, define
its $r$-th exterior bidual by
\[
  \bigcap_R^r M
  := \Hom_R\paren{\bigwedge_R^r\Hom_R(M, R), R}.
\]
The canonical map $\bigwedge_R^r M \to \bigcap_R^r M$
sends $a_1\wedge\cdots\wedge a_r$ to the functional
whose value at $\phi_1\wedge\cdots\wedge\phi_r$ is
$\det(\phi_i(a_j))_{i, j}$.
For $r = 0$, both modules are $R$.

If $R$ is a finite product of fields and $M$ is finitely
generated, this map is an isomorphism.
For a prime $\ell$, a direct factor $R$ of
$\Z_\ell[\Gal(K/k)]$, and an $R$-module $M$ finite free
over $\Z_\ell$, we therefore identify $\bigcap_R^r M$
with a lattice in
$\bigwedge_{\Q_\ell\otimes_{\Z_\ell}R}^r
(\Q_\ell\otimes_{\Z_\ell}M)$.

We will use two standard properties of exterior biduals.
For finitely generated modules over a Noetherian ring,
they commute with flat extension of coefficients, by
\cite[Tag~087R]{stacks-project} and Bourbaki
\cite[Chapter~III, \S7, no.~5, Proposition~8]
{bourbaki-algebra}.
For $R$-modules $U$ and $V$ and integers $a, b \geq 0$,
we also have a natural exterior product
\[
  \bigcap_R^a U \otimes_R\bigcap_R^b V
  \longrightarrow \bigcap_R^{a + b}(U \oplus V).
\]
We obtain this map by evaluating the two functionals on
the component of bidegree $(a, b)$ in the direct-sum
decomposition of
$\bigwedge_R^{a + b}\Hom_R(U \oplus V, R)$; see Bourbaki
\cite[Chapter~III, \S7, no.~7, Proposition~10]
{bourbaki-algebra}.
Over the rational group algebra, it agrees with the
usual exterior product.

We now return to the non-Galois cubic extension $k/k_0$
and the finite abelian extension $K/k$ in
Theorem~\ref{thm:intro-relative}(ii).
Recall that $K \subset k \cdot k_0^\ab$, that $F$ is the
quadratic resolvent, and that $L = K \cap k_0^\ab$.
As in Section~\ref{subsec:Stark-hypotheses}, we have
$K = kL$ and $\mu(L) = \mu(K)$.
The set $S_{k_0}$ contains the archimedean places,
every prime ramified in $K/k_0$, and at least two finite
places; $S$ and $S_F$ consist of the places above
$S_{k_0}$ in $k$ and $F$, respectively.
By Lemma~\ref{lem:resolvent-ramification}, $S_{k_0}$ also
contains every prime ramified in $KF/k_0$.
For the comparison below, let $r_{k_0}$ and $r$ count
all archimedean places splitting completely in
$L/k_0$ and $K/k$, respectively.

To use the factorization over $k_0$ and $F$, we pass from
$K = kL$ to $KF = k(LF)$.
Recall from Section~\ref{subsec:Stark-hypotheses} that
$KF/F$ is abelian and $KF\cap k_0^\ab = LF$.
We compare the rational Stark elements for $L/k_0$,
$KF/F$, and $K/k$.
In Corollary~\ref{cor:rationally-saturated-factorization},
applied to $KF/k$, the factor over $F$ uses the nontrivial
characters of $\Gal(KF/LF)$.
We select those characters by choosing a generator
$\tau$ of $\Gal(KF/LF)$ and setting
$e := (2 - \tau - \tau^2)/3$.
To prove integrality away from $3W_K$, fix a prime
$\ell \nmid 3W_K$ and write $R := \Z_\ell[\Gal(K/k)]$.
We will compare the modules
\begin{itemize}
  \item $U_L := \Z_\ell\otimes_\Z\calo_{L, S_{k_0}}^\times$,
  \item $U_K := \Z_\ell\otimes_\Z\calo_{K, S_{k_0}}^\times$, and
  \item $U_{KF} := \Z_\ell\otimes_\Z\calo_{KF, S_{k_0}}^\times$.
\end{itemize}
Restriction identifies $\Gal(K/k)$ with $\Gal(L/k_0)$.
Since $\tau$ fixes every root of unity in $KF$,
$e$ annihilates the torsion of $U_{KF}$.
Hence $eU_{KF}$ is free over $\Z_\ell$.

We next determine the rank over $F$ by comparing the
archimedean places. Let $r_F$ count the archimedean places
of $F$ splitting completely in $KF$.
An archimedean place $v$ of $k_0$ contributes to
$(r_{k_0}, r_F, r)$ as follows:
\[
  \begin{array}{c|c|c|c}
    v & k \otimes_{k_0}k_{0, v}
    & \text{places above }v\text{ in }L
    & (r_{k_0}, r_F, r) \\
    \hline
    \text{real} & \R^3 & \text{real} & (1, 2, 3) \\
    \text{complex} & \C^3 & \text{complex}
    & (1, 2, 3) \\
    \text{real} & \R\times\C & \text{real}
    & (1, 1, 2) \\
    \text{real} & \R\times\C & \text{complex}
    & (0, 1, 1) \\
    \text{real} & \R^3 & \text{complex}
    & (0, 0, 0).
  \end{array}
\]
Let $n_3$, $n_2$, and $n_1$ count the places contributing
$(1, 2, 3)$, $(1, 1, 2)$, and $(0, 1, 1)$,
respectively. Thus $r_{k_0} = n_3 + n_2$,
$r_F = 2n_3 + n_2 + n_1$, and
$r = 3n_3 + 2n_2 + n_1$.
In particular, $r = r_{k_0} + r_F$.
In each regulator, we list the places counted by $n_3$,
then those counted by $n_2$, then those counted by $n_1$.

Write $\eta_0$ and $\eta_F$ for the rational Stark
elements for $(L/k_0, S_{k_0})$ and $(KF/F, S_F)$,
of ranks $r_{k_0}$ and $r_F$, respectively, when they exist.
We construct the Stark element for $K/k$ from $\eta_0$
and $e\eta_F$, using $r_{k_0} + r_F = r$.

\begin{lemma}
  \label{lem:cubic-exterior-comparison}
  If $\eta_0$ and $\eta_F$ exist and
  \[
    \eta_0 \in \bigcap_R^{r_{k_0}} U_L,
    \qquad
    e\eta_F \in
    \bigcap_{e\Z_\ell[\Gal(KF/F)]}^{r_F} eU_{KF},
  \]
  then the rational rank-$r$ Stark element
  $\eta_K$ for $K/k$ exists and
  belongs to $\bigcap_R^r U_K$.
\end{lemma}

\begin{proof}
  We carry out the algebraic constructions over $\Q$
  after adjoining a primitive cube root of unity,
  and use their $\Z_\ell$-realizations to check
  integrality.
  We compare the extended regulators with the Artin
  factorization, then use injectivity to descend
  the element.

  \emph{1. The case $F \cap L = k_0$.}
  We have $\Gal(KF/k_0) \cong
  S_3\times\Gal(L/k_0)$.
  Let $s$ be the involution fixing $K$ and set
  $A := eR[\abrac{\tau}]$ and
  $\delta := \tau - \tau^2$.
  Thus $s\tau s = \tau^{-1}$ and $\delta^2 = -3e$.
  Since $3$ is invertible in $R$, we have
  \[
    U_K = U_L \oplus V,
    \qquad
    V := (eU_{KF})^s
    = \Ker(\rmn_{K/L} : U_K \to U_L).
  \]
  The map $A \otimes_R V \to eU_{KF}$,
  $a \otimes v \mapsto av$, is an isomorphism.
  Indeed, every $x \in eU_{KF}$ has a unique expression
  $x = u + \tau v$ with $u, v \in V$, given by
  $v = \delta^{-1}(x - sx)$ and $u = x - \tau v$.
  Under this isomorphism, $s$ acts on $A$ by
  $\tau \longmapsto \tau^{-1}$ and fixes $V$.
  Thus $s$ is the coefficient involution, and $A^s = R$.
  By flat extension of coefficients, we can therefore form
  $\eta_0\wedge e\eta_F$ in
  $A \otimes_R\bigcap_R^r U_K$.

  We compute the regulator one archimedean place at a time.
  Write $\alpha := \ell_{\restr{w}{L}}$ for $L/k_0$.
  The rows below follow by summing the $\tau$-conjugate
  logarithms, whose sum vanishes on $V$, and using the
  normalized absolute values.
  \begin{itemize}
    \item At a place counted by $n_3$, choose the ordered
      places $(w, sw)$ over $F$.
      Let $x$ and $y$ be the restrictions to $V$ of
      $\ell_{\restr{w}{K}}$ and
      $\ell_{\restr{\tau w}{K}}$ for $K/k$.
      The rows for $K/k$ on $U_L \oplus V$ are
      $(\alpha, x)$, $(\alpha, y)$, and
      $(\alpha, -x - y)$.
      The two rows for $KF/F$, after projection by $e$,
      are $(1 - \tau^2)x + \delta y$ and
      $(1 - \tau)x - \delta y$.
      Their determinants relative to $(\alpha, x, y)$
      and $(x, y)$ are $3$ and $-3\delta$.
    \item At a place counted by $n_2$, choose $w$ fixed
      by $s$, so that $\restr{w}{K}$ is real.
      Let $x$ be the restriction of
      $\ell_{\restr{w}{K}}$ to $V$.
      The rows for $K/k$ are $(\alpha, x)$ and
      $(2\alpha, -x)$, with determinant $-3$,
      while the row for $KF/F$ is $3x$.
    \item At a place counted by $n_1$, choose $w$ fixed
      by $sc_i$, where $c_i \in \Gal(L/k_0)$ is complex
      conjugation at the corresponding real place.
      On every component where $\eta_0$ is nonzero,
      $c_i$ acts as $-1$: otherwise that real place
      contributes an additional zero, by Tate
      \cite[Chapter~I, Proposition~3.4]
      {tate-stark-1984}.
      The trivial component also vanishes, since
      $\lvert S_{k_0}\rvert - 1 > r_{k_0}$.
      Let $x$ be the restriction of
      $\ell_{\restr{\tau w}{K}}$ to $V$.
      On those components, the rows for $K/k$ and
      $KF/F$ are $x$ and $\delta x$.
  \end{itemize}
  We multiply the determinant ratios and move all $U_L$
  rows before the $V$ rows. We obtain the required product
  of regulators by setting
  \[
    \eta_K := (-1)^{n_3 + n_2(n_2 + 1)/2}
    \delta^{n_3 + n_1}\eta_0\wedge e\eta_F.
  \]
  The correction factor is a unit of $A$.

  \emph{2. The case $F \subset L$.}
  Now $KF = K$.
  Let $J := \Gal(L/F)$, identified by restriction with
  the subgroup of $\Gal(K/k)$ fixing $F$, and choose
  $t \in \Gal(K/k)\setminus J$.
  A place counted by $n_2$ would be complex in $F$ and
  real in $L$, so $n_2 = 0$.
  Set $\calo := \Z_\ell[X]/(X^2 + X + 1)$, write
  $\zeta$ for the image of $X$, and set
  $\delta := \zeta - \zeta^{-1}$, so $\delta^2 = -3$.
  The algebra $\calo$ is quadratic \'etale over $\Z_\ell$.
  We isolate one nontrivial eigenspace of $\tau$ by setting
  \[
    P_\zeta :=
    \frac{1 + \zeta^{-1}\tau + \zeta^{-2}\tau^2}{3},
    \qquad
    W_\zeta := P_\zeta
    (\calo \otimes_{\Z_\ell}U_K).
  \]
  We identify $P_\zeta\calo[\Gal(K/F)]$ with $\calo[J]$
  by sending $\tau$ to $\zeta$.
  Since $t\tau t^{-1} = \tau^{-1}$, the map
  $h \otimes u \mapsto hu$ identifies
  \[
    \calo[\Gal(K/k)] \otimes_{\calo[J]}W_\zeta
    \cong \calo \otimes_{\Z_\ell}eU_K
    = W_\zeta \oplus tW_\zeta.
  \]
  By flat extension of coefficients, $P_\zeta\eta_F$
  therefore determines an element $\xi_\zeta$ of
  $\bigcap_{\calo[\Gal(K/k)]}^{r_F}
  (\calo \otimes_{\Z_\ell}eU_K)$.
  Since $U_K = U_L \oplus eU_K$, we can form
  $\eta_0\wedge\xi_\zeta$ in
  $\calo \otimes_{\Z_\ell}\bigcap_R^r U_K$.

  We again compute the regulator at each archimedean place.
  Let $\alpha := \ell_{\restr{w}{L}}$ for $L/k_0$,
  and let $x$ and $y$ be the restrictions to $W_\zeta$
  of $\ell_w$ and $\ell_{tw}$, with the sums taken over $J$.
  \begin{itemize}
    \item At a place counted by $n_3$, choose the ordered
      places $(w, tw)$ over $F$.
      The rows for $KF/F$, after projection by $P_\zeta$,
      are $3x$ and $3y$.
      The three rows for $K/k$ are obtained from
      $(\alpha, x, y)$ by the matrix
      \[
        \begin{pmatrix}
          1 & 1 & t \\
          1 & \zeta^{-1} & t\zeta \\
          1 & \zeta & t\zeta^{-1}
        \end{pmatrix},
      \]
      whose determinant is $3t\delta$.
    \item At a place counted by $n_1$, choose $w$ fixed
      by complex conjugation
      $c_i \in \Gal(K/k)\setminus J$.
      As in part~1, $c_i$ acts as $-1$ on every
      component where $\eta_0$ is nonzero.
      On those components, the row for $K/k$ is
      $(\zeta^{-1} + \zeta c_i)x = -\delta x$,
      while the row for $KF/F$ is $3x$.
  \end{itemize}
  Since $\delta^2 = -3$, the inverse determinant ratios
  are $-t^{-1}\delta$ and $\delta$, respectively.
  We therefore obtain the required product of regulators
  by setting
  \[
    \eta_K := (-1)^{n_3}t^{-n_3}
    \delta^{n_3 + n_1}\eta_0\wedge\xi_\zeta.
  \]
  Again, the correction factor is a unit.

  \emph{3. The Stark element and descent.}
  We verify the Stark conditions for each embedding
  of the coefficient field into $\C$.
  Let $\chi$ be a character of $\Gal(K/k)$, let
  $\chi_{k_0}$ correspond to $\chi$ under restriction,
  and set $\chi_F(g) := \chi_{k_0}(\restr{g}{L})$
  for $g \in \Gal(KF/F)$.
  For either nontrivial character $\omega_F$ of
  $\Gal(kF/F)$, inflated to $\Gal(KF/F)$,
  we apply
  Corollary~\ref{cor:rationally-saturated-factorization}
  to the inflation of $\chi$ to $\Gal(KF/k)$.
  The two factors vanish to orders at least $r_{k_0}$
  and $r_F$. Since $r_{k_0} + r_F = r$, we compare
  coefficients of $s^r$ to obtain
  \[
    \frac{L_S^{(r)}(0, \chi)}{r!}
    = \frac{L_{S_{k_0}}^{(r_{k_0})}
    (0, \chi_{k_0})}{r_{k_0}!}
    \frac{L_{S_F}^{(r_F)}
    (0, \chi_F\omega_F)}{r_F!}.
  \]
  The determinant comparisons in parts~1 and~2 show that
  $\chi^{-1}(\calr_{w_1, \ldots, w_r}(\eta_K))
  = L_S^{(r)}(0, \chi)/r!$.
  If $\ord_{s = 0}L_S(s, \chi) > r$, one factor has
  order greater than $r_{k_0}$ or $r_F$, respectively.
  The corresponding component of $\eta_0$ or $e\eta_F$
  is zero, so the $\chi^{-1}$-component of $\eta_K$
  is also zero.

  We have proved the regulator equality in the full
  extension of $\R[\Gal(K/k)]$, for both choices of
  $\omega_F$.
  Coefficient conjugation commutes with the regulator
  and fixes its value $\Theta_S^{(r)}(0)$.
  The extended regulator is injective on the components
  of exact rank $r$, by Rubin
  \cite[Lemma~2.7(iii)]{rubin-1996}.
  Since both conjugates are supported there, they coincide.
  Thus $\eta_K$ is rational.
  Its integral construction is also fixed by the
  coefficient involution. Taking invariants in
  $A \otimes_R\bigcap_R^r U_K$ in part~1, or in
  $\calo \otimes_{\Z_\ell}\bigcap_R^r U_K$ in part~2,
  recovers $\bigcap_R^r U_K$: both coefficient algebras
  have basis $1, \delta$, with $\delta \mapsto -\delta$,
  and $2$ is invertible since $\ell \nmid 3W_K$.
  Thus $\eta_K \in \bigcap_R^r U_K$.
\end{proof}

\subsection{A relative theorem on the Rubin--Stark conjecture}
\label{subsec:relative-Rubin-Stark-theorem}

We now recall the lattice in the Rubin--Stark conjecture.
We then prove Theorem~\ref{thm:intro-relative}(ii) by
showing that the element constructed in
Lemma~\ref{lem:cubic-exterior-comparison} belongs to this
lattice after inverting $3W_K$.

\subsubsection{Rubin's integral refinement}
\label{subsec:Rubin-Stark-lattice}

Rubin \cite[Conjecture~B$'$]{rubin-1996} refines the
rational Stark conjecture by requiring an element of an
integral lattice formed from $S$-units congruent to $1$
at specified finite places. We recall this lattice and
the conjecture in arbitrary rank.

Let $K/k$ be finite abelian, and let $S$ be a finite
set of places of $k$ containing the archimedean places
and every finite place ramified in $K/k$.
Choose $r \geq 0$ distinct archimedean places of $k$
splitting completely in $K$, with ordered places
$w_1, \ldots, w_r$ of $K$ above them, and assume
$\lvert S\rvert > r + 1$.
We specify the congruences by a finite set $T$ of finite
places of $k$ disjoint from $S$.
Write $U_{K, S, T}$ for the group of $S$-units of $K$
congruent to $1$ at every place above $T$, and assume
that $U_{K, S, T}$ is torsion-free.
We have imposed the hypotheses of Rubin
\cite[Hypotheses~2.1]{rubin-1996}.
The inequality $\lvert S\rvert > r + 1$ excludes the
trivial character from $e_{S, r}$, since
$\ord_{s = 0}L_S(s, \one) = \lvert S\rvert - 1$.

Since $U_{K, S, T}$ has finite index in
$\calo_{K, S}^\times$, we identify
$\Q\otimes_\Z U_{K, S, T}$ with
$U_S = \Q\otimes_\Z\calo_{K, S}^\times$.
The algebra $\Q[\Gal(K/k)]$ is a finite product of number
fields. Hence the canonical map from
$\bigwedge_{\Q[\Gal(K/k)]}^r U_S$ to
$\bigcap_{\Q[\Gal(K/k)]}^r U_S$ is an isomorphism.
Extending the integral linear maps to $U_S$ embeds
$\bigcap_{\Z[\Gal(K/k)]}^r U_{K, S, T}$ in
$\bigwedge_{\Q[\Gal(K/k)]}^r U_S$.
We require an element of this exterior bidual whose
$\chi^{-1}$-component is zero whenever
$\ord_{s = 0}L_S(s, \chi) > r$.
Following Rubin \cite[Section~2.1]{rubin-1996},
we define the lattice by these two conditions.

\begin{definition}[Rubin's lattice]
  \label{def:Rubin-lattice}
  Inside $\bigwedge_{\Q[\Gal(K/k)]}^r U_S$, define
  \[
    \Lambda_{S, T}^{(r)}
    := \bigcap_{\Z[\Gal(K/k)]}^r U_{K, S, T}
    \cap e_{S, r}\bigwedge_{\Q[\Gal(K/k)]}^r U_S.
  \]
\end{definition}

The Rubin--Stark conjecture includes the Euler factors at $T$.
Set $\delta_{K/k, T} := \prod_{v \in T}
(1 - \rmn(v)[v, K/k]^{-1})$.
We therefore multiply the rational Stark element by
$\delta_{K/k, T}$.

\begin{definition}[Rational rank-$r$ Rubin--Stark element]
  \label{def:rank-two-Rubin-Stark-element}
  Suppose that $\epsilon_{S, \Q}^{(r)}$ exists.
  The rational rank-$r$ Rubin--Stark element for
  $(K/k, S, T)$ and $(w_1, \ldots, w_r)$ is
  \[
    \epsilon_{S, T, \Q}^{(r)}
    := \delta_{K/k, T}\epsilon_{S, \Q}^{(r)}.
  \]
\end{definition}

The Rubin--Stark conjecture predicts that
$\epsilon_{S, T, \Q}^{(r)}$ exists and belongs to
$\Lambda_{S, T}^{(r)}$.

\begin{conjecture}[Rubin--Stark]
  \label{conj:integral-rank-two}
  There is an element $x \in \Lambda_{S, T}^{(r)}$ such that
  \[
    \calr_{w_1, \ldots, w_r}(x)
    = \paren{\prod_{v \in T}
    (1 - \rmn(v)[v, K/k]^{-1})}\Theta_S^{(r)}(0).
  \]
\end{conjecture}

The integral rank-$r$ Rubin--Stark conjecture for
$(K/k, S)$ is this assertion for every finite set $T$
of finite places of $k$ disjoint from $S$ such that
$U_{K, S, T}$ is torsion-free.

Rubin \cite[Proposition~3.1 and Theorem~3.5]{rubin-1996}
proves the following cases of the conjecture.

\begin{proposition}[{Rubin
  \cite[Proposition~3.1 and Theorem~3.5]{rubin-1996}}]
  \label{prop:Rubin-integral-subcases}
  Conjecture~\ref{conj:integral-rank-two} holds if
  $S$ contains more than $r$ places splitting completely
  in $K/k$, or if $K/k$ is quadratic.
\end{proposition}

\subsubsection{Euler factors and congruences at
\texorpdfstring{$T$}{T}}
\label{subsec:rank-two-congruences}

We will compare the exterior biduals of the $S$-units
and the $S$-units congruent to $1$ above $T$.
Multiplication by $\delta_{K/k, T}$ will take the first
into the second. We first prove the algebraic lemma
that accounts for one Euler factor at each place in $T$.
For example, if $p$ is a prime, $U = \Z_p^2$, and
$V = p\Z_p\oplus\Z_p$, then
$p\bigwedge_{\Z_p}^2U = \bigwedge_{\Z_p}^2V$.

\begin{lemma}
  \label{lem:congruences-exterior-square}
  Let $p$ be a prime, $G$ a finite abelian group,
  $e \in \Z_p[G]$ an idempotent, and $R := e\Z_p[G]$.
  If $V \subset U$ are $R$-modules that are finite free
  over $\Z_p$, with $U/V$ finite and embedded in
  $\bigoplus_{i = 1}^m R/(a_i)$, where each $a_i$ is
  not a zero-divisor of $R$, then, for every $r \geq 0$,
  \[
    \paren{\prod_{i = 1}^m a_i}\bigcap_R^r U
    \subseteq \bigcap_R^r V
  \]
  in $\bigwedge_{\Q_p\otimes_{\Z_p}R}^r
  (\Q_p\otimes_{\Z_p}U)$.
\end{lemma}

\begin{proof}
  The assertion for $r = 0$ is immediate.
  For $r > 0$, we work in the rational exterior power
  and test integrality by evaluating against linear maps:
  an element $x$ belongs to $\bigcap_R^r V$ precisely when
  $(\phi_1\wedge\cdots\wedge\phi_r)(x) \in R$
  for all $\phi_i \in \Hom_R(V, R)$, extended to the
  rational span.

  Write $U^* := \Hom_R(U, R)$ and $V^* := \Hom_R(V, R)$.
  Restriction identifies $U^*$ with a submodule of $V^*$.
  Since $U/V$ is finite, every map in $V^*$ extends uniquely
  to a $\Q_p$-linear map on $\Q_p\otimes_{\Z_p}U$ with
  values in $\Q_p\otimes_{\Z_p}R$. It belongs to $U^*$
  precisely when this extension maps $U$ into $R$.

  To compare $R$-linear duals with $\Z_p$-linear duals,
  define $t : \Z_p[G] \longrightarrow \Z_p$ by
  $t(\sum_{g \in G}b_g g) := b_{1_G}$.
  The pairing $(r, s) \longmapsto t(rs)$ is perfect on
  $\Z_p[G]$: the bases $(g)_{g \in G}$ and
  $(g^{-1})_{g \in G}$ are dual. The factors
  $e\Z_p[G]$ and $(1 - e)\Z_p[G]$ are orthogonal,
  so the restriction of this pairing to $R$ is also perfect.
  Consequently, for $M = U$ or $V$, the map
  \[
    \begin{tikzcd}[row sep = tiny]
      \Hom_R(M, R) \arrow[r, "\widesim{}"] &
      \Hom_{\Z_p}(M, \Z_p) \\
      \phi \arrow[r, mapsto] & t\circ\phi
    \end{tikzcd}
  \]
  is an isomorphism, as in Rubin
  \cite[Proposition~1.1, proof]{rubin-1996}.
  On the right, $R$ acts by $(r\ell)(x) := \ell(rx)$.
  We extend $t$ $\Q_p$-linearly to identify $V^*/U^*$.
  After composition with $t$ and reduction modulo $\Z_p$,
  we obtain
  \[
    \begin{tikzcd}[row sep = tiny]
      V^*/U^* \arrow[r, "\widesim{}"] &
      \Hom_{\Z_p}(U/V, \Q_p/\Z_p) \\
      \phi + U^* \arrow[r, mapsto] &
      \bigl(x + V \longmapsto
      t(\phi(x)) + \Z_p\bigr).
    \end{tikzcd}
  \]

  We next find generators of $V^*/U^*$.
  For each $i$, we use the same pairing to identify
  \[
    \begin{tikzcd}[row sep = tiny]
      a_i^{-1}R/R \arrow[r, "\widesim{}"] &
      \Hom_{\Z_p}(R/(a_i), \Q_p/\Z_p) \\
      r + R \arrow[r, mapsto] &
      \bigl(s + a_iR \longmapsto t(rs) + \Z_p\bigr).
    \end{tikzcd}
  \]
  Multiplication by $a_i$ also identifies $a_i^{-1}R/R$
  with $R/(a_i)$. Since $\Q_p/\Z_p$ is injective over $\Z_p$,
  dualizing $U/V \hookrightarrow \bigoplus_i R/(a_i)$ gives
  a surjection
  $\bigoplus_i R/(a_i) \longrightarrow V^*/U^*$.
  Choose $\phi_i \in V^*$ lifting the image of
  $e + a_iR \in R/(a_i)$ from the $i$-th summand. Then
  $V^* = U^* + \sum_i R\phi_i$ and $a_i\phi_i \in U^*$.

  To apply the criterion at the start of the proof,
  we expand an exterior product of $r$ maps in $V^*$
  using $V^* = U^* + \sum_i R\phi_i$.
  A term in which some $\phi_i$ occurs twice vanishes.
  Every remaining term lies in $\bigwedge_R^r U^*$
  after multiplication by $\prod_i a_i$, since
  $a_i\phi_i \in U^*$.
  Thus $\paren{\prod_i a_i}\bigwedge_R^r V^*
  \subseteq \bigwedge_R^r U^*$.
  For $x \in \bigcap_R^r U$, the element
  $(\prod_i a_i)x$ takes values in $R$ on
  $\bigwedge_R^r V^*$, so it belongs to $\bigcap_R^r V$.
\end{proof}

We now apply the lemma to the subgroup of $S$-units
congruent to $1$ above $T$, for a finite abelian extension
$K/k$.
Let $S$ be a finite set of places of $k$ containing the
archimedean places and every finite place ramified in
$K/k$, and let $T$ be a finite set of finite places of $k$
disjoint from $S$.
Let $p$ be a prime not dividing $W_K$,
and set $R := \Z_p[\Gal(K/k)]$,
$U := \Z_p\otimes_\Z\calo_{K, S}^\times$, and
$V := \Z_p\otimes_\Z U_{K, S, T}$.
Both modules are finite free over $\Z_p$ because
$p \nmid W_K$.

We describe $U/V$ by reducing units modulo primes above
$T$. For a finite prime $w$ of $K$, the residue field
$\calo_K/w$ has multiplicative group $(\calo_K/w)^\times$.
The kernel of the map
$a \longmapsto (a \bmod w)_{w \mid T}$ on
$\calo_{K, S}^\times$ is $U_{K, S, T}$.
Since $\Z_p$ is flat over $\Z$, this gives an injection
\[
  U/V \hookrightarrow
  \bigoplus_{v \in T}\bigoplus_{w \mid v}
  \bigl(\Z_p\otimes_\Z(\calo_K/w)^\times\bigr),
\]
as in Rubin \cite[\S1.1, Equation~1]{rubin-1996}.
In particular, $U/V$ is finite.
Each tensor product is the $p$-primary subgroup of
$(\calo_K/w)^\times$.

For each $v \in T$, set
$a_v := 1 - \rmn(v)[v, K/k]^{-1} \in R$.
The extension $K/k$ is unramified at $v$.
For a prime $w$ of $K$ above $v$, arithmetic Frobenius
$[v, K/k]$ acts on $\calo_K/w$ by
$x \longmapsto x^{\rmn(v)}$.
Since $(\calo_K/w)^\times$ is cyclic, inducing from the
decomposition group at $w$ to $\Gal(K/k)$ gives
\[
  \bigoplus_{w \mid v}
  \paren{\Z_p\otimes_\Z(\calo_K/w)^\times}
  \cong R/(a_v).
\]
For every character
$\chi : \Gal(K/k) \longrightarrow \overline{\Q}_p^\times$,
extended linearly to the group ring,
\[
  \chi(a_v) = 1 - \rmn(v)\chi([v, K/k])^{-1} \neq 0.
\]
Indeed, $\chi([v, K/k])$ is a root of unity, whereas
$\rmn(v)$ is an integer greater than $1$.
After extension to $\overline{\Q}_p$, multiplication by
$a_v$ is invertible on every character component.
Thus $a_v$ is not a zero-divisor of $R$.

We have now identified the reduction modulo primes above
$T$ with an injection
$U/V \hookrightarrow \bigoplus_{v \in T}R/(a_v)$.
We can therefore apply
Lemma~\ref{lem:congruences-exterior-square} to obtain,
for every $r \geq 0$,
\[
  \delta_{K/k, T}\bigcap_R^r U
  \subseteq \bigcap_R^r V,
\]
since $\delta_{K/k, T} = \prod_{v \in T}a_v$.

We now prove Theorem~\ref{thm:intro-relative}(ii).
We first prove the integrality of the rational Stark
element for $K/k$ away from $3W_K$.
For a finite set $T$ of finite places of $k$ disjoint
from $S$, we then multiply by the Euler factors at $T$
and use Lemma~\ref{lem:congruences-exterior-square} to prove
integrality in the exterior bidual of $S$-units congruent
to $1$ above $T$.

\begin{proof}[Proof of Theorem~\ref{thm:intro-relative}(ii)]
  Let $r := r_{k_0} + r_F$,
  and recall that $L = K \cap k_0^\ab$.
  By the hypotheses and the table before
  Lemma~\ref{lem:cubic-exterior-comparison},
  at least $r$ archimedean places of $k$ split completely
  in $K$. If there are more than $r$, we can take $0$
  for the rank-$r$ Rubin--Stark element, by Rubin
  \cite[Proposition~3.1]{rubin-1996}.
  Otherwise, the same table shows that $r_{k_0}$ and
  $r_F$ count all archimedean places splitting
  completely in $L/k_0$ and $KF/F$, respectively.
  Fix a prime $\ell \nmid 3W_K$, and set
  $R := \Z_\ell[\Gal(K/k)]$ and
  $U := \Z_\ell\otimes_\Z\calo_{K, S}^\times$.

  \emph{1. The Stark element for $L/k_0$.}
  Let $\eta_0$ be the rational rank-$r_{k_0}$ Stark element
  for $(L/k_0, S_{k_0})$.
  To recover $\eta_0$ from a Rubin--Stark element
  over $\Z_\ell$,
  we choose two finite places whose
  Euler factors are units in $\Z_\ell[\Gal(L/k_0)]$.
  Since $\ell \nmid W_L = W_K$, the group
  $\Gal(L(\mu_\ell)/L)$ is nontrivial.
  We apply the Chebotarev density theorem to a nonidentity
  element of this group to choose two primes of $k_0$
  outside $S_{k_0}$ and not above $\ell$, with distinct
  residue characteristics, that split completely in $L$
  and satisfy $\rmn(v) \not\equiv 1 \pmod\ell$.
  Let $T_0$ consist of these two primes.
  Then $U_{L, S_{k_0}, T_0}$ is torsion-free, and
  $\delta_{L/k_0, T_0} = \prod_{v \in T_0}(1 - \rmn(v))$
  is a unit of $\Z_\ell[\Gal(L/k_0)]$.
  We can therefore divide the Rubin--Stark element for
  $T_0$ by $\delta_{L/k_0, T_0}$ to conclude that
  \[
    \eta_0 \in
    \bigcap_{\Z_\ell[\Gal(L/k_0)]}^{r_{k_0}}
    (\Z_\ell\otimes_\Z\calo_{L, S_{k_0}}^\times).
  \]

  \emph{2. The Stark element for $K/k$.}
  Recall that $F$ is the quadratic resolvent of $k/k_0$.
  Let $\eta_F$ be the rational rank-$(r_F)$
  Stark element for $(KF/F, S_F)$.
  Choose a generator $\tau$ of $\Gal(KF/LF)$, which is
  the commutator subgroup of $\Gal(KF/k_0)$, and set
  $e := (2 - \tau - \tau^2)/3$.
  The group $\abrac{\tau}$ fixes
  every root of unity in $KF$,
  so $e(\Z_\ell\otimes_\Z\calo_{KF, S_F}^\times)$
  is free over $\Z_\ell$.

  Since $F(\mu_\ell)/k_0$ is abelian, $\tau$ fixes
  $KF\cap F(\mu_\ell)$.
  Hence $(\tau, 1)$ defines an automorphism of
  $KF(\mu_\ell)/F$.
  We apply the Chebotarev density theorem to $(\tau, 1)$
  to choose primes $u$ of $F$ outside $S_F$ with
  $[u, KF/F] = \tau$ and $\rmn(u) \equiv 1 \pmod\ell$.
  On the $e$-component, the Euler factor at such a prime
  reduces modulo $\ell$ to $1 - \tau^2$, which is a unit
  because $(1 - \tau)(1 - \tau^2) = 3e$.
  We choose two such primes with distinct residue
  characteristics, so the group of $S_F$-units of $KF$
  congruent to $1$ above both primes is torsion-free.
  The Euler factors are units on the $e$-component,
  so the rank-$r_F$ Rubin--Stark hypothesis
  over $\Z_\ell$ gives
  \[
    e\eta_F \in
    \bigcap_{e\Z_\ell[\Gal(KF/F)]}^{r_F}
    e(\Z_\ell\otimes_\Z\calo_{KF, S_F}^\times).
  \]
  We have verified the two integrality hypotheses of
  Lemma~\ref{lem:cubic-exterior-comparison},
  so we can construct the rational rank-$r$ Stark
  element $\eta_K \in \bigcap_R^r U$.

  \emph{3. The congruences at $T$.}
  Let $T$ be a finite set of finite places of $k$
  disjoint from $S$ such that the group $U_{K, S, T}$
  of $S$-units congruent to $1$ at every prime above $T$
  is torsion-free. Set
  $V := \Z_\ell\otimes_\Z U_{K, S, T}$.
  To pass from $U$ to the subgroup $V$ of units congruent
  to $1$ above $T$, we use the inclusion proved after
  Lemma~\ref{lem:congruences-exterior-square}:
  \[
    \delta_{K/k, T}\bigcap_R^r U
    \subseteq \bigcap_R^r V.
  \]
  Thus $\delta_{K/k, T}\eta_K \in \bigcap_R^r V$.
  By definition, $\delta_{K/k, T}\eta_K$ is the rational
  rank-$r$ Rubin--Stark element for $(K/k, S, T)$.
  Since $\ell \nmid 3W_K$ was arbitrary, this proves
  the rank-$r$ Rubin--Stark conjecture after inverting
  $3W_K$.
\end{proof}

\begin{remark}
  The proof of Theorem~\ref{thm:intro-relative}(ii)
  only uses the existence of the rational
  Stark elements $\eta_0$ and $\eta_F$
  and the two integrality conditions in
  Lemma~\ref{lem:cubic-exterior-comparison}
  for every prime $\ell \nmid 3W_K$.
  In particular, integrality over $F$
  is needed only on the $e$-component.
\end{remark}

Theorem~\ref{thm:intro-relative}(ii)
leads to new unconditional cases of the
Rubin--Stark conjecture,
such as the following example
with an $L$-function that
vanishes at $s = 0$
with order $3$.

\begin{example}
  \label{ex:relative-rank-three}
  Let $\alpha$ be a root of the
  irreducible polynomial
  $P(x) := x^3 - x^2 + 7x - 6$
  of discriminant $-1563 = -3 \cdot 521$,
  and let $\zeta_{11}$ be a primitive eleventh root of unity.
  Consider the field extensions
  \begin{align*}
    k_0 &:= \Q(i), &&k := k_0(\alpha), \\
    F_0 &:= \Q(\sqrt{-1563}),
    &&F = k_0(\sqrt{-1563}),
  \end{align*}
  in which case $k/k_0$ is non-Galois cubic
  with quadratic resolvent $F$.
  Consider the abelian extension
  $K := k(\sqrt{-1563}, \zeta_{11} + \zeta_{11}^{-1})$
  of $k$.
  With $L := K \cap k_0^\ab
  = F(\zeta_{11} + \zeta_{11}^{-1})$,
  one can check that $K = kL \subset k \cdot k_0^\ab$
  is a cyclic extension of degree $10$
  over $k$.
  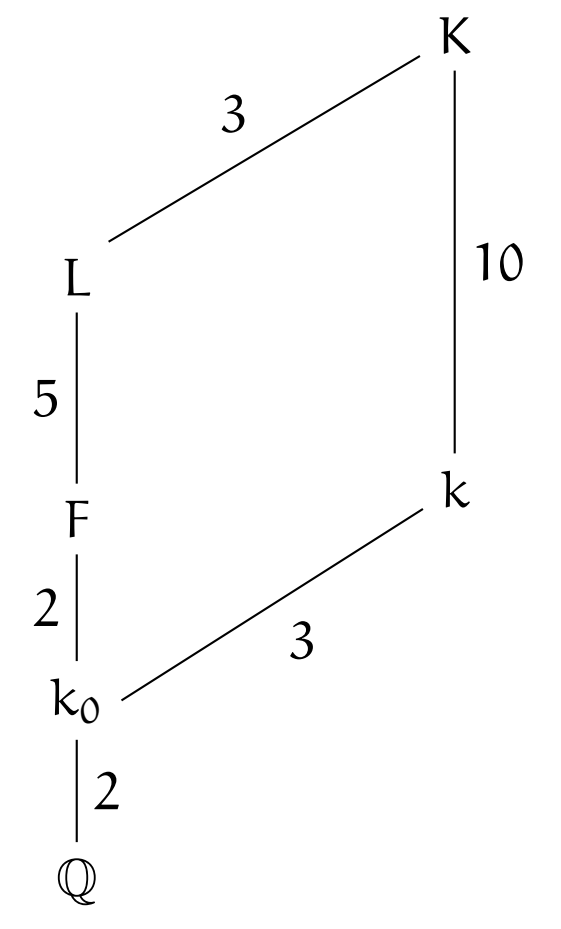
\begin{figure}[h]
    \centering
    \begin{tikzpicture}[
    x = 1cm, y = 0.8cm,
    every node/.style = {font = \small, inner sep = 3pt}
  ]
    \node (Q) at (0, -1.5) {$\Q$};
    \node (kzero) at (0, 0) {$k_0$};

    \node
      (F) at (0, 1.5)
      {$F$};
    \node
      (k) at (2.5, 1.75)
      {$k$};
    \node (Fab) at (0, 3.5) {$L$};
    \node (kcyc) at (2.5, 5.5) {$K$};

    \draw (Q) -- node[right] {$2$} (kzero);
    \draw (kzero) -- node[left] {$2$} (F);
    \draw (kzero.east) -- node[below right] {$3$}
      (k);
    \draw
      (F.north) -- node[left] (abelian)
      {\shortstack{$5$}} (Fab.south);
    \draw
      (Fab.north east) -- node[above left] (LK)
      {\shortstack{$3$}} (kcyc);
    \draw
      (k.north) -- node[right] (cyclotomic)
      {$10$} (kcyc.south);


  \end{tikzpicture}
    \caption{Field extensions in
      Example~\ref{ex:relative-rank-three}.}
    \label{figure:relative-rank-three}
  \end{figure}

  For each subfield $K'$ of $K$, let $S_{K'}$ consist
  of its archimedean places and the primes above
  $2$, $3$, $11$, and $521$, and set $S := S_k$.
  For every faithful character $\chi$ of $\Gal(K/k)$,
  \[
    \ord_{s=0} L_S(s,\chi) = 3.
  \]
  Indeed, $k$ has three complex places, and no finite
  place in $S$ splits completely in $K/k$:
  the primes $2$, $3$, and $521$ are inert in
  $\Q(\zeta_{11} + \zeta_{11}^{-1})$,
  and their residue degrees in $k$ are prime to $5$;
  the prime $11$ is totally ramified in this quintic
  field and unramified in $k$.

  The extension $K/F_0$ is abelian.
  By Theorem~\ref{thm:imaginary-quadratic-Stark}
  and Rubin~\cite[Proposition~2.5]{rubin-1996},
  the integral rank-$1$ Rubin--Stark conjectures
  hold for $(L/k_0, S_{k_0})$ and
  $(K/F_0, S_{F_0})$.
  Since $F/F_0$ is quadratic and
  $W_K = 4$ is coprime to $[K:F] = 15$,
  Popescu's base-change and comparison results
  \cite[Theorem~2.3.2 and Corollary~2.3.4(2)]
  {popescu-2004-rubin}
  give the integral rank-$2$ Rubin--Stark conjecture
  for $(K/F, S_F)$.
  Hence Theorem~\ref{thm:intro-relative}(ii), with
  $r_{k_0} = 1$ and $r_F = 2$, gives the rank-$3$
  Rubin--Stark conjecture for $(K/k, S)$ over $\Z[1/6]$.
  In particular, this includes the $5$-primary part,
  although $5$ divides $[K:k]$.
\end{example}

We now prove a relative theorem
on the Rubin--Stark conjecture for
Galois cubic extensions,
in all ranks divisible by $3$,
after inverting $3W_K$.
This extends the case $r_{k_0} = 1$,
which follows from Popescu's base-change
and comparison theorems
\cite[Theorems~3.1.1 and~5.5.1(3), (4)]
  {popescu-2002-base-change}.

\begin{theorem}
  \label{thm:Galois-Rubin-Stark}
  Let
  \begin{itemize}
    \item $k/k_0$ be a Galois cubic extension of number fields;
    \item $K/k$ be finite abelian with
      $K \subset k \cdot k_0^\ab$;
    \item $S_{k_0}$ be a finite set of places of $k_0$
      containing the archimedean places, every prime ramified
      in $K/k_0$, and at least two finite places;
    \item $S$ consist of the places above $S_{k_0}$ in $k$; and
    \item $r_{k_0}$ be the number of archimedean places of
      $k_0$ splitting completely in $K$.
  \end{itemize}
  If the rank-$r_{k_0}$ Rubin--Stark conjecture
  for $(K/k_0, S_{k_0})$ holds over $\Z[1/(3W_K)]$,
  then the rank-$3r_{k_0}$ Rubin--Stark conjecture
  for $(K/k, S)$ holds over the same ring.
  Moreover when $r_{k_0} = 1$, the same implication holds
  over $\Z[1/W_K]$.
\end{theorem}

\begin{proof}
  Since $k \subset k_0^\ab$, the extension $K/k_0$ is
  abelian and $L = K$.
  Each archimedean place of $k_0$ splits into three
  places in $k$, so $3r_{k_0}$ archimedean places of $k$
  split completely in $K$.

  \emph{1. The Stark element for $K/k_0$.}
  Fix a prime $\ell \nmid 3W_K$, and set
  $A := \Z_\ell[\Gal(K/k_0)]$,
  $R := \Z_\ell[\Gal(K/k)]$, and
  $U := \Z_\ell\otimes_\Z\calo_{K, S}^\times$.
  Since $\ell \nmid W_K$, we can choose two primes of
  $k_0$ outside $S_{k_0}$, with distinct residue
  characteristics, that split completely in $K$ and have
  norms not congruent to $1$ modulo $\ell$.
  We make this choice by the Chebotarev density theorem,
  as in part~1 of the proof of
  Theorem~\ref{thm:intro-relative}(ii).
  The Euler factors at these primes are units in $A$,
  so the Rubin--Stark hypothesis over $\Z_\ell$
  implies that the rational Stark element $\eta$ for
  $(K/k_0, S_{k_0})$ belongs to $\bigcap_A^{r_{k_0}}U$.

  \emph{2. The Stark element for $K/k$.}
  Each character of $\Gal(K/k)$ has three extensions to
  $\Gal(K/k_0)$. To separate these extensions, we enlarge
  the coefficient ring.
  Choose $t \in \Gal(K/k_0)$ generating the quotient
  $\Gal(k/k_0)$, and let $\zeta$ be a primitive cube
  root of unity. The algebra
  $R' := A\otimes_{\Z_\ell}\Z_\ell[\zeta]$ is faithfully
  flat over $R$. Since $t^3 \in \Gal(K/k)$ and
  $\ell \neq 3$, we have $R'\otimes_R A \cong (R')^3$,
  with $1\otimes t$ mapping to $q_i := t\zeta^i$
  for $i = 0, 1, 2$.
  Thus $R'\otimes_R U = U_0\oplus U_1\oplus U_2$,
  where $1\otimes t$ acts on $U_i$ as $q_i$.
  Extending coefficients carries $\eta$ to elements
  $\eta_i \in \bigcap_{R'}^{r_{k_0}}U_i$.

  If $w_j$ are the places used over $k_0$, use first
  the $w_j$, then the $tw_j$, then the $t^2w_j$ over $k$.
  On $U_i$, the equivariant logarithm for $K/k_0$,
  evaluated at $t = q_i$, is $3\ell_w$, where $\ell_w$
  is the logarithm for $K/k$; also
  $\ell_{t^a w} = q_i^{-a}\ell_w$.
  Set $D := \det(q_i^{-a})_{0 \leq a, i \leq 2}$.
  The three evaluations at $t = q_i$ of the regulator
  for $K/k_0$ have product equal to the regulator for
  $K/k$ of
  \[
    (27/D)^{r_{k_0}}
    \eta_0\wedge\eta_1\wedge\eta_2.
  \]

  By the Artin decomposition, each $L$-function over $k$
  is the product of the three $L$-functions over $k_0$
  for the extensions of its character. Their leading
  terms therefore have the required product.
  If $\ord_{s = 0}L_S(s, \chi) > 3r_{k_0}$, then one
  extension $\psi$ satisfies
  $\ord_{s = 0}L_{S_{k_0}}(s, \psi) > r_{k_0}$,
  so the corresponding component of the product is zero.

  We can form the same expression over
  $\Q[\Gal(K/k_0)]\otimes_\Q\Q(\zeta)$.
  Interchanging two $q_i$ multiplies $D^{r_{k_0}}$ and
  the exterior product by the same sign $(-1)^{r_{k_0}}$.
  Thus the expression is invariant under the permutations
  of the $q_i$ and descends to $\Q[\Gal(K/k)]$.
  We denote the resulting rational Stark element by
  $\eta_K$.
  Since $D^2 = -27t^{-6}$ is a unit, the exterior product
  is integral over $R'$. We can therefore descend this
  integrality through the faithfully flat extension
  $R'/R$ to conclude that
  $\eta_K \in \bigcap_R^{3r_{k_0}}U$.

  \emph{3. The congruences at $T$.}
  Let $T$ be a finite set of finite places of $k$ disjoint
  from $S$ such that the group $U_{K, S, T}$ of $S$-units
  congruent to $1$ above $T$ is torsion-free.
  We now use Lemma~\ref{lem:congruences-exterior-square}
  to pass to this subgroup of units:
  \[
    \delta_{K/k, T}\eta_K \in
    \bigcap_R^{3r_{k_0}}
    (\Z_\ell\otimes_\Z U_{K, S, T}).
  \]
  This proves the rank-$3r_{k_0}$ Rubin--Stark conjecture
  after inverting $3W_K$.

  If $r_{k_0} = 1$, the final assertion follows from Popescu
  \cite[Theorems~2.3.2 and~2.3.3]{popescu-2004-rubin}.
  His conjecture is equivalent to the Rubin--Stark
  conjecture in rank $1$ and, in rank $3$, at every prime
  $\ell \nmid W_K$, since $K$ has no nontrivial
  $\ell$-primary roots of unity.
\end{proof}

\subsection{The rank-\texorpdfstring{$2$}{2} case}
\label{subsec:relative-rank-two-rational}

We now construct rank-$2$ Stark elements from rank-$1$
Stark units over almost totally real cubic extensions.
In this subsection, let $k_0$ be totally real, let
$k/k_0$ be an almost totally real cubic extension with
quadratic resolvent $F$, and let $K/k$ be finite abelian
with $K \subset k \cdot k_0^\ab$.
Recall that $L = K \cap k_0^\ab$ and $K = kL$.
Write $v_0$ for the real place of $k_0$ below the complex
place of $k$. We use complex conjugation $c$ at the real
place of $k$ above $v_0$, with fixed fields $K^+$ and $L^+$.

With this construction, we will prove integrality at $2$
and $3$ under the following hypotheses.

\begin{theorem}
  \label{thm:rank-two-integral-refinement}
  Let
  \begin{itemize}
    \item $k_0$ be a totally real number field;
    \item $k/k_0$ be an almost totally real cubic
      extension, and let $v_0$ be the real place of
      $k_0$ below the complex place of $k$;
    \item $K/k$ be finite abelian with
      $K \subset k \cdot k_0^\ab$, $[K^+ : k]$ odd,
      and either $3 \nmid [K^+ : k]$ or
      $3 \nmid W_{K^+F}$;
    \item $S_{k_0}$ be a finite set of places of $k_0$
      containing the archimedean places, every prime
      ramified in $K/k_0$, and at least two finite
      places; and
    \item $S$ and $S_F$ consist of the places above
      $S_{k_0}$ in $k$ and $F$, respectively.
  \end{itemize}
  If the integral rank-$1$ abelian Stark conjecture holds
  for $(KF/F, S_F)$ at the complex place of $F$ and for
  $(L^+/k_0, S_{k_0})$ at $v_0$, then the integral
  rank-$2$ Rubin--Stark conjecture for $(K^+/k, S)$ holds
  with the real and complex places of $k$ above $v_0$.
\end{theorem}

\begin{remark}
  If $k_0 \neq \Q$,
  then we can take the rank-$2$ Rubin--Stark
  element in Theorem~\ref{thm:rank-two-integral-refinement}
  to be $0$.
  Indeed, $[L^+ : k_0] = [K^+ : k]$ is odd,
  so $L^+$ is totally real
  and every archimedean place of $k$
  splits completely in $K^+/k$.
  Since $k$ has more than two archimedean places and
  $\lvert S \rvert > 3$, we deduce from
  the order-of-vanishing formula
  \cite[Chapter~I, Proposition~3.4]{tate-stark-1984}
  that $\ord_{s = 0}L_S(s, \chi) > 2$
  for every character $\chi$.
  For $K^+/k$, this means $e_{S, 2} = 0$ and
  $\Theta_S^{(2)}(0) = 0$.
\end{remark}

We first construct the units and compute their regulators
for arbitrary $[K^+ : k]$, assuming that $K$ contains $F$.
We then remove this assumption in
Proposition~\ref{prop:relative-rank-two-regulator}.
Finally, we identify the rational rank-$2$ Stark element
in Corollary~\ref{cor:rational-rank-two-element} and prove
its integrality in Rubin's lattice.

\subsubsection{The regulator when
\texorpdfstring{$K$ contains $F$}{K contains F}}
\label{subsec:rank-two}

We now assume that $K$ contains $F$.
Then $F \subset L$, $L = L^+F$, and $K = kL = K^+F$,
as in Section~\ref{subsec:Stark-hypotheses}.
Choose a generator $\tau$ of $\Gal(K/L)$.
Let $S$, $S_{k_0}$, and $S_F$ be finite sets of places of
$k$, $k_0$, and $F$, respectively, each containing all
archimedean places of its field and satisfying
Equation~\ref{eq:local-factorization}.
Suppose also that
\begin{itemize}
  \item $S$ contains every finite place ramified in $K^+/k$;
  \item $S_{k_0}$ contains every finite place ramified in
    $L^+/k_0$ and at least two finite places; and
  \item $S_F$ contains every finite place ramified in $K/F$
    and at least two finite places.
\end{itemize}
These are the hypotheses in
Section~\ref{subsec:decomposition-partial}(ii).
Choose the real and complex places $w_\R$ and $w_\C$
of $K^+$ above one real place $w_{L^+}$ of $L^+$ over $v_0$.
Let $w_K$ be a place of $K$ above $w_\C$, and replace $\tau$
by $\tau^{-1}$ if necessary so that
\begin{equation}
  \label{eq:rank-two-place-normalization}
  c w_K = \tau w_K.
\end{equation}
Assume that there are Stark units
$\epsilon_{k_0} \in \calo_{L^+}^\times$ for
$(L^+/k_0, S_{k_0}, v_0, w_{L^+})$ and
$\epsilon_F \in \calo_K^\times$ for $(K/F, S_F, w_K)$.
The place $v_0$ splits completely in $L^+/k_0$, and
$W_{L^+} = 2$ because $L^+$ has a real embedding.
We set $u_\R := \epsilon_{k_0}$.

To express $L_{S_F}'(0, \chi_F\omega_F)$ using a unit of
$K^+$, we define
\begin{equation}
  \label{eq:rank-two-unit-definition}
  u_\C := \frac{\rmn_{K/K^+}(\epsilon_F)^3}
    {\rmn_{K/L^+}(\epsilon_F)},
\end{equation}
which lies in $K^+$ and satisfies
$\rmn_{K^+/L^+}(u_\C) = 1$.
We next express the logarithms of the conjugates of
$u_\C$ at $w_\C$ in terms of those of $\epsilon_F$ at $w_K$.

\begin{lemma}
  \label{lem:rank-two-norm-comparison}
  The unit $u_\C \in \calo_{K^+}^\times$ satisfies
  \[
    \log\verts{u_\C^\sigma}_{w_\C}
    = \log\verts{\epsilon_F^{\sigma(1 - c)}}_{w_K}
  \]
  for every $\sigma \in \Gal(K/kF)$.
\end{lemma}

\begin{proof}
  Fix $\sigma \in \Gal(K/kF)$ and set
  $a := \log\lvert\epsilon_F^\sigma\rvert_{w_K}$ and
  $b := \log\lvert\epsilon_F^{\sigma c}\rvert_{w_K}$.
  We compute the logarithm of each norm in the definition
  of $u_\C$.
  By Equation~\ref{eq:rank-two-place-normalization} and
  Lemma~\ref{lem:elliptic-conjugation-symmetry}, the
  logarithms at the three places of $K$ above $w_{L^+}$
  are $a$, $b$, and $b$.
  The normalized absolute value at $w_\C$ is the square
  of the real absolute value on $L^+$, so
  \[
    \log\verts{\rmn_{K/L^+}(\epsilon_F^\sigma)}_{w_\C}
    = 2a + 4b.
  \]
  For the norm of $K/K^+$, we similarly have
  \[
    \log\verts{\rmn_{K/K^+}(\epsilon_F^\sigma)^3}_{w_\C}
    = 3a + 3b.
  \]
  We substitute these two equalities into the definition
  of $u_\C$ in Equation~\ref{eq:rank-two-unit-definition}
  to obtain
  \begin{align*}
    \log\verts{u_\C^\sigma}_{w_\C}
      &= (3a + 3b) - (2a + 4b) \\
      &= a - b \\
      &= \log\verts{\epsilon_F^{\sigma(1 - c)}}_{w_K}.
  \end{align*}
\end{proof}

We now apply Definition~\ref{def:global-rank-one-Stark-unit}(i)
to $\epsilon_{k_0}$ and $\epsilon_F$ to evaluate
$L_S''(0, \chi)$.

\begin{proposition}
  \label{prop:rank-two-regulator}
  Under the hypotheses above,
  if $\epsilon_{k_0}$ and $\epsilon_F$ are Stark units for
  $(L^+/k_0, S_{k_0}, v_0, w_{L^+})$ and
  $(K/F, S_F, w_K)$ respectively,
  then
  \[
    L_S''(0, \chi) = \frac{1}{W_K}
    \paren{\sum_{\alpha \in \Gal(K^+/k)}\chi(\alpha)
    \log\verts{u_\R^\alpha}_{w_\R}}
    \paren{\sum_{\beta \in \Gal(K^+/k)}\chi(\beta)
    \log\verts{u_\C^\beta}_{w_\C}}
  \]
  for $u_\R := \epsilon_{k_0}$,
  $u_\C := \frac{\rmn_{K/K^+}(\epsilon_F)^3}
  {\rmn_{K/L^+}(\epsilon_F)}$,
  and every character $\chi$ of $\Gal(K^+/k)$.
\end{proposition}

\begin{proof}
  Let $\chi_{k_0}$ be the character of $\Gal(L^+/k_0)$
  corresponding to $\chi$ under restriction from $K^+$.
  Its inflation to $\Gal(L/k_0)$ is even at $v_0$.
  To compute $L_S''(0, \chi)$, we first use
  Proposition~\ref{prop:local-Artin-factorization} to write
  \[
    L_S(s, \chi)
      = L_{S_{k_0}}(s, \chi_{k_0}) L_{S_F}(s, \chi_F\omega_F).
  \]
  We will differentiate both sides twice at $s = 0$.
  Both factors vanish there, by the following applications
  of \cite[Chapter~I, Proposition~3.4]{tate-stark-1984}:
  \begin{itemize}
    \item If $\chi_{k_0}$ is nontrivial, the split place $v_0$ gives
      $L_{S_{k_0}}(0, \chi_{k_0}) = 0$.
    \item For the trivial character $\one$,
      $L_{S_{k_0}}(0, \one) = 0$ because $S_{k_0}$
      contains a finite place.
    \item The character $\chi_F\omega_F$ is nontrivial on
      $\Gal(K/L)$, and the complex place of $F$ gives
      $L_{S_F}(0, \chi_F\omega_F) = 0$.
  \end{itemize}
  We therefore have
  \[
    L_S''(0, \chi)
      = 2L_{S_{k_0}}'(0, \chi_{k_0})
        L_{S_F}'(0, \chi_F\omega_F).
  \]

  The regulator of the Stark unit $u_\R = \epsilon_{k_0}$
  for $L^+/k_0$ is
  \[
    L_{S_{k_0}}'(0, \chi_{k_0}) = -\frac{1}{2}
    \sum_{\alpha \in \Gal(K^+/k)}\chi(\alpha)
    \log\verts{u_\R^\alpha}_{w_\R}.
  \]
  We compute the derivative over $F$ as in
  Lemma~\ref{lem:elliptic-regulator-comparison}, using
  Lemma~\ref{lem:elliptic-conjugation-symmetry} to get
  \[
    L_{S_F}'\paren{0, \chi_F\omega_F}
      = -\frac{1}{W_K}
        \sum_{\beta \in \Gal(K/kF)}\chi(\beta)
        \log\verts{
        \epsilon_F^{\beta(1 - c)}}_{w_K}.
  \]
  By Lemma~\ref{lem:rank-two-norm-comparison}, we may replace
  $\log\lvert \epsilon_F^{\beta(1 - c)}\rvert_{w_K}$ by
  $\log\lvert u_\C^\beta\rvert_{w_\C}$ in the expression for
  $L_{S_F}'(0, \chi_F\omega_F)$.
  Substituting these expressions for the first derivatives
  into $L_S''(0, \chi) = 2L_{S_{k_0}}'(0, \chi_{k_0})
  L_{S_F}'(0, \chi_F\omega_F)$ gives
  \[
    L_S''(0, \chi) = \frac{1}{W_K}
    \paren{\sum_{\alpha \in \Gal(K^+/k)}\chi(\alpha)
    \log\verts{u_\R^\alpha}_{w_\R}}
    \paren{\sum_{\beta \in \Gal(K^+/k)}\chi(\beta)
    \log\verts{u_\C^\beta}_{w_\C}}.
  \]
\end{proof}

We use the abelian extension $K(\epsilon_F^{1/W_K})/F$
to find $\gamma \in \calo_{K^+}^\times$ with
$\gamma^{W_K} = u_\C$.
We will use $\gamma$ to change the denominator in the
regulator when $F \not\subset K$ and to prove integrality.

\begin{lemma}
  \label{lem:rank-two-denominators}
  There is a unit $\gamma \in \calo_{K^+}^\times$ such that
  $\gamma^{W_K} = u_\C$.
\end{lemma}

\begin{proof}
  Choose $\alpha$ with $\alpha^{W_K} = \epsilon_F$.
  Since $\epsilon_F$ is a Stark unit for $(K/F, S_F, w_K)$,
  the extension $K(\alpha)/F$ is abelian by
  Definition~\ref{def:global-rank-one-Stark-unit}(iii).
  Every root of unity in $K$ lies in $L$, so $\tau$ fixes
  $\mu(K)$.
  Choose a lift $\widetilde\tau \in \Gal(K(\alpha)/F)$ of
  $\tau$, and set $\beta := \alpha/\widetilde\tau(\alpha)$.
  Every automorphism of $K(\alpha)/K$ multiplies $\alpha$
  and $\widetilde\tau(\alpha)$ by the same root of unity,
  so $\beta \in K$. Since $\beta^{W_K} = \epsilon_F^{1 - \tau}$
  is a global unit, $\beta \in \calo_K^\times$.
  Set $\gamma := \rmn_{K/K^+}(\beta^{2 + \tau})$.
  The identity $(1 - \tau)(2 + \tau) = 2 - \tau - \tau^2$
  gives
  \begin{equation*}
    \gamma^{W_K} = \frac{\rmn_{K/K^+}(\epsilon_F)^3}
    {\rmn_{K/L^+}(\epsilon_F)} = u_\C. \qedhere
  \end{equation*}
\end{proof}

\subsubsection{Removing the assumption that
  \texorpdfstring{$K$}{K} contains \texorpdfstring{$F$}{F}}
\label{subsec:rank-two-general-extensions}

We now remove the assumption that $K$ contains $F$.
We first descend a Stark unit from $KF/F$ to
$(K^+F)/F$, then apply Proposition~\ref{prop:rank-two-regulator}
to $K' := K^+F = k(L^+F)$.
Using the root from Lemma~\ref{lem:rank-two-denominators},
we obtain a regulator with denominator $W_K$
instead of $W_{K'}$.
Let $\calr_w(a, \chi)$ denote the
character sum over $\Gal(K^+/k)$
defined in Equation~\ref{eq:calr}.

\begin{proposition}
  \label{prop:relative-rank-two-regulator}
  Let
  \begin{itemize}
    \item $k_0$ be a totally real number field;
    \item $k/k_0$ be an almost totally real cubic extension
      with quadratic resolvent $F$;
    \item $K/k$ be finite abelian with
      $K \subset k \cdot k_0^\ab$;
    \item $S_{k_0}$ be a finite set of places of $k_0$,
      including every archimedean place,
      every finite prime ramified in $K/k_0$,
      and at least two finite places;
    \item $S$ and $S_F$ consist of the places above
      $S_{k_0}$ in $k$ and $F$, respectively; and
    \item $w_\R$ and $w_\C$ be real and complex places of
      $K^+$ above the same real place of $L^+$ over $v_0$.
  \end{itemize}
  If the integral rank-$1$ Stark conjecture holds for
  $(KF/F, S_F)$ at the complex place of $F$ and for
  $(L^+/k_0, S_{k_0})$ at $v_0$,
  then there are units $u_\R \in \calo_{L^+}^\times$
  and $u_\C \in \calo_{K^+}^\times$ such that
  \[
    L_S''(0, \chi) = \frac{1}{W_K}
    \paren{\sum_{\alpha \in \Gal(K^+/k)}\chi(\alpha)
    \log\verts{u_\R^\alpha}_{w_\R}}
    \paren{\sum_{\beta \in \Gal(K^+/k)}\chi(\beta)
    \log\verts{u_\C^\beta}_{w_\C}}
  \]
  for every character $\chi$ of $\Gal(K^+/k)$.
  Furthermore, $\rmn_{K^+/L^+}(u_\C) = 1$.
\end{proposition}

\begin{proof}
  We first construct units using $K' := K^+F$ and then
  adjust the exponent of the second unit to obtain
  denominator $W_K$ rather than $W_{K'}$.

  The extensions $KF/F$ and $K'/F$ are abelian, and
  $K' \subset KF$.
  Lemma~\ref{lem:resolvent-ramification} shows that
  $S_{k_0}$ contains every finite prime ramified in $KF/k_0$.
  Since $K' \subset KF$, the sets $S$, $S_{k_0}$, and $S_F$
  contain the finite places ramified in $K'/k$,
  $(L^+F)/k_0$, and $K'/F$, respectively.
  We choose a Stark unit for $KF/F$ at a place above
  $w_\C$ and descend it to $K'/F$ by
  Lemma~\ref{lem:Stark-descent}.
  Complex conjugation fixes $K^+$ in $K'$ and $L^+$ in
  $L^+F$.
  Since $S$ and $S_F$ consist of the places above $S_{k_0}$,
  Equation~\ref{eq:local-factorization} holds for $K'/k$
  by Corollary~\ref{cor:rationally-saturated-factorization}.
  We can therefore apply Proposition~\ref{prop:rank-two-regulator}
  to the Stark units for $L^+/k_0$ and $K'/F$ to construct
  $u_\R \in \calo_{L^+}^\times$ and
  $\widetilde u_\C \in \calo_{K^+}^\times$ such that
  $L_S''(0, \chi) = W_{K'}^{-1}\calr_{w_\R}(u_\R, \chi)
  \calr_{w_\C}(\widetilde u_\C, \chi)$ for every $\chi$,
  and $\rmn_{K^+/L^+}(\widetilde u_\C) = 1$.

  To change the denominator $W_{K'}$ to $W_K$, we apply
  Lemma~\ref{lem:rank-two-denominators} to $K'/k$.
  Choose $\gamma \in \calo_{K^+}^\times$ with
  $\gamma^{W_{K'}} = \widetilde u_\C$,
  and define $u_\C := \gamma^{W_K}$.
  For every $\sigma \in \Gal(K^+/k)$, we have
  \[
    \frac{1}{W_K}\log\verts{u_\C^\sigma}_{w_\C}
    = \log\verts{\gamma^\sigma}_{w_\C}
    = \frac{1}{W_{K'}}
    \log\verts{\widetilde u_\C^\sigma}_{w_\C}.
  \]
  We substitute this equality into the expression for
  $L_S''(0, \chi)$ just obtained to replace $W_{K'}$ by $W_K$.
  Also, $\rmn_{K^+/L^+}(\gamma)^{W_{K'}} = 1$.
  Since $L^+$ has a real embedding, this norm is $1$ or
  $-1$, so $\rmn_{K^+/L^+}(u_\C) = 1$ because $W_K$ is even.
\end{proof}

\subsubsection{The determinant and the rational Stark element}
\label{subsec:rank-two-determinant}

We now interpret the regulator formulas
of Propositions~\ref{prop:rank-two-regulator}
and \ref{prop:relative-rank-two-regulator}
as determinants of logarithms.
We use the properties shared by the two
constructions: $u_\R$ belongs to $L^+$ and
$\rmn_{K^+/L^+}(u_\C) = 1$.

Let $v_\C$ be the complex place of $k$ and let $v_\R$
be the real place of $k$ above $v_0$.
Choose real and complex places $w_\R$ and $w_\C$ of $K^+$
above the same real place of $L^+$ over $v_0$.
Both $v_\R$ and $v_\C$ split completely in $K^+/k$.

Let $S$ be a finite set of places of $k$, including every
archimedean place and every finite place ramified in $K^+/k$,
with $\lvert S\rvert > 3$.
Let $u_\R \in \calo_{L^+}^\times$ and
$u_\C \in \calo_{K^+}^\times$ satisfy
$\rmn_{K^+/L^+}(u_\C) = 1$.
Set $U_S := \Q\otimes_\Z\calo_{K^+, S}^\times$.
For a place $x$ of $K^+$ and an $S$-unit $a$ of $K^+$,
define
\[
  \lambda_x(a) := \sum_{\sigma \in \Gal(K^+/k)}
    \log\verts{a^\sigma}_x\sigma^{-1}
    \in \R[\Gal(K^+/k)]
\]
and extend $\lambda_x$ linearly to $U_S$.
We use $\calr_x(a, \chi)$ for the character sum over
$\Gal(K^+/k)$, so
$\calr_x(a, \chi) = \chi^{-1}(\lambda_x(a))$.
Thus $\calr_x(a, \chi)$ is a complex number, whereas
$\lambda_x(a)$ is an element of the real group ring.
The rank-$2$ Stark regulator is the map
\[
  \begin{tikzcd}[row sep = tiny]
    \bigwedge_{\Q[\Gal(K^+/k)]}^2 U_S
    \arrow[r, "\calr_{w_\R, w_\C}"] & \R[\Gal(K^+/k)] \\
    a\wedge b \arrow[r, mapsto] &
    \lambda_{w_\R}(a)\lambda_{w_\C}(b)
    - \lambda_{w_\R}(b)\lambda_{w_\C}(a).
  \end{tikzcd}
\]

Since $u_\R$ belongs to $L^+$, we have
$\log\verts{u_\R^\sigma}_{w_\C}
= 2\log\verts{u_\R^\sigma}_{w_\R}$.
The equality $\rmn_{K^+/L^+}(u_\C) = 1$ gives
$\log\verts{u_\C^\sigma}_{w_\R}
= -\log\verts{u_\C^\sigma}_{w_\C}$.
For a character $\chi$ of $\Gal(K^+/k)$, evaluating the
regulator matrix at $\chi^{-1}$ therefore gives
\[
  \begin{pmatrix}
    \calr_{w_\R}(u_\R, \chi)
    & -\calr_{w_\C}(u_\C, \chi) \\
    2\calr_{w_\R}(u_\R, \chi)
    & \calr_{w_\C}(u_\C, \chi)
  \end{pmatrix}.
\]
The determinant is therefore
\begin{equation}
  \label{eq:actual-place-determinant}
  \chi^{-1}\paren{\calr_{w_\R, w_\C}(u_\R\wedge u_\C)}
    = 3\calr_{w_\R}(u_\R, \chi)\calr_{w_\C}(u_\C, \chi).
\end{equation}

We can now express the rational rank-$2$ Stark element
for $(K^+/k, S)$ and $(w_\R, w_\C)$ in terms of
$u_\R$ and $u_\C$.
Since $\ell_w = -\lambda_w$, the two signs cancel in
the determinant, so we can use
Equation~\ref{eq:actual-place-determinant} to compute
the regulator in
Definition~\ref{def:rational-rank-two-element}.

\begin{corollary}
  \label{cor:rational-rank-two-element}
  If $u_\R \in \calo_{L^+}^\times$ and
  $u_\C \in \calo_{K^+}^\times$ satisfy
  $\rmn_{K^+/L^+}(u_\C) = 1$ and
  \[
    L_S''(0, \chi)
    = \frac{1}{W_K}\calr_{w_\R}(u_\R, \chi)
    \calr_{w_\C}(u_\C, \chi)
  \]
  for every character $\chi$ of $\Gal(K^+/k)$,
  then
  \[
    e_{S, 2}\calr_{w_\R, w_\C}(u_\R\wedge u_\C)
    = 6W_K\Theta_S^{(2)}(0).
  \]
  In particular, the rational rank-$2$ Stark element
  exists and is given by
  \begin{equation}
    \label{eq:rational-rank-two-element}
    \epsilon_{S, \Q}^{(2)}
      = \frac{1}{6W_K}e_{S, 2}(u_\R\wedge u_\C).
  \end{equation}
\end{corollary}

\begin{proof}
  We combine the hypothesis with
  Equation~\ref{eq:actual-place-determinant} to obtain
  \[
    \chi^{-1}\paren{\calr_{w_\R, w_\C}
    (u_\R\wedge u_\C)}
    = 3W_K L_S''(0, \chi).
  \]
  If $\ord_{s = 0}L_S(s, \chi) = 2$, then
  $\chi^{-1}(6W_K\Theta_S^{(2)}(0)) = 3W_K L_S''(0, \chi)$.
  Multiplication by $e_{S, 2}$ therefore gives the asserted
  equality in the group ring. Dividing by $6W_K$ proves
  that the right side of
  Equation~\ref{eq:rational-rank-two-element} has regulator
  $\Theta_S^{(2)}(0)$.
  Since $\calr_{w_\R, w_\C}$ is injective on
  $e_{S, 2}\bigwedge_{\Q[\Gal(K^+/k)]}^2 U_S$,
  the element $(6W_K)^{-1}e_{S, 2}(u_\R\wedge u_\C)$
  is independent of the choice of $u_\R$ and $u_\C$
  satisfying the hypotheses.
\end{proof}

\subsubsection{The rank-\texorpdfstring{$2$}{2}
Rubin--Stark element}

We now study when the rational rank-$2$ Stark element
belongs to Rubin's lattice.
Recall that $S$ contains the archimedean places of $k$
and every finite place ramified in $K^+/k$, with
$\lvert S\rvert > 3$.
Let $T$ be a finite set of finite places of $k$ disjoint
from $S$, and write $U_{S, T} := U_{K^+, S, T}$ for the
group of $S$-units of $K^+$ congruent to $1$ above $T$.
Assume that $U_{S, T}$ is torsion-free.
We write $\Lambda_{S, T}^{(2)}$ and
$\epsilon_{S, T, \Q}^{(2)}$ for the lattice and rational
Rubin--Stark element of Definitions~\ref{def:Rubin-lattice}
and \ref{def:rank-two-Rubin-Stark-element} for $K^+/k$,
using the ordered places $(w_\R, w_\C)$.

We will prove integrality for the element constructed
from the units $u_\R$ and $u_\C$ of
Proposition~\ref{prop:relative-rank-two-regulator}.
By Corollary~\ref{cor:rational-rank-two-element},
$\epsilon_{S, \Q}^{(2)} =
(6W_K)^{-1}e_{S, 2}(u_\R\wedge u_\C)$.
Since $\calr_{w_\R, w_\C}$ is $\Q[\Gal(K^+/k)]$-linear,
we have
\begin{align*}
  \calr_{w_\R, w_\C}(\epsilon_{S, T, \Q}^{(2)})
    &= \paren{\prod_{v \in T}
      (1 - \rmn(v)[v, K^+/k]^{-1})}
      \calr_{w_\R, w_\C}(\epsilon_{S, \Q}^{(2)}) \\
    &= \paren{\prod_{v \in T}
      (1 - \rmn(v)[v, K^+/k]^{-1})}\Theta_S^{(2)}(0).
\end{align*}
This is the regulator required in
Conjecture~\ref{conj:integral-rank-two},
so the Rubin--Stark conjecture in our setting
reduces to the assertion that
$\epsilon_{S, T, \Q}^{(2)} \in \Lambda_{S, T}^{(2)}$.
We will prove that
$\epsilon_{S, T, \Q}^{(2)} \in
\Z[1/6] \otimes_\Z\Lambda_{S, T}^{(2)}$
and give additional conditions for integrality at $2$ and $3$,
proving Theorem~\ref{thm:rank-two-integral-refinement}.

For $u_\R = \epsilon_{k_0}$ and $u_\C$ in
Equation~\ref{eq:rank-two-unit-definition},
Lemma~\ref{lem:relative-archimedean-units}(ii) will show
that $e_\one(u_\R\wedge u_\C) = 0$.
This is why we impose $\lvert S\rvert > 3$:
when $\lvert S\rvert = 3$, the trivial character has
order of vanishing $2$ and requires a nonzero component
of the Stark element.
We will also use
Proposition~\ref{prop:Rubin-integral-subcases} to treat a
character of order $2$ in
Proposition~\ref{prop:rank-two-integrality}(iii).

\subsubsection{Denominators of the Stark element}
\label{subsec:rank-two-exterior-formulas}

We now return to the construction when $K$ contains $F$.
We use the fields, sets of places, and Stark units of
Proposition~\ref{prop:rank-two-regulator}.
In particular,
$u_\R = \epsilon_{k_0}$ and
$u_\C = \rmn_{K/K^+}(\epsilon_F)^3/\rmn_{K/L^+}(\epsilon_F)$.
Thus the root in Lemma~\ref{lem:rank-two-denominators}
satisfies $\gamma^{W_K} = u_\C$.

We now express $\epsilon_{S, \Q}^{(2)}$ as exterior products
of global units with denominators $6$ and $2W_K$.
We will use the denominator $6$ to prove integrality
at primes $p > 3$, and the denominator $2W_K$ to prove
integrality at $p = 3$ when $3 \nmid W_K$.
To obtain these expressions, we first use the archimedean
absolute values of $u_\R$ and $u_\C$ to remove $e_{S, 2}$
from Equation~\ref{eq:rational-rank-two-element}.

At the archimedean places of $L^+$ and $K$ not above $v_0$,
the Stark units $u_\R$ and $\epsilon_F$, respectively, have
absolute value $1$
(Definition~\ref{def:global-rank-one-Stark-unit}(ii)).
Taking norms in Equation~\ref{eq:rank-two-unit-definition}
shows that $u_\C$ has absolute value $1$ at every
archimedean place of $K^+$ not above $v_0$.
For a number field $K' \supset k_0$, let $U_{v_0}(K')$ be
the group of global units of $K'$ with absolute value $1$
at every archimedean place not above $v_0$.
The logarithms at places above $v_0$ bound the dimensions
of the character components of
$\C\otimes_\Z U_{v_0}(L^+)$ and
$\C\otimes_\Z U_{v_0}(K^+)$.

\begin{lemma}
  \label{lem:relative-archimedean-units}
  The groups $U_{v_0}(L^+)$ and $U_{v_0}(K^+)$ satisfy:
  \begin{enumerate}[(i)]
    \item Every $a \in U_{v_0}(L^+)$ has norm $1$ or $-1$
      over $k_0$. Furthermore,
      \[
        \bigwedge_{\Q[\Gal(K^+/k)]}^2
        (\Q\otimes_\Z U_{v_0}(L^+)) = 0.
      \]
    \item The regulator $\calr_{w_\R, w_\C}$ is injective
      on $\R\otimes_\Z\bigwedge_{\Z[\Gal(K^+/k)]}^2 U_{v_0}(K^+)$.
      Furthermore, the trivial character component of this
      exterior square is zero.
  \end{enumerate}
\end{lemma}

\begin{proof}
  For $K' = L^+$ or $K^+$, Dirichlet's unit theorem
  (as in \cite[Chapter~I, Theorem~7.3]{neukirch}
  for instance)
  says that the $\R$-linear extension of
  \[
    a \longmapsto (\log\lvert a\rvert_w)_{w \mid \infty}
  \]
  identifies $\R\otimes_\Z\calo_{K'}^\times$ with
  the hyperplane $\sum_{w \mid \infty}x_w = 0$.
  On $U_{v_0}(K')$, the coordinates at places not above
  $v_0$ vanish.
  Thus the logarithms at the places above $v_0$ give
  an injection on $\R\otimes_\Z U_{v_0}(K')$.
  For $a \in U_{v_0}(L^+)$, the logarithm
  $\log\lvert a\rvert_w$ is zero whenever $w$ does not lie
  above $v_0$. The places of $L^+$ above $v_0$ form one
  regular $\Gal(K^+/k)$-orbit, and
  $\sum_{w \mid v_0}\log\lvert a\rvert_w = 0$
  by the product formula. Hence each nontrivial character
  component of $\C\otimes_\Z U_{v_0}(L^+)$ has dimension
  at most $1$, and the trivial component is zero.
  Taking exterior squares proves the second assertion
  in part~(i). The norm $\rmn_{L^+/k_0}(a)$ has absolute
  value $1$ at the other real places of $k_0$. The product
  formula gives absolute value $1$ at $v_0$ as well.
  Since $k_0$ is totally real, this norm is $1$ or $-1$.

  The places of $K^+$ above $v_\R$ and $v_\C$ form two
  regular $\Gal(K^+/k)$-orbits. For $a \in U_{v_0}(K^+)$, the
  coordinates $\log\lvert a\rvert_w$ at archimedean places
  $w$ outside these two orbits are zero.
  The maps $\lambda_{w_\R}$ and $\lambda_{w_\C}$ in
  Section~\ref{subsec:rank-two-determinant} record the
  logarithms on these orbits in $\R[\Gal(K^+/k)]$.
  For a character $\chi$ of $\Gal(K^+/k)$, their evaluations under
  $\chi^{-1}$ therefore give an injection
  from $e_{\chi^{-1}}(\C\otimes_\Z U_{v_0}(K^+))$ into
  $\C^2$.
  Taking second exterior powers gives an injection
  into $\bigwedge_\C^2\C^2 = \C$.
  This is the $\chi^{-1}$-component of
  $\calr_{w_\R, w_\C}$, so the regulator is injective.
  The product formula bounds the dimension of the
  trivial component of $\C\otimes_\Z U_{v_0}(K^+)$ by $1$.
  Its second exterior power is therefore zero.
  This proves part~(ii).
\end{proof}

Let $\gamma \in \calo_{K^+}^\times$ be the root constructed
in the proof of Lemma~\ref{lem:rank-two-denominators}.
Since $\gamma^{W_K} = u_\C$, we have
$\gamma \in U_{v_0}(K^+)$.
We use Lemma~\ref{lem:relative-archimedean-units} to express
$\epsilon_{S, \Q}^{(2)}$ with denominators $6$ and $2W_K$.

\begin{corollary}
  \label{cor:rank-two-exterior-formulas}
  Let $\gamma \in \calo_{K^+}^\times$ satisfy
  $\gamma^{W_K} = u_\C$. Then $\epsilon_{S, \Q}^{(2)}$
  has the following expressions in
  $\bigwedge_{\Q[\Gal(K^+/k)]}^2 U_S$:
  \begin{enumerate}[(i)]
    \item $\epsilon_{S, \Q}^{(2)} = \frac{1}{6}u_\R\wedge\gamma$;
    \item $\epsilon_{S, \Q}^{(2)} = \frac{1}{2W_K}u_\R\wedge
      \rmn_{K/K^+}(\epsilon_F)$.
  \end{enumerate}
\end{corollary}

\begin{proof}
  The units $u_\R$ and $u_\C$ belong to $U_{v_0}(K^+)$.
  For a character $\chi$ with
  $\ord_{s = 0}L_S(s, \chi) > 2$, Equation
  \ref{eq:actual-place-determinant} and
  Proposition~\ref{prop:rank-two-regulator} show that the regulator
  of $e_{\chi^{-1}}(u_\R\wedge u_\C)$ vanishes.
  Lemma~\ref{lem:relative-archimedean-units}(ii), after
  extending scalars to $\C$, gives
  $e_{\chi^{-1}}(u_\R\wedge u_\C) = 0$.
  The trivial component of $u_\R\wedge u_\C$ is also zero.
  Hence $e_{S, 2}(u_\R\wedge u_\C) = u_\R\wedge u_\C$.
  Since $\gamma^{W_K} = u_\C$, multilinearity gives
  $u_\R\wedge u_\C = W_K(u_\R\wedge\gamma)$.
  Substituting in
  Equation~\ref{eq:rational-rank-two-element} proves (i).

  Lemma~\ref{lem:relative-archimedean-units}(i) applies to
  $u_\R$ and $\rmn_{K/L^+}(\epsilon_F)$ in $U_{v_0}(L^+)$.
  Hence
  $u_\R\wedge\rmn_{K/L^+}(\epsilon_F) = 0$.
  Equation~\ref{eq:rank-two-unit-definition} and
  multilinearity now give
  $u_\R\wedge u_\C
  = 3u_\R\wedge\rmn_{K/K^+}(\epsilon_F)$.
  Dividing by $6W_K$ proves (ii).
\end{proof}

We now apply these two expressions to prove integrality
in Rubin's lattice.
They settle integrality
at each prime $p > 3$
and at the prime $p = 3$ when $3 \nmid W_K$.
For the remaining cases,
we use the abelian extension
condition for the Stark units
to prove the additional
divisibility at $2$ and $3$.

\begin{proposition}
  \label{prop:rank-two-integrality}
  Consider the fields, sets of places, and Stark units of
  Proposition~\ref{prop:rank-two-regulator},
  and let $d := [K^+ : k]$.
  Let
  \begin{itemize}
    \item $p$ be a prime, and
    \item $T$ be a finite set of finite places of $k$
      disjoint from $S$ such that the group $U_{S, T}$
      of $S$-units of $K^+$ congruent to $1$ at all places
      above $T$ is torsion-free.
  \end{itemize}
  The inclusion
  $\epsilon_{S, T, \Q}^{(2)}
  \in \Z_p\otimes_\Z\Lambda_{S, T}^{(2)}$
  holds in each of the following cases:
  \begin{enumerate}[(i)]
    \item $p > 3$;
    \item $p = 2$ and $d$ is odd;
    \item $p = 3$ and either $3 \nmid W_K$ or $3 \nmid d$.
  \end{enumerate}
\end{proposition}

\begin{proof}
  \emph{Primes $p > 3$, and $p = 3$ when $3 \nmid W_K$.}
  For $p > 3$, the denominator $6$ in
  Corollary~\ref{cor:rank-two-exterior-formulas}(i) is a
  $p$-adic unit, so we have
  $\epsilon_{S, \Q}^{(2)} \in
  \bigwedge_{\Z_p[\Gal(K^+/k)]}^2
  (\Z_p\otimes_\Z\calo_{K^+, S}^\times)$.
  The same inclusion holds for $p = 3$ with $3 \nmid W_K$,
  since the denominator $2W_K$ in
  Corollary~\ref{cor:rank-two-exterior-formulas}(ii)
  is a $3$-adic unit.
  To pass to $S$-units congruent to $1$ above $T$,
  we multiply by the Euler factors at $T$
  (as in Section~\ref{subsec:rank-two-congruences})
  to obtain
  \[
    \epsilon_{S, T, \Q}^{(2)} \in
    \bigcap_{\Z_p[\Gal(K^+/k)]}^2
    (\Z_p \otimes_\Z U_{S, T}).
  \]
  Since $e_{S, 2}\epsilon_{S, T, \Q}^{(2)}
  = \epsilon_{S, T, \Q}^{(2)}$,
  this means that
  $\epsilon_{S, T, \Q}^{(2)}
  \in \Z_p \otimes_\Z \Lambda_{S, T}^{(2)}$.

  \emph{The prime $2$.} Assume that $d$ is odd.
  If $d = 1$, Lemma~\ref{lem:relative-archimedean-units}(ii)
  gives $\epsilon_{S, T, \Q}^{(2)} = 0$, so assume $d > 1$.
  We first prove that $u_\R$ is divisible by $2$ in
  $\Z_2\otimes_\Z(\calo_{L^+}^\times/\set{\pm 1})$.
  This will remove the factor $2$ from the denominator in
  $\epsilon_{S, \Q}^{(2)} = (u_\R\wedge\gamma)/6$.

  Choose $a$ with $a^2 = u_\R$.
  Since $u_\R$ is a Stark unit, $L^+(a)/k_0$ is abelian
  by Definition~\ref{def:global-rank-one-Stark-unit}(iii).
  For $h \in \Gal(L^+/k_0)$, choose a lift
  $\widetilde h \in \Gal(L^+(a)/k_0)$ and set
  $b_h := \widetilde h(a)/a$.
  If $\rho \in \Gal(L^+(a)/L^+)$, then
  $\rho(a) = \pm a$ and $\rho\widetilde h = \widetilde h\rho$.
  Hence $\rho(b_h) = b_h$, so $b_h \in L^+$.
  Moreover, $b_h^2 = h(u_\R)/u_\R$ is a global unit,
  which implies $b_h \in \calo_{L^+}^\times$.
  In $\calo_{L^+}^\times/\set{\pm 1}$, written additively,
  we therefore have $(h - 1)u_\R = 2b_h$.

  We now sum over $h \in \Gal(L^+/k_0)$.
  Since $\rmn_{L^+/k_0}(u_\R) = \pm 1$ by
  Lemma~\ref{lem:relative-archimedean-units}(i),
  its image in the quotient by $\set{\pm 1}$ is zero,
  so we have
  \[
    -d u_\R
      = \sum_{h \in \Gal(L^+/k_0)}(h - 1)u_\R
      = 2\sum_{h \in \Gal(L^+/k_0)}b_h.
  \]
  Since $d$ is a unit in $\Z_2$, this proves
  $u_\R/2 \in \Z_2\otimes_\Z
  (\calo_{L^+}^\times/\set{\pm 1})$.

  To apply Lemma~\ref{lem:congruences-exterior-square}
  at $2$, we remove the torsion $\set{\pm 1}$ by using
  the nontrivial character components.
  Define the projection for the trivial character of $\Gal(K^+/k)$,
  \[
    e_1 := d^{-1}\sum_{h \in \Gal(K^+/k)} h,
  \]
  and $R := (1 - e_1)\Z_2[\Gal(K^+/k)]$.
  The modules
  \begin{align*}
    U &:= (1 - e_1)
    (\Z_2\otimes_\Z\calo_{K^+, S}^\times) \\
    V &:= (1 - e_1)(\Z_2\otimes_\Z U_{S, T})
  \end{align*}
  are free over $\Z_2$.
  Since $u_\R \in L^+$ and
  $\rmn_{L^+/k_0}(u_\R) = \pm 1$, we have
  $e_1(1\otimes u_\R) = 0$ in
  $\Q\otimes_\Z\calo_{K^+, S}^\times$.
  Together with the divisibility by $2$ proved above,
  we conclude that $u_\R/2 \in U$. We can therefore write
  $\epsilon_{S, \Q}^{(2)}
  = \frac{1}{3}(u_\R/2)\wedge(1 - e_1)\gamma
  \in \bigwedge_R^2U$ by
  Corollary~\ref{cor:rank-two-exterior-formulas}(i).
  To compare $U$ and $V$, we reduce modulo the primes
  above $T$ and project by $1 - e_1$:
  \[
    U/V \hookrightarrow (1 - e_1)
    \bigoplus_{v \in T}\bigoplus_{w \mid v}
    \bigl(\Z_2\otimes_\Z(\calo_{K^+}/w)^\times\bigr).
  \]
  The summand for $v$ is isomorphic to
  $R/(1 - \rmn(v)[v, K^+/k]^{-1})$, as in
  Section~\ref{subsec:rank-two-congruences}.
  Lemma~\ref{lem:congruences-exterior-square} shows that
  multiplication by the Euler factors at $T$ sends
  $\epsilon_{S, \Q}^{(2)}$ into $\bigcap_R^2V$.
  Since $e_1\epsilon_{S, T, \Q}^{(2)} = 0$,
  the inclusion
  in $\bigcap_R^2 V$ proves
  $\epsilon_{S, T, \Q}^{(2)}
  \in \Z_2\otimes_\Z\Lambda_{S, T}^{(2)}$.

  \emph{The remaining case at $3$.}
  Assume that $3 \mid W_K$ and $3 \nmid d$.
  Since $3 \nmid d$, the ring $\Z_3[\Gal(K^+/k)]$ is a
  product of rings of integers of unramified extensions
  of $\Q_3$. The module
  $\Z_3\otimes_\Z\calo_{K^+, S}^\times$ is free on each
  factor, so its exterior square has no $\Z_3$-torsion.
  We may therefore regard this exterior square as a
  submodule of the exterior square over $\Q_3$.

  We use $\epsilon_{S, \Q}^{(2)} = (u_\R\wedge\gamma)/6$
  and seek divisibility by $3$.
  Via $\Gal(K/kF) \cong \Gal(K^+/k)$, let
  $\kappa : \Gal(K^+/k) \longrightarrow \set{1, -1}$
  describe the action on $\mu_3$.
  The element
  \[
    e_\kappa := \frac{1}{d}
    \sum_{h \in \Gal(K^+/k)}\kappa(h)h^{-1}
  \]
  belongs to $\Z_3[\Gal(K^+/k)]$.
  We first prove that $(1 - e_\kappa)(u_\R\wedge\gamma)$
  is divisible by $3$ in the exterior square of
  $\Z_3\otimes_\Z\calo_{K^+, S}^\times$ over
  $\Z_3[\Gal(K^+/k)]$.
  After applying the Euler factors at $T$, this will give
  integrality of $(1 - e_\kappa)\epsilon_{S, T, \Q}^{(2)}$.
  We then treat $e_\kappa\epsilon_{S, T, \Q}^{(2)}$
  using Rubin's theorem for quadratic extensions.

  To obtain divisibility by $3$, we recall the construction
  of $\gamma$ in Lemma~\ref{lem:rank-two-denominators}.
  Choose $\alpha$ with $\alpha^{W_K} = \epsilon_F$ and
  a lift $\widetilde\tau \in \Gal(K(\alpha)/F)$ of $\tau$.
  Then $\beta := \alpha/\widetilde\tau(\alpha)$ satisfies
  $\gamma = \rmn_{K/K^+}(\beta^{2 + \tau})$.
  For $h \in \Gal(K/kF)$, choose an integer $\kappa_h$
  such that $h(\xi) = \xi^{\kappa_h}$ for every
  $\xi \in \mu(K)$, and choose a lift
  $\widetilde h \in \Gal(K(\alpha)/F)$ of $h$.
  Every automorphism over $K$ multiplies
  $\widetilde h(\alpha)$ and $\alpha^{\kappa_h}$ by the
  same root of unity. Hence
  $x_h := \widetilde h(\alpha)/\alpha^{\kappa_h}$ lies in $K$.
  The equality $x_h^{W_K} = \epsilon_F^{h - \kappa_h}$
  shows that $x_h \in U_{v_0}(K)$.
  Since the lifts of $h$ and $\tau$ commute, we have
  $\beta^{h - \kappa_h} = x_h^{1 - \tau}$.
  Using $(1 - \tau)(2 + \tau) = 3 - (1 + \tau + \tau^2)$
  and taking norms to $K^+$ gives, in additive notation,
  \[
    (h - \kappa_h)\gamma
      = 3\rmn_{K/K^+}(x_h) - \rmn_{K/L^+}(x_h).
  \]
  The last term lies in $U_{v_0}(L^+)$, so
  $u_\R\wedge\rmn_{K/L^+}(x_h) = 0$ by
  Lemma~\ref{lem:relative-archimedean-units}(i).
  We therefore have
  \[
    (h - \kappa_h)(u_\R\wedge\gamma)
      = 3u_\R\wedge\rmn_{K/K^+}(x_h).
  \]
  Since $\kappa_h \equiv \kappa(h) \pmod 3$,
  $(h - \kappa(h))(u_\R\wedge\gamma)$ is divisible by $3$
  in the exterior square of $\Z_3\otimes_\Z
  \calo_{K^+, S}^\times$ over $\Z_3[\Gal(K^+/k)]$.
  We multiply by $\kappa(h)$ and sum over $h$. Since
  $\kappa(h)^2 = 1$ and
  $\sum_h\kappa(h)h = d e_\kappa$, we find that
  $d(e_\kappa - 1)(u_\R\wedge\gamma)$ is divisible by $3$.
  Since $3 \nmid d$, we conclude that
  $(1 - e_\kappa)\epsilon_{S, \Q}^{(2)}
  = (1 - e_\kappa)(u_\R\wedge\gamma)/6$ belongs to
  $\bigwedge_{\Z_3[\Gal(K^+/k)]}^2
  (\Z_3\otimes_\Z\calo_{K^+, S}^\times)$.
  By Lemma~\ref{lem:congruences-exterior-square},
  multiplying by the Euler factors at $T$ gives an
  element of $\bigcap_{\Z_3[\Gal(K^+/k)]}^2
  (\Z_3\otimes_\Z U_{S, T})$.
  Here $U_{S, T}$ consists of the $S$-units of $K^+$
  congruent to $1$ at every prime above $T$.
  Thus $(1 - e_\kappa)\epsilon_{S, T, \Q}^{(2)}
  \in \Z_3 \otimes_\Z \Lambda_{S, T}^{(2)}$.

  \emph{The component for $\kappa$.}
  It remains to prove
  $e_\kappa\epsilon_{S, T, \Q}^{(2)}
  \in \Z_3\otimes_\Z\Lambda_{S, T}^{(2)}$.
  If $\kappa$ is trivial, this component is zero by
  Lemma~\ref{lem:relative-archimedean-units}(ii), and the
  inclusion already proved for
  $(1 - e_\kappa)\epsilon_{S, T, \Q}^{(2)}$ finishes the proof.
  Suppose that $\kappa$ is nontrivial.
  Let $H := \ker\kappa$ and $M := (K^+)^H$, so $M/k$ is
  quadratic and $\lvert H \rvert = d/2$ is coprime to $3$.
  Let $U_{M, S, T}$ be the group of $S$-units of $M$
  congruent to $1$ at all primes above $T$.
  With the places $\restr{w_\R}{M}$ and $\restr{w_\C}{M}$,
  write $\epsilon_{M, S, T, \Q}^{(2)}$ for the rational
  rank-$2$ Rubin--Stark element for $(M/k, S, T)$.
  Rubin \cite[Theorem~3.5 and the following remark]{rubin-1996}
  proves that
  \[
    \epsilon_{M, S, T, \Q}^{(2)} \in
    \bigwedge_{\Z[\Gal(M/k)]}^2U_{M, S, T},
  \]
  where the exterior square is viewed in its rational span.

  To compare $\epsilon_{M, S, T, \Q}^{(2)}$ with
  $\epsilon_{S, T, \Q}^{(2)}$, extend $\rmn_{K^+/M}$
  $\Q$-linearly to the rational spaces of $S$-units.
  Define $U_{M, S} := \Q\otimes_\Z\calo_{M, S}^\times$.
  The norm induces the map
  \[
    \begin{tikzcd}[row sep = tiny]
      \bigwedge_{\Q[\Gal(K^+/k)]}^2U_S
        \arrow[r, "\nu_M"] &
        \bigwedge_{\Q[\Gal(M/k)]}^2U_{M, S} \\
      a\wedge b \arrow[r, mapsto] &
        \rmn_{K^+/M}(a)\wedge\rmn_{K^+/M}(b).
    \end{tikzcd}
  \]
  On group rings, let $\pi_M$ be the $\R$-linear map
  \[
    \begin{tikzcd}[row sep = tiny]
      \R[\Gal(K^+/k)] \arrow[r, "\pi_M"] &
      \R[\Gal(M/k)] \\
      \sigma \arrow[r, mapsto] & \restr{\sigma}{M}.
    \end{tikzcd}
  \]
  The places $w_\R$ and $w_\C$ have local degree $1$
  over their restrictions to $M$, so, for either place $w$
  and $a \in (K^+)^\times$,
  \[
    \log\verts{\rmn_{K^+/M}(a)}_{\restr{w}{M}}
      = \sum_{h \in H}\log\verts{a^h}_w.
  \]
  Substituting in the sums defining $\lambda_{w_\R}$
  and $\lambda_{w_\C}$ and taking the determinant gives,
  for every $x$ in the source of $\nu_M$,
  \[
    \calr_{\restr{w_\R}{M}, \restr{w_\C}{M}}(\nu_M(x))
      = \pi_M\paren{\calr_{w_\R, w_\C}(x)}.
  \]

  We now apply this regulator identity to
  $\epsilon_{S, T, \Q}^{(2)}$.
  For a character $\psi$ of $\Gal(M/k)$, its inflation
  to $\Gal(K^+/k)$ is
  $\sigma \longmapsto \psi(\restr{\sigma}{M})$.
  The inflated character has the same $L$-function with
  the Euler factors at $S$ omitted as $\psi$.
  Moreover, $\pi_M$ sends each factor
  $1 - \rmn(v)[v, K^+/k]^{-1}$ for $v \in T$ to
  $1 - \rmn(v)[v, M/k]^{-1}$.
  Thus $\nu_M(\epsilon_{S, T, \Q}^{(2)})$ and
  $\epsilon_{M, S, T, \Q}^{(2)}$ have the same regulator
  and vanish on every character component whose
  $L$-function has order other than $2$.
  By uniqueness of the rational Rubin--Stark element for
  $(M/k, S, T)$, we have
  \[
    \nu_M(\epsilon_{S, T, \Q}^{(2)})
      = \epsilon_{M, S, T, \Q}^{(2)}.
  \]

  To recover $e_\kappa\epsilon_{S, T, \Q}^{(2)}$
  integrally from $\epsilon_{M, S, T, \Q}^{(2)}$,
  define $e_H := \lvert H \rvert^{-1} \sum_{h \in H} h$.
  Since $\lvert H\rvert$ is a unit in $\Z_3$, the
  idempotent $e_H$ belongs to $\Z_3[\Gal(K^+/k)]$ and
  \[
    e_H(\Z_3\otimes_\Z U_{S, T})
      = \Z_3\otimes_\Z U_{M, S, T}.
  \]
  Here we include units of $M$ in $K^+$.
  Restriction identifies $e_H\Z_3[\Gal(K^+/k)]$
  with $\Z_3[\Gal(M/k)]$.
  Let $j_M$ be the induced inclusion of the second
  exterior power of $\Z_3\otimes_\Z U_{M, S, T}$ into
  that of $\Z_3\otimes_\Z U_{S, T}$, and extend $j_M$
  $\Q_3$-linearly.
  The norm is multiplication by $\lvert H\rvert$ on
  units of $M$, so $j_M\nu_M$ is multiplication by
  $\lvert H\rvert^2$ on the $H$-invariant exterior square.
  Since $e_He_\kappa = e_\kappa$, we obtain
  \[
    e_\kappa\epsilon_{S, T, \Q}^{(2)}
    = \frac{1}{\lvert H\rvert^2}
    e_\kappa j_M(\epsilon_{M, S, T, \Q}^{(2)}).
  \]
  Rubin's integrality assertion, the integrality of
  $e_\kappa$, and $3 \nmid\lvert H\rvert$ place the
  right side in
  $\bigwedge_{\Z_3[\Gal(K^+/k)]}^2
  (\Z_3\otimes_\Z U_{S, T})$.
  Its image under the canonical map to the exterior
  bidual has the same rational image.
  Since $e_{S, 2}e_\kappa\epsilon_{S, T, \Q}^{(2)}
  = e_\kappa\epsilon_{S, T, \Q}^{(2)}$, we conclude that
  $e_\kappa\epsilon_{S, T, \Q}^{(2)}
  \in \Z_3\otimes_\Z\Lambda_{S, T}^{(2)}$.
  Combining this with the inclusion for
  $(1 - e_\kappa)\epsilon_{S, T, \Q}^{(2)}$ proves
  $\epsilon_{S, T, \Q}^{(2)}
  \in \Z_3\otimes_\Z\Lambda_{S, T}^{(2)}$.
\end{proof}

\subsubsection{The relative Rubin--Stark theorem}
\label{subsec:relative-rank-two-integral-descent}

We now prove
the rank-$2$ Rubin--Stark conjecture
over $\Z[1/6]$,
without assuming that $K$ contains $F$.

We work with $K' := K^+F$ and $L' := L^+F$ so that
$K'$ contains $F$ and $(K')^{\abrac{c}} = K^+$.
We can then use Proposition~\ref{prop:rank-two-regulator}
to construct the rational rank-$2$ Stark element for
$K^+/k$ from rank-$1$ Stark units for $K'/F$ and $L^+/k_0$,
and Proposition~\ref{prop:rank-two-integrality} to prove
its integrality in Rubin's lattice.
By uniqueness, this is the element already constructed
from the units in
Proposition~\ref{prop:relative-rank-two-regulator}.

\begin{proposition}
  \label{thm:relative-rank-two-integrality}
  Assume the hypotheses of
  Proposition~\ref{prop:relative-rank-two-regulator},
  and let
  \begin{itemize}
    \item $u_\R \in \calo_{L^+}^\times$ and
      $u_\C \in \calo_{K^+}^\times$ satisfy its conclusion;
    \item $\epsilon_{S, \Q}^{(2)} =
      (6W_K)^{-1}e_{S, 2}(u_\R\wedge u_\C)$ be the rational
      Stark element given by
      Corollary~\ref{cor:rational-rank-two-element}; and
    \item $T$ be a finite set of finite places of $k$
      disjoint from $S$ such that the group $U_{S, T}$ of
      $S$-units of $K^+$ congruent to $1$ above $T$ is
      torsion-free.
  \end{itemize}
  The rational rank-$2$ Rubin--Stark element
  $\epsilon_{S, T, \Q}^{(2)}$ is in
  $\Z[1/6]\otimes_\Z\Lambda_{S, T}^{(2)}$.
  Furthermore,
  $\epsilon_{S, T, \Q}^{(2)} \in
  \Z_2\otimes_\Z\Lambda_{S, T}^{(2)}$
  if $[K^+ : k]$ is odd, and
  $\epsilon_{S, T, \Q}^{(2)} \in
  \Z_3\otimes_\Z\Lambda_{S, T}^{(2)}$
  if either $3 \nmid [K^+ : k]$ or $3 \nmid W_{K^+F}$.
\end{proposition}

\begin{proof}
  Since $K' := K^+F$ and $L' := L^+F$,
  we have that $K' = kL'$.
  We first check the hypotheses of
  Proposition~\ref{prop:rank-two-integrality} for $K'/k$.
  By Lemma~\ref{lem:resolvent-ramification}, $S_{k_0}$
  contains every finite prime ramified in $KF/k_0$.
  Since $K' \subset KF$, the sets $S$, $S_{k_0}$, and $S_F$
  contain the finite places ramified in $K^+/k$,
  $L^+/k_0$, and $K'/F$, respectively.
  Equation~\ref{eq:local-factorization} holds for $K'/k$
  by Corollary~\ref{cor:rationally-saturated-factorization}.
  We descend the assumed Stark unit for $(KF/F, S_F)$
  to a Stark unit for $(K'/F, S_F)$ by
  Lemma~\ref{lem:Stark-descent}.
  Together with the assumed Stark unit for $L^+/k_0$,
  this verifies the hypotheses of
  Proposition~\ref{prop:rank-two-integrality}.

  By Corollary~\ref{cor:rational-rank-two-element},
  the units constructed for $K'/k$ give the rational
  rank-$2$ Stark element for $K^+/k$ with ordered
  places $(w_\R, w_\C)$.
  By uniqueness, this element is
  $\epsilon_{S, \Q}^{(2)}$.
  Multiplying by the Euler factors at $T$ therefore
  gives $\epsilon_{S, T, \Q}^{(2)}$.

  Since $(K')^+ = K^+$, the unit group and Rubin's
  lattice in
  Proposition~\ref{prop:rank-two-integrality}
  are $U_{S, T}$ and $\Lambda_{S, T}^{(2)}$.
  We apply parts~(i)--(iii) of that proposition with
  $d = [K^+ : k]$ and $W_{K'} = W_{K^+F}$
  to obtain inclusions over
  $\Z[1/6]$, $\Z_2$, and $\Z_3$.
\end{proof}

We can now prove
Theorem~\ref{thm:rank-two-integral-refinement}.

\begin{proof}[Proof of
  Theorem~\ref{thm:rank-two-integral-refinement}]
  Let $T$ be any finite set of finite places of $k$
  disjoint from $S$ such that the group $U_{S, T}$ of
  $S$-units of $K^+$ congruent to $1$ above $T$ is
  torsion-free.
  We construct the units by
  Proposition~\ref{prop:relative-rank-two-regulator} and
  use Proposition~\ref{thm:relative-rank-two-integrality}
  to prove
  $\epsilon_{S, T, \Q}^{(2)} \in
  \Z[1/6]\otimes_\Z\Lambda_{S, T}^{(2)}$.
  Since $[K^+ : k]$ is odd and either
  $3 \nmid [K^+ : k]$ or $3 \nmid W_{K^+F}$, we also have
  $\epsilon_{S, T, \Q}^{(2)}$ in
  $\Z_2\otimes_\Z\Lambda_{S, T}^{(2)}$ and
  $\Z_3\otimes_\Z\Lambda_{S, T}^{(2)}$.
  Since $\Lambda_{S, T}^{(2)}$ is a $\Z$-lattice,
  these inclusions imply
  $\epsilon_{S, T, \Q}^{(2)} \in \Lambda_{S, T}^{(2)}$.
\end{proof}

\subsection{Specialization to \texorpdfstring{$k_0 = \Q$}{k0 = Q}}
\label{subsec:complex-cubic-rank-two}

When $k_0 = \Q$, the field $k$ is complex cubic.
We use the classical Stark theorems recalled in
Section~\ref{subsec:classical-Stark} to construct the units
required by Proposition~\ref{prop:relative-rank-two-regulator}.
We first construct the units in
Theorem~\ref{thm:intro-rank-two} when $S$ consists
of all places above a set of rational places containing
at least two finite primes,
and then include the additional finite places
allowed in Theorem~\ref{thm:intro-rank-two}.
We use Theorem~\ref{thm:intro-relative}(ii) and
Theorem~\ref{thm:rank-two-integral-refinement} for the
Rubin--Stark assertions.

\begin{proof}[Proof of Theorem~\ref{thm:intro-rank-two}]
  We first construct units for a subset $S_0$ of $S$
  consisting of all places above a set of rational places.
  Let $\ell_1$ and $\ell_2$ be the distinct rational primes
  in the hypothesis, and let $S_\Q$ consist of $\infty$,
  every rational prime ramified in $K/\Q$, and
  $\ell_1, \ell_2$.
  Let $S_0$ and $S_F$ consist of the places above $S_\Q$
  in $k$ and in its quadratic resolvent $F$, respectively.
  By hypothesis, $S_0 \subseteq S$.
  By Theorems~\ref{thm:rational-Stark} and
  \ref{thm:imaginary-quadratic-Stark}, we have Stark units
  for $(L^+/\Q, S_\Q)$ and $(KF/F, S_F)$.
  Applying Proposition~\ref{prop:relative-rank-two-regulator}
  with $k_0 = \Q$ and $S = S_0$, we obtain
  $u_{\R, 0} \in \calo_{L^+}^\times$ and
  $u_{\C, 0} \in \calo_{K^+}^\times$ with the required
  regulator for $S_0$ and norm
  $\rmn_{K^+/L^+}(u_{\C, 0}) = 1$.

  Let $T$ be a finite set of finite places of $k$
  disjoint from $S$ such that $U_{S, T}$ is torsion-free.
  By Rubin \cite[Proposition~2.5]{rubin-1996}, we also have
  the integral rank-$1$ Rubin--Stark conjectures for
  $(L^+/\Q, S_\Q)$ and $(K^+F/F, S_F)$.
  We apply Theorem~\ref{thm:intro-relative}(ii) to $K^+/k$;
  since $W_{K^+} = 2$,
  it implies integrality of a rational Stark element
  after inverting $6$.
  By uniqueness of the rational Stark element,
  we therefore have that
  $\epsilon_{S_0, T, \Q}^{(2)} \in
  \Z[1/6]\otimes_\Z\Lambda_{S_0, T}^{(2)}$.

  Set $K' := K^+F$ and $d := [K^+ : k]$.
  If $d$ is odd and $F \neq \Q(\sqrt{-3})$, then
  $3 \nmid W_{K'}$.
  Indeed, for $K' = k(L^+F)$,
  we have $\mu(K') = \mu(L^+F)$
  (like in Section~\ref{subsec:Stark-hypotheses}).
  The group $\Gal(L^+F/\Q) \cong
  \Gal(F/\Q)\times\Gal(L^+/\Q)$ has a unique quotient
  of order $2$, since $[L^+ : \Q] = d$ is odd.
  Hence $F$ is the unique quadratic subfield of $L^+F$.
  If $3$ divided $W_{K'}$, then $F$ would be
  $\Q(\sqrt{-3})$.
  We can therefore apply
  Theorem~\ref{thm:rank-two-integral-refinement} whenever
  $d$ is odd and either $3 \nmid d$ or
  $F \neq \Q(\sqrt{-3})$.
  By uniqueness, the rational Stark element is unchanged,
  so $\epsilon_{S_0, T, \Q}^{(2)} \in
  \Lambda_{S_0, T}^{(2)}$.
  The condition $F \neq \Q(\sqrt{-3})$ is equivalent to
  $\disc(k) \notin -3(\Q^\times)^2$.

  We now include the places in $S\setminus S_0$.
  Define
  \[
    a := \prod_{v \in S\setminus S_0}
      (1 - [v, K^+/k]^{-1}) \in \Z[\Gal(K^+/k)],
  \]
  and set $u_\R := u_{\R, 0}$ and $u_\C := u_{\C, 0}^a$.
  Every place in $S\setminus S_0$ is finite and unramified
  in $K^+/k$.
  Each $L_{S_0}(s, \chi)$ vanishes to order at least $2$
  at $0$, because both archimedean places split completely
  and $\lvert S_0\rvert > 3$.
  Removing the Euler factors at $S\setminus S_0$ therefore gives
  \[
    L_S''(0, \chi) = \chi^{-1}(a)L_{S_0}''(0, \chi).
  \]
  Reindexing the character sum for $u_{\C, 0}^a$ multiplies
  it by $\chi^{-1}(a)$, so $u_\R$ and $u_\C$ satisfy the
  asserted equality for $L_S''(0, \chi)$.
  The norm commutes with the action of $\Gal(K^+/k)$,
  so $\rmn_{K^+/L^+}(u_\C) = 1$ as well.

  The rational Stark elements satisfy
  $\epsilon_{S, T, \Q}^{(2)}
  = a \epsilon_{S_0, T, \Q}^{(2)}$.
  Rubin \cite[Proposition~3.6 and its proof]{rubin-1996}
  shows that multiplication by $a$ maps
  $\Lambda_{S_0, T}^{(2)}$ into $\Lambda_{S, T}^{(2)}$.
  Thus the same inclusions hold for
  $\epsilon_{S, T, \Q}^{(2)}$.
  Since $T$ was arbitrary, this proves the Rubin--Stark
  assertions.
\end{proof}

We record the integrality of the Rubin--Stark element
constructed from $u_\R$ and $u_\C$ in the proof of
Theorem~\ref{thm:intro-rank-two}.

\begin{corollary}
  \label{cor:rank-two-integral-odd-degree}
  Let $K/k$, $S$, $u_\R$, $u_\C$, $w_\R$, and $w_\C$ be as in
  Theorem~\ref{thm:intro-rank-two},
  and let $T$ be a finite set of finite places of $k$
  disjoint from $S$,
  such that the group $U_{S, T}$
  of $S$-units of $K^+$ congruent to $1$ at all
  places above $T$ is torsion free.
  The element $\epsilon_{S, T, \Q}^{(2)}$
  constructed from $u_\R$ and $u_\C$ satisfies
  \[
    \epsilon_{S, T, \Q}^{(2)}
    \in \Z[1/6] \otimes_\Z\Lambda_{S, T}^{(2)}.
  \]
  Furthermore, if $[K^+ : k]$ is odd and either
  $3 \nmid [K^+ : k]$ or $F \neq \Q(\sqrt{-3})$, then
  $\epsilon_{S, T, \Q}^{(2)} \in \Lambda_{S, T}^{(2)}$.
\end{corollary}

\begin{proof}
  In the proof of Theorem~\ref{thm:intro-rank-two}, we proved
  both inclusions for every finite set $T$ of finite places
  of $k$ disjoint from $S$ such that the group $U_{S, T}$
  of $S$-units of $K^+$ congruent to $1$ at every place
  above $T$ is torsion-free.
\end{proof}

\begin{remark}
  When $K$ contains $F$, we can relax the hypotheses on
  the sets of places in
  Corollary~\ref{cor:rank-two-integral-odd-degree}.
  It suffices that $S$, $S_\Q$, and $S_F$ satisfy
  Section~\ref{subsec:decomposition-partial}(ii),
  with $\lvert S \rvert > 3$.
  Let $T$ be a finite set of finite places of $k$ disjoint
  from $S$, with $U_{S, T}$ torsion-free.
  Then $\epsilon_{S, T, \Q}^{(2)} \in
  \Z[1/6]\otimes_\Z\Lambda_{S, T}^{(2)}$.
  If $[K^+ : k]$ is odd and either $3 \nmid [K^+ : k]$ or
  $F \neq \Q(\sqrt{-3})$, then
  $\epsilon_{S, T, \Q}^{(2)} \in \Lambda_{S, T}^{(2)}$.
  Indeed, Theorems~\ref{thm:rational-Stark} and
  \ref{thm:imaginary-quadratic-Stark} give Stark units for
  $L^+/\Q$ and $K/F$, so we can apply
  Proposition~\ref{prop:rank-two-integrality}.
  To check the assertion at $3$ when $[K^+ : k]$ is odd,
  note that $L = L^+F$ has exactly one quadratic subfield,
  namely $F$: the group
  $\Gal(F/\Q)\times\Gal(L^+/\Q)$ has a unique quotient of
  order $2$.
  Since $\mu(K) = \mu(L)$, the divisibility
  $3 \mid W_K$ would imply $\Q(\sqrt{-3}) \subset L$ and
  hence $F = \Q(\sqrt{-3})$. Thus $3 \nmid W_K$ when
  $F \neq \Q(\sqrt{-3})$.
\end{remark}

\subsection{Explicit class field theory over complex cubic fields}
\label{subsec:rank-two-explicit-cft}

We now ask when the units of
Theorem~\ref{thm:intro-rank-two} generate
$L^+/\Q$ and $K^+/k$.
In particular, we give instances of the approach
of Stark \cite[p.~63]{stark-1975} and
Tate \cite[p.~95]{tate-stark-1984} to explicit class field theory
when $K^+/k$ is cyclic.

\begin{corollary}
  \label{cor:rank-two-field-generation}
  Let $K/k$ and $S$ be as in
  Theorem~\ref{thm:intro-rank-two}.
  If $K^+/k$ is cyclic
  and no finite place in $S$ splits completely in $K^+/k$,
  then $L^+ = \Q(u_\R)$ and $K^+ = k(u_\C)$, where
  $u_\R \in \calo_{L^+}^\times$ and
  $u_\C \in \calo_{K^+}^\times$ are the units constructed
  in the proof of Theorem~\ref{thm:intro-rank-two}.
\end{corollary}

\begin{proof}
  Let $\chi$ be a faithful character of $\Gal(K^+/k)$.
  By Tate \cite[Chapter~I, Proposition~3.4]{tate-stark-1984},
  the order of $L_S(s, \chi)$ at $0$ is $2$.
  Thus $L_S''(0, \chi) \neq 0$.
  The sums for $u_\R$ and $u_\C$ in
  Theorem~\ref{thm:intro-rank-two} are both nonzero because
  their product is $W_KL_S''(0, \chi)$.
  If $h \in \Gal(K^+/k)$ fixes $u_\R$ or $u_\C$,
  reindexing the sum for that unit gives
  $\chi(h) = 1$. Thus $h = 1$.
  Restriction identifies
  $\Gal(K^+/k)$ with $\Gal(L^+/\Q)$,
  proving both assertions.
\end{proof}

\begin{remark}
  In Corollary~\ref{cor:rank-two-field-generation},
  we assume that no finite place in $S$ splits completely
  in $K^+/k$; otherwise $L_S''(0, \chi) = 0$
  for every faithful character $\chi$.
\end{remark}


\section{Applications to elliptic gamma functions}
\label{sec:applications}

In this section, we apply our study of Stark units over
complex cubic fields to the conjecture of
Bergeron--Charollois--Garc\'{i}a
\cite[Section~2.3]{bergeron-charollois-garcia}
on elliptic gamma functions.
In Section~\ref{subsec:elliptic-gamma}, we recall their
setting and state the independence of the elliptic
gamma values from the auxiliary element $h$.
We then compare their analogue of the second Kronecker
limit formula with the regulator of a Stark unit for
a narrow ray class field $K/k$.
For cyclotomic $K/k$, we obtain this unit from
Theorem~\ref{thm:intro-rank-one} and prove
Corollary~\ref{cor:intro-elliptic-gamma}.
We also show that the remaining choice in the conjecture
is a single root of unity, in
Proposition~\ref{prop:elliptic-gamma-phase}.
In Section~\ref{sec:elliptic-gamma-independence},
we prove the independence of
elliptic gamma values attached to
arbitrary totally imaginary narrow ray class fields
from $h$ in full generality.

\subsection{Elliptic gamma functions
and the Bergeron--Charollois--Garc\'{i}a conjecture}
\label{subsec:elliptic-gamma}

Stark \cite{stark-1980}
showed that Stark units over imaginary quadratic fields
are classical elliptic units,
which are obtained from special values
of theta functions by CM theory.
Bergeron--Charollois--Garc\'{i}a
\cite[Sections~1 and~2]{bergeron-charollois-garcia} propose
that the elliptic gamma function, defined in
Equation~\ref{eq:elliptic-gamma-function},
plays a similar role for complex cubic fields.
Felder--Varchenko
\cite[Theorem~4.1 and Section~7]{felder-varchenko} prove
its transformation laws involving $\SL_3(\Z)$, analogously
to the $\SL_2(\Z)$ transformation laws for theta functions.

We begin with the arithmetic choices that determine
the elliptic gamma values in the conjecture.

\subsubsection*{The elliptic gamma values}
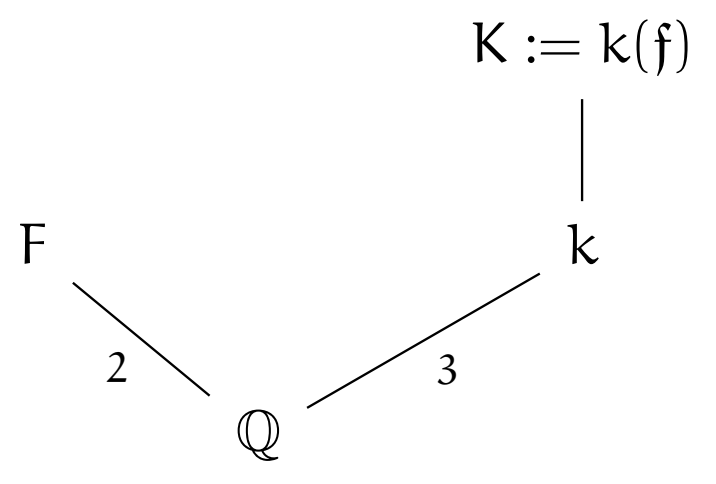
\begin{figure}[ht]
  \[
    \begin{tikzcd}
      & & K := k(\mff) \arrow[d, dash] \\
      F \arrow[dr, dash, "2"'] & &
      k \arrow[dl, dash, "3"] \\
      & \Q &
    \end{tikzcd}
  \]
  \caption{The field extensions in
  \cite[Section~2.3]{bergeron-charollois-garcia}.}
  \label{figure:elliptic-gamma-fields}
\end{figure}

Let $k$ be a complex cubic field
with quadratic resolvent $F$,
and let $K := k(\mff)$ be a totally imaginary narrow
ray class field for an ideal $\mff \subset \calo_k$.
We will attach an elliptic gamma value to each fractional
ideal $\scl := \mff\mfb^{-1}$,
where $\mfb \subset \calo_k$ is coprime to $\mff$.
The definition also uses a nonzero ideal
$\mfa \subset \calo_k$ coprime to $6\rmn(\mff)$
whose norm $\rmn(\mfa)$ is an odd rational prime.
To compare these values with units in $K$, fix an embedding
$\iota : K \longrightarrow \C$ with complex restriction
$x$ to $k$, and let $v$ be its place.
We summarize the field extensions in
Figure~\ref{figure:elliptic-gamma-fields}.

The value $\Gamma_{\mfa, h}(\scl)$ is obtained by choosing
an admissible $h$, complex parameters $\tau_1$ and $\tau_2$
for $\Gamma(z, \tau_1, \tau_2)$, and a finite set of lattice
representatives. We recall these choices from
Bergeron--Charollois--Garc\'{i}a
\cite[Sections~2.1--2.2]{bergeron-charollois-garcia}.

\emph{The admissible element.}
Let $q$ be the least positive integer in $\mff \cap \Z$.
Choose $h = \lambda/q$ with $\lambda \in \scl$ such that
$h - 1 \in \scl$ and $\lambda/\rmn(\mfa)$ belongs to
$\mfa^{-1}\scl$ and generates $\mfa^{-1}\scl/\scl$.
In particular, $\mfa^{-1}\scl = \scl + \Z\lambda/\rmn(\mfa)$.

\emph{The parameters $\tau_1$ and $\tau_2$.}
We choose the second and third arguments of
$\Gamma(z, \tau_1, \tau_2)$ from $\scl$ and a unit of $k$.
Let $x_\R$ be the real embedding of $k$, and orient
$\R\otimes_\Q k$ by $(x_\R, \re(x), \Im(x))$.
Let $\eta$ generate the group of units $u \in \calo_k^\times$
with $u - 1 \in \mff$ and $x_\R(u) > 0$, choosing the
generator with $0 < x_\R(\eta) < 1$.
Choose a primitive $\Z$-linear form
$a : \scl \longrightarrow \Z$, meaning that $a(\scl) = \Z$,
with $a(\lambda) = a(\eta^{-1}\lambda) = 0$.
Define $b(t) := a(\eta^{-1}t)$ for $t \in \scl$.
Extend $a$ and $b$ $\Q$-linearly to $k$. Their kernels
are two-dimensional $\Q$-subspaces of $k$ with
intersection $\Q\lambda$.

To obtain $\tau_1$ and $\tau_2$, choose $\alpha \in \scl$
in $\ker b$ and $\beta \in \scl$ in $\ker a$, with
$a(\alpha) > 0$ and $b(\beta) > 0$.
Let $\gamma$ generate $\scl \cap \Q\lambda$,
with $(\alpha, \beta, \gamma)$ positively oriented.
Choose the sign of $a$ and these elements so that
$\tau_1 := -x(\alpha)/x(\gamma)$ and
$\tau_2 := x(\beta)/x(\gamma)$ have positive imaginary parts.
The existence of these choices follows from
Bergeron--Charollois--Garc\'{i}a
\cite[Lemma~2.1 and Equation~2.6]{bergeron-charollois-garcia}.

\emph{The factors.}
The elements $\alpha$, $\beta$, and $\gamma$ form a
$\Q$-basis of $k$, so
$\scl/(\Z\alpha + \Z\beta + \Z\gamma)$ is finite.
Choose a set $\cald$ of representatives for this quotient
satisfying $0 \leq a(\delta) < a(\alpha)$ and
$0 \leq b(\delta) < b(\beta)$ for every $\delta \in \cald$.
These representatives index the factors in the following
definition.

\begin{definition}
  \label{def:arithmetic-elliptic-gamma}
  For these choices, define
  $z_\delta := x(h + \delta - \alpha)/x(\gamma)$ and
  \[
    \Gamma_{\mfa, h}(\scl)
      := \prod_{\delta \in \cald}
      \frac{\Gamma(z_\delta, \tau_1, \tau_2)^{\rmn(\mfa)}}
      {\Gamma(\rmn(\mfa)z_\delta, \rmn(\mfa)\tau_1,
      \rmn(\mfa)\tau_2)}.
  \]
\end{definition}

We check that Definition~\ref{def:arithmetic-elliptic-gamma}
agrees with Bergeron--Charollois--Garc\'{i}a
\cite[Definition~2.3]{bergeron-charollois-garcia},
so that we can apply their analogue of the
second Kronecker limit
formula in \cite[Theorem~3.2]{bergeron-charollois-garcia}.
For the function $\Gamma_{a, b}$ defined in
\cite[Equation~2.7]{bergeron-charollois-garcia},
we express the factors with both parameters in the
upper half-plane by using
$\Gamma(z, -\tau_1, \tau_2)
= \Gamma(z + \tau_1, \tau_1, \tau_2)^{-1}$.
Since $\mfa^{-1}\scl = \scl + \Z\gamma/\rmn(\mfa)$ by
\cite[Lemma~2.2, proof]{bergeron-charollois-garcia},
we can use the same representatives $\delta \in \cald$
for $\mfa^{-1}\scl$, replacing $\gamma$ by
$\gamma/\rmn(\mfa)$. We obtain
\begin{align}
  \label{eq:elliptic-gamma-lattice-products}
  \Gamma_{a, b}(x(h), x; \scl)
  &= \prod_{\delta \in \cald}
  \Gamma(z_\delta, \tau_1, \tau_2)^{-1}, \\
  \Gamma_{a, b}(x(h), x; \mfa^{-1}\scl)
  &= \prod_{\delta \in \cald}
  \Gamma(\rmn(\mfa)z_\delta, \rmn(\mfa)\tau_1,
  \rmn(\mfa)\tau_2)^{-1}. \nonumber
\end{align}
Bergeron--Charollois--Garc\'{i}a
\cite[Equation~2.13 and Definition~2.3]
{bergeron-charollois-garcia} define
$\Gamma_{\mfa, h}(\scl)$ as the quotient of
$\Gamma_{a, b}(x(h), x; \mfa^{-1}\scl)$ by
$\Gamma_{a, b}(x(h), x; \scl)^{\rmn(\mfa)}$.
Using Equation~\ref{eq:elliptic-gamma-lattice-products}
to evaluate the numerator and denominator, we recover
Definition~\ref{def:arithmetic-elliptic-gamma}.

Bergeron--Charollois--Garc\'{i}a
\cite[Sections~2.1--2.2]{bergeron-charollois-garcia}
show that $\Gamma_{\mfa, h}(\scl)$ is finite, nonzero,
and independent of $\alpha$, $\beta$, and $\cald$.
The following conjecture predicts that the values
$\Gamma_{\mfa, h}(\scl)$ are images of units in $K$
and describes the Galois conjugates of those units.

\begin{conjecture}[{Bergeron--Charollois--Garc\'{i}a
  \cite[Section~2.3]{bergeron-charollois-garcia}}]
  \label{conj:elliptic-gamma}
  Let
  \begin{itemize}
    \item $k$ be a complex cubic field;
    \item $\mff$ be a nonzero integral ideal of $k$ such that
      $K := k(\mff)$ is totally imaginary;
    \item $\mfa$ be an ideal of $k$ coprime to $6\rmn(\mff)$
      whose norm $\rmn(\mfa)$ is an odd rational prime; and
    \item $\iota : K \longrightarrow \C$
      be a fixed embedding that extends
      the complex embedding $x : k \longrightarrow \C$ used to
      define $\Gamma_{\mfa, h}(\scl)$.
  \end{itemize}
  Then there are units $u_{\scl, \mfa}' \in \calo_K^\times$,
  indexed by $\scl = \mff\mfb^{-1}$ for integral ideals
  $\mfb$ of $k$ coprime to $\mff$, with the following properties:
  \begin{enumerate}[(i)]
    \item The value $\Gamma_{\mfa, h}(\scl)$ is independent of
      the choice of admissible $h$ and satisfies
      \[
        \Gamma_{\mfa, h}(\scl) = \iota(u_{\scl, \mfa}')
      \]
      for every admissible $h$ for $\scl$.
    \item We have $\lvert u_{\scl, \mfa}'\rvert_{w'} = 1$
      at every place $w'$ of $K$ above the real place of $k$.
    \item For every integral ideal $\mfc$ of $k$ coprime to
      $\mff$, we have
      $\sigma_\mfc(u_{\scl, \mfa}')
      = u_{\mfc^{-1}\scl, \mfa}'$.
  \end{enumerate}
\end{conjecture}

First, we prove the independence assertion in
Conjecture~\ref{conj:elliptic-gamma}(i) in general.

\begin{proposition}
  \label{prop:elliptic-gamma-independence}
  Let
  \begin{itemize}
    \item $k$ be a complex cubic field;
    \item $\mff$ be a nonzero integral ideal of $k$ such that
      $K := k(\mff)$ is totally imaginary;
    \item $\mfa$ be an ideal of $k$ coprime to $6\rmn(\mff)$
      whose norm $\rmn(\mfa)$ is an odd rational prime;
    \item $\mfb$ be an ideal of $k$ coprime to $\mff$;
      and
    \item $x : k \longrightarrow \C$ be a complex embedding.
  \end{itemize}
  Then for any two admissible elements $h, h' \in k$
  and $\scl := \mff\mfb^{-1}$,
  \[
    \Gamma_{\mfa, h}(\scl) = \Gamma_{\mfa, h'}(\scl).
  \]
\end{proposition}

\begin{remark}
  Gray \cite[Theorem~A]{gray} independently proves
  Proposition~\ref{prop:elliptic-gamma-independence}
  when $\rmn(\mfa)$ does not divide the discriminant of the
  characteristic polynomial of multiplication
  by $\eta$ on $k$.
  Gray \cite[Sections~6.2--6.3]{gray} also associates loops
  to admissible elements: choices of $h$ congruent modulo
  $\scl$ give loops representing the same homology class
  and the same evaluation maps on equivariant cohomology.
\end{remark}

We prove Proposition~\ref{prop:elliptic-gamma-independence}
by combining a product identity of Morain
\cite[Corollary~1]{morain-gamma-II} with linear algebra
over finite fields.
Our argument applies to every ideal
$\mfa$ satisfying the hypotheses of the proposition
and does not use our results on Stark units,
so we postpone its proof until
Section~\ref{sec:elliptic-gamma-independence}.

By Proposition~\ref{prop:elliptic-gamma-independence},
we may write $\Gamma_\mfa(\scl) := \Gamma_{\mfa, h}(\scl)$
for any admissible $h$.
We now turn to the other assertions
of Conjecture~\ref{conj:elliptic-gamma}.
The link between elliptic gamma values and Stark units is
the Kronecker limit type formula of
Bergeron--Charollois--Garc\'{i}a
\cite[Theorem~3.2]{bergeron-charollois-garcia},
which we now state.
Let $\Cl_k^+(\mff)$ denote the narrow ray class group modulo
$\mff$, and write $[\mfc]$ for the class of a nonzero
integral ideal $\mfc$ coprime to $\mff$.
For $\re(s) > 1$, define
\[
  \zeta_\mff(\mfb, s)
    := \sum_{\substack{
      0 \neq \mfc \subset \calo_k \\
      (\mfc, \mff) = 1 \\
      [\mfc] = [\mfb]}}
    \rmn(\mfc)^{-s}.
\]
We use its analytic continuation at $s = 0$.
Artin reciprocity gives an isomorphism
$\Cl_k^+(\mff) \longrightarrow \Gal(K/k)$ sending
$[\mfc]$ to $[\mfc, K/k] = \sigma_\mfc$.
Let $S_\mff$ consist of the archimedean places of $k$
and the finite places dividing $\mff$.
Since $[\mfc] = [\mfb]$ if and only if
$[\mfc, K/k] = \sigma_\mfb$,
the defining sums for
$\zeta_\mff(\mfb, s)$ and $L_{S_\mff}(s, \sigma_\mfb)$
have the same terms by Equation~\ref{eq:partial-L-series}
with $S = S_\mff$.
Thus $\zeta_\mff(\mfb, s) = L_{S_\mff}(s, \sigma_\mfb)$
for $\re(s) > 1$, and hence as meromorphic functions.
Writing $\lvert z \rvert_\C := z\overline{z}$,
Bergeron--Charollois--Garc\'{i}a
\cite[Theorem~3.2]{bergeron-charollois-garcia} prove that
\begin{equation}
  \label{eq:elliptic-gamma-limit-formula}
  \log\verts{\Gamma_\mfa(\scl)}_\C
    = \zeta_\mff'(\mfa\mfb, 0)
      - \rmn(\mfa)\zeta_\mff'(\mfb, 0).
\end{equation}

\subsubsection*{Units with the prescribed absolute values}
We first assume that a Stark unit $\epsilon$ for
$(K/k, S_\mff, v)$ is given
and compare elliptic gamma values with $\epsilon$;
we will then use Theorem~\ref{thm:intro-rank-one}
to obtain the Stark unit
when proving Corollary~\ref{cor:intro-elliptic-gamma}.
For this comparison,
Equation~\ref{eq:elliptic-gamma-limit-formula} leads us to
look for a unit $u_{\mff, \mfa}$ satisfying
$u_{\mff, \mfa}^{W_K} = \epsilon^{\rmn(\mfa) - \sigma_\mfa}$.

\begin{proposition}
  \label{prop:elliptic-gamma-Stark-model}
  If the integral rank-$1$ Stark conjecture for
  $(K/k, S_\mff, v)$ holds and
  $\epsilon \in \calo_K^\times$ is the corresponding Stark unit,
  then there is a unit
  $u_{\mff, \mfa} \in \calo_K^\times$ such that
  \begin{itemize}
    \item $u_{\mff, \mfa}^{W_K}
      = \epsilon^{\rmn(\mfa) - \sigma_\mfa}$;
    \item for $\scl := \mff\mfb^{-1}$, the conjugates
      $u_{\scl, \mfa} := u_{\mff, \mfa}^{\sigma_\mfb}$ satisfy
      \[
        \verts{\Gamma_\mfa(\scl)}_\C
        = \verts{u_{\scl, \mfa}}_v;
      \]
    \item each $u_{\scl, \mfa}$ has absolute value $1$
      above the real place of $k$ and satisfies
      $\sigma_\mfc(u_{\scl, \mfa})
      = u_{\mfc^{-1}\scl, \mfa}$
      for every ideal $\mfc$ coprime to $\mff$.
  \end{itemize}
\end{proposition}

\begin{proof}
  The extension $K(\epsilon^{1/W_K})/k$ is abelian,
  and $\sigma_\mfa(\zeta) = \zeta^{\rmn(\mfa)}$
  for every $\zeta \in \mu(K)$.
  By Stark \cite[Section~4, Lemma~6]{stark-1980},
  we can choose $u_{\mff, \mfa} \in \calo_K^\times$ with
  $u_{\mff, \mfa}^{W_K}
  = \epsilon^{\rmn(\mfa) - \sigma_\mfa}$.

  We now compare the absolute values of
  $\Gamma_\mfa(\scl)$ and $u_{\scl, \mfa}$.
  Since $\zeta_\mff(\mfb, s) = L_{S_\mff}(s, \sigma_\mfb)$,
  we can substitute the logarithms of $\epsilon$ from
  Definition~\ref{def:global-rank-one-Stark-unit}(i)
  into Equation~\ref{eq:elliptic-gamma-limit-formula}
  and get
  \[
    \log\verts{\Gamma_\mfa(\scl)}_\C
      = \frac{1}{W_K}\log\verts{
    \epsilon^{(\rmn(\mfa) - \sigma_\mfa)\sigma_\mfb}}_v
      = \log\verts{u_{\scl, \mfa}}_v.
  \]

  At every place of $K$ above the real place of $k$,
  the conjugates of $\epsilon$ have absolute value $1$
  by Definition~\ref{def:global-rank-one-Stark-unit}(ii).
  Since $u_{\mff, \mfa}^{W_K} = \epsilon^{\rmn(\mfa) -
  \sigma_\mfa}$,
  we also have $\lvert u_{\scl, \mfa}\rvert_{w'} = 1$
  at every such place $w'$.
  Since $u_{\scl, \mfa} = u_{\mff, \mfa}^{\sigma_\mfb}$
  and $\sigma_\mfc\sigma_\mfb = \sigma_{\mfb\mfc}$,
  we obtain
  $\sigma_\mfc(u_{\scl, \mfa}) = u_{\mfc^{-1}\scl, \mfa}$.

  In particular, the square of the usual absolute value of
  $\Gamma_\mfa(\scl)$ is the algebraic unit
  $\iota(u_{\scl, \mfa})\overline{\iota(u_{\scl, \mfa})}$.
  Thus
  $(\iota(u_{\scl, \mfa})\overline{\iota(u_{\scl, \mfa})})^{1/2}$
  and $(\iota(u_{\scl, \mfa})\overline{\iota(u_{\scl, \mfa})})^{-1/2}$
  are algebraic integers,
  so $\lvert \Gamma_\mfa(\scl)\rvert$
  is also an algebraic unit.
\end{proof}

We now use the existence of Stark units
in Theorem~\ref{thm:intro-rank-one} to prove
Corollary~\ref{cor:intro-elliptic-gamma}.

\begin{proof}[Proof of Corollary~\ref{cor:intro-elliptic-gamma}]
  The field $K = k(\mff)$ is cyclotomic over $k$, and
  $S_\mff$ contains every place above a finite set
  $\Sigma_\Q$ containing $\infty$, every prime ramified
  in $K/\Q$, and at least two finite primes.
  We apply Theorem~\ref{thm:intro-rank-one} with
  $S = S_\mff$ to obtain a Stark unit
  $\epsilon \in \calo_K^\times$ for $(K/k, S_\mff, v)$.
  We now use $\epsilon$ in
  Proposition~\ref{prop:elliptic-gamma-Stark-model}
  to obtain units with the required absolute values
  and Artin action, as well as the algebraicity of
  $\lvert \Gamma_\mfa(\scl) \rvert$.
\end{proof}

\subsubsection*{Comparison of complex values}
We now compare the complex values $\Gamma_\mfa(\scl)$
with the images of units in $K$.
For cyclotomic $K/k$, let $u_{\scl, \mfa}' \in \calo_K^\times$
satisfy $\Gamma_\mfa(\scl) = \iota(u_{\scl, \mfa}')$
and have absolute value $1$ at every place above the real
place of $k$, as required by
Conjecture~\ref{conj:elliptic-gamma}(i) and (ii).
We will prove $u_{\scl, \mfa}'/u_{\scl, \mfa} \in \mu(K)$.
The reciprocity law in part~(iii) then requires a single
$\xi \in \mu(K)$ such that
$u_{\scl, \mfa}' = (\xi u_{\mff, \mfa})^{\sigma_\mfb}$
for every ideal $\mfb$ coprime to $\mff$.
Recall that $c$ is complex conjugation at an embedding
of $K$ that is real on $k$.
\begin{lemma}
  \label{lem:elliptic-gamma-uniqueness}
  Let $K/k$ be cyclotomic. If $r \in \calo_K^\times$
  satisfies $c(r) = r^{-1}$ and $\lvert \iota(r)\rvert = 1$,
  then $r \in \mu(K)$.
\end{lemma}

\begin{proof}
  Recall that $L = K \cap \Q^{\ab}$,
  so $KF = k(LF)$ is Galois over $\Q$.
  We first show that $r \in LF$ by comparing two complex
  conjugations on $KF$. Continue to write $c$ for complex
  conjugation on $KF$ at an embedding real on $k$.
  Extend $\iota$ to $KF$, and let $d \in \Gal(KF/\Q)$
  satisfy $\iota(dz) = \overline{\iota(z)}$ for every $z \in KF$.
  Since $\lvert \iota(r)\rvert = 1$, we have
  $d(r) = r^{-1}$ and hence $dc(r) = r$.

  Both $c$ and $d$ restrict to complex conjugation on the
  abelian field $LF$. The restrictions of $c$ and $d$ to
  $kF$ are distinct transpositions: $c$ fixes $k$,
  whereas $d$ does not. Thus $dc$ generates $\Gal(KF/LF)$,
  and $r \in LF$. Since $LF/\Q$ is abelian,
  $c(r) = r^{-1}$ implies that every conjugate of $r$
  has absolute value $1$.
  As $r \in \calo_K^\times$, it follows that $r \in \mu(K)$.
\end{proof}

If a unit $r \in K$ has absolute value $1$ above the real
place of $k$, complex conjugation at an embedding real on
$k$ gives $c(r) = r^{-1}$.
Lemma~\ref{lem:elliptic-gamma-uniqueness} therefore shows
that the absolute values in
Corollary~\ref{cor:intro-elliptic-gamma} determine each
$u_{\scl, \mfa}$ up to a root of unity.
We use this uniqueness to express
Conjecture~\ref{conj:elliptic-gamma} in terms of
$u_{\mff, \mfa}$ and its conjugates by the Artin symbol
$\sigma_\mfb$.

\begin{proposition}
  \label{prop:elliptic-gamma-phase}
  Let $K/k$ be cyclotomic and let $\epsilon \in \calo_K^\times$
  be a Stark unit for $(K/k, S_\mff, v)$. Choose
  $u_{\mff, \mfa} \in \calo_K^\times$ with
  $u_{\mff, \mfa}^{W_K}
  = \epsilon^{\rmn(\mfa) - \sigma_\mfa}$.
  Then Conjecture~\ref{conj:elliptic-gamma} for $K/k$ and
  $\mfa$ is equivalent to the existence of $\xi \in \mu(K)$
  such that
  \[
    \Gamma_\mfa(\mff\mfb^{-1})
    = \iota\paren{(\xi u_{\mff, \mfa})^{\sigma_\mfb}}
  \]
  for every ideal $\mfb$ coprime to $\mff$.
\end{proposition}

\begin{proof}
  Assume that the units $u_{\scl, \mfa}'$ satisfy
  Conjecture~\ref{conj:elliptic-gamma}.
  By parts~(i) and (ii) and
  Proposition~\ref{prop:elliptic-gamma-Stark-model},
  $u_{\scl, \mfa}'$ and $u_{\scl, \mfa}$ have the same
  absolute value at $v$ and absolute value $1$ at every
  place above the real place of $k$.
  We therefore have
  $\xi := u_{\mff, \mfa}'/u_{\mff, \mfa} \in \mu(K)$
  by Lemma~\ref{lem:elliptic-gamma-uniqueness}.
  By parts~(i) and (iii) of the conjecture,
  $\Gamma_\mfa(\mff\mfb^{-1})
  = \iota((u_{\mff, \mfa}')^{\sigma_\mfb})$.
  Since $u_{\mff, \mfa}' = \xi u_{\mff, \mfa}$,
  this is the required identity.

  Conversely, suppose that the identity holds for one
  $\xi \in \mu(K)$ and every $\mfb$ coprime to $\mff$.
  For $\scl = \mff\mfb^{-1}$, define
  $u_{\scl, \mfa}' := (\xi u_{\mff, \mfa})^{\sigma_\mfb}$.
  Then $\iota(u_{\scl, \mfa}') = \Gamma_\mfa(\scl)$,
  and Proposition~\ref{prop:elliptic-gamma-Stark-model}
  shows that $u_{\scl, \mfa}'$ has absolute value $1$
  above the real place of $k$.
  For every ideal $\mfc$ coprime to $\mff$, the identity
  $\sigma_\mfc\sigma_\mfb
  = \sigma_{\mfb\mfc}$ gives
  $\sigma_\mfc(u_{\scl, \mfa}')
  = u_{\mfc^{-1}\scl, \mfa}'$.
\end{proof}

We now explain how to
choose $u_{\mff, \mfa}$ explicitly from the
construction in Proposition~\ref{prop:relative-rank-one}.
Assume that $K/k$ is cyclotomic and $F \subset K$,
and recall that $L = K \cap \Q^{\ab}$.
Let $S_\Q$ and $S_F$ be finite sets of places of $\Q$
and $F$, respectively, containing the archimedean places
and satisfying Equation~\ref{eq:local-factorization}
with $S = S_\mff$.
Suppose that $S_\Q$ contains every prime ramified in
$L/\Q$ and that $S_F$ contains every finite place
ramified in $K/F$ and at least two finite places.
Let $\epsilon = \beta^a$ be the Stark unit constructed
in the proof of Proposition~\ref{prop:relative-rank-one}.
We have $W_K = W_L$, and
$a \in \Z[\Gal(K/kF)]$ and $\beta \in \calo_K^\times$
satisfy $W_K\Theta_\Q(0) = (1 - c)a$ and
$c(\beta) = \beta^{-1}$.

To choose $u_{\mff, \mfa}$ in
Proposition~\ref{prop:elliptic-gamma-Stark-model}, we need
a $W_K$-th root of
$\epsilon^{\rmn(\mfa) - \sigma_\mfa}$ in $K$.
Since $\Gal(K/k) = \Gal(K/kF)\times\abrac{c}$, write
$\sigma_\mfa = g_\mfa c^{\delta_\mfa}$ with
$g_\mfa \in \Gal(K/kF)$ and $\delta_\mfa \in \set{0, 1}$.
Since $c(\beta) = \beta^{-1}$, it is enough to show that
$B_\mfa := (\rmn(\mfa) - (-1)^{\delta_\mfa}g_\mfa)a/W_K$
belongs to $\Z[\Gal(K/kF)]$.
For every $\zeta \in \mu(K) = \mu(L)$, we have
$\sigma_\mfa(\zeta) = \zeta^{\rmn(\mfa)}$, so by
Equation~\ref{eq:Stickelberger-action-integrality} we have
\[
  (1 - c)B_\mfa
    = (\rmn(\mfa) - \sigma_\mfa)\Theta_\Q(0)
    \in \Z[\Gal(K/k)].
\]
The cosets $\Gal(K/kF)$ and $c\Gal(K/kF)$ are disjoint,
so the coefficients of $(1 - c)B_\mfa$ on $\Gal(K/kF)$
are precisely those of $B_\mfa$.
Thus $B_\mfa \in \Z[\Gal(K/kF)]$, and we may take
$u_{\mff, \mfa} := \beta^{B_\mfa} \in \calo_K^\times$,
since
\[
  u_{\mff, \mfa}^{W_K} = \beta^{W_K B_\mfa}
    = \epsilon^{\rmn(\mfa) - \sigma_\mfa}.
\]

\subsection{Independence of elliptic gamma values}
\label{sec:elliptic-gamma-independence}

We prove Proposition~\ref{prop:elliptic-gamma-independence}
by expressing the values attached to two admissible
elements as the same product.
We first derive the required product identity from Morain
\cite[Corollary~1]{morain-gamma-II} in
Section~\ref{subsec:gamma-product-identity}.
In Section~\ref{subsec:gamma-lattices-mod-N},
we use linear
algebra over a finite field to choose linear forms
satisfying the hypotheses of
Lemma~\ref{lem:elliptic-gamma-identity}.
We then compare the values for two admissible elements
in Section~\ref{subsec:gamma-independence-proof}.

We recall the notation that we need from
Section~\ref{subsec:elliptic-gamma},
following Bergeron--Charollois--Garc\'{i}a
\cite[Sections~1--3]{bergeron-charollois-garcia}.
Let
\begin{itemize}
  \item $k$ be a complex cubic field,
    with real embedding $x_\R : k \longrightarrow \R$
    and a fixed complex embedding $x : k \longrightarrow \C$;
  \item $\mff$ be a nonzero ideal of $\calo_k$
    whose narrow ray class field $k(\mff)$ is totally imaginary;
  \item $\mfb$ be a nonzero ideal of $\calo_k$
    coprime to $\mff$;
  \item $\mfa$ be an ideal of $\calo_k$
    coprime to $6\rmn(\mff)$
    with prime norm $\rmn(\mfa)$;
  \item $\scl$ be the fractional ideal $\mff\mfb^{-1}$ of $k$;
  \item $M$ be the fractional ideal $\mfa^{-1}\scl$ of $k$;
  \item $q$ be the least positive integer in
    $\mff \cap \Z = \scl \cap \Z$.
\end{itemize}
Thus $\rmn(\mfa) \geq 5$, $[M : \scl] = \rmn(\mfa)$,
and $\rmn(\mfa)M \subset \scl$.
Since $k(\mff)$ is totally imaginary,
Bergeron--Charollois--Garc\'{i}a
\cite[Section~2.2]{bergeron-charollois-garcia} show that
$-1 \not\equiv 1 \pmod{\mff}$.
Thus $2 \notin \mff$ and $q \geq 3$.
Since $q \mid \rmn(\mff)$ and $(\rmn(\mfa), \rmn(\mff)) = 1$,
we have $(q, \rmn(\mfa)) = 1$.
In particular, $1 \notin M$:
otherwise $\rmn(\mfa) \in \scl$
and $q \mid \rmn(\mfa)$.

An element $h = \lambda/q \in k$ is \emph{admissible} if
$\lambda \in \scl$, $h \equiv 1 \pmod{\scl}$, and
$\lambda/\rmn(\mfa) \in M$ has nonzero image in $M/\scl$.
For an admissible $h$, the coset $\lambda/\rmn(\mfa) + \scl$
generates $M/\scl$, since $M/\scl$ has prime order.
Thus $M = \scl + \Z\lambda/\rmn(\mfa)$.

\subsubsection{A product identity for elliptic gamma functions}
\label{subsec:gamma-product-identity}

We use pairs of $\Z$-linear maps
$A : M \longrightarrow \Z$ to index elliptic gamma
functions on $M$ and $\scl$.
To specify which maps we use,
let $k_\R := k \otimes_\Q \R$ and extend $x_\R$ and $x$
$\R$-linearly to $k_\R$.
We give $k_\R$ the orientation for which the isomorphism
$(x_\R, \re x, \Im x) : k_\R \longrightarrow \R^3$
preserves the standard orientation of $\R^3$.
Let $e_\R \in k_\R$ correspond to $(1, 0, 0)$ under this
isomorphism.
Every $\Z$-linear map $A : M \longrightarrow \Z$ extends
uniquely to a $\Q$-linear form on $k$ and then to an
$\R$-linear form on $k_\R$.
To evaluate elliptic gamma functions on $M$ and $\scl$,
we use those linear forms that are primitive on $\scl$
and negative at $e_\R$:
\[
  \cala := \set{
    A \in \Hom_\Z(M, \Z) :
    A(\scl) = \Z,\ A(e_\R) < 0
  }.
\]
Here $A(\scl) = \Z$ means that the restriction of $A$ to
$\scl$ is primitive.
Since $\scl \subset M$, we also have $A(M) = \Z$,
so $A$ is primitive on $M$ and $\scl$.
The inequality $A(e_\R) < 0$ defines an open halfspace in
$\Hom_\R(k_\R, \R)$.
The condition $x \in U_A$ in
\cite[Equation~2.1]{bergeron-charollois-garcia} means that
$\Im(x(u_1)\overline{x(u_2)}) > 0$, where
$(u_1, u_2)$ is an $\R$-basis of $\ker A \subset k_\R$ oriented
so that $(u_1, u_2, u_3)$ is positive for every
$u_3 \in k_\R$ with $A(u_3) > 0$.
With our orientation of $k_\R$,
$\Im(x(u_1)\overline{x(u_2)})$ has sign opposite to $A(e_\R)$,
so $x \in U_A$ is equivalent to
$A(e_\R) < 0$.
If $A, B \in \cala$ are proportional, primitivity on
$\scl$ gives $B = \pm A$.
Since $A(e_\R), B(e_\R) < 0$, we have $B = A$.
Thus distinct linear forms in $\cala$ are linearly independent.

We use the elliptic gamma function from
Equation~\ref{eq:elliptic-gamma-function}.
For $\Im(\tau) < 0 < \Im(\sigma)$,
set $\Gamma(z, \tau, \sigma)
= \Gamma(z - \tau, -\tau, \sigma)^{-1}$,
as in \cite[Section~1]{bergeron-charollois-garcia}.

Fix distinct $A, B \in \cala$ and $D = \scl$ or $M$.
Choose $\alpha, \beta \in D$ such that
$B(\alpha) = A(\beta) = 0$ and
$A(\alpha), B(\beta) > 0$.
Let $\gamma_D$ generate
$\set{y \in D : A(y) = B(y) = 0}$, and choose its sign so
that $(\alpha, \beta, \gamma_D)$ is positively oriented.

Since $x \in U_A \cap U_B$, applying
\cite[Equation~2.6]{bergeron-charollois-garcia} to
$D$ and $(A, B)$ gives $\Im(x(\alpha)/x(\gamma_D)) < 0$
and $\Im(x(\beta)/x(\gamma_D)) > 0$.
The linear forms $A$ and $B$ are primitive on $D$ and
linearly independent.
To compare admissible elements $h$ with $h - 1 \in \scl$,
we will use elliptic gamma functions
indexed by a coset of $D$.
For $t \in k$, Morain
\cite[Proposition~7 and Equation~12]{morain-gamma-I}
expresses the elliptic gamma function associated with
the coset $t + D$ and the ordered pair $(A, B)$
as the following finite product
that is meromorphic in $w \in \C$:
\begin{equation}
  \label{eq:elliptic-gamma-product}
  G_{1, A, B}(t)(w, x; D) :=
  \prod_{\substack{
    \delta \in (t + D)/\Z\gamma_D \\
    0 \leq A(\delta) < A(\alpha) \\
    0 \leq B(\delta) < B(\beta)
  }}
  \Gamma\paren{
    \frac{w + x(\delta)}{x(\gamma_D)},
    \frac{x(\alpha)}{x(\gamma_D)},
    \frac{x(\beta)}{x(\gamma_D)}
  }.
\end{equation}
The product in Equation~\ref{eq:elliptic-gamma-product}
is independent of $\alpha$ and $\beta$ and depends on $t$
only through the coset $t + D$.
Changing a representative $\delta$ modulo $\Z\gamma_D$ leaves
each factor unchanged because $\Gamma$ is $1$-periodic in $z$.
At $t = 0$, we recover the elliptic gamma function
$\Gamma_{A, B}(w, x; D) := G_{1, A, B}(0)(w, x; D)$ of
Bergeron--Charollois--Garc\'{i}a
\cite[Equation~2.7]{bergeron-charollois-garcia}.

For $D = \scl$ or $M$, every admissible $h$ satisfies
$h + D = 1 + D$.
We therefore set $t = 1$ for $D = \scl$ and $D = M$,
and define
\[
  \Phi(A, B) :=
    \frac{G_{1, A, B}(1)(0, x; M)}
    {G_{1, A, B}(1)(0, x; \scl)^{\rmn(\mfa)}}.
\]
To use $\Phi$ in multiplicative identities, we first check
that its numerator and denominator are finite and nonzero.
For $D = \scl$ and $D = M$, the value
$G_{1, A, B}(1)(0, x; D)$ is finite and nonzero.
Indeed, a zero or pole of a factor indexed by $\delta$ in
Equation~\ref{eq:elliptic-gamma-product}, with $t = 1$ and $w = 0$,
would imply
$x(\delta) \in
\Z x(\gamma_D) + \Z x(\alpha) + \Z x(\beta)$.
Injectivity of $x$ on $k$ would then give $\delta \in D$,
contrary to $\delta \in 1 + D$ and $1 \notin M$.

We will also need to reverse the order of the two linear
forms in $\Phi$.
Morain
\cite[Section~2.2, Remark after the proof of Theorem~1]
{morain-gamma-I} proves that interchanging $A$ and $B$
replaces $G_{1, A, B}(t)(w, x; D)$ by its reciprocal;
applying this identity at $(t, w) = (1, 0)$ to the numerator
and denominator defining $\Phi$ gives
$\Phi(B, A) = \Phi(A, B)^{-1}$.
We will use this equality to cancel factors
in the comparison of two admissible choices of $h$.

We now derive a product identity for $\Phi$ from Morain
\cite[Corollary~1]{morain-gamma-II},
which gives a multiplicative relation for
quotients of elliptic gamma functions on two nested lattices.

\begin{lemma}
  \label{lem:elliptic-gamma-identity}
  Let $A, B, C \in \cala$ be pairwise distinct.
  If $A, B, C$ are linearly independent, assume that a
  primitive generator of each of
  $M \cap \ker A \cap \ker B$,
  $M \cap \ker B \cap \ker C$, and
  $M \cap \ker A \cap \ker C$ lies outside $\scl$.
  Then
  \[
    \Phi(A, B)\Phi(B, C) = \Phi(A, C).
  \]
\end{lemma}

\begin{proof}
  We apply Morain \cite[Corollary~1]{morain-gamma-II}
  with lattices $L = M = \mfa^{-1}\scl$ and $L' = \scl$,
  dimension $3$, and shift $1 + \scl$.
  We first choose a positively oriented basis
  $(e_1, e_2, e_3)$ of $M$ with
  $\scl = \Z \rmn(\mfa)e_1 + \Z e_2 + \Z e_3$
  as required by \cite[Corollary~1]{morain-gamma-II}.
  Since $M/\scl$ has prime order $\rmn(\mfa)$, the inclusion
  $\scl \hookrightarrow M$ has invariant factors
  $1$, $1$, and $\rmn(\mfa)$.
  We can therefore choose the required basis by Smith
  normal form, then change the sign of one basis vector
  if necessary to make it positively oriented.

  For $(P, Q) \in \set{(A, B), (B, C), (A, C)}$,
  the quotient in \cite[Equation~9]{morain-gamma-II},
  evaluated at $w = 0$, is
  \[
    \frac{G_{1, P, Q}(1)(0, x; \scl)^{\rmn(\mfa)}}
    {G_{1, P, Q}(1)(0, x; M)} = \Phi(P, Q)^{-1}.
  \]
  For the prime index $\rmn(\mfa)$,
  \cite[Corollary~1]{morain-gamma-II}
  therefore gives
  \[
    \bigl(\Phi(B, C)^{-1}\Phi(A, C)
    \Phi(A, B)^{-1}\bigr)^{\cald(\rmn(\mfa), 3)} = 1,
  \]
  where $\cald(\rmn(\mfa), 3)
  = \rmn(\mfa)^{\lfloor 3/(\rmn(\mfa) - 1)\rfloor}$,
  provided the following conditions hold:
  \begin{enumerate}[(i)]
    \item The linear forms $A$, $B$, and $C$ are integral
      on $M$ and primitive on $\scl$.
    \item If their $\Q$-span has dimension $2$,
      the zero linear form
      is not in the convex hull of $A$, $B$, and $C$.
    \item If $A$, $B$, and $C$ are linearly independent,
      then for each pair
      $(P, Q) \in \set{(A, B), (B, C), (A, C)}$ we have
      \[
        \scl \cap \ker P \cap \ker Q
        = \rmn(\mfa)\bigl(M \cap \ker P \cap \ker Q\bigr).
      \]
  \end{enumerate}
  Conditions~(i)--(iii) are the integrality and primitivity
  requirements and the conditions in Morain
  \cite[Definitions~4 and~18]{morain-gamma-II}
  specialized to a triple of
  pairwise independent linear forms.
  The restrictions of $x$ to $M$ and $\scl$ give the
  evaluations used in the definition of $\Phi$.
  We now verify (i)--(iii).

  For (i), let $P \in \set{A, B, C}$.
  Since $P \in \cala$,
  it takes integral values on $M$
  and satisfies $P(\scl) = \Z$.
  The inclusion $\scl \subset M$ therefore
  gives $\Z = P(\scl) \subset P(M) \subset \Z$.
  Thus $P(M) = P(\scl) = \Z$, proving (i).

  For (ii), we only need to check
  when $A$, $B$, and $C$
  are linearly dependent,
  since otherwise their $\Q$-span would have
  dimension $3$.
  Suppose that two distinct forms
  $P, Q \in \set{A, B, C}$ are proportional.
  Since both are primitive on $\scl$, we have $Q = \pm P$.
  The inequalities $P(e_\R) < 0$ and $Q(e_\R) < 0$
  exclude the possibility that $Q = -P$.
  But then $Q = P$, which contradicts $P$ and $Q$
  being distinct.
  Hence $A$, $B$, and $C$ are linearly dependent
  but \textit{pairwise} linearly independent.
  Consequently, the $\Q$-span of $A$, $B$, and $C$ has
  dimension exactly $2$.
  Every convex combination of $A$, $B$, and $C$ takes a
  negative value at $e_\R$, so their convex hull
  does not include the zero linear form.
  This proves (ii).

  For (iii), suppose that $A$, $B$, and $C$ are
  linearly independent.
  For each pair
  $(P, Q) \in \set{(A, B), (B, C), (A, C)}$,
  choose a primitive generator $\gamma$ of
  $M \cap \ker P \cap \ker Q$.
  Our hypothesis gives $\gamma \notin \scl$, so
  $\gamma + \scl$ has order $\rmn(\mfa)$ in $M/\scl$.
  Consequently,
  \begin{align*}
    \scl \cap \ker P \cap \ker Q
      &= \rmn(\mfa)\Z\gamma \\
      &= \rmn(\mfa)\bigl(M \cap \ker P \cap \ker Q\bigr).
  \end{align*}
  This proves (iii).

  Having verified (i)--(iii), we apply
  \cite[Corollary~1]{morain-gamma-II} to
  $\Phi(A, B)^{-1}$, $\Phi(B, C)^{-1}$,
  and $\Phi(A, C)^{-1}$.
  Since $\rmn(\mfa) \geq 5$,
  the exponent $\cald(\rmn(\mfa), 3)$
  equals $1$.
  Rearranging gives $\Phi(A, B)\Phi(B, C) = \Phi(A, C)$.
\end{proof}

\subsubsection{Lattices modulo
  \texorpdfstring{$\rmn(\mfa)$}{N(a)}}
\label{subsec:gamma-lattices-mod-N}

To apply Lemma~\ref{lem:elliptic-gamma-identity}, we need
primitive generators of $M \cap \ker A \cap \ker B$
to lie outside $\scl$. We test this condition using linear
maps on $\scl/(\rmn(\mfa)M)$.

Since $\rmn(\mfa)$ is prime,
$\F_{\rmn(\mfa)} := \Z/\rmn(\mfa)\Z$ is
the field with $\rmn(\mfa)$ elements.
Let $W := \scl/(\rmn(\mfa)M)$.
The inclusions $\rmn(\mfa)M \subset \scl \subset M$
give the exact sequence
\[
  0 \longrightarrow W \longrightarrow M/(\rmn(\mfa)M)
  \longrightarrow M/\scl \longrightarrow 0,
\]
of $\F_{\rmn(\mfa)}$-vector spaces. The last map sends $y +
\rmn(\mfa)M$ to
$y + \scl$. Since $M$ has rank $3$ over $\Z$ and
$[M : \scl] = \rmn(\mfa)$,
we have $\dim_{\F_{\rmn(\mfa)}} W = 2$.

Let $W^\vee := \Hom_{\F_{\rmn(\mfa)}}(W, \F_{\rmn(\mfa)})$.
Every $A \in \Hom_\Z(M, \Z)$ induces a linear form
$\overline{A} \in W^\vee$ by
\[
  \begin{tikzcd}[row sep = tiny]
    W \arrow[r, "\overline{A}"] & \F_{\rmn(\mfa)} \\
    y + \rmn(\mfa)M \arrow[r, mapsto] & A(y) \pmod{\rmn(\mfa)}
  \end{tikzcd}
\]
for $y \in \scl$.
The map
\[
  \begin{tikzcd}[row sep = tiny]
    \Hom_\Z(M, \Z) \arrow[r] & W^\vee \\
    A \arrow[r, mapsto] & \overline{A}
  \end{tikzcd}
\]
is surjective.
Indeed, we can extend a linear form on $W$
to $M/(\rmn(\mfa)M)$
and lift its values on the images
of a basis of $M$ to integers.
For $A \in \cala$, the equality $A(\scl) = \Z$ implies
$\overline{A} \neq 0$.

Write $[\overline{A}] \in \P(W^\vee)$ for the
line spanned by $\overline{A}$.
For $A, B \in \cala$, the condition
$[\overline{A}] \neq [\overline{B}]$ says that the map
\[
  \begin{tikzcd}[row sep = tiny]
    W \arrow[r] & \F_{\rmn(\mfa)}^2 \\
    y \arrow[r, mapsto] & (\overline{A}(y), \overline{B}(y))
  \end{tikzcd}
\]
is an isomorphism. The next lemma uses its injectivity to
show that a primitive generator of $M \cap \ker A \cap \ker B$
lies outside $\scl$. It also shows that every point of
$\P(W^\vee)$ comes from a linear form in $\cala$.

\begin{lemma}
  \label{lem:elliptic-gamma-lattices}
  Every point of $\P(W^\vee)$ is equal to
  $[\overline{A}]$ for some $A \in \cala$.
  Furthermore, if $A, B \in \cala$ satisfy
  $[\overline{A}] \neq [\overline{B}]$, then
  $M = \scl + \Z\gamma$, where $\gamma$ is a primitive
  generator of $M \cap \ker A \cap \ker B$.
\end{lemma}

\begin{proof}
  Fix a nonzero linear form $\ell \in W^\vee$ and lift
  it to $A_0 \in \Hom_\Z(M, \Z)$ using the surjection
  $\Hom_\Z(M, \Z) \longrightarrow W^\vee$.
  Write $A_0(\scl) = d\Z$ with $d > 0$.
  Since $\ell \neq 0$, we have $\rmn(\mfa) \nmid d$.
  The inclusion $\rmn(\mfa)M \subset \scl$ gives
  $d \mid \rmn(\mfa)A_0(y)$ for $y \in M$, so $A_0/d$ is
  integral on $M$ and primitive on $\scl$.

  It remains to choose the sign so that the resulting
  linear form is negative at $e_\R$.
  The trace pairing for $k/\Q$ is nondegenerate
  (see \cite[Tag~0BIL]{stacks-project}),
  so the extension of $A_0$ to $k$ has the form
  $A_0(y) = \Tr_{k/\Q}(\theta y)$ for some
  $\theta \in k^\times$.
  Hence $A_0(e_\R) = x_\R(\theta) \neq 0$.
  Choose $A = A_0/d$ or $A = -A_0/d$ so that
  $A(e_\R) < 0$. Then $A \in \cala$ and
  $[\overline{A}] = [\ell]$.

  Now suppose $[\overline{A}] \neq [\overline{B}]$.
  The map $W \longrightarrow \F_{\rmn(\mfa)}^2$ defined by
  $\overline{A}$ and $\overline{B}$ is injective.
  A primitive generator $\gamma$ of
  $M \cap \ker A \cap \ker B$
  is not contained in $\rmn(\mfa)M$.
  If $\gamma \in \scl$,
  its class in $W$ would be a nonzero element of
  $\ker\overline{A} \cap \ker\overline{B}$,
  a contradiction.
  Thus $\gamma \notin \scl$
  and its class generates $M/\scl$,
  which has prime order $\rmn(\mfa)$.
\end{proof}

Consequently, Lemma~\ref{lem:elliptic-gamma-identity}
applies to every triple in $\cala$ with pairwise distinct
classes in $\P(W^\vee)$.
Since $\P(W^\vee)$ has $\rmn(\mfa) + 1 \geq 6$ points,
Lemma~\ref{lem:elliptic-gamma-lattices} also allows us to
choose linear forms whose classes avoid any specified set
of at most four points.

\subsubsection{Independence from \texorpdfstring{$h$}{h}}
\label{subsec:gamma-independence-proof}

We now describe the linear forms associated with an
admissible element $h$ and compare the elliptic gamma
values as $h$ varies.

Let $\eta$ generate the group of units
$u \in \calo_k^\times$ satisfying
$x_\R(u) > 0$ and $u \equiv 1 \pmod{\mff}$,
choosing the generator with $0 < x_\R(\eta) < 1$.
Let $g : k \longrightarrow k$ be multiplication by $\eta$.
Since $\eta$ is a unit and $M$ and $\scl$ are fractional
ideals, $g(M) = M$ and $g(\scl) = \scl$.
Moreover, the determinant of $g$ on $k_\R$ is
$\rmn_{k/\Q}(\eta) = 1$,
so $g$ preserves the orientation of $k_\R$.

To act on the linear forms used to define $\Phi$,
let $(gA)(y) := A(\eta^{-1}y)$ for
$A \in \Hom_\Q(k, \Q)$ and $y \in k$.
The equalities $g(M) = M$ and $g(\scl) = \scl$ show
that this action preserves integrality on $M$ and
primitivity on $\scl$.
Also, $(gA)(e_\R) = x_\R(\eta)^{-1}A(e_\R)$, so
$g$ preserves the sign at $e_\R$ and hence acts on $\cala$.

To track the reductions of these linear forms, note that
$g(\rmn(\mfa)M) = \rmn(\mfa)M$. Together with $g(M) = M$ and
$g(\scl) = \scl$, this gives induced automorphisms of
$M/(\rmn(\mfa)M)$ and $W = \scl/(\rmn(\mfa)M)$.
On $W^\vee$, define $(g\ell)(y) := \ell(g^{-1}y)$
for $\ell \in W^\vee$ and $y \in W$.
By these definitions, $\overline{gA} = g\overline{A}$.
Consequently, $g$ permutes $\P(W^\vee)$ by
$[\overline{A}] \longmapsto [\overline{gA}]$.

We next prove that applying $g$ to both linear forms
leaves $\Phi$ unchanged.
Fix distinct $A, B \in \cala$ and let $D = \scl$ or $M$.
Since $\eta - 1 \in \mff \subset \scl$,
multiplication by $\eta$ preserves $1 + D$.
For the pair $(gA, gB)$ in
Equation~\ref{eq:elliptic-gamma-product},
choose $\eta\alpha$, $\eta\beta$, and $\eta\gamma_D$
in place of $\alpha$, $\beta$, and $\gamma_D$,
and index the factors by $\eta\delta$.
Multiplication by $\eta$ preserves the orientation of $k_\R$.
The inequalities selecting a factor of
Equation~\ref{eq:elliptic-gamma-product} are
$0 \leq A(\delta) < A(\alpha)$ and
$0 \leq B(\delta) < B(\beta)$.
Since $(gA)(\eta y) = A(y)$ and $(gB)(\eta y) = B(y)$,
the two sets of inequalities are equivalent to
\begin{align*}
  0 &\leq (gA)(\eta\delta) < (gA)(\eta\alpha), \\
  0 &\leq (gB)(\eta\delta) < (gB)(\eta\beta).
\end{align*}
Thus $\delta \longmapsto \eta\delta$ identifies the
indexing sets for $(A, B)$ and $(gA, gB)$.
Since $x(\eta y) = x(\eta)x(y)$, the three arguments of
$\Gamma$ in Equation~\ref{eq:elliptic-gamma-product}
at $(t, w) = (1, 0)$ are unchanged.
Hence $G_{1, gA, gB}(1)(0, x; D) = G_{1, A, B}(1)(0, x; D)$.
Taking the quotient for $D = M$ and $D = \scl$ gives
$\Phi(gA, gB) = \Phi(A, B)$.

To express $\Gamma_{\mfa, h}(\scl)$ in terms of $\Phi$,
we associate a linear form to each admissible
$h = \lambda/q$.
Bergeron--Charollois--Garc\'{i}a
\cite[Lemma~2.1]{bergeron-charollois-garcia} prove that
there is a unique primitive $\Z$-linear map
\[
  A_\lambda: \scl \longrightarrow \Z
\]
such that
$\ker A_\lambda \cap \ker(gA_\lambda) = \Q\lambda$
and $A_\lambda(e_\R) < 0$.
Here the kernels are taken in $k$ after extending
$A_\lambda$ $\Q$-linearly, and the inequality uses its
$\R$-linear extension to $k_\R$.
Thus $A_\lambda$ and $gA_\lambda$ are distinct and
both vanish on $\lambda$.
Since $M = \scl + \Z\lambda/\rmn(\mfa)$,
both take integral values on $M$.
Multiplication by $\eta$ preserves $\scl$,
so $gA_\lambda$ is also primitive on $\scl$.
Moreover,
\[
  (gA_\lambda)(e_\R) = x_\R(\eta)^{-1}A_\lambda(e_\R) < 0.
\]
Therefore $A_\lambda$ and $gA_\lambda$ belong to $\cala$.
Bergeron--Charollois--Garc\'{i}a
\cite[Equation~2.13 and Definition~2.3]
{bergeron-charollois-garcia},
with linear forms $(A_\lambda, gA_\lambda)$ and parameters
$(w, z) = (x(h), x)$ in their notation,
express the value in
Definition~\ref{def:arithmetic-elliptic-gamma} as
\[
  \Gamma_{\mfa, h}(\scl) =
    \frac{\Gamma_{A_\lambda, gA_\lambda}(x(h), x; M)}
    {\Gamma_{A_\lambda, gA_\lambda}
    (x(h), x; \scl)^{\rmn(\mfa)}}.
\]

We now choose two more linear forms so that the values
for two admissible elements can be expressed as the same
product. We choose their reductions modulo $\rmn(\mfa)$ using
Lemma~\ref{lem:elliptic-gamma-lattices}.

\begin{proof}[Proof of
  Proposition~\ref{prop:elliptic-gamma-independence}]
  Let $h = \lambda/q$ and $h' = \lambda'/q$ be admissible,
  and set $A := A_\lambda$ and $C := A_{\lambda'}$.
  We first express $\Gamma_{\mfa, h}(\scl)$ using $\Phi$.
  Since $A(h) = (gA)(h) = 0$, replacing $\delta$ by
  $\delta + h$ preserves the bounds in
  Equation~\ref{eq:elliptic-gamma-product}.
  For $D = \scl$ or $M$, we therefore have
  \begin{align*}
    \Gamma_{A, gA}(x(h), x; D)
    &= G_{1, A, gA}(h)(0, x; D) \\
    &= G_{1, A, gA}(1)(0, x; D),
  \end{align*}
  where the second equality uses $h + D = 1 + D$.
  Taking the quotient for $D = M$ and $D = \scl$, we obtain
  $\Gamma_{\mfa, h}(\scl) = \Phi(A, gA)$.
  The same calculation with $h'$ and $C$ proves
  $\Gamma_{\mfa, h'}(\scl) = \Phi(C, gC)$.
  Moreover, $\lambda/\rmn(\mfa)$ lies in
  $M \cap \ker A \cap \ker(gA)$ but not in $\scl$, so a
  primitive generator of $M \cap \ker A \cap \ker(gA)$
  lies outside $\scl$.
  Since $\lambda'/\rmn(\mfa) \notin \scl$, a primitive
  generator of $M \cap \ker C \cap \ker(gC)$ also lies
  outside $\scl$.
  Thus admissibility gives the required condition for the
  pairs $(A, gA)$ and $(C, gC)$, even when the two reductions
  in either pair span the same line in $\P(W^\vee)$.
  In particular, this argument also applies when
  $\eta \equiv 1 \pmod{\rmn(\mfa)\calo_k}$.

  We choose two auxiliary linear forms to express
  both values by the same product. By
  Lemma~\ref{lem:elliptic-gamma-lattices}, choose
  $R, S \in \cala$ successively so that
  \begin{align*}
    [\overline{R}] & \notin \set{
      [\overline{A}], [\overline{C}],
      [\overline{g^{-1}A}],
      [\overline{g^{-1}C}]
    }, \\
    [\overline{S}] & \notin \set{
      [\overline{A}], [\overline{C}],
      [\overline{R}], [\overline{gR}]
    }.
  \end{align*}
  Each choice excludes at most four of the
  $\rmn(\mfa) + 1 \geq 6$ points of $\P(W^\vee)$.

  For $X \in \set{A, C}$,
  we will use
  Lemma~\ref{lem:elliptic-gamma-identity} to replace
  $\Phi(X, gX)$ by a product involving $R$ and $S$.
  The triples we need are $(X, gR, gX)$, $(X, S, gR)$,
  and $(R, S, X)$.
  For each pair $(P, Q)$ in these triples other than
  $(X, gX)$, we have
  $[\overline{P}] \neq [\overline{Q}]$ in $\P(W^\vee)$.
  Thus a primitive generator of $M \cap \ker P \cap \ker Q$
  lies outside $\scl$ by
  Lemma~\ref{lem:elliptic-gamma-lattices}.
  For $X = A$ and $X = C$, we proved at the start that a
  primitive generator of $M \cap \ker X \cap \ker(gX)$
  lies outside $\scl$.
  Using $\Phi(gR, gX) = \Phi(R, X)$ and
  $\Phi(S, X) = \Phi(X, S)^{-1}$,
  we obtain
  \begin{align*}
    \Phi(X, gX)
      &= \Phi(X, gR)\Phi(R, X) \\
      &= \Phi(X, S)\Phi(S, gR) \Phi(R, S)\Phi(S, X) \\
      &= \Phi(R, S)\Phi(S, gR).
  \end{align*}
  The right side is the same for $X = A$ and $X = C$,
  so $\Gamma_{\mfa, h}(\scl) = \Gamma_{\mfa, h'}(\scl)$.
\end{proof}


\section{Cyclotomic Iwasawa theory over complex cubic fields}
\label{sec:Iwasawa}

We prove Theorem~\ref{thm:intro-determinant-main}.
We first show that the Stark units
can be chosen in a system that is compatible under norms,
and then compare their inverse limits with ideal classes.
We decompose the modules of units and ideal classes
into components over $\Q$ and $F$, then apply the
Iwasawa main conjectures over these fields.

\subsection{The tower and norm compatibility}
\label{subsec:finite-norm-relations}

We show that the units in Theorems~\ref{thm:intro-rank-one}
and \ref{thm:intro-rank-two} can be chosen compatibly under
norms in a cyclotomic tower.
We first compute the groups of roots of unity at each level $n$,
and then compare regulators to choose the
cyclotomic and elliptic units satisfying compatible norm relations
that are preserved by the constructions
in Propositions~\ref{prop:rank-one-regulator}
and \ref{prop:rank-two-regulator}.

Fix a complex cubic field $k$, its Galois closure $E$,
and its quadratic resolvent $F$.
Let $L_0/\Q$ be finite abelian with $F \subset L_0$,
and let $p > 3$ be unramified in $K_0 := EL_0$, with
$p \nmid [L_0 : \Q]$.

\begin{figure}[ht]
  \[
    \begin{tikzcd}[row sep = 1em, column sep = 2.5em]
      & K_n \arrow[dd, dash, "3"] \\
      K_n^+ \arrow[ur, dash, "2"] \arrow[dd, dash, "3"']
      & \\
      & L_n \arrow[dd, dash, "p^n"] \\
      L_n^+ \arrow[ur, dash, "2"] \arrow[dd, dash, "p^n"']
      & \\
      & L_0 \\
      L_0^+ \arrow[ur, dash, "2"] &
    \end{tikzcd}
  \]
  \caption{The field extensions at level $n$.}
  \label{figure:Iwasawa-tower-fields}
\end{figure}
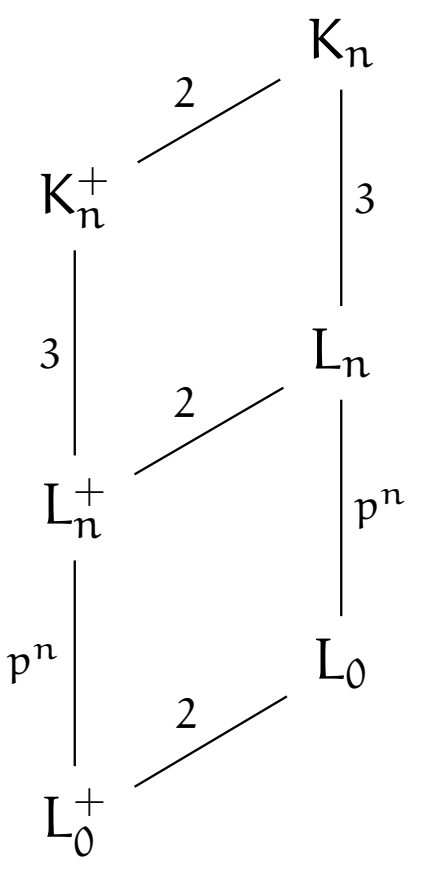

Let $\Q_\infty/\Q$ be the cyclotomic $\Z_p$-extension
contained in $\bigcup_{m \geq 1}\Q(\mu_{p^m})$,
and let $\Q_n$ be its layer of degree $p^n$.
Let $L_n := L_0\Q_n$ and $K_n := kL_n = EL_n$.
Choose compatible embeddings of $K_n$ into $\C$ that are
real on $k$, and let $c$ denote complex conjugation.
Let $L_n^+ := L_n^{\abrac{c}}$ and
$K_n^+ := K_n^{\abrac{c}}$
be the fixed fields of $c$.
Since $E \cap L_0 = F$, the degree
$[K_0 : \Q] = 3[L_0 : \Q]$ is coprime to $p$.
Hence $K_0 \cap \Q_\infty = \Q$.
Since $k \cap L_n = \Q$, restriction gives
$\Gal(K_n/k) \cong \Gal(L_0/\Q)\times\Gal(\Q_n/\Q)$.
Since $\Q_n$ is totally real, $L_n^+ = L_0^+\Q_n$,
and $K_n^+ = kL_n^+$.
We summarize the field extensions at level $n$
in Figure~\ref{figure:Iwasawa-tower-fields}.

Choose a finite set $\Sigma_\Q$ of places of $\Q$,
including $\infty$, $p$, and every prime ramified in $K_0/\Q$.
Let $\Sigma_k$ and $\Sigma_F$ consist of the places above
$\Sigma_\Q$ in $k$ and $F$, respectively.
Both $\Sigma_\Q$ and $\Sigma_F$ contain at least two
finite places: by Minkowski's theorem, some prime
divides the discriminant of $k$, and this prime differs
from $p$.
For every $n$, the extension $K_n/\Q$ is unramified outside
$\Sigma_\Q$. Thus $\Sigma_k$, $\Sigma_\Q$, and $\Sigma_F$
contain the finite places ramified in $K_n/k$, $L_n/\Q$,
and $K_n/F$, respectively.
Because $\Sigma_k$ and $\Sigma_F$ consist of all places
above $\Sigma_\Q$, we can use
Equation~\ref{eq:local-factorization} at every level by
Corollary~\ref{cor:rationally-saturated-factorization},
with $S = \Sigma_k$, $S_\Q = \Sigma_\Q$, and
$S_F = \Sigma_F$.
Choose compatible places $w_n$ above the complex place
of $k$, and choose a generator $\tau$ of
$\Gal(K_n/L_n) \cong \Gal(E/F)$ so that
$c w_n = \tau w_n$.

The absolute values in
Definition~\ref{def:global-rank-one-Stark-unit} determine
a Stark unit up to a root of unity. We compute the groups
$\mu(K_n)$ and $\mu(L_n)$ to choose units with exact norm
relations.
Since $K_n^+$ and $L_n^+$ have real embeddings,
$\mu(K_n^+) = \mu(L_n^+) = \set{\pm 1}$.

\begin{lemma}
  \label{lem:cyclotomic-roots-of-unity}
  The groups of roots of unity in $K_n$ and $L_n$ satisfy:
  \begin{enumerate}[(i)]
    \item $\mu(K_n) = \mu(K_0)$ and $\mu(L_n) = \mu(L_0)$
      for every $n \geq 0$;
    \item $p \nmid W_{K_0}$ and $p \nmid W_{L_0}$.
  \end{enumerate}
\end{lemma}

\begin{proof}
  Choose a prime $\mfp_n$ of $K_n$ above $p$,
  and let $\mfp_0$ be the prime below it in $K_0$.
  The local extension $(K_0)_{\mfp_0}/\Q_p$
  is unramified. The completion of $\Q_n$ at its unique
  prime above $p$ is totally ramified of degree $p^n$
  over $\Q_p$. Therefore
  $(K_n)_{\mfp_n}/(K_0)_{\mfp_0}$
  is totally ramified of degree $p^n$.

  Let $\xi \in \mu(K_n)$ have order coprime to $p$.
  The extension $(K_0)_{\mfp_0}(\xi)/
  (K_0)_{\mfp_0}$ is unramified, so
  $\xi \in (K_0)_{\mfp_0}$.
  The decomposition group of $\mfp_n$ in
  $K_n/K_0$ is all of $\Gal(K_n/K_0)$, since the local
  degree equals $[K_n : K_0] = p^n$.
  Thus every automorphism of $K_n/K_0$ fixes $\xi$,
  and $\xi \in K_0$.

  If a primitive $p$-th root $\xi_p$ belonged to $K_n$,
  then $\Q_p(\xi_p)/\Q_p$ would be a totally ramified
  subextension of $(K_n)_{\mfp_n}/\Q_p$
  of degree $p - 1$. This is impossible because
  $(K_n)_{\mfp_n}/\Q_p$ has ramification index $p^n$.
  Hence $\mu(K_n) = \mu(K_0)$ and $p \nmid W_{K_0}$.
  By Section~\ref{subsec:Stark-hypotheses},
  $\mu(K_n) = \mu(L_n)$ at every level, so
  $\mu(L_n) = \mu(L_0)$ and $p \nmid W_{L_0}$.
\end{proof}

To adjust the norm of a Stark unit from level $n + 1$
to level $n$, we use the norm maps on $\mu(K_n)$,
$\mu(L_n)$, $\mu(K_n^+)$, and $\mu(L_n^+)$.

\begin{corollary}
  \label{cor:norms-roots-of-unity}
  For every $n \geq 0$, norms induce isomorphisms
  \[
    \begin{tikzcd}[row sep = tiny]
      \mu(K_{n + 1}) \arrow[r, "\widesim{}"] & \mu(K_n), \\
      \mu(L_{n + 1}) \arrow[r, "\widesim{}"] & \mu(L_n).
    \end{tikzcd}
  \]
  Moreover, the norms $\rmn_{K_{n + 1}^+/K_n^+}$ and
  $\rmn_{L_{n + 1}^+/L_n^+}$ restrict to the identity on
  $\mu(K_{n + 1}^+) = \mu(K_n^+) = \set{\pm 1}$ and
  $\mu(L_{n + 1}^+) = \mu(L_n^+) = \set{\pm 1}$,
  respectively.
\end{corollary}

\begin{proof}
  For $\xi \in \mu(K_0)$, we have
  $\rmn_{K_{n + 1}/K_n}(\xi) = \xi^p$ because
  $\Gal(K_{n + 1}/K_n)$ fixes $\xi$.
  By Lemma~\ref{lem:cyclotomic-roots-of-unity},
  $\mu(K_{n + 1}) = \mu(K_n) = \mu(K_0)$ and
  $p \nmid W_{K_0}$, so this norm is an automorphism.
  Likewise, $\rmn_{L_{n + 1}/L_n}(\xi) = \xi^p$
  for $\xi \in \mu(L_0)$, and $p \nmid W_{L_0}$.

  Since $\mu(K_n^+) = \mu(L_n^+) = \set{\pm 1}$,
  the norms in the two towers of fixed fields act
  by $\xi \longmapsto \xi^p = \xi$ for
  $\xi \in \set{\pm 1}$, because $p$ is odd.
\end{proof}

Write $W_K := W_{K_0}$ and $W_L := W_{L_0}$.
Let $\Theta_{\Q, n}(0)$ be the Stickelberger element
for $(L_n/\Q, \Sigma_\Q)$.
At each level $n$,
choose $\epsilon_{F, n} \in \calo_{K_n}^\times$
from Theorem~\ref{thm:imaginary-quadratic-Stark}
applied to $(K_n/F, \Sigma_F, w_n)$, and
$\epsilon_{\Q, n} \in \calo_{L_n^+}^\times$ from
Theorem~\ref{thm:rational-Stark} with $S_\Q = \Sigma_\Q$.
For the cyclotomic regulators, use the real places
$\restr{w_n}{L_n^+}$.
We compare the regulators of $\epsilon_{F, n}$ and
$\epsilon_{\Q, n}$ at successive levels and
can adjust them by roots of unity
to obtain the following norm relations.

\begin{lemma}
  \label{lem:source-unit-norm-relations}
  The cyclotomic units $\epsilon_{\Q, n}$
  and elliptic units $\epsilon_{F, n}$ can be chosen so that,
  for every $n \geq 0$,
  \begin{align*}
    \rmn_{L_{n + 1}^+/L_n^+}(\epsilon_{\Q, n + 1})
      &= \epsilon_{\Q, n}, \\
    \rmn_{K_{n + 1}/K_n}(\epsilon_{F, n + 1})
      &= \epsilon_{F, n}.
  \end{align*}
\end{lemma}

\begin{proof}
  Let $\chi$ be a character of $\Gal(K_n/F)$.
  Its inflation to $\Gal(K_{n + 1}/F)$ is given by
  $\widetilde\chi(g) = \chi(\restr{g}{K_n})$.
  Inflation gives
  $L_{\Sigma_F}(s, \widetilde\chi) = L_{\Sigma_F}(s, \chi)$.
  By Lemma~\ref{lem:cyclotomic-roots-of-unity},
  $W_{K_n} = W_{K_{n + 1}} = W_K$.
  To compare $\epsilon_{F, n}$ with the norm of
  $\epsilon_{F, n + 1}$,
  we group the terms in
  Theorem~\ref{thm:imaginary-quadratic-Stark}
  according to the restriction $h := \restr{g}{K_n}$:
  \begin{align*}
    \sum_{h \in \Gal(K_n/F)}\chi(h)
        \log\verts{\rmn_{K_{n + 1}/K_n}
        (\epsilon_{F, n + 1})^h}_{w_n}
      &= \sum_{g \in \Gal(K_{n + 1}/F)}
        \widetilde\chi(g)
        \log\verts{\epsilon_{F, n + 1}^g}_{w_{n + 1}} \\
      &= -W_K L_{\Sigma_F}'(0, \chi) \\
      &= \sum_{h \in \Gal(K_n/F)}\chi(h)
        \log\verts{\epsilon_{F, n}^h}_{w_n}.
  \end{align*}
  By Fourier inversion,
  $\rmn_{K_{n + 1}/K_n}(\epsilon_{F, n + 1})$ and
  $\epsilon_{F, n}$
  have the same absolute value
  at every archimedean place of $K_n$.
  Since both elements are global units, their quotient
  $\rmn_{K_{n + 1}/K_n}(\epsilon_{F, n + 1}) /
  \epsilon_{F, n}$
  belongs to $\mu(K_n)$.

  For a character $\chi$ of $\Gal(L_n^+/\Q)$,
  apply Theorem~\ref{thm:rational-Stark} to
  the inflation of $\chi$ to $\Gal(L_{n + 1}^+/\Q)$;
  here $W_{L_n^+} = W_{L_{n + 1}^+} = 2$.
  Grouping the summands over the lifts of each
  $h \in \Gal(L_n^+/\Q)$ again gives equal absolute values
  for the norm of $\epsilon_{\Q, n + 1}$ and
  $\epsilon_{\Q, n}$ at every real place of $L_n^+$.
  Thus $\rmn_{L_{n + 1}^+/L_n^+}(\epsilon_{\Q, n + 1})/
  \epsilon_{\Q, n} \in \mu(L_n^+) = \set{\pm 1}$.

  By Corollary~\ref{cor:norms-roots-of-unity},
  $\rmn_{K_{n + 1}/K_n}$ restricts to an automorphism of
  $\mu(K_{n + 1})$.
  Multiplying $\epsilon_{F, n + 1}$ by a root of unity
  therefore gives norm $\epsilon_{F, n}$.
  Applying $\rmn_{L_{n + 1}^+/L_n^+}$ to
  $\mu(L_{n + 1}^+) = \set{\pm 1}$ allows the same adjustment
  for $\epsilon_{\Q, n + 1}$.
  These adjustments also preserve
  Definition~\ref{def:global-rank-one-Stark-unit}(iii).
  Indeed, a $W_K$-th root of a correcting factor in
  $\mu(K_{n + 1})$ belongs to $\mu_{W_K^2}$, and a square
  root of a correcting sign belongs to $\mu_4$.
  The fields
  $K_{n + 1}(\epsilon_{F, n + 1}^{1/W_K}, \mu_{W_K^2})$
  and $L_{n + 1}^+(\sqrt{\epsilon_{\Q, n + 1}}, \mu_4)$
  are abelian over $F$ and $\Q$, respectively.
  They contain the roots required for the adjusted
  units, so the adjusted units still satisfy (iii).
  Induction constructs $(\epsilon_{F, n})_n$ and
  $(\epsilon_{\Q, n})_n$ with the stated norms.
\end{proof}

From now on, we choose
$\epsilon_{F, n}$ and $\epsilon_{\Q, n}$
as in Lemma~\ref{lem:source-unit-norm-relations}.
For each $n \geq 0$,
let $u_{\R, n} := \epsilon_{\Q, n}$,
$u_n := \epsilon_{F, n}^{W_L\Theta_{\Q, n}(0)}$, and
$u_{\C, n} :=
\rmn_{K_n/K_n^+}(\epsilon_{F, n})^3/
\rmn_{K_n/L_n^+}(\epsilon_{F, n})$.
These are the units in
Equations~\ref{eq:rank-one-unit-definition} and
\ref{eq:rank-two-unit-definition},
formed from $\epsilon_{F, n}$ and $\epsilon_{\Q, n}$.
The group-ring actions and field norms defining $u_n$,
$u_{\R, n}$, and $u_{\C, n}$ commute with the norms
down the cyclotomic tower.
We therefore obtain the
following relations from
Lemma~\ref{lem:source-unit-norm-relations}.

\begin{theorem}
  \label{thm:finite-norm-relations}
  For every $n \geq 0$,
  the units $u_n$ satisfy
  \[
    \rmn_{K_{n + 1}/K_n}(u_{n + 1}) = u_n.
  \]
  Furthermore, for every $n \geq 0$, the units
  $u_{\R, n}$ and $u_{\C, n}$ satisfy
  \begin{align*}
    \rmn_{K_{n + 1}^+/K_n^+}(u_{\R, n + 1})
    &= u_{\R, n}, \\
    \rmn_{K_{n + 1}^+/K_n^+}(u_{\C, n + 1})
    &= u_{\C, n}.
  \end{align*}

\end{theorem}

\begin{proof}
  For a character $\chi$ of $\Gal(L_n/\Q)$, let
  $\widetilde\chi$ be its inflation to
  $\Gal(L_{n + 1}/\Q)$.
  Inflation gives $L_{\Sigma_\Q}(s, \widetilde\chi)
  = L_{\Sigma_\Q}(s, \chi)$.
  The homomorphism
  $\Z[\Gal(L_{n + 1}/\Q)] \longrightarrow
  \Z[\Gal(L_n/\Q)]$ induced by restriction to $L_n$
  therefore sends $W_L\Theta_{\Q, n + 1}(0)$ to
  $W_L\Theta_{\Q, n}(0)$.
  The integral group-ring actions commute with
  $\rmn_{K_{n + 1}/K_n}$ after restricting automorphisms
  to $K_n$.
  Since $u_n = \epsilon_{F, n}^{W_L\Theta_{\Q, n}(0)}$,
  Lemma~\ref{lem:source-unit-norm-relations} proves
  the first identity.

  The extensions $L_{n + 1}^+/L_n^+$ and $K_n/L_n^+$
  have degrees $p$ and $6$, respectively, so they are
  linearly disjoint.
  Consequently, the norms defining
  $u_{\C, n} = \frac{\rmn_{K_n/K_n^+}(\epsilon_{F, n})^3}
  {\rmn_{K_n/L_n^+}(\epsilon_{F, n})}$
  commute with the norms down the tower. This proves
  the assertion for $u_{\C, n}$.

  For $u_{\R, n} = \epsilon_{\Q, n}$, the restriction of
  $\rmn_{K_{n + 1}^+/K_n^+}$ to $L_{n + 1}^+$ is
  $\rmn_{L_{n + 1}^+/L_n^+}$.
  The required equality is therefore the norm relation
  for $\epsilon_{\Q, n}$ in
  Lemma~\ref{lem:source-unit-norm-relations}.
\end{proof}

Write $\epsilon_{n, \Q}^{(2)}$ for the element in
Equation~\ref{eq:rational-rank-two-element} for $K_n^+/k$
with $S = \Sigma_k$.
To compare successive levels, set $G_n := \Gal(K_n^+/k)$
and $V_n := \Q\otimes_\Z\calo_{K_n^+, \Sigma_k}^\times$.
Let $\rmn_n : V_{n + 1} \longrightarrow V_n$ be the
$\Q$-linear extension of $\rmn_{K_{n + 1}^+/K_n^+}$.
The norm induces the map
\[
  \begin{tikzcd}[row sep = tiny]
    \bigwedge_{\Q[G_{n + 1}]}^2 V_{n + 1}
    \arrow[r, "\rmn_n \wedge \rmn_n"] &
    \bigwedge_{\Q[G_n]}^2 V_n \\
    x\wedge y \arrow[r, mapsto] & \rmn_n(x)\wedge \rmn_n(y).
  \end{tikzcd}
\]
The scalars act through the homomorphism
\[
  \begin{tikzcd}[row sep = tiny]
    \Q[G_{n + 1}] \arrow[r, "\pi_n"] & \Q[G_n] \\
    \sigma \arrow[r, mapsto] & \restr{\sigma}{K_n^+},
  \end{tikzcd}
\]
extended $\Q$-linearly. Thus
$\rmn_n(ax) = \pi_n(a)\rmn_n(x)$ for
$a \in \Q[G_{n + 1}]$ and $x \in V_{n + 1}$.

\begin{corollary}
  \label{cor:rank-two-element-norm}
  For every $n \geq 0$, we have
  \[
    (\rmn_n\wedge \rmn_n)(\epsilon_{n + 1, \Q}^{(2)})
    = \epsilon_{n, \Q}^{(2)}.
  \]
\end{corollary}

\begin{proof}
  By Lemma~\ref{lem:cyclotomic-roots-of-unity},
  $W_{K_n} = W_K$ for every $n$.
  Hence the factor $1/(6W_K)$ does not change with $n$.

  Write $e_{\Sigma_k, 2}^{(n)}$ for the idempotent
  $e_{\Sigma_k, 2}$ in $\Q[G_n]$.
  For a character $\chi$ of $G_n$, define its inflation by
  $\widetilde\chi(\sigma) := \chi(\restr{\sigma}{K_n^+})$
  for $\sigma \in G_{n + 1}$.
  Then $L_{\Sigma_k}(s, \widetilde\chi)
  = L_{\Sigma_k}(s, \chi)$.
  Extend $\pi_n$ to $\C[G_{n + 1}] \longrightarrow \C[G_n]$.
  The definitions of the character idempotents give
  $\pi_n(e_{\widetilde\chi^{-1}}) = e_{\chi^{-1}}$.
  For a character $\psi$ of $G_{n + 1}$ that is nontrivial
  on $\Gal(K_{n + 1}^+/K_n^+)$, the sum of $\psi$ over this
  subgroup is zero, so $\pi_n(e_{\psi^{-1}}) = 0$.
  Since inflation preserves the order of vanishing,
  we can sum these identities over the characters of
  order of vanishing $2$ to obtain
  $\pi_n(e_{\Sigma_k, 2}^{(n + 1)}) = e_{\Sigma_k, 2}^{(n)}$.
  We now use the norm relations in
  Theorem~\ref{thm:finite-norm-relations} to evaluate
  Equation~\ref{eq:rational-rank-two-element} under
  $\rmn_n\wedge \rmn_n$:
  \begin{align*}
    (\rmn_n\wedge \rmn_n)(\epsilon_{n + 1, \Q}^{(2)})
    &= \frac{1}{6W_K}e_{\Sigma_k, 2}^{(n)}
    \bigl(\rmn_n(u_{\R, n + 1})\wedge
    \rmn_n(u_{\C, n + 1})\bigr) \\
    &= \frac{1}{6W_K}e_{\Sigma_k, 2}^{(n)}
    (u_{\R, n}\wedge u_{\C, n}) \\
    &= \epsilon_{n, \Q}^{(2)}. \qedhere
  \end{align*}
\end{proof}

The norm relations also hold for the Stark units
of Theorem~\ref{thm:intro-rank-one}.

\begin{corollary}
  \label{cor:Stark-unit-norm}
  For the tower and $\Sigma_k$ fixed above, the
  Stark units $\epsilon_n \in \calo_{K_n}^\times$ for
  $(K_n/k, \Sigma_k, w_n)$ can be chosen so that
  $\rmn_{K_{n + 1}/K_n}(\epsilon_{n + 1}) = \epsilon_n$
  for every $n \geq 0$.
\end{corollary}

\begin{proof}
  For each $n$, choose a Stark unit $\epsilon_n$ using
  Theorem~\ref{thm:intro-rank-one}.
  For a character $\chi$ of $\Gal(K_n/k)$, inflate
  $\chi$ to $\Gal(K_{n + 1}/k)$ and compare the regulators
  as in the proof of
  Lemma~\ref{lem:source-unit-norm-relations}.
  Fourier inversion gives equal absolute values for
  $\rmn_{K_{n + 1}/K_n}(\epsilon_{n + 1})$ and $\epsilon_n$
  at every place above the complex place of $k$.
  Above the real place of $k$, both units have absolute
  value $1$ by
  Definition~\ref{def:global-rank-one-Stark-unit}(ii).
  Thus $\rmn_{K_{n + 1}/K_n}(\epsilon_{n + 1})/\epsilon_n$
  is a global unit with absolute value $1$ at every
  archimedean place of $K_n$, so it belongs to $\mu(K_n)$.
  Corollary~\ref{cor:norms-roots-of-unity} allows an
  inductive choice with the stated norms.

  For $\xi \in \mu(K_n)$, multiplying $\epsilon_n$ by $\xi$
  preserves its logarithms. Moreover,
  $K_n((\xi\epsilon_n)^{1/W_K})
  \subset K_n(\epsilon_n^{1/W_K}, \mu_{W_K^2})$.
  The field on the right is the compositum of
  $K_n(\epsilon_n^{1/W_K})$ and $k(\mu_{W_K^2})$.
  Both extensions are abelian over $k$, so
  $K_n((\xi\epsilon_n)^{1/W_K})/k$ is abelian.
\end{proof}

\subsection{Units and ideal classes}
\label{subsec:Iwasawa-modules}

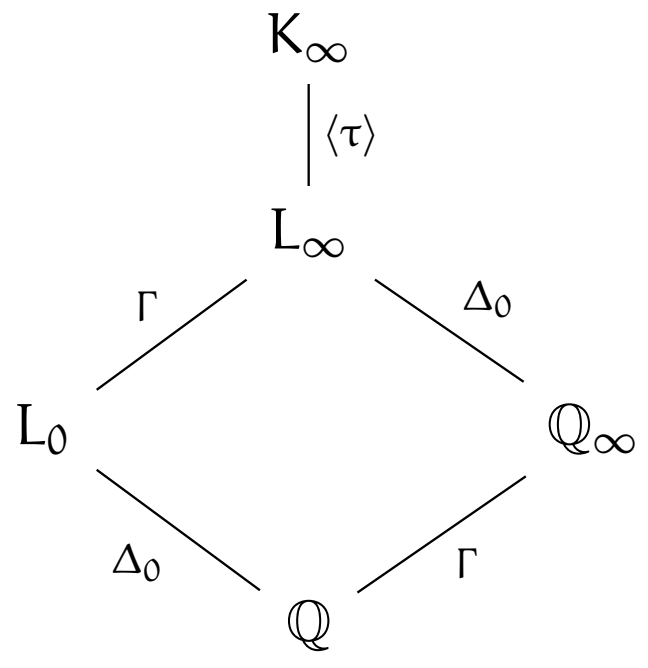
\begin{figure}[ht]
  \[
    \begin{tikzcd}[row sep = 1.8em, column sep = 2.5em]
      & K_\infty \arrow[d, dash, "\abrac{\tau}"] & \\
      & L_\infty
      \arrow[dl, dash, "\Gamma"']
      \arrow[dr, dash, "\Delta_0"] & \\
      L_0 \arrow[dr, dash, "\Delta_0"'] &&
      \Q_\infty \arrow[dl, dash, "\Gamma"] \\
      & \Q &
    \end{tikzcd}
  \]
  \caption{The extensions with Galois groups $\Gamma$,
    $\Delta_0$, and $\abrac{\tau}$.}
  \label{figure:Iwasawa-infinite-fields}
\end{figure}

We now compare the inverse limits of units and ideal
classes under norms in the tower $(K_n)_n$.
We continue to use the places above $\Sigma_\Q$ in $k$
and $F$, denoted by $\Sigma_k$ and $\Sigma_F$.
Set $L_\infty := \bigcup_n L_n$ and
$K_\infty := \bigcup_n K_n$.
To describe the Galois action on these inverse limits,
let $\Delta_0 := \Gal(L_0/\Q)$ and
$\Gamma := \Gal(\Q_\infty/\Q)$.
Since $K_0 \cap \Q_\infty = \Q$, restriction to
$\Q_\infty$ identifies $\Gal(L_\infty/L_0)$ with
$\Gamma$, and restriction to $L_0$ identifies
$\Gal(L_\infty/\Q_\infty)$ with $\Delta_0$.
Restriction to $L_\infty$ gives
$\Gal(K_\infty/k) \cong \Gal(L_\infty/\Q)$.
The map $g \longmapsto
(\restr{g}{L_0}, \restr{g}{\Q_\infty})$ then identifies
$\Gal(L_\infty/\Q)$ with $\Delta_0\times\Gamma$.
Note that $\tau$ generates $\Gal(K_\infty/L_\infty)$.
Figure~\ref{figure:Iwasawa-infinite-fields} records these
extensions.

Fix an embedding $\overline{\Q} \hookrightarrow \overline{\Q}_p$,
and use it to regard algebraic character values and
algebraic $L$-values as $p$-adic numbers.
Let $\zeta := \omega_F(\tau)$, where $\omega_F$ is the
nontrivial character of $\Gal(E/F)$ fixed in
Section~\ref{subsec:decomposition-Artin}.
Let $\calo$ be the ring of integers of a finite extension
of $\Q_p$ containing the character values of $\Delta_0$
and $\zeta$.
Set $\Lambda_n := \calo[\Gal(\Q_n/\Q)]$ and
$\Lambda := \varprojlim_n\Lambda_n = \calo[[\Gamma]]$.
The map $\Lambda_{n + 1} \longrightarrow \Lambda_n$
fixes $\calo$ and sends $g \in \Gal(\Q_{n + 1}/\Q)$
to $\restr{g}{\Q_n}$, extending $\calo$-linearly.
The norm maps on units and ideal classes commute
with these Galois actions,
so $\Lambda$ acts on their inverse limits under norms.
Choose a topological generator $\gamma$ of $\Gamma$.
Then $\Lambda$ is a power series ring over $\calo$ in
$\gamma - 1$, hence a regular local ring of dimension $2$.

Fix a character $\vartheta : \Delta_0 \longrightarrow
\calo^\times$.
Since $p \nmid \verts{\Delta_0}$, the idempotent
\[
  e_\vartheta := \frac{1}{\verts{\Delta_0}}
  \sum_{\delta \in \Delta_0}\vartheta(\delta)\delta^{-1}
\]
belongs to $\calo[\Delta_0]$.
For any $\calo[\Delta_0]$-module $M$, the map
$m \longmapsto e_\vartheta m$ is a projection from $M$ onto
$e_\vartheta M = \set{m \in M :
\delta m = \vartheta(\delta)m
\text{ for every } \delta \in \Delta_0}$.
We use the $\Delta_0$-actions induced by the restriction
isomorphisms in Section~\ref{subsec:finite-norm-relations}.
To compare units with ideal classes, form the inverse limits
\begin{align*}
  \calu_\infty &:= \varprojlim_n
    (\calo_{K_n, \Sigma_k}^\times\otimes_\Z\calo), \\
  \cala_\infty &:= \varprojlim_n
    (\Cl(\calo_{K_n, \Sigma_k})
    \otimes_\Z\calo),
\end{align*}

under norms,
and define $\cale_\Sigma := e_\vartheta\calu_\infty$
and $A_\Sigma := e_\vartheta\cala_\infty$.
The transition maps defining $\calu_\infty$ and
$\cala_\infty$ are the field norms and the norms on
ideal classes, respectively:
\[
  \begin{tikzcd}[column sep = 2.8em, row sep = tiny]
    \cdots \arrow[r] &
    \calo_{K_{n + 1}, \Sigma_k}^\times\otimes_\Z\calo
    \arrow[r, "\rmn"] &
    \calo_{K_n, \Sigma_k}^\times\otimes_\Z\calo
    \arrow[r, "\rmn"] & \cdots \\
    \cdots \arrow[r] &
    \Cl(\calo_{K_{n + 1}, \Sigma_k})\otimes_\Z\calo
    \arrow[r, "\rmn"] &
    \Cl(\calo_{K_n, \Sigma_k})\otimes_\Z\calo
    \arrow[r, "\rmn"] & \cdots
  \end{tikzcd}
\]

Tensoring these finitely generated abelian groups with
$\calo$ amounts to taking their $p$-adic completions
and extending scalars from $\Z_p$ to $\calo$.
We write the tensor products of units additively,
so group-ring exponents become scalar multiplication.
Here $\Cl(\calo_{K_n, \Sigma_k})$ is the ideal class
group of $K_n$ modulo the classes of primes above
$\Sigma_k$.

To account for the finite Euler factors omitted at
$\Sigma_\Q$, let $Y_\Sigma$ be the free $\calo$-module
on the finite primes of $K_\infty$ above $\Sigma_\Q$,
with $\Gal(K_\infty/\Q)$ acting by permutation of the primes.
Define the $\calo$-linear augmentation map $\alpha_\Sigma$ by
\[
  \begin{tikzcd}[row sep = tiny]
    Y_\Sigma \arrow[r, "\alpha_\Sigma"] & \calo \\
    \sum_w a_w w \arrow[r, mapsto] & \sum_w a_w;
  \end{tikzcd}
\]
here $w$ runs over the finite primes of $K_\infty$ above
$\Sigma_\Q$, and $a_w \in \calo$.
Thus each basis element $w$ maps to $1$.
Permuting these primes leaves
the sum of the coefficients unchanged,
so $\alpha_\Sigma$ is Galois-equivariant when
$\calo$ has trivial Galois action.
Set $\calb_\infty := \ker\alpha_\Sigma$ and
$B_\Sigma := e_\vartheta\calb_\infty$.
The primes above $p$ are totally ramified in
$K_\infty/K_0$.
At every finite prime of $K_0$ not above $p$,
Frobenius has infinite order in $\Gamma$:
the residue cardinality is an integer greater than $1$,
hence not a $p$-adic root of unity.
Thus every prime of $K_0$ above $\Sigma_\Q$
has an open decomposition subgroup in $\Gamma$.
Consequently, $Y_\Sigma$ is finite free over $\calo$
and torsion over $\Lambda$.

To compute $\rank_\Lambda\cale_\Sigma$, we also need
the archimedean places. Let $X_n$ be the group of formal sums
$\sum_{w \mid \Sigma_\Q} a_w w$ with $a_w \in \Z$
and $\sum_w a_w = 0$, where $w$ runs over all places
of $K_n$ above $\Sigma_\Q$.
The transition map $X_{n + 1} \longrightarrow X_n$
sends $\sum_w a_w w$ to $\sum_w a_w\restr{w}{K_n}$,
where $w$ runs over the places of $K_{n + 1}$
above $\Sigma_\Q$.
The sets $\Sigma_k$ and $\Sigma_F$
give the same places of $K_n$ as $\Sigma_\Q$.
The relation $\sum_w a_w = 0$ reflects the product
formula for $\Sigma_k$-units, recalled in
Section~\ref{subsec:decomposition-Artin}.

Set $r_\vartheta := 1$ if $\vartheta(c) = -1$ and
$r_\vartheta := 2$ if $\vartheta(c) = 1$.
The archimedean places determine this $\Lambda$-rank:
the complex place of $k$ contributes one for every
$\vartheta$, and the real place contributes one
precisely when $\vartheta(c) = 1$.

\begin{lemma}
  \label{lem:Iwasawa-freeness}
  The module $\cale_\Sigma$ is free of rank $r_\vartheta$
  over $\Lambda$. Moreover, $A_\Sigma$ and $B_\Sigma$
  are finitely generated torsion $\Lambda$-modules.
\end{lemma}

\begin{proof}
  The assertion for $A_\Sigma$ follows from Iwasawa's
  \cite{iwasawa-1973} torsion theorem for class groups;
  see also Washington \cite[Chapter~13]{washington}.
  Since $B_\Sigma \subset Y_\Sigma$, the module $B_\Sigma$
  is also finitely generated and torsion over $\Lambda$.

  To prove freeness, we apply Bullach--Hofer
  \cite[Proposition~3.2(b)--(c) and Remark~3.3]{bullach-hofer}
  to $K_\infty/k$, with $\Sigma = \Sigma_k$ and
  $T = \varnothing$ in their notation.
  In general, they use $\Sigma_k$-units congruent to $1$
  at primes above a finite set $T$ of finite places of $k$
  outside $\Sigma_k$, and require the tensor product of
  this group with $\Z_p$ to be torsion-free at every level.
  For $T = \varnothing$, the group at level $n$ is
  $\calo_{K_n, \Sigma_k}^\times$. Lemma
  \ref{lem:cyclotomic-roots-of-unity} gives $p \nmid W_{K_n}$,
  so $\calo_{K_n, \Sigma_k}^\times\otimes_\Z\Z_p$
  is $\Z_p$-torsion-free. Thus $T = \varnothing$ is allowed;
  their ray class groups then reduce to the
  $\Sigma_k$-class groups used to define $A_\Sigma$.

  We use the presentation of Bullach--Hofer
  \cite[Proposition~3.2(b)--(c) and Remark~3.3]{bullach-hofer}
  to express the units and ideal classes
  through a common map as follows:
  after extending scalars to $\calo$ and applying
  $e_\vartheta$,
  we obtain a $\Lambda$-linear map
  $f : \Lambda^d \longrightarrow \Lambda^d$
  for some $d$, with kernel $\cale_\Sigma$
  and an exact sequence
  \[
    0 \longrightarrow A_\Sigma
    \longrightarrow \coker f
    \longrightarrow
    e_\vartheta\varprojlim_n(X_n\otimes_\Z\calo)
    \longrightarrow 0.
  \]

  Consider the free $\calo$-modules on the archimedean
  places of $K_n$.
  The $\vartheta$-component of their
  inverse limit is free of rank $r_\vartheta$ over $\Lambda$.
  Choose a basis to identify this component with
  $\Lambda^{r_\vartheta}$.
  Retaining the archimedean coordinates defines a map
  from $e_\vartheta\varprojlim_n(X_n\otimes_\Z\calo)$
  to $\Lambda^{r_\vartheta}$.
  The kernel consists of compatible formal sums supported
  at finite places whose coefficients sum to $0$.
  It is therefore $e_\vartheta\ker\alpha_\Sigma = B_\Sigma$.

  To prove surjectivity, let $(a_n)_n$ be a compatible
  system of archimedean sums in the
  $\vartheta$-component.
  Their augmentations have a common value $b \in \calo$,
  since restriction of places preserves augmentation.
  Choose compatible primes $\mfp_n$ of $K_n$ above $p$.
  Then $(a_n - b\mfp_n)_n$ belongs to
  $\varprojlim_n(X_n\otimes_\Z\calo)$, and
  $\bigl(e_\vartheta(a_n - b\mfp_n)\bigr)_n$
  has archimedean coordinates $(a_n)_n$.
  This proves surjectivity.

  Let $q : \Lambda^d \longrightarrow
  \Lambda^{r_\vartheta}$ be the composite surjection
  \[
    \Lambda^d
      \longrightarrow \coker f
    \longrightarrow
      e_\vartheta\varprojlim_n(X_n\otimes_\Z\calo)
    \longrightarrow \Lambda^{r_\vartheta}.
  \]
  Since $\Lambda^{r_\vartheta}$ is free, $q$ splits.
  Thus $\ker q$ is a direct summand of $\Lambda^d$;
  as $\Lambda$ is local, $\ker q$ is free of rank
  $d - r_\vartheta$.
  Moreover, $qf = 0$, so $f$ takes values in $\ker q$.
  Define the kernel $Z :=
  \ker(\coker f \longrightarrow \Lambda^{r_\vartheta})$.
  Restricting the target of $f$ to $\ker q$ and choosing
  a basis gives exact sequences
  \[
    \begin{gathered}
      0 \longrightarrow \cale_\Sigma
      \longrightarrow \Lambda^d
      \xrightarrow{f} \Lambda^{d - r_\vartheta}
      \longrightarrow Z \longrightarrow 0, \\
      0 \longrightarrow A_\Sigma
      \longrightarrow Z \longrightarrow B_\Sigma
      \longrightarrow 0.
    \end{gathered}
  \]
  The module $Z$ is torsion because
  $A_\Sigma$ and $B_\Sigma$ are, so
  $\rank_\Lambda\cale_\Sigma = r_\vartheta$.
  The image of $f$ is torsion-free, hence has
  depth at least $1$.
  By the depth lemma, $\cale_\Sigma$ has depth $2$.
  Since $\Lambda$ is regular local of dimension $2$,
  we conclude from the Auslander--Buchsbaum formula that
  $\cale_\Sigma$ is free
  (see the Stacks Project
  \cite[Tags~00LX and 090V]{stacks-project}).
\end{proof}

To compare quotients of $\cale_\Sigma$ with $A_\Sigma$
and $B_\Sigma$, we use characteristic ideals.
For a finitely generated torsion $\Lambda$-module $N$,
define $\Char_\Lambda(N)$ to be the principal ideal of
$\Lambda$ whose exponent at each height-$1$ prime
$\mfq$ is
$\length_{\Lambda_{\mfq}}(N_{\mfq})$.
Characteristic ideals are multiplicative in short exact
sequences. For $0 \neq f \in \Lambda$, we have
$\Char_\Lambda(\Lambda/(f)) = (f)$.
These facts also follow from the structure
theorem in Washington \cite[Section~13.2]{washington}.

\subsection{Stark units and the Iwasawa main conjecture}
\label{subsec:Iwasawa-Stark-refinement}

We use the norm relations in
Theorem~\ref{thm:finite-norm-relations} to define
a submodule of $\cale_\Sigma$.
Let $\mathbf{u}$, $\mathbf{u}_\R$, and $\mathbf{u}_\C$
be the elements of $\calu_\infty$ represented by
$(u_n)_n$, $(u_{\R, n})_n$, and $(u_{\C, n})_n$,
respectively. We regard $u_{\R, n}$ and $u_{\C, n}$ as
units of $K_n$ through $K_n^+ \subset K_n$.
Define $\mathbf{u}_\vartheta := e_\vartheta\mathbf{u}$ and
$\mathbf{u}_{v, \vartheta} := e_\vartheta\mathbf{u}_v$
for $v \in \set{\R, \C}$, and define
\[
  \calc_\vartheta^\Sigma :=
  \begin{cases}
    \Lambda\mathbf{u}_\vartheta
    & \text{if $\vartheta(c) = -1$}, \\
    \Lambda\mathbf{u}_{\R, \vartheta}
    + \Lambda\mathbf{u}_{\C, \vartheta}
    & \text{if $\vartheta(c) = 1$}.
  \end{cases}
\]
Since $p \nmid W_K$, every root of unity in $K_n$ has
zero image in $\calo_{K_n, \Sigma_k}^\times\otimes_\Z\calo$.
Thus $\calc_\vartheta^\Sigma$ is independent of the choices
of $\epsilon_{F, n}$ and $\epsilon_{\Q, n}$ in
Lemma~\ref{lem:source-unit-norm-relations}.

We compare $(u_n)_n$ with the norm-compatible Stark units
$(\epsilon_n)_n$ of Corollary~\ref{cor:Stark-unit-norm}.
For odd $\vartheta$, define
$\boldsymbol\epsilon_\vartheta :=
e_\vartheta(\epsilon_n\otimes 1)_n \in \cale_\Sigma$.

\begin{lemma}
  \label{lem:Iwasawa-Stark-generators}
  If $\vartheta$ is odd, then
  $\mathbf{u}_\vartheta = W_L\boldsymbol\epsilon_\vartheta$
  and
  $\Lambda\mathbf{u}_\vartheta = \Lambda\boldsymbol\epsilon_\vartheta$.
\end{lemma}

\begin{proof}
  At each level, we have $u_n = \xi_n\epsilon_n^{W_L}$
  with $\xi_n \in \mu(K_n)$ by
  Proposition~\ref{prop:relative-rank-one}.
  Since $W_{K_n} = W_K$ and $p \nmid W_K$ by
  Lemma~\ref{lem:cyclotomic-roots-of-unity},
  $\xi_n\otimes 1 = 0$ in
  $\calo_{K_n, \Sigma_k}^\times\otimes_\Z\calo$.
  Thus $\mathbf{u}_\vartheta
  = W_L\boldsymbol\epsilon_\vartheta$
  in the inverse limit. Since $W_L \in \calo^\times$,
  we have
  $\Lambda\mathbf{u}_\vartheta
  = \Lambda\boldsymbol\epsilon_\vartheta$.
\end{proof}

The Iwasawa main conjecture over $k$ compares
$\cale_\Sigma/\calc_\vartheta^\Sigma$ with $A_\Sigma$
and $B_\Sigma$.

\begin{theorem}
  \label{thm:explicit-Iwasawa-Stark}
  The submodule $\calc_\vartheta^\Sigma$ is free of rank
  $r_\vartheta$ over $\Lambda$, and the quotient
  $\cale_\Sigma/\calc_\vartheta^\Sigma$ is torsion.
  Moreover,
  \[
    \Char_\Lambda(\cale_\Sigma/\calc_\vartheta^\Sigma)
    = \Char_\Lambda(A_\Sigma)\Char_\Lambda(B_\Sigma).
  \]
\end{theorem}

\begin{remark}
  If $\vartheta$ is nontrivial on the stabilizer in $\Delta_0$
  of each finite prime of $K_\infty$ above $\Sigma_\Q$,
  then $e_\vartheta Y_\Sigma = 0$ and hence $B_\Sigma = 0$.
  In this case,
  $\Char_\Lambda(\cale_\Sigma/\calc_\vartheta^\Sigma)
  = \Char_\Lambda(A_\Sigma)$.
\end{remark}

To prove Theorem~\ref{thm:explicit-Iwasawa-Stark}, we
first decompose $\cale_\Sigma$, $A_\Sigma$, and $B_\Sigma$
into components over $\Q$ and $F$.
In Lemma~\ref{lem:elliptic-unit-conjugation}, we identify
the generators of $\calc_\vartheta^\Sigma$ in those components.
We then apply the Iwasawa main conjectures over $\Q$
and $F$ to obtain the required characteristic ideal identity.

\subsection{The cyclotomic and elliptic components}
\label{subsec:integral-Artin-comparison}

We decompose $\cale_\Sigma$, $A_\Sigma$, and $B_\Sigma$
into contributions over $\Q$ and $F$ by using the action
of $\tau$ on $\calu_\infty$, $\cala_\infty$, and
$\calb_\infty$. Since $\tau$ generates
$\Gal(K_\infty/L_\infty)$, we will identify the
$\tau$-invariant submodules with modules for $L_\infty/\Q$.
The nontrivial $\tau$-eigenspaces give the contributions
from $K_\infty/F$, as for the characters in
Section~\ref{subsec:decomposition-Artin}.

To carry out this decomposition on the $\vartheta$-components,
we separate the action of $c$ from that of $\Gal(L_0/F)$.
Let $\Delta^+ := \Gal(L_0^+/\Q)$, identified by restriction
with $\Gal(L_0/F) \subset \Delta_0$.
Since $L_0 = FL_0^+$, we have
$\Delta_0 = \abrac{c}\times\Delta^+$.
Thus $\vartheta$ is determined by $\vartheta(c)$ and its
restriction $\vartheta^+$ to $\Delta^+$.
We select the $\vartheta^+$-component using
$e_{\vartheta^+} := \verts{\Delta^+}^{-1}
\sum_{\delta \in \Delta^+}\vartheta^+(\delta)\delta^{-1}
\in \calo[\Delta^+]$.
Since $E \cap L_0^+ = \Q$, we also have
$\Gal(K_\infty/\Q) \cong S_3\times\Delta^+\times\Gamma$.
In particular, $c$ and $\tau$ commute with the actions of
$\Delta^+$ and $\Gamma$.

We first describe the action of $c$ and $\tau$ on an
arbitrary $\calo[S_3]$-module $V$.
Since $\tau^3 = 1$ and $3 \in \calo^\times$, the module
$V$ splits into eigenspaces on which $\tau$ acts by $1$,
$\zeta$, or $\zeta^{-1}$.
The elements $P_0 := (1 + \tau + \tau^2)/3$ and
$P_\zeta := (1 + \zeta^{-1}\tau + \zeta^{-2}\tau^2)/3$
act on $V$ as projections onto its $\tau$-invariant
submodule and its $\zeta$-eigenspace, respectively.
The relation $c\tau c = \tau^{-1}$ shows that
$cP_\zeta V$ is the $\zeta^{-1}$-eigenspace. Thus
\[
  V = P_0V \oplus P_\zeta V \oplus cP_\zeta V.
\]

Within the sum $(1 - P_0)V$ of the two nontrivial
$\tau$-eigenspaces, we now prescribe the action of $c$.
For $s \in \set{1, -1}$,
define $e_s := (1 + s c)/2 \in \calo[S_3]$.
The elements of $e_s(1 - P_0)V$ are exactly
$v + s c v$ with $v \in P_\zeta V$.
The following maps are mutually inverse:
\[
  \begin{tikzcd}[column sep = large, row sep = small]
    P_\zeta V \arrow[r, shift left = 1ex, "j_s"] &
    e_s(1 - P_0)V
    \arrow[l, shift left = 1ex, "P_\zeta"] \\
    v \arrow[r, mapsto] & v + s c v.
  \end{tikzcd}
\]
Indeed, $P_\zeta(v + s c v) = v$, and for
$x \in e_s(1 - P_0)V$ we have
$x = P_\zeta x + s cP_\zeta x$.

Because $e_{\vartheta^+}$ commutes with $c$ and $\tau$,
this calculation applies to
$V = e_{\vartheta^+}\calu_\infty$,
$V = e_{\vartheta^+}\cala_\infty$, and
$V = e_{\vartheta^+}\calb_\infty$.
Take $s := \vartheta(c)$.
The identity $e_\vartheta = e_s e_{\vartheta^+}$ then
identifies $e_sV$ with $\cale_\Sigma$, $A_\Sigma$, and
$B_\Sigma$, respectively.

For the $\tau$-invariant contribution, apply $P_0e_\vartheta$
to $\calu_\infty$, $\cala_\infty$, and
$\calb_\infty$, and denote the images by
$\cale_\Q$, $A_\Q$, and $B_\Q$, respectively.
We identify these with the $\vartheta$-components for
$L_\infty/\Q$. On units and ideal classes, extension
from $L_n$ to $K_n$ followed by $\rmn_{K_n/L_n}$ is
multiplication by $3$; the reverse composite is
$1 + \tau + \tau^2$.
Since $3 \in \calo^\times$, extension identifies
$\varprojlim_n(\calo_{L_n, \Sigma_\Q}^\times\otimes_\Z\calo)$
with $P_0\calu_\infty$, and
$\varprojlim_n(\Cl(\calo_{L_n, \Sigma_\Q})
\otimes_\Z\calo)$ with $P_0\cala_\infty$.
For the free modules on finite primes above $\Sigma_\Q$,
identify a prime of $L_\infty$ with the average of the
primes of $K_\infty$ above it.
The number of such primes divides $3$, so this average
is defined over $\calo$ and preserves augmentation.

For the $\zeta$-eigenspace, apply $P_\zeta e_{\vartheta^+}$
to $\calu_\infty$, $\cala_\infty$, and
$\calb_\infty$, and denote the images by
$\cale_F$, $A_F$, and $B_F$, respectively.
These are the components for $K_\infty/F$ on which
$\tau$ acts by $\zeta$ and $\Delta^+$ acts by $\vartheta^+$.
Since $\tau$ acts trivially on the target $\calo$ of
$\alpha_\Sigma: Y_\Sigma \longrightarrow \calo$, we have
$P_\zeta\calo = 0$. Applying $P_\zeta e_{\vartheta^+}$
to $\calb_\infty = \ker\alpha_\Sigma$ therefore gives
$B_F = P_\zeta e_{\vartheta^+}Y_\Sigma$.

The identity on the $P_0$-component and $j_s$ on the
$P_\zeta$-component give the following decompositions of
$\cale_\Sigma$, $A_\Sigma$, and $B_\Sigma$.
These separate the contributions over $\Q$ and $F$ to the
Iwasawa main conjecture, just as the decomposition of
partial $L$-functions in
Proposition~\ref{prop:local-Artin-factorization} separates the
contributions to the Stark regulators.

\begin{proposition}
  \label{prop:derived-Artin-decomposition}
  The following three maps are isomorphisms of $\Lambda$-modules:
  \[
    \begin{tikzcd}[column sep = 9pt, row sep = tiny]
      \cale_\Q\oplus\cale_F \arrow[r, "\widesim{}"]
      & \cale_\Sigma
      & A_\Q\oplus A_F \arrow[r, "\widesim{}"] & A_\Sigma
      & B_\Q\oplus B_F \arrow[r, "\widesim{}"] & B_\Sigma
      \\
      (a, v) \arrow[r, mapsto] & a + j_s(v)
      & (a, v) \arrow[r, mapsto] & a + j_s(v)
      & (a, v) \arrow[r, mapsto] & a + j_s(v).
    \end{tikzcd}
  \]
\end{proposition}

\begin{proof}
  Recall that $s = \vartheta(c)$ and $e_s = (1 + s c)/2$.
  For any $\calo[S_3]$-module $V$, the idempotent $P_0$
  commutes with $e_s$, so
  \[
    e_sV = P_0e_sV \oplus e_s(1 - P_0)V.
  \]
  For $v \in P_\zeta V$, we have
  $P_\zeta j_s(v) = v$ and $P_0j_s(v) = 0$.
  For $x \in e_sV$, the equality $cx = s x$ gives
  \[
    x = P_0x + P_\zeta x + s cP_\zeta x
    = P_0x + j_s(P_\zeta x).
  \]
  Hence $(a, v) \longmapsto a + j_s(v)$ is an isomorphism
  from $P_0e_sV\oplus P_\zeta V$ to $e_sV$, with inverse
  $x \longmapsto (P_0x, P_\zeta x)$.

  We apply this isomorphism to
  $V = e_{\vartheta^+}\calu_\infty$,
  $V = e_{\vartheta^+}\cala_\infty$, and
  $V = e_{\vartheta^+}\calb_\infty$.
  Since $e_\vartheta = e_s e_{\vartheta^+}$,
  the three sources are
  $\cale_\Q\oplus\cale_F$, $A_\Q\oplus A_F$, and
  $B_\Q\oplus B_F$, respectively, and the targets are
  $\cale_\Sigma$, $A_\Sigma$, and $B_\Sigma$.
  The operators $P_0$, $P_\zeta$, and $j_s$ commute
  with $\Gamma$, so the isomorphisms are $\Lambda$-linear.
\end{proof}

For the characteristic ideal comparisons over $\Q$ and
$F$, we need the ranks of $\cale_\Q$ and $\cale_F$.
The real place of $\Q$ has decomposition group
$\abrac{c}$ in $L_\infty/\Q$, whereas the complex place
of $F$ splits completely in $K_\infty/F$.
The rank calculation in the proof of
Lemma~\ref{lem:Iwasawa-freeness}
therefore gives $\cale_\Q = 0$ for odd $\vartheta$,
while $\cale_\Q$ is free of rank $1$ for even $\vartheta$.
The module $\cale_F$ is free of rank $1$ in both cases.

We next express the generators of $\calc_\vartheta^\Sigma$
in the decomposition of
Proposition~\ref{prop:derived-Artin-decomposition}.
Let $\mathbf{e}_F$ be the image in $\calu_\infty$ of
$(\epsilon_{F, n})_n$ from
Lemma~\ref{lem:source-unit-norm-relations}, and set
$\mathbf{v} := P_\zeta e_{\vartheta^+}\mathbf{e}_F$.
Evaluation at $\vartheta$ identifies
$e_\vartheta\calo[\Delta_0]$ with $\calo$.
For odd $\vartheta$, the elements
$e_\vartheta\Theta_{\Q, n}(0)$ are compatible under
restriction by the proof of
Theorem~\ref{thm:finite-norm-relations}.
Write $t_\vartheta := \varprojlim_n
e_\vartheta\Theta_{\Q, n}(0) \in \Lambda$
for their inverse limit.
The element $t_\vartheta$ belongs to $\Lambda$ because
$W_L\Theta_{\Q, n}(0) \in \Z[\Gal(L_n/\Q)]$ and
$W_L \in \calo^\times$.
For a finite-order character
$\rho: \Gamma \longrightarrow \overline{\Q}_p^\times$,
evaluation at $\rho$ sends the topological generator
$\gamma$ chosen in Section~\ref{subsec:Iwasawa-modules}
to $\rho(\gamma)$ and defines
$\rho(t_\vartheta) \in \overline{\Q}_p$.
By the definition of $\Theta_{\Q, n}(0)$,
\begin{equation}
  \label{eq:Stickelberger-specialization}
  \rho(t_\vartheta)
  = L_{\Sigma_\Q}(0, (\vartheta\rho)^{-1}).
\end{equation}

We compare $\mathbf{u}_\vartheta$ and
$\mathbf{u}_{\C, \vartheta}$ with $\mathbf{v}$.
Multiplying a generator by an element of $\calo^\times$
does not change the generated submodule.

\begin{lemma}
  \label{lem:elliptic-unit-conjugation}
  Fix a character
  $\vartheta : \Delta_0 \longrightarrow \calo^\times$.
  \begin{enumerate}[(i)]
    \item If $\vartheta$ is odd, then
      $\Lambda\mathbf{u}_\vartheta
      = t_\vartheta j_{-1}(\Lambda\mathbf{v})$.
    \item If $\vartheta$ is even, then
      $\mathbf{u}_{\R, \vartheta} \in \cale_\Q$ and
      $\Lambda\mathbf{u}_{\C, \vartheta}
      = j_1(\Lambda\mathbf{v})$.
  \end{enumerate}
\end{lemma}

\begin{proof}
  Since $\mu(K_n)\otimes_\Z\calo = 0$, we have
  $c\mathbf{e}_F = \tau\mathbf{e}_F$ by
  Equation~\ref{eq:finite-elliptic-unit-conjugation}.
  It follows that
  $(1 - P_0)e_{\vartheta^+}\mathbf{e}_F
  = \mathbf{v} + \zeta c\mathbf{v}$.
  Also, $cP_0\mathbf{e}_F = P_0\tau\mathbf{e}_F
  = P_0\mathbf{e}_F$.

  If $\vartheta$ is odd, we use
  Equation~\ref{eq:rank-one-unit-definition} to express
  $\mathbf{u}_\vartheta$ in terms of $\mathbf{v}$:
  \[
    \mathbf{u}_\vartheta
    = \frac{W_L(1 - \zeta)}{2}
    t_\vartheta j_{-1}(\mathbf{v}).
  \]

  If $\vartheta$ is even, then $u_{\R, n} \in L_n^+$ and
  $\mathbf{u}_\C = 3(1 + c)(1 - P_0)\mathbf{e}_F$ by
  Equation~\ref{eq:rank-two-unit-definition}.
  Hence
  $\mathbf{u}_{\C, \vartheta} = 3(1 + \zeta)j_1(\mathbf{v})$.
  The factors $W_L(1 - \zeta)/2$ and $3(1 + \zeta)$
  are units in $\calo$, since
  $(1 - \zeta)(1 - \zeta^{-1}) = 3$ and
  $1 + \zeta = -\zeta^2$.
\end{proof}

\subsection{The Iwasawa main conjecture over
\texorpdfstring{$\Q$}{Q}}
\label{subsec:rational-source-determinant}

The cyclotomic Iwasawa main conjecture of Mazur--Wiles
\cite{mazur-wiles}, in the formulations of Huber--Kings
\cite[Theorem~4.4.1 and Corollary~5.1.3]{huber-kings},
compares ideal classes with norm-compatible cyclotomic
units in even parity and with the $p$-adic $L$-function
in odd parity.
We apply these results with $p > 3$
and omit the finite Euler factors at $\Sigma_\Q$
to prove the following identities.

\begin{proposition}
  \label{prop:rational-source-determinant}
  Fix a character
  $\vartheta : \Delta_0 \longrightarrow \calo^\times$.
  \begin{enumerate}[(i)]
    \item If $\vartheta$ is odd, then $t_\vartheta \neq 0$.
      Moreover,
      \[
        \Char_\Lambda(A_\Q)\Char_\Lambda(B_\Q)
        = (t_\vartheta).
      \]
    \item If $\vartheta$ is even, then
      $\mathbf{u}_{\R, \vartheta} \neq 0$. Moreover,
      \[
        \Char_\Lambda(A_\Q)\Char_\Lambda(B_\Q)
        = \Char_\Lambda
        (\cale_\Q/\Lambda\mathbf{u}_{\R, \vartheta}).
      \]
  \end{enumerate}
\end{proposition}

\begin{proof}
  We apply Huber--Kings
  \cite[Theorem~4.4.1 and Corollary~5.1.3]{huber-kings}.
  Our hypothesis $p > 3$ satisfies the restriction to odd
  primes in both results.
  We first identify the cohomology modules and fix the
  character convention. We then check the finite Euler
  factors at $\Sigma_\Q$. Finally, we prove (i) by comparing
  the $p$-adic $L$-function with $t_\vartheta$, and prove
  (ii) by comparing the cyclotomic element with
  $\mathbf{u}_{\R, \vartheta}$.

  \emph{1. Cohomology and characters.}
  Write $\calo(1) := \calo\otimes_{\Z_p}
  \varprojlim_m\mu_{p^m}$, with $p$-th powers as transition
  maps and Galois acting on the roots of unity.
  For $i \in \set{1, 2}$, let
  $H^i(\calo_{L_n, \Sigma_\Q}, \calo(1))$ denote the
  continuous cohomology, with coefficients in $\calo(1)$,
  of the Galois group of the maximal extension of $L_n$
  unramified outside the places above $\Sigma_\Q$.

  By Kummer theory, we can identify $\cale_\Q$ with
  $e_\vartheta\varprojlim_n
  H^1(\calo_{L_n, \Sigma_\Q}, \calo(1))$,
  with corestriction as transition maps.
  To account for ideal classes, we also use degree $2$
  and set $Z := e_\vartheta\varprojlim_n
  H^2(\calo_{L_n, \Sigma_\Q}, \calo(1))$.
  We relate $Z$ to $A_\Q$ and $B_\Q$ by taking inverse
  limits of the Kummer and Brauer exact sequences in Milne
  \cite[Chapter~II, Proposition~2.1 and
  Remark~2.2(a)]{milne-2006}:
  \[
    0 \longrightarrow A_\Q \longrightarrow Z
    \longrightarrow B_\Q \longrightarrow 0.
  \]
  Thus $\Char_\Lambda(Z)
  = \Char_\Lambda(A_\Q)\Char_\Lambda(B_\Q)$.
  It remains to identify $\Char_\Lambda(Z)$ with
  $(t_\vartheta)$ in (i) and with
  $\Char_\Lambda(\cale_\Q/\Lambda\mathbf{u}_{\R, \vartheta})$
  in (ii).

  To apply the results of Huber--Kings
  \cite[Theorem~4.4.1 and Corollary~5.1.3]{huber-kings}
  to $\cale_\Q$ and $Z$, we match the action of $\vartheta$
  with the convention for Dirichlet characters.
  Let $\chi_\vartheta$ be the Dirichlet character obtained
  by evaluating $\vartheta$ on arithmetic Artin symbols.
  Huber--Kings \cite[Section~1.1]{huber-kings} attach to
  $\chi_\vartheta$ a free rank-$1$ $\calo$-module with
  Galois action through $\vartheta^{-1}$.

  \emph{2. Finite Euler factors.}
  Before treating the two parities separately, we check
  how the characteristic ideal comparisons depend on the
  finite primes in $\Sigma_\Q$.
  Let $\ell \notin \Sigma_\Q$ be a rational prime; we compare
  $\Sigma_\Q$ with $\Sigma_\Q\cup\set{\ell}$.
  Write the arithmetic Frobenius of $\ell$ in
  $\Gal(L_\infty/\Q) = \Delta_0 \times \Gamma$ as
  $(\delta_\ell, \gamma_\ell)$, and set
  $P_\ell := 1 - \vartheta(\delta_\ell)^{-1}
  \gamma_\ell^{-1} \in \Lambda$.
  For every finite-order character $\rho$ of $\Gamma$,
  \[
    \rho(P_\ell)
    = 1 - (\vartheta\rho)^{-1}(\Frob_\ell).
  \]
  Thus $P_\ell$ represents the Euler factor omitted at
  $\ell$ at $s = 0$.

  To compare the cohomology modules for the two sets of
  places, let $Y_\ell$ be the $\vartheta$-component of the
  free $\calo$-module on the primes of $L_\infty$ above
  $\ell$. The permutation action gives
  $Y_\ell \cong \Lambda/(P_\ell)$: arithmetic Frobenius
  fixes a chosen prime, and its relation on the
  $\vartheta$-component is
  $\vartheta(\delta_\ell)\gamma_\ell - 1$,
  which differs from $P_\ell$ by a unit.
  Let $Z'$ be defined as $Z$ with
  $\Sigma_\Q\cup\set{\ell}$ in place of $\Sigma_\Q$.
  The localization sequence used by Huber--Kings
  \cite[proof of Proposition~1.2.9(b)]{huber-kings},
  together with duality, gives
  \[
    0 \longrightarrow Z \longrightarrow Z'
    \longrightarrow Y_\ell \longrightarrow 0.
  \]
  Adjoining $\ell$ to $\Sigma_\Q$ does not change
  $\cale_\Q$: compatibility under norms forces every
  valuation above $\ell$ to be divisible by arbitrarily
  large powers of $p$, and hence to vanish.
  Therefore $\Char_\Lambda(Z') = (P_\ell)\Char_\Lambda(Z)$.
  The interpolation formula for the $p$-adic $L$-function
  and the norm relations for cyclotomic elements in
  Huber--Kings \cite[Sections~3.1 and~4.2]{huber-kings}
  show that adjoining $\ell$ multiplies $t_\vartheta$ by
  $P_\ell$ for odd $\vartheta$ and
  $\mathbf{u}_{\R, \vartheta}$ by $P_\ell$ for even $\vartheta$.
  Since $\Char_\Lambda(Z') = (P_\ell)\Char_\Lambda(Z)$,
  the two characteristic ideal equalities change by the
  same factor $(P_\ell)$.

  \emph{3. Odd characters.}
  Suppose that $\vartheta(c) = -1$.
  The interpolation formula of Huber--Kings
  \cite[Definition~4.2.3]{huber-kings}, with the finite
  Euler factors at $\Sigma_\Q$ omitted, is
  Equation~\ref{eq:Stickelberger-specialization}.
  Hence their $p$-adic $L$-function is $t_\vartheta$.
  Since $\vartheta$ is odd, $\cale_\Q = 0$.
  We can therefore use Huber--Kings
  \cite[Corollary~5.1.3]{huber-kings}, with
  $\chi = \chi_\vartheta^{-1}$ and $r = 0$, in the form
  \[
    \Char_\Lambda(Z) = (t_\vartheta).
  \]
  But we have the equality
  \[
    \Char_\Lambda(A_\Q)\Char_\Lambda(B_\Q)
      = \Char_\Lambda(Z)
  \]
  from the exact sequence in step~1,
  and $t_\vartheta \neq 0$
  since $Z$ is torsion.
  This proves (i).

  \emph{4. Even characters.}
  Suppose that $\vartheta(c) = 1$.
  To prove (ii), we first compare the cyclotomic units
  obtained from roots of unity with $\epsilon_{\Q, n}$.
  We can then use the cyclotomic Iwasawa main conjecture
  to compute
  $\Char_\Lambda(\cale_\Q/\Lambda\mathbf{u}_{\R, \vartheta})$.

  Let $f$ be the conductor of $L_0$ and set
  $m_n := fp^{n + 1}$. Since $p$ is unramified in $L_0$,
  we have $p \nmid f$ and $L_n \subset \Q(\mu_{m_n})$.
  Choose compatible embeddings of $\Q(\mu_{m_n})$ into
  $\C$ extending the real places fixed for $L_n^+$,
  and let $\xi_n$ map to $e^{2\pi i/m_n}$.
  Define
  \[
    \eta_n :=
    \rmn_{\Q(\mu_{m_n})/L_n}(1 - \xi_n)^{
    \prod_{\substack{\ell \in \Sigma_\Q \\
    \ell \nmid fp}}(1 - [\ell, L_n/\Q]^{-1})},
  \]
  where the product ranges over finite primes.
  Huber--Kings \cite[Definitions~3.1.2, 3.1.3,
  and 4.1.7]{huber-kings} construct their cyclotomic element
  by corestriction of the Kummer classes of $1 - \xi_n$.
  The norm relations in \cite[Section~3.1]{huber-kings}
  account for the Euler factors in $\eta_n$.
  Restricting the resulting class from $\Q_n$ to $L_n$ gives
  \[
    \sum_{\delta \in \Delta_0}
    \vartheta(\delta)^{-1}\delta(\eta_n\otimes 1)
    = \verts{\Delta_0}e_\vartheta(\eta_n\otimes 1).
  \]
  Since $p \nmid \lvert\Delta_0\rvert$, the elements
  $\lvert\Delta_0\rvert e_\vartheta(\eta_n\otimes 1)$ and
  $e_\vartheta(\eta_n\otimes 1)$ generate the same
  $\Lambda_n$-submodule.

  To replace $\eta_n$ by $\epsilon_{\Q, n}$, we compare the
  regulators of $\epsilon_{\Q, n}$ and
  $\rmn_{L_n/L_n^+}(\eta_n)$.
  Pairing the terms indexed by $a$ and $m_n - a$ in the
  Hurwitz zeta formula gives the regulator in
  Theorem~\ref{thm:rational-Stark} for
  $\rmn_{L_n/L_n^+}(\eta_n)$, with set $\Sigma_\Q$.
  Thus $\epsilon_{\Q, n}/\rmn_{L_n/L_n^+}(\eta_n)
  \in \set{\pm 1}$, and the sign has zero image after
  tensoring with $\calo$.
  After including the units of $L_n^+$ in $L_n$,
  $\rmn_{L_n/L_n^+}$ acts by $1 + c$.
  Since $\vartheta(c) = 1$, we obtain
  \[
    e_\vartheta(\epsilon_{\Q, n}\otimes 1)
    = 2e_\vartheta(\eta_n\otimes 1).
  \]
  Consequently, restricting the cyclotomic class from
  $\Q_n$ to $L_n$ gives $\lvert\Delta_0\rvert/2$ times
  $e_\vartheta(\epsilon_{\Q, n}\otimes 1)$.

  Since the elements $\eta_n$ are compatible under norms,
  $(e_\vartheta(\eta_n\otimes 1))_n$ is
  contained in $\cale_\Q$.
  With $\chi = \chi_\vartheta$ and twist $1$, we use
  a result of Huber--Kings
  \cite[Proposition~4.2.1(a) and Theorem~4.4.1]{huber-kings}
  in the following form:
  $(e_\vartheta(\eta_n\otimes 1))_n \neq 0$ and
  \[
    \Char_\Lambda(Z) = \Char_\Lambda\bigl(
    \cale_\Q/\Lambda(e_\vartheta(\eta_n\otimes 1))_n
    \bigr).
  \]
  The factor $\lvert\Delta_0\rvert$ in the restriction
  formula does not change the generated submodule,
  since $p \nmid \lvert\Delta_0\rvert$.

  We have proved that
  $\mathbf{u}_{\R, \vartheta}
  = 2(e_\vartheta(\eta_n\otimes 1))_n$.
  Since $2 \in \calo^\times$, we conclude that
  $\mathbf{u}_{\R, \vartheta} \neq 0$ and
  $\Char_\Lambda(Z) = \Char_\Lambda
  (\cale_\Q/\Lambda\mathbf{u}_{\R, \vartheta})$.
  The exact sequence
  $0 \longrightarrow A_\Q \longrightarrow Z
  \longrightarrow B_\Q \longrightarrow 0$
  now proves (ii).

\end{proof}

\subsection{The Iwasawa main conjecture over
\texorpdfstring{$F$}{F}}
\label{subsec:elliptic-Rubin-Stark}

Rubin \cite{rubin-1991} proved Iwasawa main conjectures
relating elliptic units to ideal classes over imaginary
quadratic fields.
We prove an analogous comparison for
$K_\infty/F$ in the rank-$1$ formulation of
Burns--Kurihara--Sano
\cite[Conjecture~3.14]{burns-kurihara-sano}.
Recall that $\mathbf{v} = P_\zeta e_{\vartheta^+}\mathbf{e}_F$,
where $\mathbf{e}_F \in \calu_\infty$ is represented by
$(\epsilon_{F, n})_n$. We compare
$\Char_\Lambda(\cale_F/\Lambda\mathbf{v})$ with
$\Char_\Lambda(A_F)\Char_\Lambda(B_F)$.

We need a relation between elliptic units and ideal classes
that retains the Galois action.
Bullach--Hofer \cite[Theorem~5.8]{bullach-hofer}
prove the $p$-part of the equivariant Tamagawa number
conjecture at $s = 0$
(an equivariant refinement of the
analytic class number formula)
for the extensions $K_n/F$ that we consider here.
We combine Bullach--Hofer
\cite[Theorem~5.8]{bullach-hofer} with the
determinant formula of Burns--Kurihara--Sano
\cite[Section~2D and Proposition~2.6(iii)]
{burns-kurihara-sano} to compare determinant generators
with rank-$1$ Rubin--Stark elements for $K_n/F$.
We identify the Rubin--Stark elements with
$W_K^{-1}\epsilon_{F, n}$ and use the norm relations
to obtain the equality of characteristic ideals over
$\Lambda$.

We recall the rank-$1$ Rubin--Stark element for $K_n/F$
at the complex place of $F$.
In rank $1$, Rubin~\cite[Proposition~2.5]{rubin-1996}
identifies the integral rank-$1$ Rubin--Stark conjecture
with the integral rank-$1$ abelian Stark conjecture,
which Stark~\cite[Theorem~2]{stark-1980}
proved over imaginary quadratic fields
(as recalled in Theorem~\ref{thm:imaginary-quadratic-Stark}).
We compare the Rubin--Stark element for $K_n/F$ with
$W_K^{-1}\epsilon_{F, n}$ using the logarithmic isomorphism
from Section~\ref{subsec:decomposition-Artin}:
\[
  \begin{tikzcd}[row sep = tiny, column sep = large]
    \R\otimes_\Z\calo_{K_n, \Sigma_F}^\times
    \arrow[r, "\lambda_n"] & \R\otimes_\Z X_n \\
    a \arrow[r, mapsto]
    & -\sum_{w \mid \Sigma_F}\log\verts{a}_w\,w
  \end{tikzcd}
\]
We extend $\lambda_n$ $\R$-linearly.
Here,
$X_n = \{\sum_{w \mid \Sigma_F}a_w w \mid
a_w \in \Z,\ \sum_w a_w = 0\}$, as in
Section~\ref{subsec:Iwasawa-modules}.
In the notation of Bullach--Hofer
\cite[Section~2]{bullach-hofer},
we take $T = \varnothing$.

Choose a finite place $w'_n$ of $K_n$ above $\Sigma_F$.
For the complex place $w_n$ fixed in
Section~\ref{subsec:finite-norm-relations}, Bullach--Hofer
\cite[Definition~2.1]{bullach-hofer} define
$\epsilon_{F, n}^{\rs} \in
\R\otimes_\Z\calo_{K_n, \Sigma_F}^\times$ by
\[
  \lambda_n(\epsilon_{F, n}^{\rs})
  = \Theta'_{K_n/F, \Sigma_F}(0)(w_n - w'_n).
\]
Injectivity of $\lambda_n$ makes this element unique.
We use the Rubin--Stark elements $\epsilon_{F, n}^{\rs}$
to prove the characteristic ideal formula for
$\cale_F/\Lambda\mathbf{v}$.

\begin{proposition}
  \label{prop:elliptic-source-identification}
  The element $\mathbf{v}$ is nonzero. Moreover,
  \[
    \Char_\Lambda(\cale_F/\Lambda\mathbf{v})
    = \Char_\Lambda(A_F)\Char_\Lambda(B_F).
  \]
\end{proposition}

\begin{proof}
  We first identify the Rubin--Stark elements with
  $W_K^{-1}\epsilon_{F, n}$. We then construct a nonzero
  element $z \in \cale_F$ from the maximal minors of a
  matrix presenting the units and ideal classes.
  We compare $z$ with $W_K^{-1}\mathbf{v}$ at each layer,
  then use the norm relations to prove
  $\Lambda z = \Lambda\mathbf{v}$.
  Finally, we compute the characteristic ideal of
  $\cale_F/\Lambda z$ from the same matrix.

  \emph{1. Identifying the Rubin--Stark elements.}
  Fix $n \geq 0$. We prove the equality
  $\epsilon_{F, n}^{\rs} = W_K^{-1}\epsilon_{F, n}$ in
  $\R\otimes_\Z\calo_{K_n, \Sigma_F}^\times$;
  the right side is the image of the rational element
  $(1/W_K) \otimes \epsilon_{F, n}$.
  If a nontrivial character $\psi$ of $\Gal(K_n/F)$ occurs
  in the free module on finite places above $\Sigma_F$,
  then $\ord_{s = 0}L_{\Sigma_F}(s, \psi^{-1}) \geq 2$
  by Tate \cite[Chapter~I, Proposition~3.4]{tate-stark-1984}.
  Also, $L'_{\Sigma_F}(0, 1) = 0$ because $\Sigma_F$
  contains at least two finite places.
  Consequently, $\Theta'_{K_n/F, \Sigma_F}(0)$ annihilates
  the free module on these finite places, including $w'_n$.
  By Theorem~\ref{thm:imaginary-quadratic-Stark},
  \[
    \lambda_n(W_K^{-1}\epsilon_{F, n})
    = \Theta'_{K_n/F, \Sigma_F}(0)w_n
    = \lambda_n(\epsilon_{F, n}^{\rs}).
  \]
  The finite coordinates on the left vanish because
  $\epsilon_{F, n}$ is a global unit.
  Since $\lambda_n$ is injective, we conclude that
  $\epsilon_{F, n}^{\rs} = W_K^{-1}\epsilon_{F, n}$,
  independently of $w'_n$.
  We may now regard
  $\epsilon_{F, n}^{\rs}$ as the element
  $W_K^{-1}(\epsilon_{F, n} \otimes 1)
  \in \calo_{K_n, \Sigma_F}^\times \otimes_\Z \calo$
  because $W_K \in \calo^\times$.

  \emph{2. A matrix for the units and ideal classes.}
  As in the proof of Lemma~\ref{lem:Iwasawa-freeness},
  Lemma~\ref{lem:cyclotomic-roots-of-unity} allows us
  to take $T = \varnothing$ in Bullach--Hofer
  \cite[Proposition~3.2(b)--(c) and Remark~3.3(a)]
  {bullach-hofer}.
  We apply their cohomology presentation to $K_\infty/F$;
  after extending scalars to $\calo$ and applying
  $P_\zeta e_{\vartheta^+}$, we obtain
  a map $f : \Lambda^d \longrightarrow \Lambda^d$
  for some integer $d \geq 1$
  with $\ker f = \cale_F$ and an exact sequence
  \[
    0 \longrightarrow A_F \longrightarrow \coker f
      \longrightarrow P_\zeta e_{\vartheta^+}
      \varprojlim_n(X_n\otimes_\Z\calo)
      \longrightarrow 0.
  \]

  To reduce the rank of the target of $f$ by one, we use
  the coordinates at complex places of $K_n$.
  The $P_\zeta e_{\vartheta^+}$-component of the inverse
  limit of the free $\calo$-modules on these places is
  free of rank $1$ over $\Lambda$, with basis
  $P_\zeta e_{\vartheta^+}(w_n)_n$.
  Since $P_\zeta\calo = 0$, projection onto the complex
  coordinates is surjective after applying
  $P_\zeta e_{\vartheta^+}$. We obtain a composite surjection
  \[
    \coker f \longrightarrow P_\zeta e_{\vartheta^+}
    \varprojlim_n(X_n\otimes_\Z\calo)
    \longrightarrow \Lambda,
  \]
  where the last map retains the complex coordinates and
  sends $P_\zeta e_{\vartheta^+}(w_n)_n$ to $1$.
  Denote the composite by $\pi : \coker f \longrightarrow \Lambda$.
  Let $q : \Lambda^d \longrightarrow \Lambda$ be its
  composite with the quotient map to $\coker f$.
  Choose a $\Lambda$-linear splitting of $q$.
  The kernel $\ker q$ is free of rank $d - 1$.
  Since $qf = 0$, the image of $f$ lies in $\ker q$.
  Choose bases of the source $\Lambda^d$ of $f$ and of
  $\ker q$, and let $D$ be the matrix of
  $f : \Lambda^d \longrightarrow \ker q$.
  Writing $Z := \ker\pi$, we obtain
  \[
    \begin{gathered}
      0 \longrightarrow \cale_F \longrightarrow \Lambda^d
      \xrightarrow{D} \Lambda^{d - 1}
      \longrightarrow Z \longrightarrow 0, \\
      0 \longrightarrow A_F \longrightarrow Z
      \longrightarrow B_F \longrightarrow 0.
    \end{gathered}
  \]
  Here $Z \longrightarrow B_F$ records the finite
  coordinates, and its kernel is $A_F$.

  We now choose $z \in \cale_F$ so that the minors of $D$
  will compute $\Char_\Lambda(\cale_F/\Lambda z)$.
  For $1 \leq i \leq d$, let $D^{(i)}$ be the matrix obtained
  from $D$ by deleting column $i$, and set
  \[
    z := \bigl((-1)^{i + 1}\det D^{(i)}\bigr)_{i = 1}^d.
  \]
  By expanding determinants, we obtain $Dz = 0$.
  Since $\cale_F = \ker f = \ker D \subset \Lambda^d$,
  we may regard $z$ as an element of $\cale_F$.
  Since $A_F$ and $B_F$ are torsion, $Z$ is torsion.
  Thus $D$ has rank $d - 1$, and $z \neq 0$.

  \emph{3. Comparing the elements at level $n$.}
  Fix $n \geq 0$. We will prove that the image $z_n$ of $z$
  and $P_\zeta e_{\vartheta^+}\epsilon_{F, n}^{\rs}$
  differ by multiplication by a unit of $\Lambda_n$.
  Recall that $\gamma$ is the topological generator of
  $\Gamma$ chosen in Section~\ref{subsec:Iwasawa-modules}
  and $\Lambda_n := \Lambda/(\gamma^{p^n} - 1)$.
  Let $D_n$ be the image of $D$ under
  $\Lambda \longrightarrow \Lambda_n$.
  The splitting of $q$ is defined over $\Lambda$,
  so it also splits the target of
  $f\otimes_\Lambda\Lambda_n$.
  Bullach--Hofer
  \cite[Proposition~3.2(c) and Remark~3.3(a)]{bullach-hofer}
  therefore identify
  \[
    \ker(D_n : \Lambda_n^d \longrightarrow \Lambda_n^{d - 1})
    = P_\zeta e_{\vartheta^+}
    (\calo_{K_n, \Sigma_F}^\times\otimes_\Z\calo).
  \]
  In particular, $z_n$ is an element of this module of units.

  To compare $z_n$ with the Rubin--Stark element, we
  describe the map from the determinant module to units.
  Reduce the chosen splitting of $q$ to level $n$.
  Take the image of $1$ under the reduced splitting as
  the first target basis vector, followed by the reduced
  basis of $\ker q$.
  The first vector maps to $1$ under $q\otimes 1$,
  hence to $P_\zeta e_{\vartheta^+}w_n$ in the free module
  on complex places.
  We also reduce the source basis
  of $f$ to level $n$.

  We regard $[\Lambda_n^d\xrightarrow{f\otimes 1}
  \Lambda_n^d]$ as a complex in degrees $0$ and $1$.
  Its determinant module is the free rank-$1$
  $\Lambda_n$-module
  \[
    \Det_{\Lambda_n}
    [\Lambda_n^d\xrightarrow{f\otimes 1}\Lambda_n^d]
    = \bigwedge_{\Lambda_n}^d\Lambda_n^d
    \otimes_{\Lambda_n}
    \Hom_{\Lambda_n}\bigl(
    \bigwedge_{\Lambda_n}^d\Lambda_n^d, \Lambda_n\bigr).
  \]
  The first exterior power uses the source of
  $f\otimes 1$, and the exterior power inside $\Hom$
  uses its target.
  Define the $\Lambda_n$-linear map
  \[
    \Phi_n : \Det_{\Lambda_n}
    [\Lambda_n^d\xrightarrow{f\otimes 1}\Lambda_n^d]
    \longrightarrow \ker D_n
  \]
  by sending the determinant generator from the ordered
  source and target bases to $(-1)^{d - 1}z_n$.
  To identify $\Phi_n$ in the notation of
  Burns--Kurihara--Sano
  \cite[Section~2D and Proposition~2.6(iii)]
  {burns-kurihara-sano}, we use their map
  $\pi^V_{K_n/F, \Sigma_F, \varnothing}$,
  where $V := \set{\restr{w_n}{F}}$ consists of the
  complex place of $F$, with chosen lift $w_n$.
  After extending scalars from $\Z_p$ to $\calo$ and
  applying $P_\zeta e_{\vartheta^+}$, its source is
  $\Det_{\Lambda_n}
  [\Lambda_n^d\xrightarrow{f\otimes 1}\Lambda_n^d]$
  and its target in rank $1$ is
  \[
    P_\zeta e_{\vartheta^+}
    (\calo_{K_n, \Sigma_F}^\times\otimes_\Z\calo)
    = \ker D_n.
  \]

  After applying $P_\zeta e_{\vartheta^+}$, their formula
  uses a first target basis vector mapping to
  $P_\zeta e_{\vartheta^+}(w_n - w'_n)$.
  To meet this condition, we use the natural surjection
  \[
    \Lambda_n^d \longrightarrow
    P_\zeta e_{\vartheta^+}(X_n\otimes_\Z\calo)
  \]
  to choose a lift of
  $P_\zeta e_{\vartheta^+}(w_n - w'_n)$.
  This lift and our first target basis vector both
  map to $1$ under $q\otimes 1$, so their difference
  belongs to $\ker(q\otimes 1)$.
  Replacing the first target basis vector by this lift
  leaves the other target basis vectors unchanged.
  The change of basis has determinant $1$ and does
  not change $D_n$, so it preserves the determinant
  generator and $z_n$.
  We now apply the formula of Burns--Kurihara--Sano
  \cite[Proposition~2.6(iii)]{burns-kurihara-sano}
  with $r = 1$ to obtain the image $(-1)^{d - 1}z_n$.
  Thus $\Phi_n$ is the extension of
  $\pi^V_{K_n/F, \Sigma_F, \varnothing}$ to $\calo$,
  followed by passage to the
  $P_\zeta e_{\vartheta^+}$-component.

  Let $F_\infty^{(2)}$ be the compositum of the
  $\Z_p$-extensions of $F$.
  Combining Bullach--Hofer
  \cite[Theorem~5.8]{bullach-hofer} with
  Burns--Kurihara--Sano
  \cite[Sections~2C--2D]{burns-kurihara-sano},
  we obtain the following consequence for $K_n/F$:
  if the prime $p$ splits in $F$ or if
  $\Gal(K_nF_\infty^{(2)}/F)$ has no element of order $p$,
  then the determinant module
  $\Det_{\Lambda_n}[\Lambda_n^d\xrightarrow{f\otimes 1}
  \Lambda_n^d]$ has a generator whose image under
  $\Phi_n$ is
  $P_\zeta e_{\vartheta^+}\epsilon_{F, n}^{\rs}$.

  In our tower, $F\Q_n \subset F_\infty^{(2)}$ and
  $[K_0 : F] = 3[L_0 : F]$ is prime to $p$.
  Since $F_\infty^{(2)}/F$ is pro-$p$ and
  $\Gal(F_\infty^{(2)}/F) \cong \Z_p^2$, we have
  $K_0 \cap F_\infty^{(2)} = F$ and
  \[
    \Gal(K_nF_\infty^{(2)}/F)
    \cong \Gal(K_0/F)\times\Z_p^2.
  \]
  Since $p \nmid [K_0 : F]$, the group
  $\Gal(K_0/F)\times\Z_p^2$ has no element of order $p$.
  Thus $(-1)^{d - 1}z_n$ and
  $P_\zeta e_{\vartheta^+}\epsilon_{F, n}^{\rs}$
  are images under $\Phi_n$ of two generators of
  $\Det_{\Lambda_n}[\Lambda_n^d\xrightarrow{f\otimes 1}
  \Lambda_n^d]$. The determinant module is free of rank
  $1$ over $\Lambda_n$, so its generators differ by a
  unit. Absorbing the sign into the unit, we obtain
  \[
    a_n z_n
      = P_\zeta e_{\vartheta^+}\epsilon_{F, n}^{\rs}
      \qquad\text{for some }a_n \in \Lambda_n^\times.
  \]

  \emph{4. Passing to the inverse limit.}
  We prove that a single $a \in \Lambda^\times$ satisfies
  $az = W_K^{-1}\mathbf{v}$.
  Bullach--Hofer \cite[Remark~3.3(b)]{bullach-hofer}
  identify the map on kernels obtained by applying
  $\Lambda_{n + 1} \longrightarrow \Lambda_n$ to every
  coordinate with $\rmn_{K_{n + 1}/K_n}\otimes 1$ on
  $P_\zeta e_{\vartheta^+}
  (\calo_{K_{n + 1}, \Sigma_F}^\times\otimes_\Z\calo)$.
  Since reduction sends each signed minor of $D_{n + 1}$
  to the corresponding signed minor of $D_n$, we have
  $(\rmn_{K_{n + 1}/K_n}\otimes 1)(z_{n + 1}) = z_n$.
  Since $\epsilon_{F, m}^{\rs} = W_K^{-1}\epsilon_{F, m}$
  for every $m$ by step~1,
  we can use the norm relation for $\epsilon_{F, n}$ in
  Lemma~\ref{lem:source-unit-norm-relations} to obtain
  \[
    (\rmn_{K_{n + 1}/K_n}\otimes 1)
    (P_\zeta e_{\vartheta^+}\epsilon_{F, n + 1}^{\rs})
      = P_\zeta e_{\vartheta^+}\epsilon_{F, n}^{\rs}.
  \]

  For each $n$, let $C_n$ be the set of $a \in \Lambda^\times$
  whose image in $\Lambda_n$ satisfies
  $a z_n = P_\zeta e_{\vartheta^+}\epsilon_{F, n}^{\rs}$.
  Each $C_n$ is closed and nonempty: the unit $a_n$ from
  step~3 lifts through the surjection of local rings
  $\Lambda \longrightarrow \Lambda_n$.
  Because both $z_n$ and
  $P_\zeta e_{\vartheta^+}\epsilon_{F, n}^{\rs}$ satisfy
  the norm relations, every $a \in C_{n + 1}$ also satisfies
  $a z_n = P_\zeta e_{\vartheta^+}\epsilon_{F, n}^{\rs}$.
  Thus $C_{n + 1} \subseteq C_n$.
  By compactness of $\Lambda^\times$, we can choose
  $a \in \bigcap_{n \geq 0} C_n$; for this $a$, the equalities
  at every level combine to $az = W_K^{-1}\mathbf{v}$.
  Since $a$ and $W_K$ are units and $z \neq 0$,
  $\Lambda z = \Lambda\mathbf{v}$ and $\mathbf{v} \neq 0$.

  \emph{5. Characteristic ideals.}
  It remains to compute $\Char_\Lambda(\cale_F/\Lambda z)$.
  Let $\mfq$ be a height-$1$ prime of $\Lambda$.
  The ring $\Lambda_\mfq$ is a discrete valuation ring.
  We use Smith normal form over $\Lambda_\mfq$ to choose
  bases in which the matrix of
  $D : \Lambda_\mfq^d \longrightarrow
  \Lambda_\mfq^{d - 1}$ is
  \[
    \bigl(\operatorname{diag}(s_1, \ldots, s_{d - 1})
    \quad 0\bigr),
  \]
  with each $s_i \neq 0$.
  Thus $Z_\mfq \cong
  \bigoplus_{i = 1}^{d - 1}\Lambda_\mfq/(s_i)$.
  The kernel $(\cale_F)_\mfq$ is generated by the last
  vector of the new source basis.
  The changes of bases have unit determinants, so $z$
  is a unit times $\prod_{i = 1}^{d - 1}s_i$ times that
  last basis vector. Consequently,
  \[
    \length_{\Lambda_\mfq}
    ((\cale_F/\Lambda z)_\mfq)
    = \sum_{i = 1}^{d - 1}
    \length_{\Lambda_\mfq}(\Lambda_\mfq/(s_i))
    = \length_{\Lambda_\mfq}(Z_\mfq).
  \]
  Since this holds at every height-$1$ prime,
  $\Char_\Lambda(\cale_F/\Lambda z) = \Char_\Lambda(Z)$.
  Using $\Lambda z = \Lambda\mathbf{v}$ and the exact
  sequence $0 \longrightarrow A_F \longrightarrow Z
  \longrightarrow B_F \longrightarrow 0$, we conclude
  \[
    \Char_\Lambda(\cale_F/\Lambda\mathbf{v})
    = \Char_\Lambda(A_F)\Char_\Lambda(B_F). \qedhere
  \]
\end{proof}

\subsection{Proof of the Iwasawa main conjecture over \texorpdfstring{$k$}{k}}
\label{subsec:proof-Iwasawa}

We combine the characteristic ideal identities in
Propositions~\ref{prop:rational-source-determinant} and
\ref{prop:elliptic-source-identification}, using the
decompositions of Proposition~\ref{prop:derived-Artin-decomposition}.

\begin{proof}[Proof of Theorem~\ref{thm:explicit-Iwasawa-Stark}]
  Suppose first that $\vartheta$ is odd. Then
  $\cale_\Q = 0$, and
  Lemma~\ref{lem:elliptic-unit-conjugation}(i)
  identifies $\cale_\Sigma/\calc_\vartheta^\Sigma$ with
  $\cale_F/t_\vartheta\Lambda\mathbf{v}$.
  Since $t_\vartheta$ and $\mathbf{v}$ are nonzero,
  $\calc_\vartheta^\Sigma$ is free of rank $1$, and
  $\cale_\Sigma/\calc_\vartheta^\Sigma$ is torsion.
  To compare $\cale_F/t_\vartheta\Lambda\mathbf{v}$ with
  $\cale_F/\Lambda\mathbf{v}$, we use the exact sequence
  \[
    \begin{tikzcd}[column sep = small]
      0 \arrow[r] & \Lambda/(t_\vartheta) \arrow[r] &
      \cale_F/t_\vartheta\Lambda\mathbf{v} \arrow[r] &
      \cale_F/\Lambda\mathbf{v} \arrow[r] & 0.
    \end{tikzcd}
  \]
  The first nonzero map sends the class of $b \in \Lambda$
  to the class of $b\mathbf{v}$; it is injective because
  $\mathbf{v}$ is nonzero in the free rank-$1$ module
  $\cale_F$. The next map is the natural quotient map.
  Multiplicativity of characteristic ideals gives
  \[
    \Char_\Lambda(\cale_\Sigma/\calc_\vartheta^\Sigma)
    = (t_\vartheta)
    \Char_\Lambda(\cale_F/\Lambda\mathbf{v}).
  \]
  We now express $(t_\vartheta)$ and
  $\Char_\Lambda(\cale_F/\Lambda\mathbf{v})$ in terms of
  ideal classes by
  Propositions~\ref{prop:rational-source-determinant}(i)
  and \ref{prop:elliptic-source-identification}.
  Since $A_\Sigma \cong A_\Q\oplus A_F$ and
  $B_\Sigma \cong B_\Q\oplus B_F$,
  we conclude that
  \[
    \Char_\Lambda(\cale_\Sigma/\calc_\vartheta^\Sigma)
      = \Char_\Lambda(A_\Sigma)\Char_\Lambda(B_\Sigma)
  \]
  for odd $\vartheta$.

  If $\vartheta$ is even, we use
  Lemma~\ref{lem:elliptic-unit-conjugation}(ii) to identify
  the generators in the two components of
  Proposition~\ref{prop:derived-Artin-decomposition}:
  \[
    \cale_\Sigma/\calc_\vartheta^\Sigma
    \cong (\cale_\Q/\Lambda\mathbf{u}_{\R, \vartheta})
    \oplus (\cale_F/\Lambda\mathbf{v}).
  \]
  The elements $\mathbf{u}_{\R, \vartheta}$ and
  $\mathbf{v}$ are nonzero in the free rank-$1$ modules
  $\cale_\Q$ and $\cale_F$. Thus $\calc_\vartheta^\Sigma$
  is free of rank $2$, and
  $\cale_\Sigma/\calc_\vartheta^\Sigma$ is torsion.
  We can therefore compute the characteristic ideal of
  this direct sum from
  Propositions~\ref{prop:rational-source-determinant}(ii)
  and \ref{prop:elliptic-source-identification}, obtaining
  \[
    \Char_\Lambda(\cale_\Sigma/\calc_\vartheta^\Sigma)
      = \Char_\Lambda(A_\Sigma)\Char_\Lambda(B_\Sigma)
  \]
  for even $\vartheta$.
\end{proof}

To finish Theorem~\ref{thm:intro-determinant-main}, we
compare $e_\vartheta(\epsilon_n\otimes 1)_n$,
$e_\vartheta(u_{\R, n}\otimes 1)_n$, and
$e_\vartheta(u_{\C, n}\otimes 1)_n$ with
$\mathbf{u}_\vartheta$, $\mathbf{u}_{\R, \vartheta}$,
and $\mathbf{u}_{\C, \vartheta}$, respectively.

\begin{proof}[Proof of
  Theorem~\ref{thm:intro-determinant-main}]
  
  Recall that we have the submodule
  \[
    \calc_\vartheta^\Sigma :=
    \begin{cases}
      \Lambda e_\vartheta(\epsilon_n\otimes 1)_n
      & \text{if $\vartheta(c) = -1$}, \\
      \Lambda e_\vartheta(u_{\R, n}\otimes 1)_n
      + \Lambda e_\vartheta(u_{\C, n}\otimes 1)_n
      & \text{if $\vartheta(c) = 1$}.
    \end{cases}
  \]
  For odd $\vartheta$, we have
  $\Lambda e_\vartheta(\epsilon_n\otimes 1)_n
  = \Lambda\mathbf{u}_\vartheta$ by
  Lemma~\ref{lem:Iwasawa-Stark-generators}.
  For even $\vartheta$, the generators are
  $\mathbf{u}_{\R, \vartheta}
  = e_\vartheta(u_{\R, n}\otimes 1)_n$ and
  $\mathbf{u}_{\C, \vartheta}
  = e_\vartheta(u_{\C, n}\otimes 1)_n$.
  Thus $\calc_\vartheta^\Sigma$ is the submodule in
  Theorem~\ref{thm:explicit-Iwasawa-Stark} in both cases.

  By Lemma~\ref{lem:Iwasawa-freeness},
  $\cale_\Sigma$ is free of rank $r_\vartheta$,
  and $A_\Sigma$ and $B_\Sigma$ are
  finitely generated torsion $\Lambda$-modules.
  In Theorem~\ref{thm:explicit-Iwasawa-Stark},
  we proved that
  $\calc_\vartheta^\Sigma$ is free of rank $r_\vartheta$
  and $\cale_\Sigma/\calc_\vartheta^\Sigma$ is torsion,
  with
  \[
    \Char_\Lambda(\cale_\Sigma/\calc_\vartheta^\Sigma)
    = \Char_\Lambda(A_\Sigma)\Char_\Lambda(B_\Sigma). \qedhere
  \]
\end{proof}

\appendix

\renewcommand{\thesection}{\arabic{section}}

\renewcommand{\theHsection}{appendix.\arabic{section}}

\section{Dasgupta's example}
\label{app:dasgupta}

We verify the Stark unit predicted by Dasgupta
\cite[Proposition~5.5.3]{dasgupta1999}
for a cyclic sextic extension of a complex cubic field.
Using Proposition~\ref{prop:local-correction-factors},
we construct a product of conjugates of an elliptic unit
whose logarithms are $18$ times those required of the
predicted Stark unit.
We then use the Mahler measure bound of Mossinghoff--Rhin--Wu
\cite[Theorem~1.1]{mossinghoff-rhin-wu}
to identify this unit with the $18$-th power of Dasgupta's
predicted unit up to a root of unity.\footnote{
  Alternatively, one can use the height bound of \cite{voutier}
  instead of the Mahler measure bound of Mossinghoff--Rhin--Wu
  \cite{mossinghoff-rhin-wu}
  (which an earlier version of this paper did),
  but that proof is slightly longer.
}
Finally, we apply
Proposition~\ref{prop:elliptic-gamma-Stark-model} to obtain
units with the absolute values of the corresponding
elliptic gamma values.

Let $\beta$ be a root of $P = x^3 - x^2 + 5x + 1$.
The corresponding complex cubic field $k = \Q(\beta)$
has class number $3$, discriminant $-588$,
and quadratic resolvent $\Q(\sqrt{-588}) = \Q(\sqrt{-3})$.
Its Hilbert class field is $k(\zeta_7 + \zeta_7^{-1})$,
where $\zeta_m$ denotes a primitive $m$-th root of unity.
Adjoining $\zeta_3$ gives the cyclotomic extension
$K = k(\zeta_3, \zeta_7 + \zeta_7^{-1})
= \Q(\beta, \zeta_3, \zeta_7 + \zeta_7^{-1})$
of degree $6$ over $k$.
Dasgupta \cite[Proposition~5.5.3]{dasgupta1999}
predicts that the Stark unit
(up to a root of unity)
for the cyclic sextic extension
$K := k(\zeta_3, \zeta_7 + \zeta_7^{-1})$
and the set $S := \set{\infty_\R, \infty_\C, \mfp_1}$,
where $\mfp_1 := (3, \beta + 1)$,
is a root of
\begin{align}
  \label{eq:Dasgupta-polynomial}
  f(X) := X^6
  & -(\beta^2 - 4\beta + 1)X^5
  + (-7\beta^2 + 10\beta + 5)X^4
  - (22\beta^2 + 17\beta + 6)X^3 \\
  & + (-7\beta^2 + 10\beta + 5)X^2
  - (\beta^2 - 4\beta + 1)X + 1.
  \nonumber
\end{align}

In our notation from Section~\ref{sec:relative-cubic},
we have that
\begin{itemize}
  \item $k = \Q(\beta)$ is a complex cubic field with
    quadratic resolvent $F = \Q(\sqrt{-3})$;
  \item $E = kF = \Q(\beta, \sqrt{-3})$;
  \item $L = \Q(\zeta_3, \zeta_7 + \zeta_7^{-1})$
    is the maximal abelian subextension of $K/\Q$;
  \item $K = kL = \Q(\beta, \zeta_3, \zeta_7 + \zeta_7^{-1})$
    is a cyclotomic extension of $k$;
  \item $L^+ = \Q(\zeta_7 + \zeta_7^{-1})$;
  \item $K^+ = kL^+ = k(\zeta_7 + \zeta_7^{-1})$;
  \item $W_K = W_L = 6$; and
  \item $S = \set{\infty_\R, \infty_\C, (3, \beta + 1)}$.
\end{itemize}

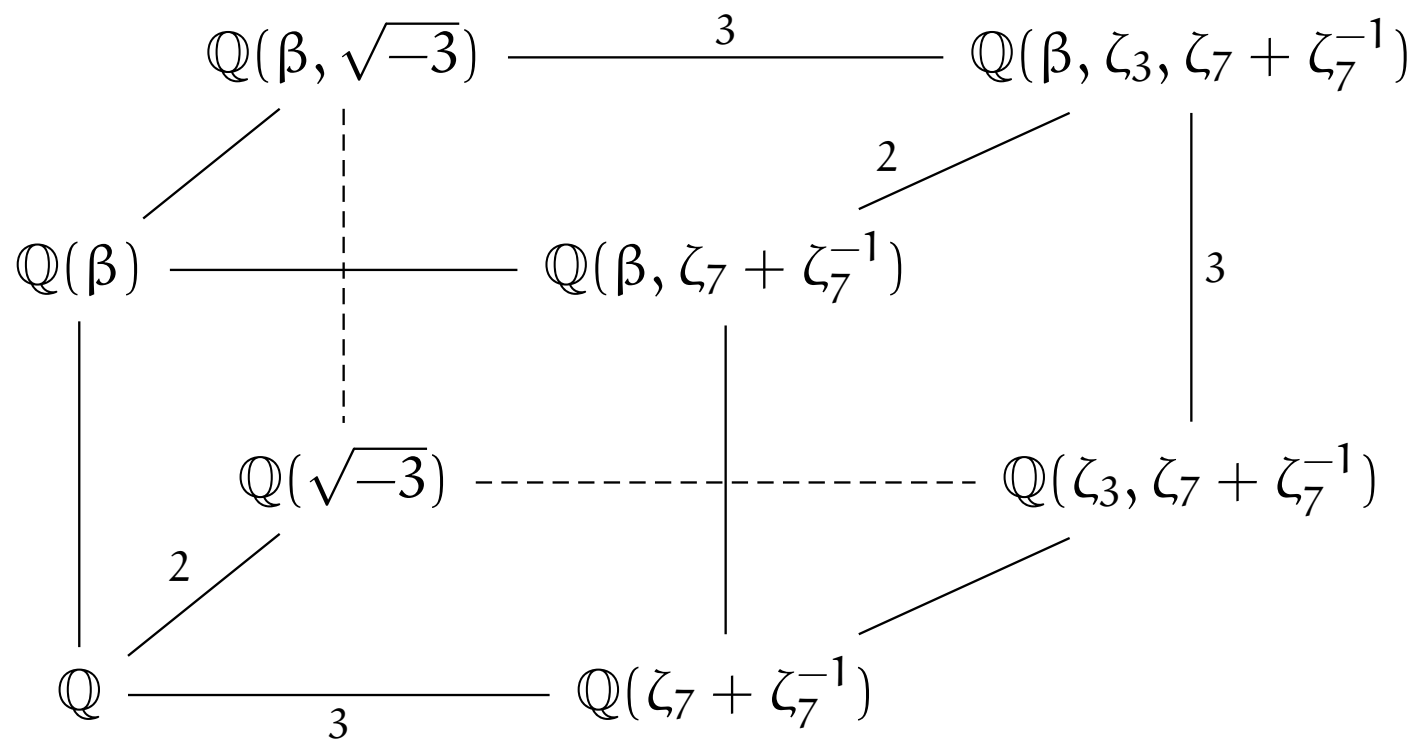
\begin{figure}[ht]
  \[
    \begin{tikzcd}[row sep = 1.7em, column sep = 0.2em]
      & \Q(\beta, \sqrt{-3}) \arrow[rr, dash, "3"]
      \arrow[dd, dashed, no head]
      && \Q(\beta, \zeta_3, \zeta_7 + \zeta_7^{-1}) \arrow[dd, dash, "3"] \\
      \Q(\beta) \arrow[ur, dash] \arrow[rr, dash]
      \arrow[dd, dash]
      && \Q(\beta, \zeta_7 + \zeta_7^{-1}) \arrow[ur, dash, "2", near start]
      \arrow[dd, dash] \\
      & \Q(\sqrt{-3}) \arrow[rr, dashed, no head] && \Q(\zeta_3, \zeta_7 + \zeta_7^{-1}) \\
      \Q \arrow[rr, dash, "3"'] \arrow[ur, dash, "2"]
      && \Q(\zeta_7 + \zeta_7^{-1}) \arrow[ur, dash]
    \end{tikzcd}
  \]
  \caption{The field extensions in Dasgupta's example.}
  \label{figure:Dasgupta-fields}
\end{figure}

We summarize the field extensions in
Figure~\ref{figure:Dasgupta-fields},
which is a specialization of
Figure~\ref{figure:field-extensions-cube}.

By Dasgupta \cite[Proposition~5.1.1]{dasgupta1999},
$\mfp_1$ is the only finite prime ramified in $K/k$.
Thus $S$ contains all places required by the integral
rank-$1$ abelian Stark conjecture. The stronger condition
on $S$ in Theorem~\ref{thm:intro-rank-one} requires every
place above $2$, $3$, and $7$, since these primes divide
$\disc(k) = -2^2\cdot 3\cdot 7^2$.
Dasgupta's set contains only $\mfp_1$ among the finite
places, so that theorem does not apply directly.
We instead use Proposition~\ref{prop:local-correction-factors}
to compare the Euler factors for this $S$.
Let $S_F$ consist of the complex place of $F$ and all
finite places above $2$ and $7$.
The set $S_F$ contains all places ramified in $K/F$.
Fix a place $w$ of $K$ above $\infty_\C$,
and let
$\epsilon_F \in \calo_K^\times$ be the elliptic unit of
Theorem~\ref{thm:imaginary-quadratic-Stark} for
$(K/F, S_F, w)$.

Set $\gamma := [(2, \beta - 1), K/k]$. The image of
$\gamma$ has order $3$ on the Hilbert class field
$k(\zeta_7 + \zeta_7^{-1})$ and order $2$ on $kF$, since
$2$ is inert in $F$. Thus $\gamma$ generates
$\Gal(K/k)$ and $c = \gamma^3$ is complex conjugation.
We next express $L_S'(0, \sigma)$ in terms of the
logarithms of conjugates of $\epsilon_F$.

\begin{samepage}
\begin{proposition}
  \label{prop:Dasgupta-exact-unit}
  The unit
  $u_{18} := \epsilon_F^{(7 + 5\gamma^{-1} + \gamma^{-2})(1 - c)}
  \in \calo_K^\times$
  satisfies
  \[
    L_S'(0, \sigma)
      = -\frac{1}{108} \log\verts{u_{18}^\sigma}_w
  \]
  for every $\sigma \in \Gal(K/k)$. Furthermore, $u_{18}$ has
  absolute value $1$ at every place above the real place
  of $k$.
\end{proposition}
\end{samepage}

\begin{proof}
  We use the character sums
  $\calr_w(x, \chi) := \sum_{\sigma \in \Gal(K/k)}\chi(\sigma)
  \log \lvert x^\sigma \rvert_w$
  from Section~\ref{subsec:rank-one} to prove that
  $\calr_w(u_{18}, \chi) = -108L_S'(0, \chi)$
  for odd characters $\chi$ of $\Gal(K/k)$,
  and then apply Fourier inversion.

  We first compute the Euler factors at $7$ needed for the
  decomposition of partial $L$-functions in
  Proposition~\ref{prop:local-correction-factors}.
  Let $\mfp_2 := (3, \beta - 1)$,
  $\mfp_3 := (2, \beta - 1)$, and
  $\mfp_4 := (7, \beta + 2)$. The factorizations
  $(\beta - 1) = \mfp_2\mfp_3$ and
  $(\beta + 2) = \mfp_2\mfp_4$ in Dasgupta
  \cite[Section~5.1]{dasgupta1999} show that
  $\mfp_3$ and $\mfp_4$ have the same ideal class.
  Their Artin symbols agree on the Hilbert class field.
  Since $7$ splits in $F$, whereas $2$ is inert, this
  gives $[\mfp_4, K/k] = \gamma^4$.

  For $j \in \set{1, 3, 5}$, let $\chi_j$ be the
  character of $\Gal(K/k)$ determined by
  $\chi_j(\gamma) = \zeta_6^j$, and set
  $\psi_j := \chi_{j, F}\omega_F$. Inertia at
  $\mfp_1$ is generated by $c$, so the odd character
  $\chi_j$ has trivial Euler factor there. Hence
  $L_S(s, \chi_j) = L(s, \chi_j)$. At $2$, the Euler
  factors over $k$ and $\Q$ agree, so the product of
  the Euler factors of $L(s, \psi_j)$ above $2$ is $1$.
  At $7$, the Euler factor over $k$ is
  $(1 - \chi_j(\gamma^4)7^{-s})^{-1}$.
  For $j \in \set{1, 5}$, the character $\chi_{j, \Q}$ has
  conductor $21$ and hence trivial Euler factor at $7$.
  For $j = 3$, it has conductor $3$ and value $1$
  at $7$, so its Euler factor agrees with the Euler factor
  over $k$. Therefore
  $L_{S_F}(s, \psi_j) = E_j(s)L(s, \psi_j)$,
  where $E_3(s) = 1$ and
  $E_j(s) = 1 - \chi_j(\gamma^4)7^{-s}$ for $j \in \set{1, 5}$.
  Differentiating the decomposition of Artin $L$-functions
  in Equation~\ref{eq:complete-L-factorization} gives
  \begin{equation}
    \label{eq:Dasgupta-factorization}
    L_S'(0, \chi_j)
      = \frac{L(0, \chi_{j, \Q})}{E_j(0)}
        L_{S_F}'(0, \psi_j).
  \end{equation}

  We compute the remaining factors using
  $L(0, \chi) = -f_\chi^{-1}\sum_{a = 1}^{f_\chi}a\chi(a)$
  for a primitive Dirichlet character of conductor $f_\chi$.
  Here,
  \begin{align*}
    L(0, \chi_{1, \Q}) &= 2\zeta_6, \\
    L(0, \chi_{3, \Q}) &= \frac{1}{3}, \\
    L(0, \chi_{5, \Q}) &= 2\zeta_6^{-1},
    \intertext{so for $j \in \{1, 3, 5\}$:}
    \frac{L(0, \chi_{j, \Q})}{E_j(0)}
      &= \frac{7 + 5\zeta_6^j + \zeta_6^{2j}}{9}.
  \end{align*}
  The value of $(7 + 5\gamma^{-1} + \gamma^{-2})(1 - c)$
  under $\chi_j^{-1}$ is $18L(0, \chi_{j, \Q})/E_j(0)$.
  Choose $\tau$ so that $c w = \tau w$. Since $S_F$ is
  invariant under complex conjugation, we have
  $\calr_w(\epsilon_F, \chi_j) = -6L_{S_F}'(0, \psi_j)$
  by Lemma~\ref{lem:elliptic-regulator-comparison}.
  We now substitute into
  Equation~\ref{eq:character-logarithm-action} and use
  Equation~\ref{eq:Dasgupta-factorization} to obtain
  $\calr_w(u_{18}, \chi_j) = -108L_S'(0, \chi_j)$.
  For even characters both sides vanish, so by Fourier
  inversion we conclude that
  \[
    L_S'(0, \sigma) = -(1/108)\log\verts{u_{18}^\sigma}_w.
  \]

  Since $c(u_{18}) = u_{18}^{-1}$ and $c$ preserves the
  absolute value at every place
  above the real place of $k$,
  the unit $u_{18}$ has absolute value $1$ at each
  such place.
\end{proof}

Dasgupta \cite[Section~5.5]{dasgupta1999} identifies, for
the chosen place $w$,
a distinguished root $\epsilon_w \in K$ of
the polynomial $f$ from
Equation~\ref{eq:Dasgupta-polynomial}.
Dasgupta \cite[Propositions~5.5.2 and 5.5.3]{dasgupta1999}
proves that $\epsilon_w$ is a unit with absolute value $1$
above the real place of $k$ and that
$K(\epsilon_w^{1/6})/k$ is abelian. Moreover,
in \cite[Proposition~5.5.3]{dasgupta1999}
he proves the bound
\begin{equation}
  \label{eq:Dasgupta-certified-bound}
  \verts{\log\verts{\epsilon_w^\sigma}_w
    + 6L_S'(0, \sigma)}
    < 10^{-25}
\end{equation}
for every $\sigma \in \Gal(K/k)$.
The equality
$L_S'(0, \sigma)
= -\frac{1}{6}\log\verts{\epsilon_w^\sigma}_w$
is the assertion left open in
Dasgupta \cite[Proposition~5.5.3]{dasgupta1999}.
We verify Dasgupta's prediction
by comparing $\epsilon_w^{18}$ with
our $u_{18}$ from Proposition~\ref{prop:Dasgupta-exact-unit}.

\begin{theorem}
  \label{thm:Dasgupta-exact}
  The distinguished root $\epsilon_w \in \calo_K^\times$
  of the polynomial $f$ in Equation~\ref{eq:Dasgupta-polynomial}
  is a Stark unit for $(K/k, S, w)$.
  In particular,
  \[
    L_S'(0, \sigma)
      = -\frac{1}{6} \log\lvert \epsilon_w^\sigma\rvert_w
  \]
  for every $\sigma \in \Gal(K/k)$.
\end{theorem}

\begin{proof}
  To prove the regulator identity, we show that
  $\delta := u_{18}\epsilon_w^{-18}$ is a root of unity.
  We compare the logarithms of $u_{18}$ in
  Proposition~\ref{prop:Dasgupta-exact-unit}
  with Dasgupta's bound for $\epsilon_w$
  in Equation~\ref{eq:Dasgupta-certified-bound}
  to obtain
  $\lvert\log\lvert\delta^\sigma\rvert_w\rvert
  < 18\cdot10^{-25}$.
  The global unit $\delta$ has absolute value $1$ above
  the real place of $k$. Since $[K : \Q] = 18$, its
  absolute logarithmic height satisfies
  \[
    18h(\delta) < 6\cdot18\cdot10^{-25}
    = 1.08\cdot10^{-23}.
  \]

  Now suppose that $\delta$ is not a root of unity.
  In particular, its minimal polynomial
  has degree $d := [\Q(\delta) : \Q] \leq 18$
  and Mahler measure $\exp(d h(\delta)) > 1$.
  But Mossinghoff--Rhin--Wu
  \cite[Theorem~1.1]{mossinghoff-rhin-wu}
  prove that the Mahler measure of an integer polynomial
  of degree less than $56$
  is either equal to $1$ or
  is at least Lehmer's number $1.1762808\ldots$;
  so in this case
  \begin{align*}
    d h(\delta) &\geq \log(1.1762808\ldots) > 0.16.
    \intertext{But this contradicts}
    dh(\delta) &\leq 18h(\delta) < 1.08 \cdot 10^{-23}.
  \end{align*}
  Therefore $\delta \in \mu(K)$,
  so $u_{18} = \zeta\epsilon_w^{18}$ for some
  $\zeta \in \mu(K)$.
  Since $\log\verts{u_{18}^\sigma}_w
  = 18\log\verts{\epsilon_w^\sigma}_w$,
  we conclude from
  Proposition~\ref{prop:Dasgupta-exact-unit} that
  $L_S'(0, \sigma) = -(1/6)\log\verts{\epsilon_w^\sigma}_w$.
  Dasgupta \cite[Propositions~5.5.2 and 5.5.3]{dasgupta1999}
  proves that $\epsilon_w$ has absolute value $1$ above the
  real place of $k$ and that $K(\epsilon_w^{1/6})/k$ is abelian.
\end{proof}

Finally, we use $\epsilon_w$ in this example
to construct units with the
same absolute value
as Bergeron--Charollois--Garc\'{i}a's
elliptic gamma values.
Bergeron--Charollois--Garc\'{i}a
\cite[Section~2.4]{bergeron-charollois-garcia} identify
$K$ with the narrow ray class field $k(\mfp_1)$.
Set $\mff := \mfp_1$,
and let $\mfa$ be an ideal of $k$ coprime to $6\rmn(\mff)$
whose norm $\rmn(\mfa)$ is an odd rational prime.
In Theorem~\ref{thm:Dasgupta-exact}, we proved that
$\epsilon_w$ is a Stark unit for
$S_{\mfp_1} = \set{\infty_\R, \infty_\C, \mfp_1}$.
Applying Proposition~\ref{prop:elliptic-gamma-Stark-model},
we obtain $u_{\mfp_1, \mfa} \in \calo_K^\times$ with
$u_{\mfp_1, \mfa}^6
= \epsilon_w^{\rmn(\mfa) - \sigma_\mfa}$
and
\[
  \verts{\Gamma_\mfa(\mfp_1\mfb^{-1})}_\C
  = \verts{u_{\mfp_1, \mfa}^{\sigma_\mfb}}_w
\]
for every ideal $\mfb$ coprime to $\mfp_1$.
Consequently, $\lvert\Gamma_\mfa(\mfp_1\mfb^{-1})\rvert$
is an algebraic unit for every $\mfb$ coprime to $\mfp_1$.
Choose an embedding $\iota : K \longrightarrow \C$ inducing $w$.
It remains to prove the equality
$\Gamma_\mfa(\mfp_1\mfb^{-1})
= \iota((\xi u_{\mfp_1, \mfa})^{\sigma_\mfb})$
for one $\xi \in \mu(K)$ and every $\mfb$ coprime to $\mfp_1$,
as in Proposition~\ref{prop:elliptic-gamma-phase}.


\bibliographystyle{alpha}
\bibliography{bibliography-complex-cubic}

\end{document}